\documentclass[10pt]{article}

\usepackage[left=1in, right=1in, bottom=1in, top=1in, footskip=0.3in, 
]{geometry}

\usepackage{xcolor}
\definecolor{LinkColor}{rgb}{0,0,1}
\definecolor{LinkColor2}{rgb}{1,0,0}
\definecolor{lbcolor}{rgb}{0.85,0.85,0.85}
\definecolor{FrameColor}{rgb}{0.85,0.85,0.85}

\usepackage[ngerman, english]{babel}
\usepackage[utf8]{inputenc}
\usepackage{enumerate}
\usepackage[normalem]{ulem}
\usepackage{setspace}

\def\pskip{\\[-3mm]}

\usepackage{multicol}
\usepackage[babel,french=guillemets,german=swiss]{csquotes}

\usepackage{array}
\usepackage{booktabs}
\usepackage{framed}
\usepackage{rotating}
\usepackage{longtable}
\usepackage{multirow}
\usepackage{tabularx}
\newcolumntype{L}[1]{>{\raggedright\arraybackslash}p{#1}} 
\newcolumntype{C}[1]{>{\centering\arraybackslash}p{#1}} 
\newcolumntype{R}[1]{>{\raggedleft\arraybackslash}p{#1}} 

\usepackage{graphicx,import}
\usepackage{caption}

\usepackage{amssymb}
\usepackage{dsfont}
\usepackage{nicefrac}
\usepackage{empheq}
\allowdisplaybreaks
\usepackage{MnSymbol}

\usepackage[%
pdftitle={Titel},%
pdfauthor={Autor},%
pdfcreator={LaTeX, LaTeX with hyperref and KOMA-Script},
pdfsubject={Betreff}, 
pdfkeywords={Keywords}
]{hyperref} 

\hypersetup{%
	colorlinks=true,
	linkcolor=LinkColor,%
	anchorcolor=LinkColor,%
	citecolor=LinkColor2,%
	filecolor=LinkColor,%
	menucolor=LinkColor,%
	urlcolor=LinkColor,%
}

\usepackage[savemem]{listings}
\usepackage{paralist}
{\begin{list}{$\diamondsuit$}{}}%
	{\end{list}}

\usepackage{amsthm}
\usepackage{upgreek}

\newtheoremstyle{tstyle}
{15pt}	
{5pt}	
{\itshape}	
{}	
{\bfseries}	
{.}	
{0.5em}	
{}	

\theoremstyle{tstyle}

\newtheorem{theorem}{Theorem}[section]
\newtheorem{lemma}[theorem]{Lemma}

\newtheorem{corollary}[theorem]{Corollary}
\newtheorem{definition}[theorem]{Definition}
\newtheorem{remark}[theorem]{Remark}

\newtheoremstyle{cstyle}
{15pt}	
{5pt}	
{}	
{}	
{\bfseries}	
{}	
{0.2222em}	
{}	

\theoremstyle{cstyle}

\makeatletter
\g@addto@macro{\thm@space@setup}{\thm@headpunct{}}
\renewenvironment{proof}[1][\proofname.]{\par
	\pushQED{\qed}%
	\normalfont \topsep6\p@\@plus6\p@\relax
	\trivlist
	\item[\hskip\labelsep
	\bfseries
	#1\@addpunct{\,}]\ignorespaces
}{%
	\popQED\endtrivlist\@endpefalse
}
\makeatother

\makeatletter
\g@addto@macro{\thm@space@setup}{\thm@headpunct{}}
\newenvironment{sketch-proof}[1][Sketch of the proof]{\par
	\pushQED{\qed}%
	\normalfont \topsep6\p@\@plus6\p@\relax
	\trivlist
	\item[\hskip\labelsep
	\bfseries
	#1\@addpunct{\,}]\ignorespaces
}{%
	\popQED\endtrivlist\@endpefalse
}
\makeatother

\newcounter{subeq}
\renewcommand{\thesubeq}{\theequation\alph{subeq}}
\newcommand{\newsubeqblock}{\setcounter{subeq}{0}\refstepcounter{equation}}
\newcommand{\mysubeq}{\refstepcounter{subeq}\tag{\thesubeq}}

\usepackage{mathtools}  
\usepackage{amsmath}
\usepackage{fancyhdr}
\usepackage{tikz}
\usepackage{pdflscape}
\usepackage{subfig}
\usepackage{float}
\usepackage{hyperref}
\usepackage{tikz}
\usepackage{pgfplots}
\pgfplotsset{compat=1.16} 
\usepackage{algorithmic}

\usepackage[section]{placeins}         

\usepackage[backend=biber,style=numeric,maxbibnames=99]{biblatex} 

\numberwithin{equation}{section}

\makeatletter
\newcommand{\mylabel}[2]{%
   \protected@write \@auxout {}{\string \newlabel {#1}{{#2}{\thepage}{#2}{#1}{}}}%
   \hypertarget{#1}{#2}%
}
\newcommand{\mylabelHIDE}[2]{%
   \protected@write \@auxout {}{\string \newlabel {#1}{{#2}{\thepage}{#2}{#1}{}}}%
   \hypertarget{#1}{}%
}
\makeatother

\renewcommand{\rho}{\varrho}
\renewcommand{\phi}{\varphi}

\newcommand{\pmbn}{\pmb{\mathrm{n}}}
\newcommand{\pmbP}{\pmb{\mathrm{P}}}

\newcommand{\bJ}{\ensuremath{\mathbf{J}}}
\newcommand{\bn}{\ensuremath{\mathbf{n}}}
\newcommand{\bu}{\ensuremath{\mathbf{u}}}
\newcommand{\bv}{\ensuremath{\mathbf{v}}}
\newcommand{\bw}{\ensuremath{\mathbf{w}}}
\newcommand{\bz}{\ensuremath{\mathbf{z}}}

\newcommand{\bbA}{\ensuremath{\mathbb{A}}}
\newcommand{\bbB}{\ensuremath{\mathbb{B}}}
\newcommand{\bbC}{\ensuremath{\mathbb{C}}}
\newcommand{\bbD}{\ensuremath{\mathbb{D}}}
\newcommand{\bbF}{\ensuremath{\mathbb{F}}}
\newcommand{\bbG}{\ensuremath{\mathbb{G}}}
\newcommand{\bbI}{\ensuremath{\mathbb{I}}}
\newcommand{\bbN}{\ensuremath{\mathbb{N}}}
\newcommand{\bbO}{\ensuremath{\mathbb{O}}}
\newcommand{\bbR}{\ensuremath{\mathbb{R}}}
\newcommand{\bbT}{\ensuremath{\mathbb{T}}}

\newcommand{\D}{\ensuremath{\mathbb{D}}}

\newcommand{\mycal}[1]{\ensuremath{\mathcal{#1}}}
\newcommand{\calA}{\mycal{A}}
\newcommand{\calB}{\mycal{B}}
\newcommand{\calC}{\mycal{C}}
\newcommand{\calD}{\ensuremath{\mathcal{D}}}
\newcommand{\calE}{\mycal{E}}
\newcommand{\calF}{\ensuremath{\mathcal{F}}}
\newcommand{\calG}{\ensuremath{\mathcal{G}}}
\newcommand{\calH}{\ensuremath{\mathcal{H}}}
\newcommand{\calI}{\mycal{I}}

\newcommand{\calN}{\ensuremath{\mathcal{N}}}

\newcommand{\calQ}{\mycal{Q}}

\newcommand{\calP}{\ensuremath{\mathcal{P}}}
\newcommand{\calS}{\ensuremath{\mathcal{S}}}
\newcommand{\calT}{\ensuremath{\mathcal{T}}}

\newcommand{\calV}{\ensuremath{\mathcal{V}}}
\newcommand{\calW}{\ensuremath{\mathcal{W}}}

\newcommand{\divergenz}[1]{\operatorname{div}\left({#1}\right)}
\DeclareMathOperator{\Div}{div}
\DeclareMathOperator{\diam}{diam}

\DeclareMathOperator{\trace}{Tr}

\newcommand{\abs}[1]{\left\lvert{#1}\right\rvert}
\newcommand{\norm}[1]{\|{#1}\|}
\newcommand{\nnorm}[1]{\left\|{#1}\right\|}
\newcommand{\dualp}[2]{\left<{#1\, ,\, #2}\right>}
\newcommand{\skp}[2]{\left({#1\, ,\, #2}\right)}

\newcommand{\skpp}[2]{\left\llangle {#1\, , \, #2}\right\rrangle }  

\newcommand{\dv}[1]{\,{\mathrm d}#1}
\newcommand{\dx}{\dv{x}}

\newcommand{\dt}{\dv{t}}

\bibliography{main}

\usepackage[draft]{todonotes}   

\begin{document}
%
%
	
\begin{center}	
	\LARGE{Approximation and existence of a viscoelastic phase-field model for tumour growth in two and three dimensions}
\end{center}
\bigskip

\begin{center}	
	\normalsize{Harald Garcke and Dennis Trautwein}\\[1mm]
	\textit{Fakult\"at f\"ur Mathematik, Universit\"at Regensburg, 93053 Regensburg, Germany}\\[1mm]
	\texttt{\{Harald.Garcke, Dennis.Trautwein\}@ur.de}
\end{center}



\begin{abstract}~
\footnotesize
In this work, we present a phase-field model for tumour growth, where a diffuse interface separates a tumour from the surrounding host tissue. In our model, we consider transport processes by an internal, non-solenoidal velocity field. We include viscoelastic effects with the help of a general Oldroyd-B type description with relaxation and possible stress generation by growth.
The elastic energy density is coupled to the phase-field variable which allows to model invasive growth towards areas with less mechanical resistance.

The main analytical result is the existence of weak solutions in two and three space dimensions in the case of additional stress diffusion. The idea behind the proof is to use a numerical approximation with a fully-practical, stable and (subsequence) converging finite element scheme. The physical properties of the model are preserved with the help of a regularization technique, uniform estimates and a limit passage on the fully-discrete level.
Finally, we illustrate the practicability of the discrete scheme with the help of numerical simulations in two and three dimensions.
\pskip

\pskip
	
\end{abstract}



%
%
%

\section{Introduction}

\begin{figure}
\centering
\includegraphics[width=0.5\textwidth]{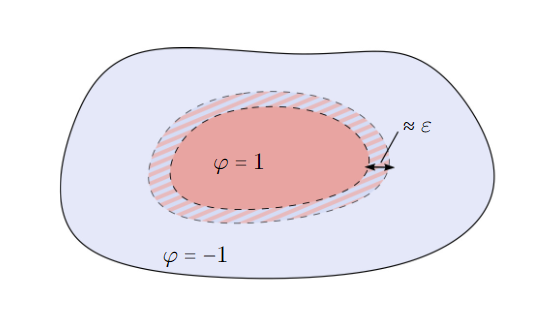}
\caption{Typical setting for phase-field models, where a smooth interface with width related to $\varepsilon>0$ separates the pure phases. }
\label{fig:phases}
\end{figure}

Diffuse interface models for tumour growth have become quite popular in recent years. 
The so-called order parameter $\phi\in[-1,1]$ is given by the difference of volume fractions of two components such that $\phi=1$ in the pure tumour phase and $\phi=-1$ in the pure host phase, which usually consists of a healthy tissue or in experiments of an agarose gel. The evolution of $\phi$ is typically described by a convected diffusion equation with source terms:
\begin{align}
\label{eq:0_phi}
    \partial_t\phi + \Div(\phi\bv) + \Div\bJ_\phi = \Gamma_\phi.
\end{align}
Here, $\bv$ is an internal velocity field and $\Gamma_\phi$ decribes possible sources. 
In the phase-field setting of Cahn--Hilliard type, the diffusive flux $\bJ_\phi$ satisfies Fick's law \cite{AbelsGG_2012} and is given by
\begin{align}
\label{eq:0_Jphi}
    \bJ_\phi = - m_\phi \nabla \mu,
    \quad \mu = \frac{\partial \calE}{\partial \phi} - \Div \frac{\partial \calE}{\partial \nabla\phi},
\end{align}
where $m_\phi$ is a non-negative mobility function and $\mu$ denotes the chemical potential for the order parameter $\phi$. Typically, the energy density $\calE$ is split additively into the Ginzburg--Landau free energy density $\calE_{\mathrm{GL}}(\phi,\nabla\phi) = \frac\beta\varepsilon \psi(\phi) + \frac{\beta\varepsilon}{2} \abs{\nabla\phi}^2$, which is used to account for adhesive forces between the different types of cells or tissues, and other contributions that may depend on further variables \cite{GarckeLSS_2016}. 
The parameter $\beta>0$ is related to surface tension and $\varepsilon>0$ is a small constant that is proportional to the width of a smooth interface between the pure regions where $\phi=\pm1$, see also Figure~\ref{fig:phases}.
Moreover, $\psi(\cdot)$ is a nonlinear potential with global minima at $\phi=\pm1$. Common choices are the polynomial double-well potential, the non-smooth double obstacle potential and the logarithmic potential:
\begin{align*}
    \psi_{\mathrm{dw}}(\phi) &= \frac14 (1-\phi^2)^2 \quad \forall \, \phi\in\bbR ,
    \qquad\qquad 
    \psi_{\mathrm{do}}(\phi) = 
    \begin{cases} 
    \frac12 (1-\phi^2) &\text{for }\phi\in[-1,1],
    \\
    +\infty &\text{else},
    \end{cases}
    \\
    \psi_{\mathrm{log}}(\phi) &= \frac{\theta}{2}  \big(
    (1+\phi) \log(1+\phi) + (1-\phi) \log(1-\phi) \big)
    + \frac{\theta_c}{2} (1-\phi^2) \quad \forall\, \phi\in [-1,1],
\end{align*}
where $0 < \theta < \theta_c$ are positive constants related to temperature.

The simplest mathematical approaches for tumour growth that are based on \eqref{eq:0_phi} neglect the influence of the internal velocity field $\bv$ and have been studied in, e.g., \cite{agosti_ACGV_2017, frigeri_grasselli_rocca_2015, frigeri_2017_tumour_degenerate, garcke_lam_2017, hawkins_2010}.
More general models include velocity effects with a Darcy law \cite{agosti_CGAC_2018, Chen_Lowengrub_2014, GarckeLSS_2016, GLNS_2018_multiphase_tumour_necrosis, wise_lowengrub_2008} or a Stokes/Brinkman law \cite{ebenbeck_garcke_nurnberg_2020, knopf_2022_multiphase} for the fluid velocity.
%
%
%
%
%
%
Other multiphase approaches with a (visco-)elastic description have been addressed in, e.g., \cite{garcke_2022_viscoelastic, garcke_2023_viscoelastic, ambrosi_2009, garcke_lam_signori_2021, lucci_agosti_2022, bresch_2009, GKT_2022_viscoelastic}. The motivation for (visco-)elastic models is that solid stress can be generated by growth mechanisms and also affect growth itself. The most common reason is that elastic stresses built up due to local overcrowding or oriented growth \cite{mao_2013}. It is also widely known that solid stresses increase with the size of the tumour resulting in a highly nonlinear stress-strain response. On the other side, interstitial fluid pressure is always present in living tissues but not as dominant as solid stresses or growth induced stresses. The latter is usually not as dominant as solid stresses, but it is not negligible in general \cite{nia_2016, northcott_2018, voutouri_2014}.
To account for these effects, we assume a momentum balance law for a viscoelastic fluid with negligible Reynolds number
\begin{align}  
\label{eq:0_momentum}
    - \Div \big( \eta(\phi) (\nabla \bv + (\nabla\bv)^\top) + \bbT_e \big) + \nabla p = \mathbf{f},
\end{align}
with a scalar pressure $p$ and external forces $\mathbf{f}$. 
The full stress tensor is composed of elastic stress $\bbT_e$ and Newtonian stress $\eta(\phi) (\nabla \bv + (\nabla\bv)^\top)$, where $\eta(\phi)$ is the shear viscosity.
%
%
%
In many cases, the force on the right-hand side of \eqref{eq:0_momentum} is given by $\mathbf{f} = - \Div \big( \nabla\phi \otimes \frac{\partial\calE}{\partial\nabla\phi} \big)$ and describes capillary forces, but also other related forces are often used in order to simplify the existence analysis \cite{ebenbeck_garcke_nurnberg_2020, GarckeLSS_2016, GKT_2022_viscoelastic}.
As the pure tumour tissue ($\phi=1$) and the pure host tissue ($\phi=-1$) in general have different (constant) mass densities $\rho_t, \rho_h>0$, we consider a variable mass density 
\begin{align*}
    \rho(\phi) &= \frac12 \rho_t (1+\phi) + \frac12 \rho_h (1-\phi),
\end{align*}
for the mixture. Consequently, as local volume changes of the mixture are possible, we account for them with 
\begin{align}
\label{eq:0_div}
    \Div \bv &= \Gamma_{\bv},
\end{align}
where $\Gamma_{\bv}$ is related to the mass source $\Gamma_\phi$.

\begin{figure}
\centering
\includegraphics[width=0.75\textwidth]{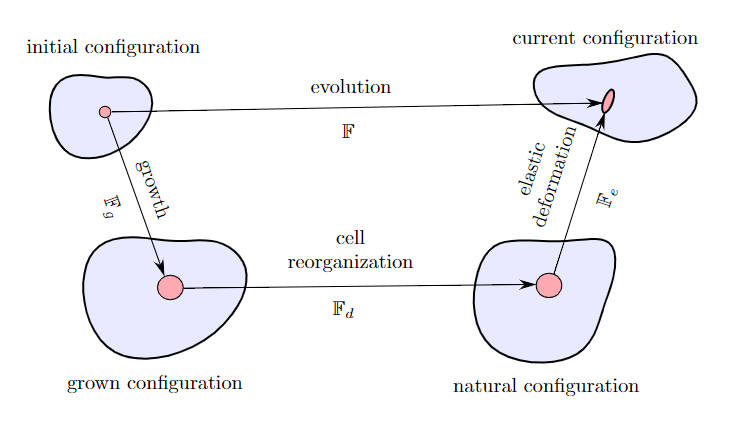}
\caption{Framework of multiple configurations.}
\label{fig:configurations}
\end{figure}

In some cases, the counterpart of stress generation --- the relaxation --- comes into play. Relaxation describes the phenomenom when solid stresses reduce or even vanish after some time, while the strain remains nearly constant. In biology, this can have several reasons like oriented cell division \cite{legoff_2016} or cell reordering \cite{forgacs_1998}. The time scale of relaxation is often short when the properties of the tissue (e.g.~an embryotic tissue) ressemble the one of a fluid. However, experimental observations suggest that the relaxation time scale is comparable to that of morphogenetic processes \cite{doubrovinski_2017, lowengrub_2021_viscoelastic}. Therefore, from the modelling point of view, we consider a viscoelastic model which can account both for stress generation and relaxation.
Here, we choose the Oldroyd-B model with relaxation time $\tau(\phi)$ and stress sources $\Gamma_{\bbB}$ related to growth:
\begin{align}
\label{eq:0_B}
\partial_t \bbB + (\bv\cdot\nabla)\bbB
    + \frac{1}{\tau(\phi)} \bbT_e
    &= 
    \nabla\bv \bbB + \bbB (\nabla\bv)^\top
    - \Gamma_{\bbB} \bbB,
\end{align}
where $\bbB$ denotes the left Cauchy--Green tensor associated with the elastic part of the total mechanical response of the viscoelastic fluid, and the relation between the Cauchy--Green tensor $\bbB$ and the elastic stress tensor is given by $\bbT_e = 2 \frac{\partial\calE}{\partial\bbB} \bbB$.
The Oldroyd-B model with relaxation and stress source can be derived within a framework of multiple configurations \cite{ambrosi_2009, GKT_2022_viscoelastic, malek_prusa_2018, RAJAGOPAL2000207}, see also Figure \ref{fig:configurations} for a schematic sketch. The main idea is a virtual multiplicative decomposition of the full deformation tensor $\bbF = \bbF_e \bbF_d \bbF_g$ into one part that describes deformation by growth ($\bbF_g$), one part that accounts for dissipative processes like cell reorganization ($\bbF_d$) and the elastic part of the total mechanical response ($\bbF_e$). Then, the sought quantity is the left Cauchy--Green tensor $\bbB \coloneqq \bbF_e \bbF_e^\top$ associated with $\bbF_e$, whose evolution in Eulerian coordinates is described with \eqref{eq:0_B}. For viscoelastic fluids, the quantities $\bbF_d$ and $\bbF_g$ do not explicitly enter the system of equations, but they are strongly related to the terms $\frac{1}{\tau(\phi)} \bbT_e$ and $\Gamma_{\bbB} \bbB$ from \eqref{eq:0_B} due to constitutive assumptions, see also \cite{GKT_2022_viscoelastic, malek_prusa_2018, RAJAGOPAL2000207}.
This viscoelastic description allows a (formal) interpolation of material laws between the phases by selecting different values for the shear viscosity $\eta(\phi)$ and the relaxation time $\tau(\phi)$ in each phase.
For example, by neglecting the influence of growth ($\Gamma_{\bbB}=0$), an elastic tumour and a fully-viscous host tissue can be reached with the choices $\tau(1)\to +\infty$, $\eta(1)\to 0$ and $\tau(-1)\to 0$, respectively. In particular, we then have 
\begin{align*}
    - \Div \bbT_e  + \nabla p = \mathbf{f},
    \qquad 
    \partial_t \bbB + (\bv\cdot\nabla)\bbB
    = \nabla\bv \bbB + \bbB (\nabla\bv)^\top ,
    \qquad \Div \bv = \Gamma_{\bv},
\end{align*}
in the pure tumour phase ($\phi=1$), where $\bbT_e = 2 \frac{\partial\calE}{\partial\bbB} \bbB$. In the pure host phase $(\phi=-1)$ it is
\begin{align*}
    - \eta(-1) \Div \big(  \nabla \bv + (\nabla\bv)^\top \big) + \nabla p = \mathbf{f},
    \qquad 
    \bbT_e = \mathbb{O}, 
    \qquad \Div \bv = \Gamma_{\bv},
\end{align*}
where $\mathbb{O}\in \bbR^{d\times d}$, $d\in\{2,3\}$, denotes the zero matrix. Note that the pressure serves as a Lagrange multiplier for \eqref{eq:0_div}. 


%
%
%
Possible consequences of present stresses in growing tissues have been known for a long time. The cell division mechanism (proliferation) of tumour cells can be inhibited by compression \cite{voutouri_2014}, while the programmed cell death (apoptosis) and necrosis can increase \cite{delarue_2014, helmlinger_1997, northcott_2018}.
This motivates the following example for the mass source:
\begin{align}
\label{eq:0_source_phi}
    \Gamma_\phi(\phi,\sigma,\bbB) = \frac12 (1+\phi) (\mathcal{P} \sigma f_1(\bbB)
    - \mathcal{A} f_2(\bbB) ).
\end{align}
Here, as in \cite{GarckeLSS_2016}, we assume that proliferation and apoptosis with constant rates $\calP, \calA\geq 0$, respectively, only take place in the growing phase where $\phi=1$. To account for the influence of elastic stress, it is reasonable to include a scaling with functions $f_1(\bbB), f_2(\bbB)$ that depend on $\bbB$. 
Besides, $\sigma$ denotes the concentration of a chemical species that serves as a nutrient for the tumour (e.g.~glucose or oxygen). Typically, the source terms $\Gamma_{\bv}$ and $\Gamma_{\bbB}$ have a similar structure as \eqref{eq:0_source_phi}, see \cite{ambrosi_2009, GKT_2022_viscoelastic, GarckeLSS_2016}.

Usually, small tumours have not developed an own vascular system yet, so that nutrients are distributed by diffusion. As the time scale of diffusion compared to morphogenic processes is usually very short \cite{byrne_tumour_2003}, it is reasonable to assume a quasi-static description by
\begin{align}
\label{eq:0_sig}
    - \Div(m_\sigma \nabla \sigma) = - \sigma \Gamma_\sigma,
\end{align}
where $m_\sigma$ is a non-negative mobility function and $\Gamma_\sigma$ a source or sink term related to nutrient consumption.
For example, assuming that the nutrient is primarily consumed by the tumour with a rate $\calC\geq 0$, a possible choice for $\Gamma_\sigma$ is 
\begin{align}
\label{eq:0_source_sig}
    \Gamma_\sigma(\phi,\bbB) = \frac12 \mathcal{C} (1+\phi) f_3(\bbB),
\end{align}
which can be scaled with a function $f_3(\bbB)$ that depends on $\bbB$.
Additional sources for nutrient supply can be included with Dirichlet or mixed boundary conditions.
Moreover, compression also has an indirect influence on growth. In particular, the nutrient supply can be affected due to, e.g., vascular compression \cite{voutouri_2014}. This can be included in the model by choosing mobility functions that depend on $\phi$ and $\bbB$.





Another point we want to address in our model is the invasive behaviour of tumours. It has been observed that the invasive potential of tumours can be enhanced by mechanical forces like compression, see, e.g., \cite{northcott_2018, voutouri_2014} and references therein. Also, it is widely known in medicine that the movement and growth of tumour cell aggregates can be directed by mechanical effects like local differences in stress or mechanical resistance \cite{helmlinger_1997}. This behaviour is sometime referred to as \textit{mechanotaxis}.
Moreover, a lack of nutrients, e.g., due to compression of blood vessels, can favour \textit{chemotaxis}, which in general describes the movement along the gradient of a chemical species like nutrients or signaling molecules. 
For these reasons, we consider the specific form of the energy density 
\begin{align*}
    \calE(\phi,\nabla\phi, \sigma,\bbB) 
    = \frac\beta\varepsilon \psi(\phi) + \frac{\beta\varepsilon}{2} \abs{\nabla\phi}^2
    + \chi_\phi \sigma (1-\phi)
    + W(\phi,\bbB).
\end{align*}
The constant $\chi_\phi\geq0$ can be interpreted as an sensitivity for chemotaxis \cite{GarckeLSS_2016}. Moreover, $W(\phi,\bbB)$ denotes the elastic energy density of the system for which we assume the form
\begin{align*}
    W(\phi,\bbB) \coloneqq 
    W_0(\bbB) 
    + W_1(\phi, \trace \bbB), 
\end{align*}
where $W_0(\bbB)$ is convex and frame invariant, and $W_1(\phi,\trace\bbB)$ only depends on $\trace\bbB$ and the phase-field variable $\phi$.
%
%
With this choice for the energy density, \eqref{eq:0_Jphi} translates into
\begin{align}
    \label{eq:0_Jphi1}
    \bJ_\phi &= -m_\phi \nabla\mu ,
    \quad
    \mu = \frac\beta\varepsilon \psi'(\phi) - \beta\varepsilon \Delta \phi 
    - \chi_\phi \sigma + \frac{\partial W_1}{\partial\phi}.
\end{align}
Here, the term $m_\phi \chi_\phi \nabla\sigma$ in $\bJ_\phi$
lets the tumour move along the gradient of the nutrient, which can be interpreted as the chemotaxis effect. Moreover, the term $- m_\phi \nabla \frac{\partial W_1}{\partial\phi}$ in $\bJ_\phi$ drives the tumour towards regions which are energetically more advantageous due to, e.g., less mechanical resistance. 
%
%
%
%
%
In this work, we consider
\begin{align}
    \label{eq:0_W}
    W_0(\bbB) = \frac14 \abs{\bbB}^2 - \frac12 \trace\ln\bbB,
    \qquad
    W_1(\phi,\trace\bbB) = \frac12 \kappa(\phi) \trace\bbB,
\end{align}
where $W_0$ is strictly convex, only depends on $\trace(\bbB\bbB)$ and $\det(\bbB)$, i.e., on invariants of $\bbB$, and it has the global minimizer $\bbB=\bbI$, where $\bbI\in\bbR^{d\times d}$, $d\in\{2,3\}$, denotes the identity matrix. Note that $\trace\ln\bbB = \ln \det \bbB$ for symmetric positive definite matrices \cite{barrett_boyaval_2009}. 
The specific form of $W_1(\phi,\trace\bbB) = \frac12 \kappa(\phi) \trace\bbB$ is chosen as a mechanical counterpart of the chemotaxis energy $\chi_\phi\sigma(1-\phi)$ to account for growth along the gradient of $-\kappa'(\phi) \trace\bbB$,
as the invasive potential of tumours can be intensified and directed by mechanical effects \cite{helmlinger_1997, northcott_2018, voutouri_2014}.
Here $\kappa(\phi)$ is a phase-depending material parameter that can have different values in the pure phases where $\phi=\pm 1$.
The term $- m_\phi \nabla \frac{\partial W_1}{\partial\phi} = - \frac12 m_\phi \nabla( \kappa'(\phi) \trace\bbB)$ in the diffusive flux $\bJ_\phi$ in \eqref{eq:0_Jphi1} can enforce, weaken and direct the tumour movement, depending on local differences in $\bbB$. In particular, if $\kappa(\phi)$ is an affine function so that $\kappa'(\phi)$ is constant, the tumour moves along $\nabla\trace\bbB$.
For numerical examples, see also Section \ref{sec:numeric}.
The elastic stress tensor can be written as $\bbT_e = 2 \frac{\partial\calE}{\partial\bbB} \bbB 
    = \bbB^2 - \bbI + \kappa(\phi) \bbB$,
and it depends nonlinearly on the Cauchy--Green tensor $\bbB$. 
This agrees with the observation that the stress-strain response in tumours is typically nonlinear \cite{nia_2016, voutouri_2014}.
In the literature, there are many suggestions for $W_0$ for viscoelastic fluids. Assuming that the small particles in the mixture behave like Hookean dumbbells leads to the classical Oldroyd-B model with $W_0(\bbB) = \frac12 \big( \trace\bbB - \trace\ln\bbB - d \big)$.
Here the first term accounts for stretching to infinity and the second term penalizes compression to a single point. This elastic energy density has been used, e.g., in \cite{barrett_boyaval_2009, barrett_lu_sueli_2017, constantin_kliegl_2012, GKT_2022_viscoelastic, LT_2022_arxiv}, but it has the disadvantage that long-time existence results are only possible in two space dimensions due to lack of compactness for $\bbB$, even with additional stress diffusion. 
The FENE-P model \cite{barrett_2018_fene-p} with $W_0(\bbB) = -\frac b2 \ln \big( 1 - b^{-1} \trace\bbB \big) - \frac12 \trace\ln\bbB$, where $b\gg1$, is a nonlinear variant of the classical Oldroyd-B model and it has the additional feature that particles are only finitely stretchable, i.e., $W_0(\bbB)\to \infty$ as $\trace\bbB\to b$. The classical Oldroyd-B model can be recovered in the limit $b\to \infty$.
Besides, models with quadratic contributions in $W_0$, similarly to \eqref{eq:0_W}, can be motivated by non-Hookean behaviour of the particles on a microscopic scale \cite{Lukacova_2017}. Here, existence results of fluid models have also been established in three dimensions \cite{bathory_2021_viscoelastic_3D, brunk_2022_3d}. 
For further viscoelastic models, we refer to, e.g., \cite{garcke_2022_viscoelastic, garcke_2023_viscoelastic, chupin_2018, dostalik_prusa_skrivan_2019} and references therein.
\subsection{The full mathematical model}

We now summarize the full mathematical model that was described so far. For the notation, we refer to the end of this section.

Let $T>0$ and $\Omega\subset\bbR^d$, $d\in\{2,3\}$, be a bounded Lipschitz domain with boundary $\partial\Omega$. The full system of our interest is given as follows.
For any $(\mathbf{x},t)\in \Omega_T \coloneqq \Omega\times(0,T)$, find the order parameter and its chemical potential $\phi(\mathbf{x},t), \mu(\mathbf{x},t) \in \bbR$, the nutrient concentration $\sigma(\mathbf{x},t)\in\bbR$, the velocity field and the scalar pressure $\bv(\mathbf{x},t) \in\bbR^d$, $p(\mathbf{x},t)\in\bbR$ as well as the Cauchy--Green tensor $\bbB(\mathbf{x},t)\in\bbR^{d\times d}_{\mathrm{S}}$, such that
\begin{subequations}
\begin{align}
    \label{eq:phi}
    \partial_t \phi + \Div(\phi \bv) 
    &= \Div(m_\phi(\phi, \bbB) \nabla\mu) + \Gamma_\phi(\phi,\sigma,\bbB),
    \\
    \label{eq:mu}
    \mu &= \frac\beta\varepsilon \psi'(\phi) - \beta\varepsilon \Delta\phi 
    -\chi_\phi \sigma
    + \frac12 \kappa'(\phi) \trace\bbB,
    \\
    \label{eq:sigma}
    0 &= \Div( m_\sigma(\phi, \bbB) \nabla \sigma) - \sigma \Gamma_\sigma(\phi,\bbB),
    \\
    \label{eq:div}
    \Div \bv &= \Gamma_{\bv}(\phi,\sigma,\bbB),
    \\
    \label{eq:v}
    - \divergenz{ \eta(\phi) ( \nabla\bv + (\nabla\bv)^\top)} + \nabla p
    &=  \Div \bbT_e
    + (\mu  + \chi_\phi \sigma ) \nabla\phi
    + \frac12 \kappa(\phi) \nabla \trace\bbB ,
    \\
    \label{eq:B}
    \partial_t \bbB + (\bv\cdot\nabla)\bbB
    + \frac{1}{\tau(\phi)} \bbT_e
    &= 
    \nabla\bv \bbB + \bbB (\nabla\bv)^\top
    - \Gamma_{\bbB}(\phi,\sigma,\bbB) \bbB
    + \alpha \Delta \bbB,
\end{align}
\end{subequations}
where the elastic part of the stress tensor $\bbT_e$ is denoted by
\begin{align*}
    \bbT_e &= \bbB^2 + \kappa(\phi) \bbB - \bbI,
\end{align*}
where $\bbI\in\bbR^{d\times d}$ denotes the unit matrix.
To close the system of partial differential equations, we employ the initial conditions
\begin{align*}
    \phi(0) = \phi_0, \quad \bbB(0) = \bbB_0 \quad \text{in } \Omega,
\end{align*}
and the boundary conditions
\begin{subequations}
\begin{alignat}{4}
    \label{eq:bc_phimu}
    \nabla\phi\cdot\bn &= m_\phi(\phi,\bbB) \nabla\mu\cdot\bn = 0 \quad &&\text{on } \partial\Omega,
    \\
    \label{eq:bc_sigma}
    m_\sigma(\phi, \bbB) \nabla \sigma \cdot \bn &= K(\sigma_\infty - \sigma) \quad &&\text{on } \partial\Omega,
    \\
    \label{eq:bc_B}
    (\bn\cdot\nabla)\bbB &= \bbO \quad &&\text{on } \partial\Omega,
    \\
    \label{eq:bc_T}
    \big(\eta(\phi) (\nabla\bv + (\nabla\bv)^\top) + \bbT_e - p\bbI \big) \bn 
    &= \mathbf{0} \quad &&\text{on } \partial_{\mathrm{N}}\Omega,
    \\
    \label{eq:bc_v}
    \bv&= \mathbf{0} \quad &&\text{on } \partial_{\mathrm{D}}\Omega, 
\end{alignat}
\end{subequations}
where $\bn$ denotes the outer unit normal to $\partial\Omega$.
Here, $\partial_{\mathrm{D}}\Omega \subset \partial\Omega$ is assumed to be closed and of positive surface measure $\calH^{d-1}(\partial_{\mathrm{D}}\Omega) > 0$, and we set $\partial_\mathrm{N}\Omega = \partial\Omega \setminus \partial_{\mathrm{D}}\Omega$.
Moreover, $\mathbf{0}\in\bbR^d$ and $\bbO\in\bbR^{d\times d}$ denote the zero vector and zero matrix, respectively.
Similarly to, e.g., \cite{barrett_boyaval_2009, barrett_lu_sueli_2017, bathory_2021_viscoelastic_3D, GKT_2022_viscoelastic}, we have the additional term $\alpha\Delta \bbB$ on the right-hand side of \eqref{eq:B}, which is needed for the analysis of global-in-time weak solutions.
The no-flux boundary conditions \eqref{eq:bc_phimu}, \eqref{eq:bc_B} and the no-stress boundary condition \eqref{eq:bc_T} are common choices in the literature \cite{barrett_boyaval_2009, ebenbeck_garcke_nurnberg_2020}.
The boundary condition \eqref{eq:bc_v} is needed for analytical reasons and can be motivated by, e.g., the presence of a bone along $\partial_{\mathrm{D}}\Omega$. Moreover, \eqref{eq:bc_sigma} describes possible nutrient supply on the boundary, where $K>0$ and $\sigma_\infty\colon\, \partial\Omega\times(0,T)\to\bbR$.

The external force $\mathbf{f}$ from \eqref{eq:0_momentum} has been specified in \eqref{eq:v} as $(\mu+\chi_\phi \sigma) \nabla\phi + \frac12 \kappa(\phi) \nabla\trace\bbB$, which is closely related to the capillary force $\mathbf{f} = -\Div\big( \beta\varepsilon \nabla\phi \otimes \nabla\phi\big)$ after a reformulation of the scalar pressure $p$, see also \cite{AbelsGG_2012, ebenbeck_garcke_nurnberg_2020, GKT_2022_viscoelastic, GarckeLSS_2016}. 
In particular, defining a new scalar pressure $\tilde p = p + \frac{\beta}{\varepsilon} \psi(\phi) + \frac{\beta\varepsilon}{2} \abs{\nabla\phi}^2 + \frac12 \kappa(\phi) \trace\bbB$ and using the relation
\begin{align*}
\nabla \tilde p = \nabla p + \Div\big( \beta\varepsilon \nabla\phi \otimes \nabla\phi\big) + (\mu+\chi_\phi \sigma)\nabla\phi + \frac12 \kappa(\phi) \nabla\trace\bbB,
\end{align*}
one observes that the momentum equation \eqref{eq:0_momentum} with $\mathbf{f} = -\Div\big( \beta\varepsilon \nabla\phi \otimes \nabla\phi\big)$ translates into \eqref{eq:v}, where we again write $p$ instead of $\tilde p$.
This form is chosen to simplify the existence analysis of weak solutions, but both forms are equivalent for sufficiently smooth solutions.

\subsection{Main results and key strategy}
The main results of this work are as follows.
We study the global-in-time existence of weak solutions to the system \eqref{eq:phi}--\eqref{eq:B} endowed with the boundary conditions \eqref{eq:bc_phimu}--\eqref{eq:bc_v}, see Theorem \ref{theorem:weak_solution}. The existence result is valid in two and in three space dimensions.
In addition to that, we present a fully-practical numerical approximation of the system, where existence and stability of discrete solutions are shown under a very mild constraint on the time step size, see Theorem \ref{theorem:existence_FE}. 
Under specific assumptions on some model functions, it is also possible to obtain the unconditional stability and existence of discrete solutions, see Remark \ref{remark:scheme}. 
Later, we prove subsequence convergence of discrete solutions, where the limit functions form a global-in-time weak solution. From this we obtain the existence result in Theorem \ref{theorem:weak_solution}.
To illustrate the practicability of the numerical approximation of the system, we present simulations for possible scenarios in two and in three space dimensions in Section~\ref{sec:numeric}.

First, we state the precise definition of a weak solution. For the notation, we refer to the end of this section.


\begin{definition}[Weak solution]
\label{def:weak_solution}
The tuple of functions $(\phi,\mu,\sigma,p,\bv,\bbB)$ is called a weak solution of the system \eqref{eq:phi}--\eqref{eq:B} subject to the boundary conditions \eqref{eq:bc_phimu}--\eqref{eq:bc_v}, if
\begin{align*}
    &\phi \in L^2(0,T;H^1(\Omega)) 
    \cap H^1(0,T; (H^1(\Omega))'),
    \quad
    \mu \in L^2(0,T; H^1(\Omega)),
    \\
    &\sigma \in L^2(0,T; H^1(\Omega)) ,
    \quad 
    p \in L^{4/3}(0,T; L^2(\Omega)),
    \quad
    \bv \in L^2(0,T; H^1_\mathrm{D}(\Omega;\bbR^d)),
    \\
    &\bbB \in L^2\left(0,T; H^1(\Omega;\bbR^{d\times d}_{\mathrm{S}})\right)
    \cap W^{1,\frac{4}{3}}\left(0,T; (H^1(\Omega; \bbR^{d\times d}_{\mathrm{S}}))'\right),
\end{align*}
with $\phi(0) = \phi_0$ in $L^2(\Omega)$, $\bbB(0) = \bbB_0$ in $L^2(\Omega;\bbR^{d\times d}_{\mathrm{S}})$ and with $\bbB$ positive definite a.e.~in $\Omega_T$, such that
\begin{subequations}
\begin{align}
    \label{eq:div_weak}
    \Div \bv &= \Gamma_{\bv}(\phi,\sigma,\bbB),
    \quad \text{a.e.~in } \Omega_T,
\end{align}
and, for a.e.~$t\in(0,T)$ and for all $\zeta, \rho, \xi\in H^1(\Omega)$, $\bw\in H^1_{\mathrm{D}}(\Omega;\bbR^d)$ and $\bbG\in H^1(\Omega;\bbR^{d\times d}_{\mathrm{S}})$,
\begin{align}
    \label{eq:phi_weak}
    \dualp{\partial_t\phi}{\zeta}_{H^1} 
    + \skp{\Div(\phi\bv)}{\zeta}_{L^2}
    &=  
    - \skp{ m_\phi(\phi,\bbB)\nabla\mu}{\nabla\zeta}_{L^2} 
    + \skp{\Gamma_\phi(\phi,\sigma,\bbB)}{\zeta}_{L^2}, 
    \\
    \label{eq:mu_weak} 
    \skp{\mu}{\rho}_{L^2} &= 
    \skp{\frac\beta\varepsilon \psi'(\phi) - \chi_\phi \sigma + \tfrac12 \kappa'(\phi) \trace\bbB}{\rho}_{L^2} 
    + \beta\varepsilon \skp{\nabla\phi}{\nabla\rho}_{L^2} ,
    \\
    \label{eq:sigma_weak} \nonumber
    0 &= \skp{ m_\sigma(\phi, \bbB) \nabla \sigma}{\nabla\xi}_{L^2} 
    + \skp{\sigma \Gamma_\sigma(\phi,\bbB)}{\xi}_{L^2}
    \\
    &\quad
    + K \skp{\sigma-\sigma_\infty}{\xi}_{L^2(\partial\Omega)},
    \\
    \label{eq:v_weak} \nonumber
    \skp{2 \eta(\phi) \D(\bv)}{\D(\bw)}_{L^2} 
    &=  \skp{p}{\Div \bw}_{L^2}
    - \skp{\bbT_e}{\nabla \bw}_{L^2}
    + \skp{(\mu  + \chi_\phi \sigma ) \nabla\phi}{\bw}_{L^2}
    \\
    &\quad
    + \tfrac12 \skp{\kappa(\phi) \nabla \trace\bbB }{\bw}_{L^2},
    \\
    \label{eq:B_weak} \nonumber
    \dualp{\partial_t \bbB}{\bbG}_{H^1} 
    + \skp{(\bv\cdot\nabla)\bbB}{\bbG}_{L^2}
    &= 
    - \skp{\tfrac{1}{\tau(\phi)} \bbT_e}{\bbG}_{L^2}
    + \skp{\nabla\bv \bbB + \bbB (\nabla\bv)^\top}{\bbG}_{L^2}
    \\
    &\quad
    - \skp{ \Gamma_{\bbB}(\phi,\sigma,\bbB) \bbB }{\bbG}_{L^2}
    - \alpha \skp{\nabla \bbB}{\nabla \bbG}_{L^2},
\end{align}
where $\bbT_e = \bbB^2 + \kappa(\phi) \bbB - \bbI$.
\end{subequations}
\end{definition}

For the existence result, we need the following assumptions.
\begin{itemize}	    
\item[\mylabel{A1}{$(\mathbf{A1})$}] 
The function $\Gamma_\sigma \colon \, \bbR \times \bbR^{d\times d} \to\bbR$ is non-negative, bounded and Lipschitz continuous. Moreover, the source functions $\Gamma_f\colon\, \bbR \times \bbR \times \bbR^{d\times d} \to \bbR$, $f\in\{\phi,\bv,\bbB\}$, are Lipschitz continuous and fulfill the growth condition $\abs{\Gamma_f(s,r,\bbG)} \leq C (1 + \abs{r})$ for all $s,r\in\bbR$ and $\bbG\in\bbR^{d\times d}$.

\item[\mylabel{A2}{$(\mathbf{A2})$}]
Let $\chi_\phi\geq 0$ and $\beta, \varepsilon, K, \alpha>0$ be constants. Let $m_\phi, m_\sigma \colon\, \bbR\times \bbR^{d\times d} \to \bbR$ and $\eta \colon\, \bbR\to\bbR$ be Lipschitz continuous, uniformly positive and bounded. Moreover, let $\tau\colon\, \bbR\to\bbR$ be continuous, uniformly positive and bounded, with $(1/\tau)(\cdot)$ being Lipschitz continuous. 
In particular, suppose there exist positive constants 
such that, for all $s\in\bbR, \bbG\in\bbR^{d\times d}$,
\begin{align*}
    m_0 &\leq m_\phi(s,\bbG) \leq m_1, 
    \quad 
    n_0 \leq m_\sigma(s,\bbG) \leq n_1,
    \quad
    \eta_0 \leq \eta(s) \leq \eta_1, 
    \quad
    \tau_0 \leq \tau(s) \leq \tau_1.
\end{align*}

\item[\mylabel{A3}{$(\mathbf{A3})$}]
Suppose that the function $\kappa\colon\, \bbR\to\bbR$ is continuously differentiable and satisfies the growth condition
$\abs{\kappa(s)} \leq C_\kappa (1 + \abs{s})$ $\forall \,s\in\bbR,$, where $C_\kappa>0$, and it has a bounded and Lipschitz continuous derivative.

\item[\mylabel{A4}{$(\mathbf{A4})$}] The potential $\psi\colon\, \bbR\to\bbR$ is non-negative, continuously differentiable with a Lipschitz continuous derivative and fulfills
\begin{align*}
    \psi(s) &\geq C_\psi( \abs{s}^2 - 1) ,
    \quad
    \abs{\psi^\prime(s)} \leq C_{\psi^\prime} (1+\abs s) \quad \forall \,  s\in\bbR,
\end{align*}
where $C_\psi, C_{\psi^\prime}>0$.
Moreover, we assume that $\frac\beta\varepsilon > C_\kappa^2 C_\psi^{-1}$. 

\item[\mylabel{A5}{$(\mathbf{A5})$}] 
For the initial and boundary data, we assume
\begin{align*}
    &\phi_0 \in H^1(\Omega), 
    \quad 
    \sigma_\infty \in L^\infty(0,T;L^2(\partial\Omega)),
    \quad
    \bbB_0 \in L^2(\Omega;\bbR^{d\times d}_{\mathrm{S}}),
\end{align*}
and that there exists a constant $b_0>0$ such that
\begin{align*}
    &\pmb\xi^\top \bbB_0(\mathbf{x}) \pmb\xi \geq b_0 \abs{\pmb\xi}^2 \quad \forall\, \pmb\xi\in \bbR^d, \text{ for a.e.~} \mathbf{x}\in\Omega.
\end{align*}
\end{itemize}

We now state the main existence result.

\begin{theorem}[Existence of weak solutions]
\label{theorem:weak_solution}
Let $T>0$ and $\Omega\subset\bbR^d$, $d\in\{2,3\}$, be a bounded Lipschitz domain with polygonal (or polyhedral, respectively) boundary $\partial\Omega$. Moreover, we assume that $\partial_\mathrm{D}\Omega \subset \partial\Omega$ is closed and of positive surface measure, and we set $\partial_\mathrm{N}\Omega = \partial\Omega \setminus \partial_\mathrm{D}\Omega$. In addition, we assume that a mesh exists such that $\partial_\mathrm{D}\Omega$ is matched exactly by sides of the mesh, see \ref{S} for more details.
Let \ref{A1}--\ref{A5} hold true. Then, there exists a weak solution $(\phi,\mu,\sigma,p,\bv,\bbB)$ of the system \eqref{eq:phi}--\eqref{eq:B} subject to the boundary conditions \eqref{eq:bc_phimu}--\eqref{eq:bc_v} in the sense of Definition \ref{def:weak_solution}. Moreover, it holds
\begin{align}
    \label{eq:stability_sigma}
    \nnorm{\sigma}_{L^\infty(0,T;H^1)} \leq C \norm{\sigma_\infty}_{L^\infty(0,T;L^2(\partial\Omega))},
\end{align}
and
\begin{align}
    \label{eq:stability} \nonumber
    &\nnorm{\phi}_{L^\infty(0,T;H^1)}^2
    + \nnorm{\partial_t\phi}_{L^2(0,T;(H^1)')}
    + \norm{\mu}_{L^2(0,T;H^1)}^2
    + \nnorm{\bv}_{L^2(0,T;H^1)}^2 
    \\
    \nonumber
    &\quad 
    + \nnorm{p}_{L^{4/3}(0,T;L^2)}
    + \nnorm{\bbB}_{L^\infty(0,T;L^2)} ^2
    + \nnorm{\bbB}_{L^2(0,T;H^1)}^2 
    + \nnorm{\partial_t \bbB}_{L^{4/3}(0,T;(H^1)')}
    \\
    \nonumber
    &\quad
    + \nnorm{\trace\ln\bbB}_{L^\infty(0,T;L^1)}
    + \nnorm{\trace\ln\bbB}_{L^2(0,T;H^1)}^2
    + \nnorm{\bbB^{-1}}_{L^1(\Omega_T)}^{2/3}
    \\
    &\leq C(T,c_\infty)
    \big(1 + \norm{\bbB_0}_{L^2}^2 + \abs{\ln b_0} + 
    \norm{\phi_0}_{H^1}^2 \big) ,
\end{align}
with constants $C, C(T, c_\infty)>0$, where $C(T, c_\infty)$ depends exponentially on $T$ and $c_\infty \coloneq \norm{\sigma_\infty}_{L^\infty(0,T;L^2(\partial\Omega)}$.
\end{theorem}


The fundamental idea behind the proof of Theorem \ref{theorem:weak_solution} is an approximation of \eqref{eq:phi}--\eqref{eq:B} endowed with the boundary conditions \eqref{eq:bc_phimu}--\eqref{eq:bc_v} with a fully-discrete finite element scheme in Section~\ref{sec:approximation}. For the numerical scheme, we prove stability and existence of discrete solutions in Theorem \ref{theorem:existence_FE}. In Section~\ref{sec:convergence}, we derive higher order estimates for all solutions of the discrete scheme and we use compactness results to extract a converging subsequences of discrete solutions, which converge to a weak solution in the sense of Definition \ref{def:weak_solution}. This will prove Theorem \ref{theorem:weak_solution}.

To avoid the approximation of the boundary in the numerical scheme in Section \ref{sec:approximation}, we restrict the analysis in this work to a polygonal (or polyhedral, respectively) Lipschitz domain $\Omega\subset\bbR^d$, $d\in\{2,3\}$, which however is allowed to be non-convex.
The smallness assumption on $\varepsilon$ in \ref{A4} is for technical reasons and it is not restrictive in practice, as $\varepsilon$ is usually very small.
Due to \ref{A4}, the function $\psi(\cdot)$ has at most quadratic growth. This is a technical assumption which is needed to handle the source terms in the \textit{a priori} estimates, see also \cite{garcke_lam_2017}. 
In practice, the order parameter $\phi$ always stays close to the interval $[-1,1]$. Therefore, for numerical computations, one can use more general potentials like the double-well potential $\psi_{\mathrm{dw}}(\phi)=\frac14(1-\phi^2)^2$ with (theoretical) cut-offs outside of the interval $[-M,M]$ for some large $M\gg1$.
Moreover, the example \eqref{eq:0_source_phi} can be used for the source terms $\Gamma_\phi$, $\Gamma_{\bv}$ and $\Gamma_{\bbB}$ in practice, while one has to prescribe cut-offs for the existence analysis. The same applies for $\Gamma_\sigma$, where \eqref{eq:0_source_sig} can be used. More examples for the model functions can be found in, e.g., \cite{ebenbeck_garcke_nurnberg_2020, GKT_2022_viscoelastic, GarckeLSS_2016}.

As we consider numerical integration in the fully-discrete scheme in Section \ref{sec:approximation}, we need the Lipschitz continuity of all functions in \ref{A1}--\ref{A5} for the limit passing in the discretization parameters $(h,\Delta t)\to (0,0)$. In particular, we will often use the desired error estimates \eqref{eq:interp_lipschitz}--\eqref{eq:interp_lipschitz_multivar}.

For the existence analysis, the last condition in \ref{A5} can be relaxed to $\trace\ln\bbB_0 \in L^1(\Omega)$, see also \cite{barrett_lu_sueli_2017}. 
For the explicit construction of a numerical approximation of $\bbB_0$, we have to use an initial datum that is uniformly positive definite a.e.~in $\Omega$. 
This ensures that the discrete counterpart of the logarithmic energy is bounded uniformly in the discretization parameters $h,\Delta t>0$.
We note that it follows from $\trace\ln\bbB_0\in L^1(\Omega)$ that $\bbB_0$ is (only) positive definite a.e.~in $\Omega$, i.e., 
\begin{align*}
    \pmb\xi^\top \bbB_0(\mathbf{x}) \pmb\xi > 0 \quad \forall\, \pmb\xi\in \bbR^d, \text{ for a.e.~} \mathbf{x}\in\Omega.
\end{align*}
Therefore, one can start with $\bbB_0$ such that $\trace\ln\bbB_0 \in L^1(\Omega)$, and use $\bbB_0 + b_0\bbI$ with $b_0>0$ to construct a discrete initial datum.
After the limit passage $(h,\Delta t)\to (0,0)$ and possibly after further energy estimates on the level of weak solutions, one can then send $b_0\to 0$ to recover the original initial datum $\bbB_0$.

We now explain the key strategy of the proof of Theorem \ref{theorem:weak_solution} and the outline of the paper.
For the numerical scheme in Section \ref{sec:approximation}, it is essential to understand the structure of the mathematical system \eqref{eq:phi}--\eqref{eq:B}.
In order to handle the non-homogeneous divergence equation \eqref{eq:div}, we use a splitting ansatz for the velocity field $\bv = \bu + \bz$, where $\bu$ solves \eqref{eq:div} and $\bz$ is solenoidal, i.e., $\Div \bz = 0$. 
This ansatz allows to temporarily forget about the pressure $p$ and to derive a formal energy identity which is essential for \textit{a priori} estimates. In particular, multiplying \eqref{eq:phi} with $\mu + \chi_\phi\sigma$, \eqref{eq:mu} with $\partial_t\phi$, \eqref{eq:v} with $\bz$, \eqref{eq:B} with $\frac12 \bbB + \frac12 \kappa(\phi)\bbI - \frac12 \bbB^{-1}$, one observes after integrating over $\Omega$ and by parts, that
\begin{align*}
    &\frac{d}{dt} \calF(\phi,\bbB) 
    + \calD(\phi,\mu,\bz,\bbB)
    = \mathcal{R}(\phi, \mu, \sigma, \bz, \bu, \bbB).
\end{align*}
Here, we use the energy functional
\begin{align*}
    \calF(\phi,\bbB) =
    \int_\Omega \Big( 
    \frac\beta\varepsilon \psi(\phi) + \frac{\beta\varepsilon}{2} \abs{\nabla\phi}^2
    + \frac14 \abs{\bbB}^2 - \frac12 \trace\ln\bbB + \frac12 \kappa(\phi) \trace\bbB \Big) \dx,
\end{align*}
and the dissipation functional
\begin{align*}
    \calD(\phi,\mu,\bz,\bbB) 
    & =
    \skp{m_\phi(\phi,\bbB) \nabla\mu}{\nabla\mu}_{L^2}
    + \skp{ 2 \eta(\phi) \D(\bz)}{\D(\bz)}_{L^2}
    \\ \nonumber
    &\quad
    + \skp{\tfrac{1}{2\tau(\phi)} \bbT_{e}}{\bbB + \kappa(\phi) \bbI - \bbB^{-1}}_{L^2}
    + \frac12 \alpha \norm{\nabla \bbB}_{L^2}^2
    - \frac{\alpha}{2} \skp{\nabla\bbB}{\nabla \bbB^{-1}}_{L^2},
\end{align*}
which is non-negative because of \ref{A2} and the identities $- \frac{\alpha}{2} \skp{\nabla\bbB}{\nabla \bbB^{-1}}_{L^2} \geq \frac{\alpha}{2d} \norm{\nabla\trace\ln\bbB}_{L^2}^2$ (see \cite[Lem.~3.1]{barrett_lu_sueli_2017}) and 
\begin{align*}
    \bbT_{e} : (\bbB + \kappa(\phi) \bbI - \bbB^{-1}) = \bbT_{e} : ( \bbB^{-1} \bbT_e) = \abs{\bbT_e \bbB^{-1/2} }^2 \geq 0.
\end{align*}
Here we temporarily assume that $\bbB$ is positive definite.
The mixed right-hand side is given by
\begin{align*}
    \mathcal{R}(\phi, \mu, \sigma, \bz, \bu, \bbB)
    & =
    - \chi_\phi \skp{m_\phi(\phi,\bbB) \nabla\mu}{\nabla\sigma}_{L^2}
    - \skp{\bu \cdot \nabla\phi}{\mu + \chi_\phi \sigma}_{L^2}
    \\
    &\quad
    + \skp{\Gamma_\phi(\phi,\sigma,\bbB) - \phi \Gamma_{\bv}(\phi,\sigma,\bbB)}{\mu + \chi_\phi \sigma}_{L^2}
    - \skp{2 \eta(\phi) \D(\bu)}{\D(\bz)}_{L^2}
    \\ 
    &\quad 
    - \frac12 \skp{(\bu\cdot\nabla) \bbB}{\bbB + \kappa(\phi) \bbI - \bbB^{-1}}_{L^2}
    - \frac12 \skp{\Gamma_{\bbB}(\phi,\sigma,\bbB)}{\trace \bbT_{e}}_{L^2}
    \\ 
    &\quad
    + \skp{\nabla\bu}{\trace \bbT_{e}}_{L^2}
    - \frac12 \alpha \skp{\nabla\trace\bbB}{\nabla \kappa( \phi) }_{L^2}.
\end{align*}
The key idea is to use the assumptions \ref{A1}--\ref{A5} and \textit{a priori} information for $\bu$ and $\sigma$, which can be computed independently of the other variables, to control the mixed terms with clever estimates, using Hölder's and Young's inequalities and a Gronwall argument. Note that an additional testing procedure for the chemical potential $\mu$ is required due to the presence of the source terms $\Gamma_\phi$ and $\Gamma_{\bv}$. This is where the growth assumptions on the potential $\psi(\cdot)$ from \ref{A4} are used, see also \cite{garcke_lam_2017, ebenbeck_2019_analysis}. 
In the end, the pressure $p$ can be reconstructed with standard arguments for saddle point problems \cite{girault_raviart_1986}.

To guarantee the positive definiteness of the Cauchy--Green tensor $\bbB$, we first introduce a regularization of the logarithmic function on the fully-discrete level. For technical reasons, we also have to regularize some other terms arising from the quadratic part of the elastic energy density.
To mimic the formal testing procedure from above on the regularized, fully-discrete level, we make use of \textit{discrete chain rules}. One very important aspect is a careful and non-trivial approximation of the convective term 
$\skp{(\bv\cdot\nabla)\bbB}{\bbG}_{L^2}$ from \eqref{eq:B_weak} and the term $\frac12 \skp{\kappa(\phi) \nabla\trace\bbB}{\bw}_{L^2}$ from \eqref{eq:v_weak}, respectively,
for which we refer to Remark \ref{remark:scheme0} and Section \ref{sec:regularizations}.
We note that similar strategies have already been used in \cite{barrett_boyaval_2009, barrett_2018_fene-p}. 
However, in our case, the approximation of the convective term is more involved, as we have to take into account additional terms that result from the quadratic part of the elastic energy density. For more details, we refer to Section \ref{sec:regularizations}. 

After using an analogue of the formal testing procedure from above on the regularized fully-discrete level, we derive \textit{a priori} estimates which are uniform in the regularization parameter $\delta>0$ and in the spatial and temporal discretization parameters $h, \Delta t>0$, respectively. 
To justify the existence of such a discrete solution, we combine the \textit{a priori} estimates and a fixed-point argument on the finite dimensional level.
After passing to the limit in $\delta\to 0$ and using subsequence convergence, we find a discrete solution with a positive definite Cauchy--Green tensor.
The precise existence result of discrete solutions is summarized in Theorem \ref{theorem:existence_FE}.

Section \ref{sec:convergence} is devoted to the proof of Theorem \ref{theorem:weak_solution}. First, higher order estimates are derived for all solutions of the discrete scheme from Section \ref{sec:approximation}, so that one can apply compactness results and find subsequences of discrete solutions, which converge to some limit functions. The most involved part is to show that the limit functions satisfy the weak formulation \eqref{eq:div_weak}--\eqref{eq:B_weak} from Definition \ref{def:weak_solution}. Here, we need to pass to the limit in the numerical scheme, which turns out to be quite technical, as we include numerical integration in the discrete system. 
After tackling these difficulties, we pass to the limit in the discrete system which will finally prove Theorem \ref{theorem:weak_solution}. 

Afterwards, in Section \ref{sec:numeric}, we present several numerical results in two and three space dimensions in order to demonstrate the practicability of the discrete scheme from Section \ref{sec:approximation}.

\subsection{Notation}
We close this section by introducing the used notation.
Vector or matrix valued quantities are usually represented with a bold or blackboard bold font, respectively.
Let $d\in\{2,3\}$. The Euclidean scalar product of two vectors $\bv,\bw\in\bbR^d$ and the Frobenius inner product of two matrices $\bbA,\bbB\in\bbR^{d\times d}$ are defined as $\bv\cdot\bw = \bv^\top \bw = \bw^\top \bv$ and $\bbA:\bbB = \trace(\bbA^\top \bbB) = \trace(\bbB^\top \bbA)$, respectively, where $\trace(\bbA)$ is the trace of a matrix $\bbA\in\bbR^{d\times d}$.
Moreover, by $\bbR^{d\times d}_{\mathrm{S}}$ we mean the set of symmetric real-valued $d\times d$ matrices.
The Euclidean norm for vectors and the Frobenius norm for matrices are both denoted by $\abs{\cdot}$.
We use the notation $(\nabla\bv)_{i,j}=\partial_{x_j} \bv_i$, $i,j\in\{1,...,d\}$, and $(\nabla\bbB)_{i,j,k} = \partial_{x_k} \bbB_{i,j}$, $i,j,k\in\{1,...,d\}$, for the gradient of any vector $\bv\in\bbR^d$ and of any matrix $\bbB\in\bbR^{d\times d}$, respectively. 
The divergence of a vector $\bv\in\bbR^d$ and the divergence of a matrix $\bbB\in\bbR^{d\times d}$, respectively, are defined by $\Div \bv = \sum_{j=1}^d \partial_{x_j} \bv_j$ and $(\Div \bbB)_i = \sum_{j=1}^d \partial_{x_j} \bbB_{i,j}$, $i\in\{1,...,d\}$.
The symmetrized gradient of a vector $\bv\in\bbR^d$ is defined as $\D(\bv)\coloneqq \frac12 \big( \nabla\bv + (\nabla\bv)^\top \big)$. Besides, given a vector $\bv\in\bbR$ and a matrix $\bbB\in\bbR^{d\times d}$, we denote the corresponding convective derivative by $(\bv\cdot\nabla)\bbB = \sum_{k=1}^d \bv_k \partial_{x_k} \bbB$. 
Also, we write $\nabla\bbA : \nabla\bbB = \sum_{k=1}^d (\partial_{x_k}\bbA) : (\partial_{x_k} \bbB)$ for the inner product of gradients of matrices.
Given a real Banach space $X$, we denote by $\norm{\cdot}_{X}$ its norm, by $X'$ its dual space, and by $\dualp{\cdot}{\cdot}_{X}$ the duality pairing between $X$ and $X'$. If $X=H$ is a Hilbert space, we denote its inner product by $\skp{\cdot}{\cdot}_H$.
For $p\in[1,\infty]$, an integer $m\geq 0$ and an open set $U\subset\bbR^n$, $n\in\bbN$, we use the standard notation for Lebesgue and Sobolev spaces with values in X, i.e., $L^p(U;X)$ and $W^{m,p}(U;X)$, respectively. If $p=2$ and if $X$ is a Hilbert space, we also write $H^m(U;X)=W^{m,2}(U;X)$.
In the case $U=(0,T)$, where $T>0$, we write $W^{m,p}(0,T;X)=W^{m,p}((0,T);X)$. 
If $X=\bbR$ and $U=\Omega\subset\bbR^d$ is a bounded domain, we use the notation $W^{m,p}(\Omega) = W^{m,p}(\Omega;\bbR)$.
Also, if $X\in\{\bbR, \bbR^d, \bbR^{d\times d}_{\mathrm{S}}\}$, we sometimes write $W^{m,p} = W^{m,p}(\Omega) = W^{m,p}(\Omega;X)$ if the choice of $X\in\{\bbR, \bbR^d, \bbR^{d\times d}_{\mathrm{S}}\}$ is clear from the context. For the $W^{m,p}$ semi-norms we write $\abs{\cdot}_{W^{m,p}}$ where $p\in[1,\infty]$ and $m\geq 0$ is an integer.
We adopt the above notation for $L^p(U;X)$ and $H^m(U;X)$ in a natural way.
Moreover, by $C(\overline\Omega; X)$ we mean the set of continuous functions on the closure of $\Omega$ with values in a Banach space $X$. By $C^\infty(\overline\Omega;X)$ we denote the set of smooth functions on $\overline\Omega$ with values in $X$, and $C_0^\infty(\Omega;X)$ is the set of smooth functions with compact support and with values in $X$.
If $\Omega$ is like in \eqref{eq:phi}--\eqref{eq:bc_v}, i.e., with Lipschitz boundary $\partial\Omega$ and with $\partial_{\mathrm{D}}\Omega \subset\partial\Omega$ being closed and of positive surface measure, we define 
$H^1_\mathrm{D}(\Omega;\bbR^d)$ and $W^{1,p}_{\mathrm{D}}(\Omega;\bbR^d)$, $p\in[2,\infty]$, as the closures of $\{ \bv\in C^\infty(\overline\Omega;\bbR^d) \mid \bv|_{\partial_{\mathrm{D}}\Omega}=\mathbf{0} \}$ with respect to the $H^1$ and $W^{1,p}$ norms, respectively.

\section{Approximation with a fully-practical numerical scheme}
\label{sec:approximation}

For numerical computations, it is of high interest to investigate a well-posed discrete problem in the sense that, on the one hand, there exists at least one solution, and, on the other hand, all solutions are stable. Both results should be valid without a CFL constraint, which means, that the time step size $\Delta t$ and the mesh parameter $h$ can be chosen independently of each other.
One common strategy is to discretize the mathematical model in a way such that energy estimates on the formal level translate to discrete analogues on the approximate level. 
Before we present the numerical approximation of the model \eqref{eq:phi}--\eqref{eq:B}, we recall some definitions and useful results concerning finite elements.

\subsection{Discrete setting}
For future reference, we state the following assumption.

\begin{itemize}

\item[\mylabel{S}{$(\mathbf{S})$}] 
Let $T>0$ and let $\Omega\subset\bbR^d$, $d\in\{2,3\}$, be a 
bounded Lipschitz domain with polygonal (or polyhedral, respectively) boundary $\partial\Omega$. 
Moreover, we assume that $\partial_\mathrm{D}\Omega \subset \partial\Omega$ is closed and of positive surface measure, and we set $\partial_\mathrm{N}\Omega = \partial\Omega \setminus \partial_\mathrm{D}\Omega$.
We split the time interval $[0,T)$ into equidistant subintervals $[t^{n-1},t^n)$ with $t^n = n \Delta t$ and $t^{N_T}=T$, where $\Delta t = \frac{T}{N_T}$, $N_T\in\bbN$ and $n \in \{0, ..., N_T\}$.
We require $\{\calT_h\}_{h>0}$ to be a family of conforming partitionings of $\Omega$ into disjoint open simplices $K$ with $h_K = \diam(K)$ and $h=\max_{K\in\calT_h} h_K$, such that $\overline\Omega = \bigcup_{K\in\calT_h} \overline K$. 
The set of closed edges of triangles ($d=2$) in the mesh $\calT_h$ or closed facets of tetrahedra ($d=3$), respectively, is denoted by $\partial\calT_h$. 
We always assume that $\partial_{\mathrm{D}}\Omega$ and $\calT_h$ are given such that $\partial_{\mathrm{D}}\Omega$ is matched exactly by sides in $\calT_h$, i.e., $\partial_{\mathrm{D}}\Omega = \bigcup \{ S_h \mid S_h\in \partial\calT_h, \, S_h \subset \partial_{\mathrm{D}}\Omega \}$.
The set of all the vertices of $\calT_h$ is denoted by $\calN_h$. 
Moreover, we assume that $\{\calT_h\}_{h>0}$ is shape regular (or non-degenerate), i.e., it holds 
\begin{align*}
    \sup_{K\in\calT_h} h_K \rho_K^{-1} \leq C,
\end{align*}
where $\rho_K$ denotes the diameter of the largest inscribed ball in the simplex $K\in\calT_h$.
In addition, we assume that the family of meshes $\{\calT_h\}_{h>0}$ consists only of non-obtuse simplices, i.e., all dihedral
angles of any simplex in $\calT_h$ are less than or equal to $\frac\pi2$.
\end{itemize}


Sometimes in the literature, the terms ``locally quasi-uniform'' 
or ``quasi-uniform'' 
are used instead of shape regular. Here, we use the notation of, e.g., \cite{bartels_2016}, where ``quasi-uniform'' is used for the condition $h_K \geq Ch$ $\forall\, K\in\calT_h$. Note that in some other works, 
a family of triangulations is called ``quasi-uniform'' if it is shape regular and also the condition $h_K \geq Ch$ $\forall\, K\in\calT_h$ holds true.

Let $\hat K$ denote the standard open reference simplex in $\bbR^d$. Given a simplex $K\in\calT_h$, we denote the affine transformation from $\hat K$ to $K\in\calT_h$ by
\begin{align}
    \label{eq:trafo_K}
    \calB_K \colon\, \hat K \to K, \quad \hat{\mathbf{x}}\mapsto P_0^K + \calA_K \hat{\mathbf{x}},
\end{align}
where $\calA_K \in \bbR^{d\times d}$ is a non-singular matrix. 
Under the assumption \ref{S}, in particular, if the family of triangulations is shape regular, it holds
\begin{align}
    \label{eq:shape_regular}
    \abs{\calA_K^\top} \abs{(\calA_K^\top)^{-1}} \leq C \quad \forall\, K\in\calT_h.
\end{align}
Moreover, a consequence of \ref{S}, i.e., if the simplices $K\in\calT_h$ are non-obtuse, is the following inequality, which will be used in the proof of \eqref{eq:nabla_ln}.
Let $K\in\calT_h$ with local vertices $P_0^K, ..., P_d^K$ and let $\eta_0^K(\cdot), ..., \eta_d^K(\cdot) \in \calP_1(K)$ such that $\eta_i^K(P_j^K) = \delta_{ij}$ for all $i,j\in\{0,...,d\}$, where $\delta_{ij}$ denotes the Kronecker delta. Under the assumption \ref{S}, i.e., if the simplex $K$ is non-obtuse, it holds
\begin{align}
\label{eq:non-obtuse}
    \nabla \eta_i^K \cdot \nabla \eta_j^K \leq 0 
    \quad \text{on } K, \, 
    \forall\, i, j\in\{0,...,d\} \text{ with } i\not= j.
\end{align}

\bigskip
\subsubsection{Discrete function spaces}
For the approximation of the system \eqref{eq:phi}--\eqref{eq:B} with a finite dimensional system, we introduce the following discrete function spaces: 
\begin{subequations}
\begin{align*}
    \calS_h &\coloneqq \left\{ q_h \in C(\overline\Omega) \mid q_h|_{K} \in \calP_1(K) \ \forall \,  K\in\calT_h \right\} 
    \subset W^{1,\infty}(\Omega),
    \\
    \calW_h &\coloneqq \left\{ \bbB_h \in C(\overline\Omega;\bbR^{d\times d}_{\mathrm{S}}) \mid \bbB_h|_{K} \in \calP_1(K; \bbR^{d\times d}_{\mathrm{S}}) \ \forall \,  K\in\calT_h \right\} 
    \subset W^{1,\infty}(\Omega;\bbR^{d\times d}_{\mathrm{S}}),
    \\
    \calV_h &\coloneqq \left\{ \bv_h \in C(\overline\Omega;\bbR^d) \cap H^1_{\mathrm{D}}(\Omega;\bbR^d) \mid \bv_{h}|_{K} \in \calP_2(K;\bbR^d) \ \forall \,  K\in\calT_h \right\} \subset W_{\mathrm{D}}^{1,\infty}(\Omega;\bbR^d).
\end{align*}
Here, $\calP_j(K;X)$ denotes the set of polynomials of order $j\in\{1,2\}$ on $K\in\calT_h$ with values in $X\in\{\bbR, \bbR^d, \bbR^{d\times d}_{\mathrm{S}} \}$, and we use the convention $\calP_j(K) = \calP_j(K;\bbR)$.
The function space $\calV_h \times \calS_h$ is also referred to as the $\calP_2$/$\calP_1$-Taylor--Hood element \cite{girault_raviart_1986} for the discrete velocity field and pressure.
Later, we will also need
\begin{align*}
    \calV_{h,\mathrm{div}} &\coloneqq \left\{ \bv_h \in \calV_h \mid  \skp{\Div\bv_h}{q_h}_{L^2} = 0 \ \forall \,  q_h \in \calS_h \right\}.
\end{align*}
\end{subequations}
For any $r, r'\in(1,\infty)$ with $\frac1r + \frac{1}{r'} = 1$, it is well-known that the $\calP_2$/$\calP_1$-Taylor--Hood element $\calV_h \times \calS_h$ satisfies the stability condition
\begin{align}
    \label{eq:LBB}
    \sup_{\bw_h\in \calV_h \setminus\{0\}}
    \frac{\skp{\Div\bw_h}{q_h}_{L^2} }{ \norm{\bw_h}_{W^{1,r}}  }
    \geq C_{\mathrm{stab}} \norm{q_h}_{L^{r'}} \quad \forall\, q_h \in \calS_h,
\end{align}
where $C_{\mathrm{stab}}>0$ denotes a constant that is independent of the mesh parameter $h>0$, see, e.g., \cite[Chap.~4.2.5]{ern_guermond_2004}. 
The inequality \eqref{eq:LBB} is often referred to as inf--sup inequality or Ladyzhenskaya--Babu{\v s}ka--Brezzi (LBB) condition and it holds true under the assumption \ref{S}, more precisely, if the family of meshes $\{\calT_h \}_{h>0}$ is shape-regular.
Sometimes, an additional (but very mild) condition on the triangulation is posed due to technical reasons, but usually it can be dropped if the triangulation $\calT_h$ is fine enough, see \cite[Chap.~8.8]{BBF2013_fem}. 
The classical inf--sup stability condition with $r=r'=2$ is a special case of \eqref{eq:LBB}. 
The results of this work remain valid if other discrete function spaces $\widetilde\calV_h$ for the velocity field are used, supposed that the discrete inf--sup stability condition \eqref{eq:LBB} is fulfilled for $\widetilde\calV_h \times \calS_h$. Two possible examples are the $\calP_1$-bubble/$\calP_1$ element (\textit{mini-element}) and the $\calP_1$-iso-$\calP_2$/$\calP_1$ element, see \cite[Chap.~4.2.4]{ern_guermond_2004} and \cite[Chap.~4.2.6]{ern_guermond_2004}, respectively.


For future reference, we note the following two results which are consequences of the inf--sup stability condition \eqref{eq:LBB}, see, e.g., 
\cite[Lem.~4.1]{girault_raviart_1986}: 
For any linear functional $\ell_h \colon\,\calS_h\to \bbR$, there exists a velocity field 
$\bu_h \in  \calV_h$ (which is unique in a subset $(\calV_{h,\mathrm{div}})^\perp \subset \calV_h$, where $(\calV_{h,\mathrm{div}})^\perp$ denotes the orthogonal of $\calV_{h,\mathrm{div}}$ in $H^1_{\mathrm{D}}(\Omega;\bbR^d)$ for the scalar product associated with $\norm{\nabla\cdot}_{L^2}$) such that, for any $r, r'\in(1,\infty)$ with $\frac1r + \frac{1}{r'} = 1$,
\begin{subequations}
\begin{align}
    \label{eq:LBB_u}
    \skp{\Div\bu_h}{q_h}_{L^2} &= \ell_h(q_h) \quad \forall\, q_h\in\calS_h,
    \\
    \label{eq:LBB_u_bound}
    \norm{\bu_h}_{W^{1,r}} &\leq C_{\mathrm{stab}}^{-1} 
    \sup_{q_h\in\calS_h \setminus \{0\}} \frac{\ell_h(q_h)}{\norm{q_h}_{L^{r'}}},
\end{align}
\end{subequations}
where $C_{\mathrm{stab}}>0$ is the constant from \eqref{eq:LBB}.
%
%
%
%
Moreover, for any linear functional $\calF_h \colon\,\calV_h\to \bbR$ with $\calF_h|_{\calV_{h,\mathrm{div}}}=0$,
there exists a unique scalar pressure $p_h\in\calS_h$
such that, for any $r, r'\in(1,\infty)$ with $\frac1r + \frac{1}{r'} = 1$,
\begin{subequations}
\begin{align}
    \label{eq:LBB_p}
    \skp{\Div\bw_h}{p_h}_{L^2} &= \calF_h(\bw_h) \quad \forall\, \bw_h \in \calV_h,
    \\
    \label{eq:LBB_p_bound}
    \norm{p_h}_{L^{r'}} &\leq C_{\mathrm{stab}}^{-1} 
    \sup_{\bw_h\in\calV_h \setminus \{0\}} \frac{\calF_h(\bw_h)}{\norm{\bw_h}_{W^{1,r}}} .
\end{align}
\end{subequations}

\subsubsection{Results concerning interpolation and projections}
In this work, we will frequently use the standard nodal interpolation operator $\calI_h\colon\, C(\overline\Omega)\to \calS_h$, given by $(\calI_h \eta)(P) \coloneqq \eta(P)$ $\forall\, P \in \calN_h$.
As we use numerical integration (mostly in terms of \textit{mass lumping}), we introduce the following semi-inner product and the induced semi-norm on $C(\overline\Omega)$ by
\begin{alignat*}{5}
    &\skp{\eta_1}{\eta_2}_h &&\coloneqq 
    \int_\Omega \calI_h \big[ \eta_1 \eta_2 \big] \dx,
    \quad\quad
    &&\norm{\eta_1}_h &&\coloneqq \sqrt{\skp{\eta_1}{\eta_1}_h},
    \quad\quad
    &&\forall \,  \eta_1, \eta_2 \in C(\overline\Omega).
\end{alignat*}
Moreover, we introduce the \textit{lumped} $L^2$-projector $\calQ_h\colon\, L^2(\Omega) \to \calS_h$, where, for $\eta\in L^2(\Omega)$, its \textit{lumped} $L^2$-projection $\calQ_h\eta \in \calS_h$ is characterized by
\begin{align}
    \label{eq:proj_def}
    \skp{\calQ_h \eta}{\zeta_h}_h = \skp{\eta}{\zeta_h}_{L^2}
    \quad \forall\, \zeta_h \in \calS_h.
\end{align}
Now, we recall the following well-known results concerning the finite element space $\calS_h$ and the operators $\calI_h$ and $\calQ_h$, which hold true under the assumption \ref{S}, i.e., if the family of triangulations is shape regular:
\begin{alignat}{2}
    \label{eq:interp_H2}
    \norm{\eta - \calI_h \eta}_{L^p(K)}
    + h_K \norm{\nabla(\eta-\calI_h\eta)}_{L^p(K)} 
    &\leq C h_K^2 \abs{\eta}_{W^{2,p}(K)}       \quad &&\forall\, \eta\in W^{2,p}(K), \ p\in[2,\infty], \ \forall\, K\in\calT_h,
    \\
    \label{eq:interp_stability}
    \norm{\nabla \calI_h \eta}_{L^\infty(K)}
    &\leq C \norm{\nabla \eta}_{L^\infty(K)}
    \quad
    &&\forall\, \eta\in W^{1,\infty}(K), \ \forall\, K\in\calT_h,
    \\
    \label{eq:interp_continuous}
    \lim_{h\to 0} \norm{\eta - \calI_h\eta}_{L^\infty}
    &= 0
    \quad
    &&\forall \,  \eta \in C(\overline\Omega),
    \\
    \label{eq:proj_error}
    \norm{\eta - \calQ_h \eta}_{L^2(K)}
    + h_K \norm{\nabla(\eta - \calQ_h\eta)}_{L^2(K)}
    &\leq C h_K \norm{\nabla\eta}_{L^2(K)}    
    \quad 
    &&\forall\, \eta\in H^1(K), \ \forall\, K\in\calT_h.
\end{alignat}
For proofs, we refer to, e.g., \cite{barrett_nurnberg_styles_2004, barrett_suli_2011, GKT_2022_viscoelastic} or the books \cite{bartels_2016, brenner_scott_2008} and references therein. We note that similar results hold true for the matrix valued finite element space $\calW_h$ in a natural way with $\calI_h, \calQ_h$ naturally extended to matrix valued functions.
Also, it holds
\begin{align}
    \label{eq:norm_equiv}
    \norm{q_h}_{L^2(K)}^2
    \leq 
    \int_K \calI_h \big[ \abs{q_h}^2 \big] \dx
    \leq C \norm{q_h}_{L^2(K)}^2
    \quad 
    \forall\, q_h\in \calS_h, \ \forall\, K\in\calT_h,
\end{align}
which, after summation over $K\in\calT_h$, implies norm equivalence of $\norm{\cdot}_{L^2}$ and $\norm{\cdot}_h$ on $\calS_h$ independently of $h>0$. 
Besides, the following local inverse estimate holds true for any $u_h \in \calS_h$, $K\in\calT_h$, $s,m\in\{0,1\}$ with $s\leq m$ and $1\leq r \leq p \leq \infty$:
\begin{align}
    \label{eq:inverse_estimate}
    \abs{u_h}_{W^{m,p}(K)} 
    &\leq 
    C h_K^{s - m + \frac{d}{p} - \frac{d}{r}} \abs{u_h}_{W^{s,r}(K)}.
\end{align}
%
%
%
%
Note that \eqref{eq:inverse_estimate} is also valid if $u_h \in \calV_h$ or $u_h\in \calW_h$.
If, in addition to \ref{S}, the family of triangulations $\{\calT_h\}_{h>0}$ is quasi-uniform, i.e., $h_K \geq Ch$ $\forall\, K\in\calT_h$, then \eqref{eq:inverse_estimate} holds true globally, i.e., with $h_K$ replaced by $h$ and $K$ replaced by $\Omega$ (after summation over all $K\in\calT_h$).

For future reference, we also note the following local error estimate 
\begin{align}
    \label{eq:interp_error_Sh_Sh}
    \norm{\calI_h[ q_h\zeta_h] - q_h \zeta_h}_{L^p(K)}
    \leq C h_K^{1+m} \abs{q_h}_{W^{m,\infty}(K)} \norm{\nabla\zeta_h}_{L^p(K)}
    \quad \forall\, q_h, \zeta_h \in \calS_h, \ m\in\{0,1\}, \ p\in[1,\infty], \ \forall K \in \calT_h,
\end{align}
which follows from \eqref{eq:interp_H2}, \eqref{eq:inverse_estimate}, Hölder's inequality and the fact that $q_h, \zeta_h\in\calS_h$ are affine functions on $K$.
We also recall the well-known mass lumping error estimate, which is a direct consequence of \eqref{eq:interp_error_Sh_Sh} with $p=1$ and after summation over $K\in\calT_h$:
\begin{align}
    \label{eq:lump_Sh_Sh}
    \abs{ \skp{q_h}{\zeta_h}_h 
    - \skp{q_h}{\zeta_h}_{L^2} } 
    &\leq
    C h^{1+m} \abs{q_h}_{H^m} \norm{\nabla \zeta_h}_{L^2}
    \quad
    \forall \,   q_h, \zeta_h \in \calS_h, \ m\in\{0,1\}.
\end{align}

Moreover, for any Lipschitz continuous functions $g\colon\, \bbR^m\to\bbR$, $\tilde g\colon\,\bbR^m \times \bbR^{d\times d} \to \bbR$, $m\in\bbN$, with Lipschitz constants $L_g, L_{\tilde g}>0$, respectively, it holds for any $\mathbf{q}_h = (q_{h,1},...,q_{h,m}) \in (\calS_h)^m$, $\bbG_h\in\calW_h$, $p\in[1,\infty]$ and $K\in \calT_h$, that
\begin{subequations}
\begin{alignat}{2}
    \label{eq:interp_lipschitz}
    \norm{ \calI_h [g(\mathbf{q}_h)] - g(\mathbf{q}_h)}_{L^p(K)}
    &\leq C L_g h_K \sum_{i=1}^m \norm{\nabla q_{h,i}}_{L^p(K)},
    \\
    \label{eq:interp_lipschitz_multivar}
    \norm{ \calI_h [\tilde g(\mathbf{q}_h, \bbG_h)] - \tilde g(\mathbf{q}_h, \bbG_h)}_{L^p(K)}
    &\leq C L_{\tilde g} h_K \Big( \sum_{i=1}^m \norm{\nabla q_{h,i}}_{L^p(K)} + \norm{\nabla \bbG_h}_{L^p(K)} \Big),
\end{alignat}
\end{subequations}
see \cite[Lem.~6.8]{barrett_2018_fene-p} for the first inequality with $p=2$ and $m=1$. 
The first inequality with arbitrary $p\in[1,\infty]$, $m\in\bbN$ and the second inequality can be shown in the same way after minor adaptions.

\subsection{The fully-practical numerical scheme}

For the well-posedness of the numerical scheme, we need the following assumption which helps us to approximate the nonlinear model functions $\psi(\cdot)$ and $\kappa(\cdot)$.

\begin{itemize}	    
\item[\mylabel{A6}{$(\mathbf{A6})$}] 
In addition to \ref{A3}--\ref{A4}, we define $\kappa_h'(a,b) \coloneqq \frac{\kappa(a)-\kappa(b)}{a-b}$ $\forall\,a,b \in\bbR$ with $a\not=b$, and $\kappa_h'(a,a) \coloneqq \kappa'(a)$ $\forall\, a\in\bbR$.
Moreover, we assume that $\psi_h'(\cdot,\cdot)$ is a Lipschitz continuous function and fulfills $\forall\, a,b\in\bbR$:
\begin{align*}
    \psi_h'(a,a) & = \psi'(a), \\
    \psi_h'(a,b) (a-b) &\geq \psi(a) - \psi(b),  \\
    \abs{\psi_h'(a,b)} &\leq C (1 + \abs{a} + \abs{b}).
\end{align*}
\end{itemize}
Possible examples for $\psi_h'(\cdot,\cdot)$ are the classical convex-concave splitting 
\begin{align*}
    \psi_h'(\phi_h^n,\phi_h^{n-1}) 
    &\coloneqq \psi_1'(\phi_h^n) + \psi_2'(\phi_h^{n-1}) \quad \text{where } \psi=\psi_1+\psi_2 \text{ with } \psi_1 \text{ convex and } \psi_2 \text{ concave,}
\end{align*}
or an approximation with a difference quotient like for $\kappa_h'(\cdot,\cdot)$.
Note that a convex-concave splitting for $\psi$ is possible due to the quadratic growth assumption \ref{A4}.

We now present the numerical scheme.

\subsubsection*{Problem \ref{P_FE}:}
\mylabelHIDE{P_FE}{$(\pmbP_{h}^{\Delta t})$}
Let $n\in\{1,...,N_T\}$ and suppose that $\phi_h^{n-1} \in \calS_h$, $\sigma_{\infty,h}^n\in\calS_h$ and $\bbB_h^{n-1} \in\calW_h$ 
are given. Then, the goal is to find a solution tuple 
\begin{align*}
    (\phi_h^{n}, \mu_h^n, \sigma_{h}^{n}, p_h^n, \bv_{h}^{n}, \bbB_{h}^{n}) \in (\calS_h)^4 \times \calV_h  \times \calW_{h}
    \quad \text{with } \bbB_h^n \text{ positive definite},
\end{align*}
which satisfies for any test function tuple $(\zeta_h, \rho_h, \xi_h, q_h, \bw_h, \bbG_h) \in (\calS_h)^4 \times \calV_h \times \calW_h$:
\begin{subequations}
\begin{align}
    \label{eq:phi_FE} \nonumber
    0 &= 
    \frac{1}{\Delta t} \skp{\phi_h^n - \phi_h^{n-1}}{\zeta_h}_h
    + \skp{\calI_h [m_{\phi,h}^{n-1}]  \nabla\mu_h^n}{\nabla\zeta_h}_{L^2}
    + \skp{\bv_h^n \cdot \nabla\phi_h^{n-1}}{\zeta_h}_{L^2}
    \\
    &\quad
    + \skp{\phi_h^{n-1} \Gamma_{\bv,h}^{n-1} - \Gamma_{\phi,h}^{n-1}}{\zeta_h}_h,
    \\[5pt]
    \label{eq:mu_FE}
    0 &= 
    \skp{- \mu_h^n + \frac\beta\varepsilon \psi_h'(\phi_h^n,\phi_h^{n-1}) - \chi_\phi \sigma_h^n + \frac12 \kappa_h'(\phi_h^n, \phi_h^{n-1}) \trace\bbB_h^n} {\rho_h}_h
    + \beta\varepsilon \skp{ \nabla\phi_h^n}{\nabla\rho_h}_{L^2},
    \\[5pt]
    \label{eq:sig_FE}
    0 &= \skp{\calI_h [m_{\sigma,h}^{n-1}] \nabla\sigma_h^n}{\nabla\xi_h}_{L^2}
    + \skp{\sigma_h^n \Gamma_{\sigma,h}^{n-1}}{\xi_h}_h
    + K \skp{\sigma_h^n - \sigma_{\infty,h}^n}{\xi_h}_{L^2(\partial\Omega)},
    \\[5pt]
    \label{eq:div_FE}
    0 &= \skp{\Div\bv_h^n}{q_h}_{L^2} - \skp{\Gamma_{\bv,h}^{n-1}}{q_h}_{h},
    \\[5pt]
    \nonumber
    \label{eq:v_FE}
    0 &= \skp{2 \calI_h[\eta(\phi_h^{n-1})] \D(\bv_h^n)}{\D(\bw_h)}_{L^2}
    - \skp{p_h^n}{\Div\bw_h}_{L^2}
    + \skp{\calI_h \bbT_{e,h}^n}{\nabla\bw_h}_{L^2}
    \\
    &\quad \nonumber
    - \skp{(\mu_h^n + \chi_\phi \sigma_h^n) \nabla\phi_h^{n-1}}{\bw_h}_{L^2}
    - \frac12 \skp{\nabla \calI_h\left[ \kappa(\phi_h^{n-1}) \trace\bbB_h^n \right]}{\bw_h}_{L^2}
    \\
    &\quad
    + \frac12 \sum_{i,j=1}^d \skp{\partial_{x_j} \calI_h\left[ \kappa(\phi_h^{n-1}) \right]  \trace \mathbf{\Lambda}_{i,j}(\bbB_h^n)}{(\bw_h)_i}_{L^2},
    \\[5pt]
    \label{eq:B_FE}
    \nonumber
    0 &= 
    \frac{1}{\Delta t} \skp{\bbB_h^n - \bbB_h^{n-1}}{\bbG_h}_h
    + \skp{\bv_h^n}{ \nabla\calI_h\left[\bbB_h^n : \bbG_h \right]}_{L^2}
    - \sum_{i,j=1}^d \skp{(\bv_h^n)_i \mathbf{\Lambda}_{i,j}(\bbB_h^n)}{\partial_{x_j} \bbG_h}_{L^2}
    \\
    &\quad 
     + \skp{\tfrac{1}{\tau(\phi_h^{n-1})} \bbT_{e,h}^n }{\bbG_h}_h 
    + \skp{\Gamma_{\bbB,h}^{n-1} \bbB_h^n}{\bbG_h}_h
    - \skp{2 \nabla\bv_h^n}{\calI_h\left[ \bbG_h \bbB_h^n \right]}_{L^2}
    + \alpha \skp{\nabla\bbB_h^n}{\nabla\bbG_h}_{L^2},
\end{align}
\end{subequations}
where $\bbD(\bw_h) \coloneqq \frac12 \big(\nabla\bw_h + (\nabla\bw_h)^\top \big)$ denotes the symmetrized gradient for $\bw_h\in\calV_h$. Moreover, for the elastic stress tensor, the source terms and the mobility functions, we use the notation
\begin{align*}
    \bbT_{e,h}^n 
    &\coloneqq  (\bbB_h^n)^2 + \kappa(\phi_h^{n-1}) \bbB_h^n - \bbI  ,
    \\
    \Gamma_{f,h}^{n-1} &\coloneqq \Gamma_{f}(\phi_h^{n-1}, \sigma_h^n, \bbB_h^{n-1}) , \quad \text{for } f\in\{\phi,\bv,\bbB\},
    \\
    \Gamma_{\sigma,h}^{n-1} &\coloneqq \Gamma_\sigma(\phi_h^{n-1}, \bbB_h^{n-1}),
    \\
    m_{f,h}^{n-1} &\coloneqq m_f(\phi_h^{n-1}, \bbB_h^{n-1}), \quad \text{for } f\in\{\phi,\sigma\}.
\end{align*}


\begin{remark}~
\begin{enumerate}[(i)]
\label{remark:scheme0}



\item 
The structure of the numerical scheme \ref{P_FE} is based on a discretization of the weak formulation from Definition \ref{def:weak_solution} with a semi-implicit Euler scheme in time and a finite element ansatz in space. 
In \eqref{eq:v_FE} and \eqref{eq:B_FE}, we made an approximation of the terms $\kappa(\phi) \nabla \trace\bbB \cdot \bw$ and $(\bv\cdot\nabla)\bbB : \bbG$ which appear in \eqref{eq:v_weak} and \eqref{eq:B_weak}, respectively. More precisely, we approximate them via
\begin{align*}
    -\kappa(\phi) \nabla \trace\bbB \cdot \bw
    &= - \nabla\left[ \kappa(\phi) \trace\bbB \right] \cdot \bw 
    + \trace\bbB \nabla\kappa(\phi) \cdot \bw
    \\
    &\approx 
    -\nabla \calI_h\left[ \kappa(\phi_h^{n-1}) \trace\bbB_h^n \right] \cdot \bw_h
    + \sum_{i,j=1}^d \partial_{x_j} \calI_h\left[ \kappa(\phi_h^{n-1}) \right]  \trace \mathbf{\Lambda}_{i,j}(\bbB_h^n) (\bw_h)_i
\end{align*}
and
\begin{align*}
    (\bv\cdot\nabla)\bbB : \bbG
    &= \bv\cdot\nabla \left[\bbB : \bbG\right]
    - (\bv\cdot\nabla)\bbG : \bbB
    \\
    &\approx
    \bv_h^n \cdot \nabla\calI_h\left[\bbB_h^n : \bbG_h \right]
    - \sum_{i,j=1}^d (\bv_h^n)_i \mathbf{\Lambda}_{i,j}(\bbB_h^n) : \partial_{x_j} \bbG_h,
\end{align*}
where $\mathbf{\Lambda}_{i,j}(\bbB_h^n) \approx \bbB_h^n \, \delta_{i,j} \in \bbR^{d\times d}$ for any $i,j\in\{1,...,d\}$, where $\delta_{i,j}$ denotes the Kronecker delta. We refer to Section \ref{sec:regularizations} for more details.
The reason for that is, that in the testing procedure for energy estimates, we choose $\bbG_h$ as a finite element approximation of $\bbB + \kappa(\phi)\bbI - \bbB^{-1}$. The specific structure of the approximation of $(\bv\cdot\nabla)\bbB : \bbG$ in \eqref{eq:B_FE} allows a discrete analogue of the chain rule
\begin{align*}
    (\bv\cdot\nabla)\bbB : \left[ \bbB + \kappa(\phi)\bbI - \bbB^{-1} \right] 
    = \bv\cdot\nabla \left( \tfrac12 \abs{\bbB}^2 - \trace\ln\bbB \right) 
    + \kappa(\phi) \bv\cdot\nabla \trace\bbB,
\end{align*}
where $\bv\cdot\nabla \left( \tfrac12 \abs{\bbB}^2 - \trace\ln\bbB \right)$ can be controlled with other terms, and $\kappa(\phi) \bv\cdot\nabla \trace\bbB$ cancels out with the corresponding term in \eqref{eq:v_FE}.
This also justifies the specific structure of the approximation of $\kappa(\phi) \nabla \trace\bbB \cdot \bw$ in \eqref{eq:v_FE}.

\item 
The formulation of the numerical scheme \ref{P_FE} allows to evaluate all integrals exactly by using quadrature rules of third order on each simplex $K\in\calT_h$. If the $\calP_2$/$\calP_1$-Taylor--Hood element $\calV_h \times \calS_h$ is replaced by the inf--sup stable mini-element \cite{girault_raviart_1986}, then a quadrature rule of fifth order for $d=2$ or seventh order for $d=3$, respectively, should be used instead. All nonlinear functions are either interpolated nodewise with the operator $\calI_h$ or contained within the \textit{lumped} $L^2$-inner product $\skp{\cdot}{\cdot}_h$.
In order to solve the nonlinear system at each time step, one can use a fixed-point or (quasi-)Newton method.
We refer to Section \ref{sec:numeric}, where we present an iterative method for the nonlinear system.





\end{enumerate}
\end{remark}

We now state an existence and stability result for the numerical scheme \ref{P_FE}.

\begin{theorem}[Well-posedness of the numerical scheme]
\label{theorem:existence_FE} ~\\
Let \ref{A1}--\ref{A6} and \ref{S} hold true. Let the discrete initial data $\phi_h^0\in\calS_h$ and $\bbB_h^0\in\calW_h$ be given with $\bbB_h^0$ being  positive definite. 
Moreover, for any $n\in\{1,...,N_T\}$, let the discrete boundary data $\sigma_{\infty,h}^n \in\calS_h$ be given. Besides, assume that $\Delta t < \Delta t_*$,
where the constant $\Delta t_*>0$ depends only on the model parameters and on $c_\infty \coloneqq \max_{n\in\{1,...,N_T\}} \norm{\sigma_{\infty,h}^n}_{L^2(\partial\Omega)}$. 
Then, for all $n\in\{1,...,N_T\}$, there exists at least one solution tuple $(\phi_h^{n},\mu_h^n,\sigma_{h}^{n},p_h^n,\bv_{h}^{n},\bbB_{h}^{n}) \in (\calS_h)^4 \times \calV_h \times \calW_h$ to the discrete problem \ref{P_FE} with $\bbB_{h}^{n}$ being positive definite.
Moreover, all solutions of \ref{P_FE} are stable in the sense that
\begin{align}
    \label{eq:stability_FE_sigma}
    &\max_{n=1,...,N_T} \norm{\sigma_h^n}_{H^1}^2
    \leq C 
    \max_{n=1,...,N_T}
    \norm{\sigma_{\infty,h}^n}_{L^2(\partial\Omega)}^2,
    \\[2ex]
    \label{eq:stability_FE} \nonumber
    &  \max_{n=1,...,N_T} \Big( 
    \norm{\phi_h^n}_{H^1}^2
    +  \norm{ \bbB_h^n }_{L^2}^2
    +  \norm{\calI_h \trace \ln(\bbB_h^n)}_{L^1}
    \Big)
    \\
    \nonumber
    &\quad 
    +  \sum_{n=1}^{N_T} \Big(
    \norm{\nabla\phi_h^n - \nabla\phi_h^{n-1}}_{L^2}^2 
    +  \norm{\bbB_h^n - \bbB_h^{n-1}}_{L^2}^2 \Big)
    + \Delta t \sum_{n=1}^{N_T} \Big( 
    \norm{\mu_h^n}_{H^1}^2 
    + \norm{\bv_h^n}_{H^1}^2  
    + \norm{\bbB_h^n}_{H^1}^2  \Big)
    \\
    \nonumber
    &\quad
    + \Delta t \sum_{n=1}^{N_T} \Big( 
    \norm{\nabla\calI_h \trace \ln( \bbB_h^n)}_{L^2}^2
    + \nnorm{\bbT_{e,h}^n \, [\bbB_h^n]^{-1/2}}_h^2
    \Big)
    + \bigg( \Delta t \sum_{n=1}^{N_T} \norm{p_h^n}_{L^2}^{4/3} \bigg)^{3/4}
    \\
    &\leq
    C(T, c_\infty) \big(
    1 + \norm{\bbB_h^0}_{L^2}^2 + \norm{\calI_h\trace\ln\bbB_h^0}_{L^1}
    + \norm{\nabla\phi_h^0}_{L^2}^2
    + \norm{\calI_h \psi(\phi_h^0)}_{L^1} \big) 
\end{align}
where the constants $C, C(T, c_\infty)>0$ are independent of $h, \Delta t$, but $C(T,c_\infty)$ depends exponentially on $T$ and $c_\infty$.
\end{theorem}

We give a remark on the very mild constraint on the time step size, i.e., $\Delta t < \Delta t_*$.

\begin{remark} \label{remark:scheme}
In absence of source terms and chemotaxis, i.e., $\Gamma_\phi = \Gamma_{\bv} = \Gamma_{\bbB} = \chi_\phi = 0$, and for $\kappa(\phi) = \overline\kappa \in\bbR$, the discrete scheme has a dissipative discrete energy and is unconditionally stable in the sense that, for any $h,\Delta t>0$ and $n \in \{1,...,N_T\}$, all solutions of \ref{P_FE} fulfill \eqref{eq:stability_FE_sigma} and
\begin{align*}
    & \frac{\beta\varepsilon}{2} \norm{\nabla\phi_h^n}_{L^2}^2 
    + \frac{\beta}{\varepsilon} \skp{\psi(\phi_h^n)}{1}_h 
    + \frac{\beta\varepsilon}{2} 
    \norm{\nabla\phi_h^n - \nabla\phi_h^{n-1}}_{L^2}^2
    \\
    &\quad
    + \frac14 \norm{\bbB_h^n}_h^2
    + \frac12 \skp{\overline\kappa \trace\bbB_h^n - \trace\ln(\bbB_h^n)}{1}_h
    + \frac{1}{4} \norm{\bbB_h^n- \bbB_h^{n-1}}_{h}^2
    \\
    &\quad
    + \Delta t \skp{\calI_h[ m_{\phi,h}^{n-1}] }{ \abs{\nabla\mu_h^n}^2}_{L^2}
    + 2 \Delta t \skp{\calI_h [\eta(\phi_h^{n-1})] }{ \abs{\D(\bv_h^n)}^2}_{L^2}
    \\
    &\quad
    + \Delta t \skp{\frac{1}{2\tau(\phi_h^{n-1})} }{ \abs{\bbT_{e,h}^{n} \, [\bbB_h^n]^{-1/2}}^2 }_h 
    + \frac12 \alpha \Delta t \norm{\nabla \bbB_h^n}_{L^2}^2
    + \frac{\alpha}{2d} \Delta t \norm{\nabla \calI_h \trace \ln (\bbB_h^n)}_{L^2}^2
    \\
    &\leq 
    \frac{\beta\varepsilon}{2} \norm{\nabla\phi_h^{n-1}}_{L^2}^2 
    + \frac{\beta}{\varepsilon} \skp{\psi(\phi_h^{n-1})}{1}_h 
    + \frac14 \norm{\bbB_h^{n-1}}_h^2
    + \frac12 \skp{\overline\kappa \trace\bbB_h^{n-1} - \trace\ln(\bbB_h^{n-1})}{1}_h .
\end{align*}
The fact that this inequality holds can be shown similarly to \eqref{eq:energy_estimate_dissipative} whereby 
we implicitly assume that $\bbB_h^{n-1}, \bbB_h^n \in \calW_h$ are positive definite. 
As the existence proof of discrete solutions mainly relies on the discrete energy inequality (stability) of the discrete system and a fixed-point argument, it is possible to establish the existence of discrete solutions under the same assumptions as for the stability estimates. More precisely, if the source terms and chemotaxis are neglected and if $\kappa(\phi)$ is constant, then it is also possible to prove the unconditional existence of discrete solutions for any $h,\Delta t>0$ and $n \in \{1,...,N_T\}$.
%
%
%
%

%
Note that, in Theorem \ref{theorem:existence_FE}, the time step size has to satisfy a minor smallness constraint, i.e., $\Delta t < \Delta t_*$ for a constant $\Delta t_*>0$ which depends only on the model parameters but not on $h$. Typically in Cahn--Hilliard models with source terms, one observes the condition $\Delta t = \mathcal{O}(\varepsilon)$ where the small constant $\varepsilon>0$ is related to the interface width, see also \cite{GKT_2022_viscoelastic, trautwein_2021}.
\end{remark}

We now comment on the explicit construction of possible examples for the initial and boundary data under the assumptions \ref{A5} and \ref{S}.

\begin{remark}[Approximation of the initial and boundary data]~
\label{remark:initial_data}
\begin{enumerate}[(i)]
\item One possible choice for the discrete initial datum $\bbB_h^0\in\calW_h$ of the Cauchy--Green tensor is given as the unique solution of the projection problem
\begin{align}
    \label{eq:B0_constr}
    \skp{\bbB_h^0}{\bbG_h}_h + \Delta t \skp{\nabla \bbB_h^0}{\nabla \bbG_h}_{L^2} = \skp{\bbB_0 }{\bbG_h}_{L^2} \quad \forall\, \bbG_h\in\calW_h,
\end{align}
where $\bbB_0$ is given as in \ref{A5}.
Setting $\bbG_h=\bbB_h^0$ and using Hölder's and Young's inequalities and \eqref{eq:norm_equiv}, one obtains the stability estimate
\begin{align}
    \label{eq:B0_constr_stability}
    \norm{\bbB_h^0}_h^2 + \Delta t \norm{\nabla\bbB_h^0}_{L^2}^2 
    \leq C \norm{\bbB_0 }_{L^2}^2.
\end{align}
Here, $\bbB_0$ is uniformly positive definite a.e.~in $\Omega$. In particular, it holds $\pmb\xi^\top \bbB_0(\mathbf{x}) \pmb\xi \geq b_0 \abs{\pmb\xi}^2$ $\forall\, \pmb\xi\in \bbR^d$ and for a.e.~$\mathbf{x}\in \Omega$.
Using \cite[Lem.~5.2]{barrett_boyaval_2009}, one observes that the approximate initial datum $\bbB_h^0$ is also uniformly positive definite in all vertices $P \in \calN_h$, i.e.,
\begin{align*}
    \pmb \xi^\top \bbB_h^0(P) \pmb \xi 
    \geq b_0 \abs{\pmb\xi}^2, 
    \quad \forall\, \pmb\xi\in\bbR^d, \ \forall\, P\in\calN_h.
\end{align*}
A simple calculation shows
\begin{align*}
    \norm{\calI_h \trace\ln\bbB_h^0}_{L^1}
    \leq C \big( \norm{\trace \bbB_h^0}_{L^1}
    + \abs{\ln b_0} \big) 
    \leq C \big( 1 + \norm{\bbB_0}_{L^2}^2 
    + \abs{\ln b_0} \big).
\end{align*}
In addition, it follows from \eqref{eq:B0_constr}, \eqref{eq:B0_constr_stability}, \eqref{eq:lump_Sh_Sh} and a density argument, that $\bbB_h^0\to\bbB_0$ weakly in $L^2(\Omega;\bbR^{d\times d}_{\mathrm{S}})$, as $(h,\Delta t)\to(0,0)$. 


\item Let $\phi_0\in H^1(\Omega)$. Then, $\phi_0$ can be approximated with the help of the \textit{lumped} $L^2$-projector, i.e., $\phi_h^0 = \calQ_h \phi_0 \in \calS_h$. In this case, it holds with \eqref{eq:proj_error}, that
\begin{align*}
    \norm{\phi_h^0}_{H^1} \leq C \norm{\phi_0}_{H^1},
\end{align*}
and $\phi_h^0 \to \phi_0$ strongly in $L^2(\Omega)$, as $h\to 0$.
For more regular initial data, e.g., $\phi_0\in H^2(\Omega)$, one can approximate $\phi_0$ with the help of the nodal interpolation operator, i.e., $\phi_h^0 = \calI_h \phi_0 \in \calS_h$. Then, on noting \eqref{eq:interp_H2} and \eqref{eq:interp_continuous}, it holds
\begin{align*}
    \norm{\phi_h^0}_{H^1}
    \leq C \norm{\phi_0}_{H^2},
\end{align*}
and $\phi_h^0\to \phi_0$ strongly in $H^1(\Omega) \cap C(\overline\Omega)$,
as $h\to 0$.

\item The approximation of the boundary values $\sigma_\infty \in L^\infty(0,T;L^2(\partial\Omega))$ can be established with the following strategy. For $n\in\{1,...,N_T\}$, let $\sigma_{\infty,h}^n \in \calS_h$ be given, such that $\sigma_{\infty,h}^n|_{\partial\Omega}$ solves the projection problem
\begin{align} 
    \label{eq:sigma_infty_constr}
    \skp{\sigma_{\infty,h}^n}{\xi_h}_{L^2(\partial\Omega)}
    = \frac{1}{\Delta t} \int_{t^{n-1}}^{t^n} \skp{\sigma_\infty}{\xi_h}_{L^2(\partial\Omega)} \dt \quad \forall\, \xi_h\in\calS_h,
\end{align}
and, e.g., $\sigma_{\infty,h}^n(P) \coloneqq 0$ in all vertices $P\in\calN_h$ that lie in the interior of $\Omega$. For all $t\in(t^{n-1}, t^n]$, $n\in\{1,...,N_T\}$, we denote the piecewise constant extension by $\sigma_{\infty,h}^{\Delta t,+}(\cdot,t) = \sigma_{\infty,h}^n(\cdot)$.
It follows from \eqref{eq:sigma_infty_constr}, Hölder's and Young's inequalities, that
\begin{align}
    \label{eq:sigma_infty_constr_stability}
    \max_{n=1,...,N_T} \norm{\sigma_{\infty,h}^n}_{L^2(\partial\Omega)} 
    &= \norm{\sigma_{\infty,h}^{\Delta t,+}}_{ L^\infty(0,T; L^2(\partial\Omega))}
    \leq \norm{\sigma_\infty}_{L^\infty(0,T;L^2(\partial\Omega))}.
\end{align}
Using \eqref{eq:sigma_infty_constr}, \eqref{eq:sigma_infty_constr_stability} and a density argument, one can also show that $\sigma_{\infty,h}^{\Delta t,+}|_{\partial\Omega} \to \sigma_\infty$ strongly in $L^2(0,T;L^2(\partial\Omega))$, as $(h,\Delta t)\to (0,0)$.

\end{enumerate}
\end{remark}

The rest of Section \ref{sec:approximation} is devoted to the proof of Theorem \ref{theorem:existence_FE}. After that, in Section \ref{sec:convergence}, we show convergence (up to subsequences) of discrete solutions of the above scheme to a global-in-time weak solution in the sense of Definition \ref{def:weak_solution}.
The outline of the rest of Section \ref{sec:approximation} is as follows.
First, we state some technical preliminaries, including a regularization strategy of the elastic energy density which extends the approach of Barrett and Boyaval \cite{barrett_boyaval_2009}. We also explain the definition of the approximative convective terms in detail. Then, we introduce a regularized discrete scheme for which we show uniform stability estimates and existence of solutions. Passing to the limit $\delta\to 0$ in the regularization parameter (for $h,\Delta t$ fixed) will then prove Theorem \ref{theorem:existence_FE} with the help of converging subsequences.

\subsection{Technical preliminaries and regularizations}
\label{sec:regularizations}
\subsubsection{Some regularizations}

First we need the following identities from Barrett and Boyaval \cite{barrett_boyaval_2009}. After that, we provide extended results which are necessary for our case. Let $\delta>0$.
For any $s\in\bbR$, we define a regularization of the logarithmic function by
\begin{align*}
    g_\delta(s) \coloneqq 
    \begin{cases}
    \frac{s}{\delta} + \ln(\delta) - 1, 
    & s < \delta,
    \\
    \ln(s), 
    & s \geq \delta.
    \end{cases}
\end{align*}
Moreover, we define $\beta_\delta(s) \coloneqq g_\delta'(s)^{-1} = \max\{s,\delta\}$ for any $s\in\bbR$. These functions are extended to matrix valued functions in a natural way via eigenvalues. In particular, let $\bbB \in\bbR^{d\times d}_S$ be a symmetric matrix such that $\bbB = \mathbb{U}^\top \D \mathbb{U}$ with an orthogonal matrix $\mathbb{U}$ and a diagonal matrix $\D = \mathrm{diag}(D_{11}, ..., D_{dd})$. Then, for any scalar valued function $f:\bbR\to\bbR$, we define $f(\bbB)\coloneqq \mathbb{U}^\top f(\D) \mathbb{U}$, where $f(\D) \coloneqq \mathrm{diag}(f(D_{11}), ..., f(D_{dd}))$.


We now prove some identities concerning the function $g_\delta$.

\begin{lemma}~\\
Let $\delta\in(0,1)$ and $\kappa\in\bbR$. Then, for any $\bbB,\bbC\in\bbR^{d\times d}_{\mathrm{S}}$, it holds
\begin{subequations}
\begin{align}
    \label{eq:reg_a}
    \beta_\delta(\bbB) g_\delta'(\bbB) 
    &= g_\delta'(\bbB) \beta_\delta(\bbB) 
    = \bbI, 
    \\
    \label{eq:reg_b}
    (\bbB-\bbC) : \beta_\delta^{-1}(\bbC) 
    &\geq \trace\big( g_\delta(\bbB) - g_\delta(\bbC) \big), 
    \\
    \label{eq:reg_c}
    \frac12 \abs{\bbB}^2 - \trace g_\delta(\bbB) 
    &\geq 
    \frac14 \abs{\bbB}^2 + \frac{1}{\delta} \abs{ [\bbB]_- },
    \\
    \label{eq:reg_d}
    \left(\bbB \beta_\delta(\bbB) + \kappa \beta_\delta(\bbB) - \bbI \right) :
    \left(\bbB + \kappa \bbI - \beta_\delta^{-1}(\bbB) \right) &\geq 0,
    \\
    \label{eq:reg_e}
    \bbB : (\bbB - \beta_\delta^{-1}(\bbB)) &\geq \abs{\bbB}^2 - d,
\end{align}
\end{subequations}
where $\abs{s}_- \coloneqq \min\{s, 0\}$ for any $s\in\bbR$.
\end{lemma}
\begin{proof}
For the identities \eqref{eq:reg_a}--\eqref{eq:reg_b}, we refer to \cite[Lem.~2.1]{barrett_boyaval_2009}. It suffices to show \eqref{eq:reg_c}--\eqref{eq:reg_e}.
Let $\delta\in (0, 1)$ and $\kappa\in\bbR$. 
Then, for any $s,t\in\bbR$, it can be shown with simple arguments that the following scalar identities hold true:
\begin{subequations}
\begin{align}
    \label{eq:reg_c_scalar}
    \frac12 s^2  - g_\delta(s) 
    \geq 
    \frac14 s^2 + \frac{1}{\delta} \abs{\min\{s, 0\}},
    \\
    \label{eq:reg_d_scalar}
    (s \beta_\delta(s) + \kappa \beta_\delta(s) - 1)
    (s + \kappa - \beta_\delta^{-1}(s)) 
    = \underbrace{\beta_\delta(s)}_{=\max\{s, \delta\}} (s + \kappa - \beta_\delta^{-1}(s))^2 \geq 0,
    \\
    \label{eq:reg_e_scalar}
    s (s  - \beta_\delta^{-1}(s)) \geq s^2 - 1.
\end{align}
\end{subequations}
Now let $\bbB\in\bbR^{d\times d}_S$ with $\bbB = \mathbb U^\top \mathbb D \mathbb U$ with $\mathbb U \in \bbR^{d\times d}$ orthogonal and $\mathbb D = \mathrm{diag}(D_{11}, ..., D_{dd})$. Then, on noting \eqref{eq:reg_c_scalar}, it holds
\begin{align*}
    \frac12 \abs{\bbB}^2 - \trace g_\delta(\bbB)
    &= \trace\left( \frac12 (\mathbb U^\top \mathbb D \mathbb U) (\mathbb U^\top\mathbb D \mathbb U)
    - \mathbb U^\top g_\delta(\mathbb D) \mathbb U \right)
    = \trace \left( \frac12 \mathbb D^2 - g_\delta(\mathbb D) \right)
    = \sum_{i=1}^d  \left( \frac12 D_{ii}^2 
    - g_\delta(D_{ii}) \right) 
    \\
    &\geq \sum_{i=1}^d \left(
    \frac14 D_{ii}^2 + \frac{1}{\delta} \abs{\min\{D_{ii},0\}}
    \right)
    \geq \sum_{i=1}^d \frac14 D_{ii}^2  
    + \frac{1}{\delta} \left( \sum_{i=1}^d \abs{\min\{D_{ii},0\}}^2
    \right)^{1/2}
    \\
    &= \frac14 \abs{\bbD}^2 + \frac{1}{\delta} \abs{[\bbD]_-}
    = \frac14 \abs{\bbB}^2 + \frac{1}{\delta} \abs{[\bbB]_-}.
\end{align*}
Moreover, using \eqref{eq:reg_d_scalar}, it holds with similar arguments
\begin{align*}
    &\left(\bbB \beta_\delta(\bbB) + \kappa \beta_\delta(\bbB) - \bbI \right) :
    \left(\bbB + \kappa \bbI - \beta_\delta^{-1}(\bbB) \right) 
    = \trace \Big( 
    ( \bbD \beta_\delta(\bbD)
    + \kappa \beta_\delta(\bbD)
    - \bbI )
    ( \bbD + \kappa \bbI - \beta_\delta^{-1}(\bbD))
    \Big)
    \\
    &= \trace \left( 
    \beta_\delta(\bbD)
    ( \bbD + \kappa \bbI - \beta_\delta^{-1}(\bbD))^2
    \right)
    = \sum_{i=1}^d \underbrace{\beta_\delta(D_{ii})}_{= \max\{ D_{ii}, \delta\} > 0} 
    \underbrace{(D_{ii} + \kappa - \beta_\delta^{-1}(D_{ii}))^2 }_{\geq 0}
    \geq 0.
\end{align*}
On noting \eqref{eq:reg_e_scalar}, we have with similar arguments
\begin{align*}
    \bbB : (\bbB - \beta_\delta^{-1}(\bbB)) 
    &= 
    \trace\left( 
    \bbB^2 - \bbB \beta_\delta^{-1}(\bbB)) 
    \right)
    = \trace\left( 
    \bbD^2 - \bbD \beta_\delta^{-1}(\bbD) 
    \right)
    = \sum_{i=1}^d \left( 
    D_{ii}^2  - D_{ii} \beta_\delta^{-1}(D_{ii})
    \right)
    \\
    &\geq \sum_{i=1}^d \left( 
    D_{ii}^2 - 1
    \right)
    = \abs{\bbD}^2 - d
    = \abs{\bbB}^2 - d.
\end{align*}
This proves the lemma.
\end{proof}

%
%
%
%
%


\subsubsection{Approximation of the convective term}

To motivate the approximation of the convective term, we temporarily assume $\bbB \colon\, \overline\Omega\to\bbR^{d\times d}_{\mathrm{S}}$ to be smooth and positive definite, and $\bv \colon\, \overline\Omega\to\bbR$ to be a smooth velocity field. One key identity that we want to use is
\begin{align*}
     (\bv\cdot\nabla \bbB) : (\bbB - \bbB^{-1}) 
    = \bv\cdot\nabla \big(\frac12 \abs{\bbB}^2 - \trace \ln \bbB \big) ,
\end{align*}
On the finite element level, we would consider $\bbB_h\in\calW_h$ positive definite and $\bv_h\in\calV_h$ and an admissible approximation of the nonlinear testfunction $\bbB - \bbB^{-1}$, i.e., $\bbB_h - \calI_h [\bbB_h^{-1}] \in \calW_h$. However, there is no direct analogue of the chain rule from above, i.e.,
\begin{align*}
     (\bv_h\cdot\nabla \bbB_h) : \big( \bbB_h - \calI_h [\bbB_h^{-1}] \big) 
    \not=  \bv_h\cdot\nabla \calI_h \big[ \frac12 \abs{\bbB_h}^2 - \trace \ln \bbB_h \big] .
\end{align*}
This problem can be overcome using the concept of \textit{discrete chain rules} which has been widely used in the literature, see, e.g., \cite{barrett_boyaval_2009, barrett_2018_fene-p, barrett_nurnberg_2004, barrett_nurnberg_styles_2004, sieber_2020}. 
The basic idea is to approximate $(\bv\cdot \nabla \bbB) : \bbG$ by
\begin{align*}
    (\bv\cdot \nabla \bbB) : \bbG
    &= - \bbB : (\bv\cdot \nabla \bbG)
    + \bv\cdot \nabla \big( \bbB : \bbG \big)
    \\
    &\approx - \sum_{i,j=1}^d (\bv_h)_i \, \mathbf{\Lambda}_{i,j}(\bbB_h) : \partial_{x_j} \bbG_h 
    + \bv_h \cdot \nabla \calI_h \big[ \bbB_h : \bbG_h \big],
\end{align*}
where $\mathbf\Lambda_{i,j}(\bbB_h) \in L^\infty(\Omega;\bbR^{d\times d}_{\mathrm{S}})$,
$i,j\in\{1,...,d\}$, is a nonlinear quantity that will be introduced below, such that, on any simplex $K\in\calT_h$, it approximatively holds $\mathbf\Lambda_{i,j}(\bbB_h)|_{K} \approx \delta_{i,j} \bbB_h |_{K}$, where $\delta_{i,j}$ denotes the Kronecker symbol. 
In particular, $\mathbf\Lambda_{i,j}(\bbB_h)$ will be constructed in a specific way such that it holds
\begin{align*}
    - \sum_{j=1}^d {\mathbf\Lambda} _{i,j}(\bbB_h) : \partial_{x_j} \calI_h[ \bbB_h - \bbB_h^{-1}] 
    = \partial_{x_i} \calI_h\big[ -\tfrac12 \abs{\bbB_h}^2 + \trace \ln(\bbB_h^{-1}) \big] \, ,
    \qquad \text{on any } K\in\calT_h, \ \forall\, i\in\{1,...,d\},
\end{align*}
see below.
This allows the discrete chain rule
\begin{align*}
    - \sum_{i,j=1}^d (\bv_h)_i \, \mathbf{\Lambda}_{i,j}(\bbB_h) : \partial_{x_j} \calI_h[\bbB_h - \bbB_h^{-1}] 
    + \bv_h \cdot \nabla \calI_h \big[ \bbB_h : (\bbB_h - \bbB_h^{-1}) \big]
    = \bv_h\cdot\nabla \calI_h \big[ \frac12 \abs{\bbB_h}^2 - \trace \ln \bbB_h \big].
\end{align*}
In principle, one could think of approximating $(\bv\cdot\nabla\bbB) : \bbG$ directly with $\sum_{i,j=1}^d (\bv_h)_i \, \widetilde{\mathbf{\Lambda}}_{i,j}(\bbG_h) : \partial_{x_j} \bbB_h$, for some $\widetilde{\mathbf{\Lambda}}_{i,j}(\bbG_h) \in L^\infty(\Omega;\bbR^{d\times d}_\mathrm{S})$, $i,j\in\{1,...,d\}$.
However, this would result in a nonlinear dependence on the test function $\bbG_h \in\calW_h$. To avoid this, we use the identity $(\bv\cdot \nabla \bbB) : \bbG
= - \bbB : (\bv\cdot \nabla \bbG)
+ \bv\cdot \nabla \big( \bbB : \bbG \big)$. 
Note that, a priori, we have no information whether $\trace \ln(\bbB_h)$ is well-defined for any arbitrary $\bbB_h\in \calW_h$. For this, we use the regularization strategy of, e.g., \cite{barrett_boyaval_2009}, in order to regularize the logarithmic function as above. 
Then, the goal is to find a formal chain rule for
$- (\bv\cdot\nabla \beta_\delta(\bbB)) :\big( \bbB - \beta_\delta^{-1}(\bbB) \big) $. For that reason, we introduce another regularized function, which approximates a quadratic function, i.e.,
\begin{align*}
    f_\delta(s) = 
    \begin{cases}
        \frac12 s^2, & s\geq \delta,
        \\
        \delta s - \frac12 \delta^2, & s < \delta,
    \end{cases}
\end{align*}
such that $f_\delta'(s) = \beta_\delta(s) = \max\{s, \delta\}$ for any $s\in\bbR$. To motivate a new discrete chain rule, we note the formal identity
\begin{align*}
     (\bv\cdot\nabla \beta_\delta(\bbB)) : (\bbB - \beta_\delta^{-1}(\bbB))
    &= - \bv\cdot\nabla \big(\bbB - \beta_\delta^{-1}(\bbB) \big) : \beta_\delta(\bbB)  
    +  \bv\cdot\nabla \big( \beta_\delta(\bbB) : (\bbB - \beta_\delta^{-1}(\bbB)) \big) 
    \\
    &=  \bv\cdot\nabla \big( - \trace f_\delta(\bbB) + \trace \ln\beta_\delta^{-1}(\bbB) \big) 
    +  \bv\cdot\nabla \big( \beta_\delta(\bbB) : \bbB - d \big) .
\end{align*}

Next, we construct a regularized version ${\mathbf\Lambda} _{i,j}^\delta(\bbB_h) \in L^\infty(\Omega;\bbR^{d\times d}_{\mathrm{S}})$ for $\bbB_h\in\calW_h$, such that ${\mathbf\Lambda} _{i,j}^\delta(\bbB_h)|_{K} \approx \delta_{i,j} \, \beta_\delta(\bbB_h)|_{K}$, $i,j\in\{1,...,d\}$, on any $K\in\calT_h$, and such that a discrete analogue of the identity $\nabla \big(- \bbB + \beta_\delta^{-1}(\bbB) \big) : \beta_\delta(\bbB) = \nabla \big( - \trace f_\delta(\bbB) + \trace \ln\beta_\delta^{-1}(\bbB) \big)$ is fulfilled, i.e.,
\begin{align}
    \label{eq:Lambda_delta_kettenregel}
    \sum\limits_{j=1}^d 
    {\mathbf\Lambda} _{i,j}^\delta (\bbB_h) : \partial_{x_j} \calI_h\big[ - \bbB_h + \beta_\delta^{-1}(\bbB_h) \big]
    = \partial_{x_i} \calI_h\left[ - \trace f_\delta(\bbB_h) + \trace \ln \beta_\delta^{-1}(\bbB_h) \right] \, ,
    \quad \text{ on any } K \in \calT_h, \, \forall\, i\in\{1,...,d\}.
\end{align}

Let $\hat K$ be the standard open reference simplex in $\bbR^d$. Let $\zeta_h\in \calS_h$. We then consider the affine-linear transformation $\calB_K$
from the reference simplex $\hat K$ to any simplex $K\in\calT_h$, see \eqref{eq:trafo_K}. We introduce the notation
\begin{align*}
    \hat\bbB_h(\hat{\mathbf{x}}) \coloneqq \bbB_h(\calB_K(\hat{\mathbf{x}})),
    \quad\quad 
    (\hat\calI_h \hat\bbB_h)(\hat{\mathbf{x}}) \coloneqq (\calI_h \bbB_h)(\calB_K(\hat{\mathbf{x}})),
    \quad\quad \forall \,  \hat{\mathbf{x}}\in\hat K, \, \bbB_h\in\calW_h,
\end{align*}
and we define $\bbB_j^K \coloneqq \bbB_h(P_j^K)$, $j\in\{0,...,d\}$, for $\bbB_h\in\calW_h$ and $K\in\calT_h$, where $P_0^K, ..., P_d^K$ denote the vertices of the simplex $K$. 
Then, we define $\hat{\mathbf\Lambda} _i^\delta(\hat\bbB_h) \in \bbR^{d\times d}_{\mathrm{S}}$, $i\in\{1,...,d\}$, on the reference element $\hat K$ by
\begin{align*}
    \hat{\mathbf\Lambda} _i^\delta(\hat\bbB_h) 
    \coloneqq 
    \begin{cases}
    \beta_\delta(\bbB_i^K) & \text{ if } \beta_\delta(\bbB_i^K) = \beta_\delta(\bbB_0^K),
    \\
    \beta_\delta(\bbB_i^K) + \lambda_i^\delta(\hat\bbB_h)
    \big( \beta_\delta(\bbB_0^K) - \beta_\delta(\bbB_i^K) \big)
    & \text{ if } \beta_\delta(\bbB_i^K) \not= \beta_\delta(\bbB_0^K),
    \end{cases}
\end{align*}
where the scalar quantities $\lambda_i^\delta(\hat\bbB_h) \in \bbR$, $i\in\{1,...,d\}$, are defined as
\begin{align*}
    \lambda_i^\delta(\hat\bbB_h) 
    &\coloneqq
    \frac{ \big( -\trace f_\delta(\bbB_i^K) 
    + \trace \ln \beta_\delta^{-1}(\bbB_i^K)
    +\trace f_\delta(\bbB_0^K) 
    - \trace \ln \beta_\delta^{-1}(\bbB_0^K) \big)}
    { \big(\beta_\delta(\bbB_i^K) - \beta_\delta(\bbB_0^K) \big)  
    : \big( \bbB_i^K - \beta_\delta^{-1}(\bbB_i^K) 
    - \bbB_0^K + \beta_\delta^{-1}(\bbB_0^K) \big)   }
    \\
    &\quad + 
    \frac{\beta_\delta(\bbB_i^K) : 
    \big( \bbB_i^K - \beta_\delta^{-1}(\bbB_i^K) 
    - \bbB_0^K + \beta_\delta^{-1}(\bbB_0^K))}
    { \big(\beta_\delta(\bbB_i^K) - \beta_\delta(\bbB_0^K) \big)  
    : \big( \bbB_i^K - \beta_\delta^{-1}(\bbB_i^K) 
    - \bbB_0^K + \beta_\delta^{-1}(\bbB_0^K) \big) } \,.
\end{align*}
Finally, we define ${\mathbf\Lambda} _{i,j}^\delta(\bbB_h)|_{K} \in \bbR^{d\times d}_{\mathrm{S}}$, $i,j\in\{1,...,d\}$, on any simplex $K\in\calT_h$ by
\begin{align*}
    {\mathbf\Lambda} _{i,j}^\delta(\bbB_h) |_{K}
    \coloneqq
    \sum\limits_{m=1}^d [(\calA_K^\top)^{-1}]_{i,m} \,
    \hat{\mathbf\Lambda} _{m}^\delta(\hat\bbB_h) \, [\calA_K^\top]_{m,j} \in\bbR^{d\times d}_{\mathrm{S}}
    \qquad \text{ on any } K\in\calT_h, \ \forall\, i,j\in\{1,...,d\}.
\end{align*}
Clearly, ${\mathbf\Lambda} _{i,j}^\delta(\bbB_h)|_{K} \in \bbR^{d\times d}_{\mathrm{S}}$, $i,j\in\{1,...,d\}$, depends continuously on $\bbB_h|_{K}$ for any $K\in\calT_h$.

In the next lemma, we show that ${\mathbf\Lambda} _{i,j}^\delta(\bbB_h)$, $i,j\in\{1,...,d\}$, is well-defined. The proof uses the identity
\begin{align}
\label{eq:convex_scalar}
    (b-a) F'(a) \leq F(b)-F(a) \leq (b-a) F'(b), \quad \forall\, a,b\in\bbR,
\end{align}
which holds true for any convex function $F\in C^1(\bbR)$. An analogue identity holds true (cf.~\cite[eq.~(2.15)]{barrett_boyaval_2009}) if $a,b\in\bbR$ are replaced by symmetric matrices $\bbA,\bbB \in\bbR^{d\times d}_{\mathrm{S}}$, i.e.,
\begin{align}
\label{eq:convex_matrix}
    (\bbB-\bbA) : F'(\bbA) 
    \leq \trace F(\bbB) - \trace F(\bbA) 
    \leq (\bbB-\bbA) : F'(\bbB), \quad \forall\,\bbA,\bbB \in\bbR^{d\times d}_{\mathrm{S}}.
\end{align}
If $F(\cdot)$ is strictly convex, then the inequalities in \eqref{eq:convex_scalar}--\eqref{eq:convex_matrix} are strict for $a,b\in\bbR$ with $a\not=b$ and $\bbA,\bbB \in\bbR^{d\times d}_{\mathrm{S}}$ with $\bbA\not=\bbB$, respectively.

\begin{lemma}
Let $\delta \in (0,1)$, $K\in\calT_h$ and $\bbB_h\in \calW_h$ be given. Then, it is $\lambda_i^\delta(\hat\bbB_h) \in [0,1]$, $i\in\{1,...,d\}$, and the discrete chain rule \eqref{eq:Lambda_delta_kettenregel} holds true. In addition, suppose that \ref{S} is satisfied. Then, as the family of triangulations is shape regular, it holds
\begin{align}
    \label{eq:Lambda_delta}
    \max_{i,j\in\{1,...,d\}} 
    \norm{ {\mathbf\Lambda} _{i,j}^\delta(\bbB_h) }_{L^\infty(K)} 
    \leq C \norm{\calI_h\left[ \beta_\delta(\bbB_h) \right]}_{L^\infty(K)}
    \quad \forall\, K\in\calT_h, \, \forall \, \bbB_h\in \calW_h.
\end{align}
\end{lemma}

\begin{proof}
Let $i\in\{1,...,d\}$ be given. We first show that the denominator of $\lambda_i^\delta(\hat\bbB_h)$ is strictly positive for $\beta_\delta(\bbB_i^K) \not= \beta_\delta(\bbB_0^K)$. To check this, we rewrite the denominator as
\begin{align*}
    &\big(\beta_\delta(\bbB_i^K) - \beta_\delta(\bbB_0^K) \big)  
    : \big( \bbB_i^K - \beta_\delta^{-1}(\bbB_i^K) 
    - \bbB_0^K + \beta_\delta^{-1}(\bbB_0^K) \big)
    \\
    &=  \big(\beta_\delta(\bbB_i^K) - \beta_\delta(\bbB_0^K) \big) 
    : \big( \bbB_i^K - \bbB_0^K \big)
    + \big(- \beta_\delta(\bbB_i^K) + \beta_\delta(\bbB_0^K) \big) 
    : \big( \beta_\delta^{-1}(\bbB_i^K) -\beta_\delta^{-1}(\bbB_0^K) \big).
\end{align*}
The first term is non-negative due to \eqref{eq:convex_matrix}, because $f_\delta(\cdot)$ is convex and $f_\delta'(\cdot) = \beta_\delta(\cdot)$. 
Setting $F(\cdot) = -\ln(\cdot)$, $\bbA=\beta_\delta^{-1}(\bbB_0^K)$ and $\bbB=\beta_\delta^{-1}(\bbB_i^K)$ in \eqref{eq:convex_matrix}, we have that the second term is strictly positive, as $-\ln(\cdot)$ is strictly convex. This shows that the denominator of $\lambda_i^\delta(\hat\bbB_h)$ is strictly positive for $\beta_\delta(\bbB_i^K) \not= \beta_\delta(\bbB_0^K)$.
We use this result to prove $\lambda_i^\delta(\hat\bbB_h) \in[0,1]$. In particular, now it suffices to show
\begin{align*}
    \big( -\trace f_\delta(\bbB_i^K) 
    + \trace \ln \beta_\delta^{-1}(\bbB_i^K)
    +\trace f_\delta(\bbB_0^K) 
    - \trace \ln \beta_\delta^{-1}(\bbB_0^K) \big) 
    + \beta_\delta(\bbB_i^K) : 
    \big( \bbB_i^K - \beta_\delta^{-1}(\bbB_i^K) 
    - \bbB_0^K + \beta_\delta^{-1}(\bbB_0^K))
    \geq 0,
\end{align*}
and 
\begin{align*}
    \big( -\trace f_\delta(\bbB_i^K) 
    + \trace \ln \beta_\delta^{-1}(\bbB_i^K)
    +\trace f_\delta(\bbB_0^K) 
    - \trace \ln \beta_\delta^{-1}(\bbB_0^K) \big) 
    + \beta_\delta(\bbB_0^K) : 
    \big( \bbB_i^K - \beta_\delta^{-1}(\bbB_i^K) 
    - \bbB_0^K + \beta_\delta^{-1}(\bbB_0^K))
    \leq 0.
\end{align*}
Again, both inequalities directly follow from \eqref{eq:convex_matrix} with $F(\cdot)=f_\delta(\cdot)$, $\bbA=\bbB_0^K$, $\bbB=\bbB_i^K$ and $F(\cdot)=-\ln(\cdot)$, $\bbA=\beta_\delta^{-1}(\bbB_0^K)$, $\bbB=\beta_\delta^{-1}(\bbB_i^K)$, respectively, as $f_\delta(\cdot)$ and $-\ln(\cdot)$ are convex and strictly convex, respectively.

Suppose that \ref{S} is fulfilled. Then, we obtain \eqref{eq:Lambda_delta} with \eqref{eq:shape_regular} and the definition of ${\mathbf\Lambda} _{i,j}^\delta(\bbB_h)|_{K}$, as $\lambda_i^\delta(\hat\bbB_h) \in [0,1]$.
\end{proof}

\begin{remark}[$\delta\to 0$]
Later, we will send the regularization parameter $\delta>0$ to zero. For that, we define the unregularized version ${\mathbf\Lambda} _{i,j}(\bbB_h)$ for $\bbB_h\in\calW_h$ with $\bbB_h$ positive definite similarly to ${\mathbf\Lambda} _{i,j}^\delta(\bbB_h)$, with $\hat{\mathbf\Lambda} _i^\delta(\hat\bbB_h)$ replaced by $\hat{\mathbf\Lambda} _i(\hat\bbB_h)$, which is defined similarly with $\lambda_i^\delta(\hat\bbB_h)$, $\beta_\delta$ replaced by $\lambda_i(\hat\bbB_h)$ and the identity function $\bbR\to\bbR, \, s\mapsto s$, where $\lambda_i(\hat\bbB_h)$ is defined similarly to $\lambda_i^\delta(\hat\bbB_h)$ with $f_\delta$, $\beta_\delta$ replaced by a quadratic function and the identity function, i.e., $\bbR\to\bbR, \, s\mapsto \frac12 s^2$ and $\bbR\to\bbR, \, s\mapsto s$, respectively.

Here, we have the analogue properties to the regularized case. In particular, let $K\in\calT_h$ and $\bbB_h\in\calW_h$ with $\bbB_h$ positive definite. Then, it is $\lambda_i(\hat\bbB_h)\in[0,1]$, $i\in\{1,...,d\}$, and the discrete chain rule
\begin{align*}
    \sum\limits_{j=1}^d 
    {\mathbf\Lambda} _{i,j} (\bbB_h) : \partial_{x_j} \calI_h\big[ - \bbB_h + \bbB_h^{-1} \big]
    = \partial_{x_i} \calI_h\left[ - \frac12 \abs{\bbB_h}^2 + \trace \ln (\bbB_h^{-1}) \right] \, ,
    \quad \text{ on any } K \in \calT_h, \, \forall\, i\in\{1,...,d\},
\end{align*}
holds true.
Moreover, if \ref{S} is satisfied, i.e., if $\{\calT_h\}_{h>0}$ is shape regular, then it holds
\begin{align}
    \label{eq:Lambda}
    \max_{i,j\in\{1,...,d\}} 
    \norm{ {\mathbf\Lambda} _{i,j}(\bbB_h) }_{L^\infty(K)} 
    \leq C \norm{\bbB_h}_{L^\infty(K)}
    \quad \forall\, K\in\calT_h, \, \forall \, \bbB_h\in \calW_h \text{ with } \bbB_h \text{ positive definite}.
\end{align}
\end{remark}

We specify the approximation property of ${\mathbf\Lambda} _{i,j}(\bbB_h)$ for the case without the regularization. 
\begin{lemma}
Let \ref{S} hold true. For any $p\in[1,\infty]$, $K\in\calT_h$ and all $\bbB_h\in\calW_h$ with $\bbB_h$ positive definite, it holds
\begin{align}
\label{eq:error_Lambda}
    \max_{i,j\in\{1,...,d\}} \norm{{\mathbf\Lambda} _{i,j}(\bbB_h) - \delta_{i,j} \, \bbB_h }_{L^p(K)}
    \leq C h_K \norm{\nabla\bbB_h}_{L^p(K)}.
\end{align}
\end{lemma}
\begin{proof}
The proof is based on a simple calculation with Hölder's inequality, \eqref{eq:shape_regular}, \eqref{eq:inverse_estimate} and the fact that $\lambda_i(\hat\bbB_h)\in[0,1]$, $i\in\{1,...,d\}$. In particular, for any $i,j\in\{1,...,d\}$, it holds
\begin{align*}
    \norm{{\mathbf\Lambda} _{i,j}(\bbB_h) - \delta_{i,j} \, \bbB_h }_{L^p(K)}
    &= \nnorm{ \sum\limits_{m=1}^d [(\calA_K^\top)^{-1}]_{i,m} \,
    \big(\hat{\mathbf\Lambda}_{m}(\hat\bbB_h) - \bbB_h\big) \, [\calA_K^\top]_{m,j} }_{L^p(K)}
    \\
    &\leq C \sum\limits_{m=1}^d  
    \nnorm{ \hat{\mathbf\Lambda}_{m}(\hat\bbB_h) - \bbB_h }_{L^p(K)}
    \\
    &\leq C \abs{K}^{1/p} \max_{\ell,m\in\{0,...,d\}} \abs{\bbB_h(P_\ell^K) - \bbB_h(P_m^K)}
    \\
    &\leq C h_K \norm{\nabla\bbB_h}_{L^p(K)},
\end{align*}
where $P_\ell^K$, $\ell\in\{0,...,d\}$, denote the vertices on the simplex $K\in\calT_h$.
\end{proof}

\subsubsection{Another discrete chain rule}
Another key identity, that we want to use, is the discrete analogue of
\begin{align*}
    - \int_\Omega \nabla \bbB : \nabla \bbB^{-1} \dx
    \geq \int_\Omega \frac1d \abs{\nabla\trace\ln\bbB}^2 \dx.
\end{align*}
Here, we improve the results from the literature, where either a too weak estimate was shown \cite{barrett_boyaval_2009} or very strong assumptions on the triangulation were required \cite{sieber_2020}.

\begin{lemma}~\\
Let $\delta\in(0,1)$ and let \ref{S} hold true, in particular, let all simplices $K\in\calT_h$ be non-obtuse. Then, for all $\bbB_h\in\calW_h$ and all $K\in\calT_h$, it holds
\begin{align}
\label{eq:nabla_ln}
    - \skp{\nabla \bbB_h}{\nabla\calI_h \beta_\delta^{-1}(\bbB_h)}_{L^2(K)} 
    \geq 
    \frac1d \norm{\nabla\calI_h \trace\ln \beta_\delta(\bbB_h)}_{L^2(K)}^2
    = \frac1d  \norm{\nabla\calI_h \trace\ln \beta_\delta^{-1}(\bbB_h)}_{L^2(K)}^2.
\end{align}
\end{lemma}

\begin{proof}
Let $K\in\calT_h$ be given with local vertices $P_0^K, ..., P_d^K$ and nodal basis functions $\eta_0^K, ..., \eta_d^K \in\calS_h$ on $K$. As $K$ is non-obtuse, we can use \eqref{eq:non-obtuse}.
Now let $a_i, b_i \in\bbR$ and $\pmb a_i, \pmb b_i \in \bbR^{d\times d}_{\mathrm{S}}$, $i\in\{0,...,d\}$, be some scalar and matrix valued coefficients, respectively. Then, with a simple calculation (cf.~\cite[Lem.~5.1]{barrett_boyaval_2009}) which consists of $\sum_{j=0}^d \eta_j^K=1$ on $K$, taking the gradient on both sides and rearranging some terms, it holds 
\begin{subequations}
\begin{align}
    \label{eq:bb11_lem5.1a}
    \nabla \left( \sum_{i=0}^d a_i \eta_i^K \right)
    \cdot  \nabla \left( \sum_{j=0}^d b_j \eta_j^K \right)
    &=
    - \sum_{i=0}^d \sum_{j>i}^d \left[ 
    ( a_i - a_j) (b_i - b_j) \right]
    \nabla \eta_i^K \cdot \nabla \eta_j^K,
    \\
    \label{eq:bb11_lem5.1b}
    \nabla \left( \sum_{i=0}^d \pmb a_i \eta_i^K \right)
    : \nabla \left( \sum_{j=0}^d \pmb b_j \eta_j^K \right)
    &=
    - \sum_{i=0}^d \sum_{j>i}^d \left[ 
    ( \pmb a_i - \pmb a_j) : (\pmb b_i - \pmb b_j) \right]
    \nabla \eta_i^K \cdot \nabla \eta_j^K.
\end{align}
\end{subequations}
Now, we require the following inequality from \cite[Appendix C, (C.4)]{barrett_lu_sueli_2017},
\begin{align}
\label{eq:bls17_(c.4)}
    - ( \bbB - \bbC) : (\beta_\delta^{-1}(\bbB) - \beta_\delta^{-1}(\bbC))
    \geq \frac1d \abs{ \trace\ln \beta_\delta(\bbB) - \trace\ln \beta_\delta(\bbC) }^2 
    \quad \forall\, \bbB,\bbC \in \bbR^{d\times d}_{\mathrm{S}},
\end{align}
which holds true for any $\delta\in(0,1)$. Hence, using \eqref{eq:bb11_lem5.1b} with $\pmb a_i \coloneqq \bbB_h(P_i^K)$, $\pmb b_i \coloneqq - \beta_\delta^{-1}(\bbB_h(P_i^K))$ for all $i\in\{0,...,d\}$ and applying \eqref{eq:non-obtuse} and \eqref{eq:bls17_(c.4)}, we have locally on $K\in\calT_h$,
\begin{align*}
    - \nabla \bbB_h : \nabla\calI_h \beta_\delta^{-1}(\bbB_h)
    &= \nabla \left( \sum_{i=0}^d \pmb a_i \eta_i^K \right)
    : \nabla \left( \sum_{j=0}^d \pmb b_j \eta_j^K \right)
    \\
    &= - \sum_{i=0}^d \sum_{j>i}^d \left[ 
    ( \pmb a_i - \pmb a_j) : (\pmb b_i - \pmb b_j) \right]
    \nabla \eta_i^K \cdot \nabla \eta_j^K
    \\
    &\geq 
    - \frac1d \sum_{i=0}^d \sum_{j>i}^d 
    \abs{ \trace\ln \beta_\delta(\bbB_h(P_i^K)) - \trace\ln \beta_\delta(\bbB_h(P_j^K)) }^2 \nabla \eta_i^K \cdot \nabla \eta_j^K .
\end{align*}
Next, using \eqref{eq:bb11_lem5.1a} with $a_i \coloneqq b_i \coloneqq \trace\ln \beta_\delta(\bbB_h(P_i^K))$ for all $i\in\{0,...,d\}$, it holds locally on $K\in\calT_h$,
\begin{align*}
    - \nabla \bbB_h : \nabla\calI_h \beta_\delta^{-1}(\bbB_h)
    &\geq 
    - \frac1d \sum_{i=0}^d \sum_{j>i}^d 
    \abs{ \trace\ln \beta_\delta(\bbB_h(P_i^K)) - \trace\ln \beta_\delta(\bbB_h(P_j^K)) }^2 \nabla \eta_i^K \cdot \nabla \eta_j^K
    \\
    &= 
    - \frac1d \sum_{i=0}^d \sum_{j>i}^d \left[ 
    ( a_i - a_j) (b_i - b_j) \right]
    \nabla \eta_i^K \cdot \nabla \eta_j^K
    \\
    &= \frac1d \nabla \left( \sum_{i=0}^d a_i \eta_i^K \right)
    \cdot \nabla \left( \sum_{j=0}^d b_j \eta_j^K \right)
    \\
    &= \frac1d \abs{ \nabla \calI_h \trace\ln \beta_\delta(\bbB_h)}^2. 
\end{align*}
The claim follows by taking the integral over $K\in\calT_h$ and noting $\ln(s^{-1}) = - \ln(s)$ $\forall\, s>0$.
\end{proof}

\subsection{Proof of the existence and stability theorem}

We prove the existence result in Theorem \ref{theorem:existence_FE} with (subsequence) convergence and a limit passing in a regularized discrete scheme with cut-offs in certain terms. The difference to the unregularized scheme \ref{P_FE} is that the discrete Cauchy--Green tensor $\bbB_h^n$ does not necessarily have to be positive definite in presence of cut-offs. The key ideas are a uniform estimate in the regularization parameter $\delta$ and the limit passing in $\delta \to 0$, which retrieves the positive definiteness. Throughout this section, we assume that \ref{A1}--\ref{A6} and \ref{S} hold true.

\subsubsection{Regularized discrete scheme with cut-offs}

The $\delta$-regularized finite element scheme reads as follows.

\subsubsection*{Problem \ref{P_delta_FE}:}
\mylabelHIDE{P_delta_FE}{$(\pmbP_{\delta,h}^{\Delta t})$} 
Let $\delta\in (0, 1)$.  
Let $n\in\{1,...,N_T\}$ and suppose that $\phi_h^{n-1} \in \calS_h$, $\sigma_{\infty,h}^n\in\calS_h$ and $\bbB_h^{n-1} \in\calW_h$ 
are given. Then, the goal is to find a solution tuple 
\begin{align*}
    (\phi_h^{n}, \mu_h^n, \sigma_{h}^{n}, p_h^n, \bv_{h}^{n}, \bbB_{h}^{n}) \in (\calS_h)^4 \times \calV_h  \times \calW_{h},
\end{align*}
which satisfies for any test function tuple $(\zeta_h, \rho_h, \xi_h, q_h, \bw_h, \bbG_h) \in (\calS_h)^4 \times \calV_h \times \calW_h$:
\begin{subequations}
\begin{align}
    \label{eq:phi_reg}
    0 &= 
    \frac{1}{\Delta t} \skp{\phi_h^n - \phi_h^{n-1}}{\zeta_h}_h
    + \skp{\calI_h [m_{\phi,h}^{n-1}] \nabla\mu_h^n}{\nabla\zeta_h}_{L^2}
    + \skp{\bv_h^n \cdot \nabla\phi_h^{n-1}}{\zeta_h}_{L^2}
    + \skp{\phi_h^{n-1} \Gamma_{\bv,h}^{n-1} - \Gamma_{\phi,h}^{n-1}}{\zeta_h}_h,
    \\[5pt]
    \label{eq:mu_reg}
    0 &= 
    \skp{- \mu_h^n + \frac\beta\varepsilon \psi_h'(\phi_h^n,\phi_h^{n-1}) - \chi_\phi \sigma_h^n + \frac12 \kappa_h'(\phi_h^n, \phi_h^{n-1}) \trace\bbB_h^n} {\rho_h}_h
    + \beta\varepsilon \skp{\nabla\phi_h^n}{\nabla\rho_h}_{L^2},
    \\[5pt]
    \label{eq:sig_reg}
    0 &= \skp{\calI_h[m_{\sigma,h}^{n-1}] \nabla\sigma_h^n}{\nabla\xi_h}_{L^2}
    + \skp{\sigma_h^n \Gamma_{\sigma,h}^{n-1}}{\xi_h}_h
    + K \skp{\sigma_h^n - \sigma_{\infty,h}^n}{\xi_h}_{L^2(\partial\Omega)},
    \\[5pt]
    \label{eq:div_reg}
    0 &= \skp{\Div\bv_h^n}{q_h}_{L^2} - \skp{\Gamma_{\bv,h}^{n-1}}{q_h}_{h},
    \\[5pt]
    \nonumber
    \label{eq:v_reg}
    0 &= \skp{2 \calI_h[\eta(\phi_h^{n-1})] \D(\bv_h^n)}{\D(\bw_h)}_{L^2}
    - \skp{p_h^n}{\Div\bw_h}_{L^2}
    + \skp{\calI_h \bbT_{e,h}^{\delta,n}}{\nabla\bw_h}_{L^2}
    - \skp{(\mu_h^n + \chi_\phi \sigma_h^n) \nabla\phi_h^{n-1}}{\bw_h}_{L^2}
    \\
    &\quad
    - \frac12 \skp{\nabla \calI_h\left[ \kappa(\phi_h^{n-1}) \trace\beta_\delta(\bbB_h^n) \right]}{\bw_h}_{L^2}
    + \frac12 \sum_{i,j=1}^d \skp{\partial_{x_j} \calI_h\left[ \kappa(\phi_h^{n-1}) \right]  \trace \mathbf{\Lambda}^\delta_{i,j}(\bbB_h^n)}{(\bw_h)_i}_{L^2},
    \\[5pt]
    \label{eq:B_reg}
    \nonumber
    0 &= 
    \frac{1}{\Delta t} \skp{\bbB_h^n - \bbB_h^{n-1}}{\bbG_h}_h
    + \skp{\bv_h^n}{ \nabla\calI_h\left[\beta_\delta(\bbB_h^n) : \bbG_h \right]}_{L^2}
    - \sum_{i,j=1}^d \skp{(\bv_h^n)_i \, \mathbf{\Lambda}^\delta_{i,j}(\bbB_h^n)}{\partial_{x_j} \bbG_h}_{L^2}
    \\
    &\quad 
     + \skp{\tfrac{1}{\tau(\phi_h^{n-1})} \bbT_{e,h}^{\delta,n} }{\bbG_h}_h 
    + \skp{\Gamma_{\bbB,h}^{n-1} \beta_\delta(\bbB_h^n)}{\bbG_h}_h
    - \skp{2 \nabla\bv_h^n}{\calI_h\left[ \bbG_h \beta_\delta(\bbB_h^n) \right]}_{L^2}
    + \alpha \skp{\nabla\bbB_h^n}{\nabla\bbG_h}_{L^2},
\end{align}
\end{subequations}
where, for the regularized elastic stress tensor, we use the notation
\begin{align*}
    \bbT_{e,h}^{\delta,n} 
    &\coloneqq \bbB_h^n \beta_\delta(\bbB_h^n) + \kappa(\phi_h^{n-1}) \beta_\delta(\bbB_h^n) - \bbI  .
\end{align*}
For the source terms and the mobility functions, we adopt the notation from \ref{P_FE}.


For future reference, we note some estimates that will be frequently used in the next sections. These estimates follow with the assumptions \ref{A1}--\ref{A6} and \ref{S} and standard techniques for the interpolation operator $\calI_h$ that are based on \eqref{eq:interp_H2}, \eqref{eq:interp_stability}, \eqref{eq:norm_equiv}, \eqref{eq:inverse_estimate} and Hölder's inequality.
\begin{subequations}
\begin{align}
    \label{eq:bound_psi_h'}
    \norm{\calI_h \psi_h'(\phi_h, \tilde\phi_h)}_{L^q} 
    &\leq C ( 1 + \norm{\phi_h}_{L^q} + \norm{\tilde\phi_h}_{L^q}),
    \\
    \label{eq:bound_k_h}
    \norm{\calI_h \kappa(\phi_h)}_{L^q} 
    &\leq C ( 1 + \norm{\phi_h}_{L^q}),
    \\
    \label{eq:bound_grad_k_h}
    \norm{\nabla \calI_h \kappa(\phi_h)}_{L^2} 
    &\leq C \norm{\nabla \phi_h}_{L^2},
    \\
    \label{eq:bound_kh'_tau}
    \norm{\calI_h \kappa_h'(\phi_h, \tilde\phi_h)}_{L^q}
    + \norm{\calI_h \tfrac{1}{\tau(\phi_h)}}_{L^q}
    &\leq C,
    \\
    \nonumber \label{eq:bound_source}
    \norm{\calI_h \Gamma_\phi(\phi_h,\sigma_h,\bbB_h)}_{L^q}
    + \norm{\calI_h \Gamma_{\bv}(\phi_h,\sigma_h,\bbB_h)}_{L^q} \qquad
    & 
    \\ + \norm{\calI_h \Gamma_{\bbB}(\phi_h,\sigma_h,\bbB_h)}_{L^q}
    &\leq C (1 + \norm{\sigma_h}_{L^q}),
    \\
    \label{eq:bound_beta_delta}
    \norm{\calI_h \beta_\delta(\bbB_h)}_{L^q} 
    &\leq C (1 + \norm{\bbB_h}_{L^q}),
    \\
    \label{eq:bound_grad_beta_delta}
    \norm{\nabla \calI_h \beta_\delta(\bbB_h)}_{L^2} 
    &\leq C \norm{\nabla \bbB_h}_{L^2},
\end{align}
for any $q\in[1,\infty]$, $\phi_h,\tilde\phi_h, \sigma_h \in\calS_h$ and $\bbB_h \in\calW_h$, where, unless otherwise stated, $C>0$ always denotes a generic constant which is independent of $\delta, h, \Delta t$
\end{subequations}

Besides, we now state some well-known results that will be needed several times. Let $r\in[2,\infty)$ for $d=2$ and $r\in[2,6]$ for $d=3$, respectively. Then, the embedding $H^1(\Omega)\hookrightarrow L^r(\Omega)$ is continuous. Moreover, the following Gagliardo--Nirenberg interpolation inequality holds for any $f\in H^1(\Omega)$,
\begin{align}
\label{eq:Gagliardo}
    \norm{f}_{L^r(\Omega)} 
    \leq C \norm{f}_{L^2(\Omega)}^{1-\theta} 
    \norm{f}_{H^1(\Omega)}^\theta,
\end{align}
where $C=C(d,r,\Omega)$ and $\theta = \frac d2 \frac{(r-2)}{r}$, 
see, e.g., \cite[eq.~(II.3.19)]{galdi_2011}.
The following analogues for boundary integrals are a direct consequence of a variant of the trace theorem which can be found in \cite[Thm.~II.4.1]{galdi_2011}. As the bounded domain $\Omega\subset\bbR^d$ with $d\in\{2,3\}$ is locally Lipschitz (due to \ref{S}), the embedding $H^1(\Omega) \hookrightarrow L^s(\partial\Omega)$ is continuous with $s\in[2,\infty)$ for $d=2$ and $s\in[2,4]$ for $d=3$, respectively. Moreover, the following interpolation inequality is valid for any $f \in H^1(\Omega)$,
\begin{align}
\label{eq:Gagliardo_boundary}
    \norm{f}_{L^s(\partial\Omega)}
    \leq C \left( \norm{f}_{L^2(\Omega)}^{1-\lambda}
    \norm{f}_{H^1(\Omega)}^\lambda  
    + \norm{f}_{L^2(\Omega)}^{(1-\frac1s)(1-\lambda)} 
    \norm{f}_{H^1(\Omega)}^{\frac1s + \lambda (1-\frac1s) }
    \right),
\end{align}
where $C=C(d,s,\Omega)$ and $\lambda = \frac d2 \frac{(s-2) }{(s-1)}$. For example, on noting \eqref{eq:Gagliardo} and \eqref{eq:Gagliardo_boundary} with dimension $d=3$, $r=3$, $\theta=\frac12$, $s=2$ and $\lambda=0$, respectively, it holds for any $f\in H^1(\Omega)$,
\begin{align*}
    \norm{f}_{L^3(\Omega)} + \norm{f}_{L^2(\partial\Omega)}
    &\leq C 
    \norm{f}_{L^2(\Omega)}^{1/2} \norm{f}_{H^1(\Omega)}^{1/2}.
\end{align*}

Later, we will need the following results.
It holds with Hölder's inequality, \eqref{eq:Gagliardo}, \eqref{eq:bound_k_h}, \eqref{eq:bound_grad_k_h}, \eqref{eq:bound_beta_delta}, \eqref{eq:bound_grad_beta_delta} and standard techniques for $\calI_h$, that
\begin{align}
    \label{eq:bound_Te_delta} \nonumber
    \norm{\calI_h \bbT_{e,h}^{\delta,n}}_{L^2}
    &\leq C \left( \norm{\calI_h \beta_\delta(\bbB_h^n)}_{L^3}
    \norm{\bbB_h^n}_{L^6}
    + \norm{\calI_h \beta_\delta(\bbB_h^n)}_{L^3}
    \norm{\calI_h \kappa(\phi_h^{n-1})}_{L^6}
    + 1 \right)
    \\
    &\leq C \left( 1 + \left( 1 + \norm{\bbB_h^n}_{L^2}^{1/2} \norm{\bbB_h^n}_{H^1}^{1/2} \right)
    \left( \norm{\bbB_h^n}_{H^1} 
    + \norm{\phi_h^{n-1}}_{H^1} + 1 \right) \right).
\end{align}
Similarly, it holds
\begin{align}
    \label{eq:bound_Te} \nonumber
    \norm{\calI_h \bbT_{e,h}^{n}}_{L^2}
    &\leq C \left( \norm{\bbB_h^n}_{L^3}
    \norm{\bbB_h^n}_{L^6}
    + \norm{\bbB_h^n}_{L^3}
    \norm{\calI_h \kappa(\phi_h^{n-1})}_{L^6}
    + 1 \right)
    \\
    &\leq C \left( 1 + \left( \norm{\bbB_h^n}_{L^2}^{1/2} \norm{\bbB_h^n}_{H^1}^{1/2} \right)
    \left( \norm{\bbB_h^n}_{H^1} 
    + \norm{\phi_h^{n-1}}_{H^1} + 1 \right) \right).
\end{align}

Let $\bv_h\in\calV_h$ and $\bbB_h, \bbG_h \in \calW_h$. It holds with the Gauß theorem, that
\begin{align*}
    \skp{\bv_h}{\nabla \calI_h[ \bbB_h : \bbG_h]}_{L^2}
    = - \skp{\divergenz{\bv_h}}{\calI_h\left[ \bbB_h:\bbG_h \right]}_{L^2}
    + \skp{\bv_h\cdot\pmbn}{\calI_h\left[\bbB_h : \bbG_h \right]}_{L^2(\partial\Omega)}.
\end{align*}
From this identity, we get with Hölder's inequality, \eqref{eq:inverse_estimate} and \eqref{eq:Gagliardo}--\eqref{eq:Gagliardo_boundary}, that
%
%
%
%
\begin{align}
    \label{eq:askp} \nonumber
    \abs{ \skp{\bv_h}{\nabla \calI_h[ \bbB_h : \bbG_h]}_{L^2}} 
    &\leq C 
    \norm{\bv_h}_{H^1} \norm{\bbB_h}_{L^3} \norm{\bbG_h}_{L^6}
    + C \norm{\bv_h}_{L^4(\partial\Omega)}
    \norm{\bbB_h}_{L^2(\partial\Omega)}
    \norm{\bbG_h}_{L^4(\partial\Omega)}
    \\
    &\leq
    C \norm{\bv_h}_{H^1} 
    \norm{\bbB_h}_{L^2}^{1/2} \norm{\bbB_h}_{H^1}^{1/2}
    \norm{\bbG_h}_{H^1}.
\end{align}

In the following, we show existence of discrete solutions for the regularized discrete problem \ref{P_delta_FE} and we derive stability estimates that are uniform in $\delta>0$. 

\subsubsection{Existence of nutrients}
First, we note that \eqref{eq:sig_reg} can be solved separately, as it is decoupled from the system \eqref{eq:phi_reg}--\eqref{eq:B_reg}. The existence of a unique solution $\sigma_h^n \in\calS_h$ of \eqref{eq:sig_reg} follows from, e.g., the Lax--Milgram theorem. Moreover, we have the following identity, which follows from setting $\xi_h=\sigma_h^n$ in \eqref{eq:sig_reg}:
\begin{align*}
    \skp{\calI_h m_{\sigma,h}^{n-1}} {\abs{\nabla\sigma_h^n}^2}_{L^2} + \skp{ \Gamma_{\sigma,h}^{n-1}}{\abs{\sigma_h^n}^2}_h
    + K \norm{\sigma_h^n}_{L^2(\partial\Omega)}^2 
    = K \skp{\sigma_{\infty,h}^n}{\sigma_h^n}_{L^2(\partial\Omega)}.
\end{align*}
Due to \ref{A1}--\ref{A2}, we obtain with Hölder's, Young's and Poincaré's inequalities, 
\begin{align}
    \label{eq:bound_sigma}
    \norm{\sigma_h^n}_{H^1}^2 \leq C \norm{\sigma_{\infty,h}^n}_{L^2(\partial\Omega)}^2,
\end{align}
for a constant $C>0$ that is independent of $\delta, h, \Delta t$.


\subsubsection{Solving the divergence equation}
If we try to test \eqref{eq:div_reg}--\eqref{eq:v_reg} with $q_h = p_h^n$ and $\bw_h = \bv_h^n$, respectively, we would end up with problems. In particular, we would have to deal with the term $\skp{p_h^n}{\Div \bv_h^n}_{L^2} = \skp{p_h^n}{\Gamma_{\bv,h}^{n-1}}_h$, for which we have no control as there isn't any \textit{a priori} information about the pressure $p_h^n\in\calS_h$.
Therefore, we decompose the velocity field $\bv_h^n = \bz_h^n + \bu_h^n$ into a ``divergence-free'' part $\bz_h^n \in \calV_{h,\mathrm{div}}$ and a non-homogeneous part $\bu_h^n\in\calV_h$ which solves \eqref{eq:div_reg}. In particular, the strategy is as follows.
It follows from \eqref{eq:LBB_u}--\eqref{eq:LBB_u_bound} that there exists a velocity field $\bu_h \in \calV_h$ solving \eqref{eq:div_reg} and which satisfies the stability estimate
\begin{align}
    \label{eq:bound_uh}
    \norm{\bu_h^n}_{W^{1,r}} 
    \leq C \norm{\calI_h \Gamma_{\bv,h}^{n-1}}_{L^r}
    \leq C (1 + \norm{\sigma_h^n}_{L^r})
    \leq C (1 + \norm{\sigma_h^n}_{H^1})
    \leq C (1 + \norm{\sigma_{\infty,h}^n}_{L^2(\partial\Omega)} ),
\end{align}
where $r\in(1,\infty)$ for $d=2$, and $r\in(1,6]$ for $d=3$, such that the Sobolev embedding $H^1(\Omega)\hookrightarrow L^r(\Omega)$ is valid.
For future reference, we note that choosing $r \in (d,\frac{2d}{d-2})$ and using the Sobolev embedding $W^{1,r}(\Omega) \hookrightarrow L^\infty(\Omega)$ leads to $\norm{\bu_h^n}_{L^\infty} \leq C (1 + \norm{\sigma_{\infty,h}^n}_{L^2(\partial\Omega)} )$.

Now we consider a reduced problem, where we use the ansatz $\bz_h^n \coloneqq \bv_h^n - \bu_h^n \in \calV_{h,\mathrm{div}}$ and restrict to weakly solenoidal test functions $\bw_h \in \calV_{h,\mathrm{div}}$ in \eqref{eq:v_reg}. This allows to temporarily forget about the pressure $p_h^n\in\calS_h$, which we will reconstruct afterwards using \eqref{eq:LBB_p}--\eqref{eq:LBB_p_bound}. The reduced problem reads as follows.

The goal is to find a solution tuple 
\begin{align*}
    (\phi_h^{n}, \mu_h^n, \bz_{h}^{n}, \bbB_{h}^{n}) \in (\calS_h)^2 \times \calV_{h,\mathrm{div}}  \times \calW_{h},
\end{align*}
which satisfies for any test function tuple $(\zeta_h, \rho_h, \bw_h, \bbG_h) \in (\calS_h)^2 \times \calV_{h,\mathrm{div}} \times \calW_h$:
\begin{subequations}
\begin{align}
    \nonumber \label{eq:phi_reduced}
    0 &= 
    \frac{1}{\Delta t} \skp{\phi_h^n - \phi_h^{n-1}}{\zeta_h}_h
    + \skp{\calI_h [m_{\phi,h}^{n-1}] \nabla\mu_h^n}{\nabla\zeta_h}_{L^2}
    + \skp{(\bz_h^n + \bu_h^n) \cdot \nabla\phi_h^{n-1}}{\zeta_h}_{L^2}
    \\
    &\quad 
    + \skp{\phi_h^{n-1} \Gamma_{\bv,h}^{n-1} - \Gamma_{\phi,h}^{n-1}}{\zeta_h}_h,
    \\[5pt]
    \label{eq:mu_reduced}
    0 &= 
    \skp{- \mu_h^n + \frac\beta\varepsilon \psi_h'(\phi_h^n,\phi_h^{n-1}) - \chi_\phi \sigma_h^n + \frac12 \kappa_h'(\phi_h^n, \phi_h^{n-1}) \trace\bbB_h^n} {\rho_h}_h
    + \beta\varepsilon \skp{\nabla\phi_h^n}{\nabla\rho_h}_{L^2},
    \\[5pt]
    \nonumber
    \label{eq:v_reduced}
    0 &= \skp{2 \calI_h[\eta(\phi_h^{n-1})] \D(\bz_h^n)}{\D(\bw_h)}_{L^2}
    + \skp{2 \calI_h[\eta(\phi_h^{n-1})] \D(\bu_h^n)}{\D(\bw_h)}_{L^2}
    \\
    \nonumber
    &\quad
    + \skp{\calI_h \bbT_{e,h}^{\delta,n}}{\nabla\bw_h}_{L^2}
    - \skp{(\mu_h^n + \chi_\phi \sigma_h^n) \nabla\phi_h^{n-1}}{\bw_h}_{L^2}
    \\
    &\quad
    - \frac12 \skp{\nabla \calI_h\left[ \kappa(\phi_h^{n-1}) \trace\beta_\delta(\bbB_h^n) \right]}{\bw_h}_{L^2}
    + \frac12 \sum_{i,j=1}^d \skp{\partial_{x_j} \calI_h\left[ \kappa(\phi_h^{n-1}) \right]  \trace \mathbf{\Lambda}^\delta_{i,j}(\bbB_h^n)}{(\bw_h)_i}_{L^2},
    \\[5pt]
    \label{eq:B_reduced}
    \nonumber
    0 &= 
    \frac{1}{\Delta t} \skp{\bbB_h^n - \bbB_h^{n-1}}{\bbG_h}_h
    + \skp{\bz_h^n + \bu_h^n}{ \nabla\calI_h\left[\beta_\delta(\bbB_h^n) : \bbG_h \right]}_{L^2}
    - \sum_{i,j=1}^d \skp{(\bz_h^n + \bu_h^n)_i \, \mathbf{\Lambda}^\delta_{i,j}(\bbB_h^n)}{\partial_{x_j} \bbG_h}_{L^2}
    \\
    &\quad 
     + \skp{\tfrac{1}{\tau(\phi_h^{n-1})} \bbT_{e,h}^{\delta,n} }{\bbG_h}_h 
    + \skp{\Gamma_{\bbB,h}^{n-1} \beta_\delta(\bbB_h^n)}{\bbG_h}_h
    - \skp{2 \nabla(\bz_h^n + \bu_h^n)}{\calI_h\left[ \bbG_h \beta_\delta(\bbB_h^n) \right]}_{L^2}
    + \alpha \skp{\nabla\bbB_h^n}{\nabla\bbG_h}_{L^2}.
\end{align}
\end{subequations}

Next, we derive stability estimates for any solution $(\phi_h^{n}, \mu_h^n, \bz_{h}^{n}, \bbB_{h}^{n}) \in (\calS_h)^2 \times \calV_{h,\mathrm{div}}  \times \calW_{h}$ of the reduced system \eqref{eq:phi_reduced}--\eqref{eq:B_reduced}. After that, we use the stability estimates and a fixed-point argument to justify the existence of a solution of \eqref{eq:phi_reduced}--\eqref{eq:B_reduced}. In the end, we reconstruct the pressure $p_h^n\in\calS_h$ and show, that $\phi_h^n, \mu_h^n, \sigma_h^n, p_h^n, \bv_h^n, \bbB_h^n$ with $\bv_h^n = \bz_h^n + \bu_h^n$ forms a solution of the original regularized system \eqref{eq:phi_reg}--\eqref{eq:B_reg}.

\subsubsection{Stability of a reduced problem}
We now provide a stability estimate for the regularized scheme, supposed that there exists a solution (which will be proved after that).

For future reference, we recall the elementary identity
\begin{align}
    \label{eq:elementary_identity}
    2x(x-y) = x^2-y^2 + (x-y)^2 \quad\quad \forall \,  x,y\in\bbR.
\end{align}
Also, we note the following discrete version of Gronwall's inequality. For a proof, we refer to, e.g., \cite[pp.~401--402]{dahmen_reusken_numerik}. 
Assume that $e_n, a_n, b_n \geq 0$ for all $n\in \bbN_0$. Then 
\begin{align}
\label{eq:gronwall_discrete}
    e_n \leq a_n + \sum\limits_{i=0}^{n-1}  b_i e_i 
    \quad \forall \,  n\in \bbN_0
    \quad
    \Longrightarrow \quad
    e_n  \leq a_n \cdot 
    \exp\Big( \sum\limits_{i=0}^{n-1} b_i \Big)
    \quad \forall \,  n\in \bbN_0.
\end{align}

Due to the high complexity of the reduced system, it is quite hard to provide an explicit definition for the generic constants $C>0$, which, however, only depend on the domain $\Omega$, the space dimension $d$, the model parameters and the constants from the assumptions on the model functions. In particular, the generic constants $C$ are always independent of the discretization parameters $h,\Delta t$ and the regularization parameter $\delta>0$.


Now, testing \eqref{eq:phi_reduced} with $\xi_h = \mu_h^n + \chi_\phi \sigma_h^n$ and \eqref{eq:mu_reduced} with $\rho_h = \frac{1}{\Delta t} (\phi_h^n - \phi_h^{n-1})$. By adding the resulting equations and on noting \ref{A6}, \eqref{eq:elementary_identity}, we obtain
\begin{align*}
    &\frac{\beta\varepsilon}{2\Delta t} \Big( \norm{\nabla\phi_h^n}_{L^2}^2 - \norm{\nabla\phi_h^{n-1}}_{L^2}^2  
    +  \norm{\nabla\phi_h^n - \nabla\phi_h^{n-1}}_{L^2}^2 \Big)
    + \frac{\beta}{\varepsilon \Delta t} \skp{\psi(\phi_h^n)-\psi(\phi_h^{n-1})}{1}_h
    \\
    &\quad
    + \frac{1}{2\Delta t} \skp{\kappa(\phi_h^n) - \kappa(\phi_h^{n-1})}
    {\trace \bbB_h^n}_h
    + \skp{\calI_h m_{\phi,h}^{n-1}}{\abs{\nabla\mu_h^n}^2}_{L^2}
    \\
    &\quad 
    + \skp{(\bz_h^n + \bu_h^n) \cdot \nabla\phi_h^{n-1}}{\mu_h^n + \chi_\phi \sigma_h^n}_{L^2}
    + \skp{\phi_h^{n-1} \Gamma_{\bv,h}^{n-1} - \Gamma_{\phi,h}^{n-1}}{\mu_h^n + \chi_\phi \sigma_h^n}_h
    + \chi_\phi \skp{\calI_h[m_{\phi,h}^{n-1}] \nabla\mu_h^n}{\nabla\sigma_h^n}_{L^2}
    \\
    &\leq 0.
\end{align*}
Choosing $\bw_h = \bz_h^n \in \calV_{h,\mathrm{div}}$ in \eqref{eq:v_reduced} gives
\begin{align*}
    0 &= \skp{2 \calI_h[\eta(\phi_h^{n-1})]}{\abs{ \D(\bz_h^n)}^2}_{L^2}
    + \skp{2 \calI_h[\eta(\phi_h^{n-1})] \D(\bu_h^n)}{\D(\bz_h^n)}_{L^2}
    \\
    \nonumber
    &\quad
    + \skp{\calI_h \bbT_{e,h}^{\delta,n}}{\nabla\bz_h^n}_{L^2}
    - \skp{(\mu_h^n + \chi_\phi \sigma_h^n) \nabla\phi_h^{n-1}}{\bz_h^n}_{L^2}
    \\
    &\quad
    - \frac12 \skp{\nabla \calI_h\left[ \kappa(\phi_h^{n-1}) \trace\beta_\delta(\bbB_h^n) \right]}{\bz_h^n}_{L^2}
    + \frac12 \sum_{i,j=1}^d \skp{\partial_{x_j} \calI_h\left[ \kappa(\phi_h^{n-1}) \right]  \trace \mathbf{\Lambda}^\delta_{i,j}(\bbB_h^n)}{(\bz_h^n)_i}_{L^2}.
\end{align*}
Moreover, setting $\bbG_h = \frac12\bbB_h^n + \frac12 \calI_h \kappa(\phi_h^{n-1}) \bbI - \frac12 \calI_h \beta_\delta^{-1}(\bbB_h^n)$ in \eqref{eq:B_reduced} yields
\begin{align*}
    0 &= 
    \frac{1}{2\Delta t} \skp{\bbB_h^n - \bbB_h^{n-1}}{\bbB_h^n +  \kappa(\phi_h^{n-1}) \bbI - \beta_\delta^{-1}(\bbB_h^n)}_h
    \\
    &\quad
    + \frac12 \skp{\bz_h^n + \bu_h^n}{ \nabla\calI_h\left[\beta_\delta(\bbB_h^n) : \big( 
    \bbB_h^n + \calI_h \kappa(\phi_h^{n-1}) \bbI - \calI_h \beta_\delta^{-1}(\bbB_h^n)
    \big) \right]}_{L^2}
    \\
    &\quad
    - \frac12 \sum_{i,j=1}^d \skp{(\bz_h^n + \bu_h^n)_i \, \mathbf{\Lambda}^\delta_{i,j}(\bbB_h^n)}{\partial_{x_j} \big( \bbB_h^n + \calI_h \kappa(\phi_h^{n-1}) \bbI - \calI_h \beta_\delta^{-1}(\bbB_h^n) \big)}_{L^2}
    \\
    &\quad
    + \skp{\frac{1}{2\tau(\phi_h^{n-1})} \bbT_{e,h}^{\delta, n}}{\bbB_h^n + \kappa(\phi_h^{n-1}) \bbI - \beta_\delta^{-1}(\bbB_h^n)}_h 
    + \frac12 \skp{\Gamma_{\bbB,h}^{n-1} \beta_\delta(\bbB_h^n)}{\bbB_h^n + \kappa(\phi_h^{n-1}) \bbI - \beta_\delta^{-1}(\bbB_h^n)}_h
    \\
    &\quad
    - \skp{\nabla (\bz_h^n + \bu_h^n)}{ \calI_h \big[ \big(\bbB_h^n + \kappa(\phi_h^{n-1}) \bbI - \beta_\delta^{-1}(\bbB_h^n) \big) \beta_\delta(\bbB_h^n) \big] }_{L^2}
    \\
    &\quad
    + \frac12 \alpha \skp{\nabla\bbB_h^n}{\nabla(\bbB_h^n + \calI_h \kappa(\phi_h^{n-1}) \bbI - \calI_h \beta_\delta^{-1}(\bbB_h^n))}_{L^2}.
\end{align*}
Using \eqref{eq:elementary_identity}, \eqref{eq:nabla_ln}, \eqref{eq:Lambda_delta_kettenregel}, the definition of $\bbT_{e,h}^{\delta, n}$, \eqref{eq:reg_d} and the fact that $\bz_h^n \in \calV_{h,\mathrm{div}}$, we obtain
\begin{align*}
    &\frac{1}{4\Delta t} 
    \Big( \norm{\bbB_h^n}_h^2 - \norm{\bbB_h^{n-1}}_h^2 
    + \norm{\bbB_h^n - \bbB_h^{n-1}}_h^2 \Big) 
    + \frac1{2\Delta t} \skp{\kappa(\phi_h^{n-1})}{\trace\bbB_h^n - \trace\bbB_h^{n-1}}_h
    \\
    &\quad 
    + \frac1{2\Delta t} \skp{-\trace g_\delta(\bbB_h^n) + \trace g_\delta(\bbB_h^{n-1})}{1}_h
    \\
    &\quad
    + \frac12 \skp{\nabla \calI_h \big[ \trace f_\delta(\bbB_h^n) + \trace \ln \beta_\delta^{-1}(\bbB_h^n) + \kappa(\phi_h^{n-1}) \trace \beta_\delta(\bbB_h^n) \big]}{\bu_h^n}_{L^2}
    \\
    &\quad
    - \frac12 \sum_{i,j=1}^d \skp{(\bu_h^n)_i \, \trace \mathbf{\Lambda}^\delta_{i,j}(\bbB_h^n)}{\partial_{x_j} \calI_h \kappa(\phi_h^{n-1})}_{L^2}
    + \skp{\frac{1}{2\tau(\phi_h^{n-1})}} {\abs{\bbT_{e,h}^{\delta, n} \, \beta_\delta^{-1/2}(\bbB_h^n)}^2 }_h 
    \\
    &\quad 
    + \frac12 \skp{\Gamma_{\bbB,h}^{n-1}}{\trace \bbT_{e,h}^{\delta, n}}_h
    - \skp{\nabla (\bz_h^n + \bu_h^n)}{ \calI_h \bbT_{e,h}^{\delta, n} }_{L^2}
    \\
    &\quad
    + \frac12 \alpha \norm{\nabla\bbB_h^n}_{L^2}^2
    + \frac12 \alpha \skp{\nabla\trace\bbB_h^n}{\nabla \calI_h \kappa(\phi_h^{n-1})}_h
    + \frac{\alpha}{2d} \norm{\nabla \calI_h \trace \ln \beta_\delta(\bbB_h^n)}_{L^2}^2
    \\
    & \leq 0.
\end{align*}

Summing everything gives rise to
\begin{align}
    \label{eq:energy_estimate_dissipative}
    \nonumber
    &\frac{\beta\varepsilon}{2\Delta t} \Big( \norm{\nabla\phi_h^n}_{L^2}^2 - \norm{\nabla\phi_h^{n-1}}_{L^2}^2  
    +  \norm{\nabla\phi_h^n - \nabla\phi_h^{n-1}}_{L^2}^2 \Big)
    + \frac{\beta}{\varepsilon \Delta t} \skp{\psi(\phi_h^n)-\psi(\phi_h^{n-1})}{1}_h
    \\ \nonumber
    &\quad
    +\frac{1}{4\Delta t} \Big( \norm{\bbB_h^n}_{h}^2 - \norm{\bbB_h^{n-1}}_{h}^2 + \norm{\bbB_h^n- \bbB_h^{n-1}}_{h}^2\Big)
    + \frac{1}{2\Delta t} \skp{\kappa(\phi_h^n) \trace\bbB_h^n - \kappa(\phi_h^{n-1}) \trace\bbB_h^{n-1}}{1}_h
    \\ \nonumber
    &\quad 
    + \frac{1}{2\Delta t} \skp{- \trace g_\delta(\bbB_h^n) + \trace g_\delta(\bbB_h^{n-1})}{1}_h
    + \skp{\calI_h m_{\phi,h}^{n-1}}{\abs{\nabla\mu_h^n}^2}_{L^2}
    + \skp{2 \calI_h[\eta(\phi_h^{n-1})]}{\abs{ \D(\bz_h^n)}^2}_{L^2}
    \\ \nonumber
    &\quad
    + \skp{\frac{1}{2\tau(\phi_h^{n-1})}} {\abs{\bbT_{e,h}^{\delta, n} \, \beta_\delta^{-1/2}(\bbB_h^n)}^2 }_h 
    + \frac12 \alpha \norm{\nabla \bbB_h^n}_{L^2}^2
    + \frac{\alpha}{2d} \norm{\nabla \calI_h \trace \ln \beta_\delta(\bbB_h^n)}_{L^2}^2
    \\ \nonumber
    &\leq
    - \chi_\phi \skp{\calI_h[m_{\phi,h}^{n-1}] \nabla\mu_h^n}{\nabla\sigma_h^n}_{L^2}
    - \skp{\bu_h^n \cdot \nabla\phi_h^{n-1}}{\mu_h^n + \chi_\phi \sigma_h^n}_{L^2}
    + \skp{\Gamma_{\phi,h}^{n-1} - \phi_h^{n-1} \Gamma_{\bv,h}^{n-1}}{\mu_h^n + \chi_\phi \sigma_h^n}_h
    \\ \nonumber
    &\quad 
    - \skp{2 \calI_h[\eta(\phi_h^{n-1})] \D(\bu_h^n)}{\D(\bz_h^n)}_{L^2}
    - \frac12 \skp{\nabla \calI_h \big[ \trace f_\delta(\bbB_h^n) + \trace \ln \beta_\delta^{-1}(\bbB_h^n) + \kappa(\phi_h^{n-1}) \trace \beta_\delta(\bbB_h^n) \big]}{\bu_h^n}_{L^2}
    \\ \nonumber
    &\quad
    + \frac12 \sum_{i,j=1}^d \skp{(\bu_h^n)_i \trace \mathbf{\Lambda}^\delta_{i,j}(\bbB_h^n)}{\partial_{x_j} \calI_h\left[ \kappa(\phi_h^{n-1}) \right] }_{L^2}
    - \frac12 \skp{\Gamma_{\bbB,h}^{n-1} }{\trace \bbT_{e,h}^{\delta,n}}_h
    \\
    &\quad
    + \skp{\nabla\bu_h^n}{\calI_h \bbT_{e,h}^{\delta,n}}_{L^2}
    - \frac12 \alpha \skp{\nabla\trace\bbB_h^n}{\nabla\calI_h[ \kappa( \phi_h^{n-1}) ]}_{L^2}.
\end{align}
Here we see that the right-hand side of \eqref{eq:energy_estimate_dissipative} would vanish in absence of source terms and chemotaxis, i.e., if $\chi_\phi=\Gamma_{\bv}=\Gamma_\phi = \Gamma_{\bbB} = 0$, and if $\kappa(\phi)= \bar\kappa$ for a constant $\bar\kappa\in\bbR$, see also Remark \ref{remark:scheme}. Note that, in the case $\Gamma_{\bv}=0$, $\bu_h^n$ would not appear and we would have $\bv_h^n = \bz_h^n$.

Now we estimate the terms on the right-hand side of \eqref{eq:energy_estimate_dissipative}. On noting \ref{A2} and Hölder's and Young's inequalities, we have
\begin{align*}
    \abs{ \chi_\phi \skp{\calI_h[m_{\phi,h}^{n-1}] \nabla\mu_h^n}{\nabla\sigma_h^n}_{L^2} }
    \leq 
    \frac14 m_0 \norm{\nabla\mu_h^n}_{L^2}^2 
    + C \norm{\sigma_h^n}_{H^1}^2.
\end{align*}
Moreover, using \ref{A2}, \eqref{eq:bound_uh} and Hölder's and Young's inequalities, we calculate
\begin{align*}
    \abs{\skp{2 \calI_h[\eta(\phi_h^{n-1})] \D(\bu_h^n)}{\D(\bz_h^n)}_{L^2}}
    \leq 
    \eta_0 \norm{\D(\bz_h^n)}_{L^2}^2 
    + C \norm{\bu_h^n}_{H^1}^2
    \leq 
    \eta_0 \norm{\D(\bz_h^n)}_{L^2}^2 
    + C ( 1 + \norm{\sigma_h^n}_{H^1}^2).
\end{align*}
Similarly, we have with \eqref{eq:bound_Te_delta}, \eqref{eq:bound_source}, \eqref{eq:bound_uh}, \eqref{eq:norm_equiv} and Hölder's and Young's inequalities, that
\begin{align*}
    &\abs{\skp{\nabla\bu_h^n}{\calI_h \bbT_{e,h}^{\delta,n}}_{L^2}}
    + \abs{\frac12 \skp{\Gamma_{\bbB,h}^{n-1} }{\trace \bbT_{e,h}^{\delta,n}}_h}
    \\
    &\leq C \big( \norm{\nabla\bu_h^n}_{L^2} + \norm{\calI_h \Gamma_{\bbB,h}^{n-1}}_{L^2} \big) \norm{\calI_h \bbT_{e,h}^{\delta,n}}_{L^2}
    \\
    &\leq \frac\alpha8 \norm{\nabla\bbB_h^n}_{L^2}^2 
    + C (1+\norm{\sigma_h^n}_{H^1}^4) (1+ \norm{\bbB_h^n}_h^2) 
    + C \norm{\phi_h^{n-1}}_{H^1}^2.
\end{align*}
Besides, using \eqref{eq:bound_uh}, \eqref{eq:bound_grad_k_h}, Hölder's and Young's inequalities, it holds
\begin{align*}
    \abs{\frac12 \alpha \skp{\nabla\trace\bbB_h^n}{\nabla\calI_h[ \kappa( \phi_h^{n-1}) ]}_{L^2}}
    \leq 
    \frac\alpha8 \norm{\nabla\bbB_h^n}_{L^2}^2 
    + C \norm{\nabla\calI_h \kappa(\phi_h^{n-1})}_{L^2}^2
    \leq 
    \frac\alpha8 \norm{\nabla\bbB_h^n}_{L^2}^2 
    + C \norm{\nabla\phi_h^{n-1}}_{L^2}^2.
\end{align*}
We deduce from \eqref{eq:Lambda_delta}, \eqref{eq:inverse_estimate}, \eqref{eq:bound_beta_delta}, \eqref{eq:bound_grad_k_h} and Hölder's and Young's inequalities, that
\begin{align*}
    &\abs{\frac12 \sum_{i,j=1}^d \skp{(\bu_h^n)_i \trace \mathbf{\Lambda}^\delta_{i,j}(\bbB_h^n)}{\partial_{x_j} \calI_h\left[ \kappa(\phi_h^{n-1}) \right] }_{L^2}}
    \\
    &\leq  C \norm{\bu_h^n}_{L^\infty} (1 + \norm{\bbB_h^n}_{L^2}) \norm{\nabla\calI_h \kappa(\phi_h^{n-1})}_{L^2}
    \\
    &\leq C (1+\norm{\sigma_h^n}_{H^1}) ( 1 + \norm{\bbB_h^n}_{L^2}) \norm{\nabla \phi_h^{n-1}}_{L^2}
    \\
    &\leq C \big(1 + \norm{\sigma_h^n}_{H^1}^2)  
    (1 + \norm{\bbB_h^n}_{L^2}^2) 
    + C \norm{\nabla\phi_h^{n-1}}_{L^2}^2 .
\end{align*}
Using Hölder's and Young's inequalities, \eqref{eq:bound_uh} and the fact that $\ln(s^{-1}) = - \ln(s)$ $\forall\, s>0$, we obtain
\begin{align*}
    \abs{\frac12 \skp{\nabla \calI_h \trace \ln \beta_\delta^{-1}(\bbB_h^n)}{\bu_h^n}_{L^2} }
    &\leq \frac{\alpha}{4d} \norm{\nabla \calI_h \trace \ln \beta_\delta(\bbB_h^n)}_{L^2}^2 + C \norm{\bu_h^n}_{L^2}^2
    \\
    &\leq \frac{\alpha}{4d} \norm{\nabla \calI_h \trace \ln \beta_\delta(\bbB_h^n)}_{L^2}^2 + C ( 1 + \norm{\sigma_h^n}_{H^1}^2).
\end{align*}
On noting $f_\delta'(\cdot) = \beta_\delta(\cdot)$, \eqref{eq:bound_beta_delta}, \eqref{eq:bound_grad_beta_delta}, \eqref{eq:interp_stability}, \eqref{eq:inverse_estimate}, \eqref{eq:bound_k_h}, \eqref{eq:bound_grad_k_h}, \eqref{eq:bound_uh} and Hölder's and Young's inequalities, we have
\begin{align*}
    &\abs{ \frac12 \skp{\nabla \calI_h \big[ \trace f_\delta(\bbB_h^n) + \kappa(\phi_h^{n-1}) \trace \beta_\delta(\bbB_h^n) \big]}{\bu_h^n}_{L^2} }
    \\
    &\leq C \norm{\nabla\bbB_h^n}_{L^2} \big( 1 + \norm{\bbB_h^n}_{L^2}  \big) \norm{\bu_h^n}_{L^\infty}
    + C \norm{\nabla\phi_h^{n-1}}_{L^2} \big( 1 + \norm{\bbB_h^n}_{L^2}  \big) \norm{\bu_h^n}_{L^\infty}
    \\
    &\quad
    + C \big( 1 + \norm{\phi_h^{n-1}}_{L^2}\big) \norm{\nabla\bbB_h^n}_{L^2} \norm{\bu_h^n}_{L^\infty}
    \\
    &\leq \frac\alpha8 \norm{\nabla\bbB_h^n}_{L^2}^2 
    + C ( 1 + \norm{\sigma_h^n}_{H^1}^2) \big(1 + \norm{\bbB_h^n}_{L^2}^2 
    + \norm{\phi_h^{n-1}}_{H^1}^2  \big).
\end{align*}
%
%
%
%
%
%
%
%
For the remaining terms on the right-hand side of \eqref{eq:energy_estimate_dissipative}, we deduce from \eqref{eq:bound_source}, \eqref{eq:norm_equiv}, \eqref{eq:bound_uh}, Hölder's and Young's inequalities, that
\begin{align*}
    &\abs{\skp{\bu_h^n \cdot \nabla\phi_h^{n-1}}{\mu_h^n + \chi_\phi \sigma_h^n}_{L^2}}
    + \abs{\skp{\Gamma_{\phi,h}^{n-1} - \phi_h^{n-1} \Gamma_{\bv,h}^{n-1}}{\mu_h^n + \chi_\phi \sigma_h^n}_h}
    \\
    &\leq 
    C \norm{\bu_h^n}_{L^\infty} \norm{\nabla\phi_h^{n-1}}_{L^2} ( \norm{\mu_h^n}_{L^2} + \norm{\sigma_h^n}_{L^2})
    + C (1+\norm{\sigma_h^n}_{L^4}) (1+\norm{\phi_h^{n-1}}_{L^4}) (\norm{\mu_h^n}_{L^2} + \norm{\sigma_h^n}_{L^2})
    \\
    &\leq \frac12 \norm{\mu_h^n}_h^2 + C (1 + \norm{\phi_h^{n-1}}_{H^1}^2) ( 1 + \norm{\sigma_h^n}_{H^1}^4) .
\end{align*}
Due to the presence of the source terms $\Gamma_\phi$ and $\Gamma_{\bv}$ on the right-hand side of \eqref{eq:energy_estimate_dissipative}, we need a further estimate for $\mu_h^n$ in order to control $\norm{\mu_h^n}_h^2$. 
Therefore, we test \eqref{eq:mu_reduced} with $\rho_h=\mu_h^n$ and we get with \eqref{eq:bound_psi_h'}, \eqref{eq:bound_kh'_tau}, Hölder's and Young's inequalities, that
\begin{align*}
    \norm{\mu_h^n}_h^2 
    &= \skp{\frac\beta\varepsilon \psi_h'(\phi_h^n,\phi_h^{n-1}) - \chi_\phi \sigma_h^n + \frac12 \kappa_h'(\phi_h^n, \phi_h^{n-1}) \trace\bbB_h^n} {\mu_h^n}_h
    + \beta\varepsilon \skp{\nabla\phi_h^n}{\nabla\mu_h^n}_{L^2}
    \\
    &\leq \frac12 \norm{\mu_h^n}_h^2
    + C ( 1 + \norm{\sigma_h^n}_h^2 + \norm{\phi_h^{n-1}}_h^2
    + \norm{\phi_h^n}_h^2 + \norm{\bbB_h^n}_h^2 + \norm{\nabla\phi_h^n}_{L^2}^2
    + \frac{m_0}{8} \norm{\nabla\mu_h^n}_{L^2}^2,
\end{align*}
from which we deduce
\begin{align*}
    \norm{\mu_h^n}_h^2 
    \leq \frac{m_0}{4} \norm{\nabla\mu_h^n}_{L^2}^2
    + C ( 1 + \norm{\sigma_h^n}_{H^1}^2 
    + \norm{\phi_h^{n-1}}_h^2
    + \norm{\phi_h^n}_h^2 + \norm{\bbB_h^n}_h^2 + \norm{\nabla\phi_h^n}_{L^2}^2 ).
\end{align*}

Combining the estimates after \eqref{eq:energy_estimate_dissipative} and using \eqref{eq:norm_equiv} and \ref{A2}, we have 
\begin{align*}
    &\frac{\beta\varepsilon}{2\Delta t} \Big( \norm{\nabla\phi_h^n}_{L^2}^2 - \norm{\nabla\phi_h^{n-1}}_{L^2}^2  
    +  \norm{\nabla\phi_h^n - \nabla\phi_h^{n-1}}_{L^2}^2 \Big)
    + \frac{\beta}{\varepsilon \Delta t} \skp{\psi(\phi_h^n)-\psi(\phi_h^{n-1})}{1}_h
    \\
    &\quad
    +\frac{1}{4\Delta t} \Big( \norm{\bbB_h^n}_{h}^2 - \norm{\bbB_h^{n-1}}_{h}^2 + \norm{\bbB_h^n- \bbB_h^{n-1}}_{h}^2\Big)
    + \frac{1}{2\Delta t} \Big( \skp{\kappa(\phi_h^n)}{\trace\bbB_h^n}_h
    - \skp{\kappa(\phi_h^{n-1})}{\trace\bbB_h^{n-1}}_h \Big)
    \\
    &\quad 
    + \frac{1}{2\Delta t} \skp{- \trace g_\delta(\bbB_h^n) + \trace g_\delta(\bbB_h^{n-1})}{1}_h
    + \frac{m_0}{2} \norm{\nabla\mu_h^n}^2_{L^2}
    + \frac12 \norm{\mu_h^n}_h^2
    + \eta_0 \norm{ \D(\bz_h^n)}^2_{L^2}
    \\
    &\quad
    + \frac{1}{2\tau_1}\norm{\bbT_{e,h}^{\delta, n} \, \beta_\delta^{-1/2}(\bbB_h^n)}_h^2 
    + \frac18 \alpha \norm{\nabla \bbB_h^n}_{L^2}^2
    + \frac{\alpha}{4d} \norm{\nabla \calI_h \trace \ln \beta_\delta(\bbB_h^n)}_{L^2}^2
    \\
    &\leq
    C \big( 
    1 + \norm{\phi_h^{n-1}}_h^2 + \norm{\nabla\phi_h^{n-1}}_{L^2}^2
    + \norm{\bbB_h^n}_h^2 
    + \norm{\phi_h^n}_h^2
    + \norm{\nabla\phi_h^n}_{L^2}^2
    \big) ( 1 + \norm{\sigma_h^n}_{H^1}^4). 
\end{align*}
Then, we multiply both sides by $\Delta t$ and sum over $n\in\{1,...,m\}$, where $m\in\{1,...,N_T\}$, and we get
\begin{align}
    \label{eq:energy_estimate_dissipative2} \nonumber
    &\frac{\beta\varepsilon}{2} \norm{\nabla\phi_h^m}_{L^2}^2 
    + \frac\beta\varepsilon \skp{\psi(\phi_h^m)}{1}_h
    + \skp{\frac14 \abs{\bbB_h^m}^2 + \frac12 \kappa(\phi_h^m) \trace\bbB_h^m - \frac12 \trace g_\delta(\bbB_h^m)}{1}_h
    \\ \nonumber
    &+  C \sum_{n=1}^{m} \Big( \norm{\nabla\phi_h^n - \nabla\phi_h^{n-1}}_{L^2}^2
    + \norm{\bbB_h^n- \bbB_h^{n-1}}_{h}^2 \Big)
    + C \Delta t \sum_{n=1}^{m} \Big( \norm{\nabla\mu_h^n}^2_{L^2}
    + \norm{\mu_h^n}_h^2
    + \norm{ \D(\bz_h^n)}^2_{L^2}
    \Big)
    \\ \nonumber
    &\quad 
    + C \Delta t \sum_{n=1}^{m} \Big( 
    \norm{\bbT_{e,h}^{\delta,n} \beta_\delta^{-1/2}(\bbB_h^n) }_h^2
    + \norm{\nabla \bbB_h^n}_{L^2}^2
    + \norm{\nabla \calI_h \trace \ln \beta_\delta^{-1}(\bbB_h^n)}_{L^2}^2 \Big)
    \\ \nonumber
    &\leq
    \frac{\beta\varepsilon}{2} \norm{\nabla\phi_h^0}_{L^2}^2  
    + \frac\beta\varepsilon \skp{\psi(\phi_h^0)}{1}_h
    + \skp{\frac14 \abs{\bbB_h^0}^2 + \frac12 \kappa(\phi_h^0) \trace\bbB_h^0 - \frac12 \trace g_\delta(\bbB_h^0)}{1}_h
    \\ 
    &\quad
    + C \Delta t \sum_{n=0}^{m} \Big( 
    1 + \norm{\bbB_h^n}_h^2
    + \norm{\phi_h^n}_h^2
    + \norm{\nabla\phi_h^n}_{L^2}^2
    \Big) ( 1 + c_\infty^4),
    \\ \nonumber
\end{align}
where we used $\max_{n\in\{1,...,N_T\}} \norm{\sigma_h^n}_{H^1} \leq c_\infty \coloneqq \max_{n\in\{1,...,N_T\}} \norm{\sigma_{\infty,h}^n}_{L^2(\partial\Omega)}$.
To apply a Gronwall argument, the terms on the left-hand side of \eqref{eq:energy_estimate_dissipative2} have to be non-negative, so we need to control the term $\kappa(\phi_h^m) \trace\bbB_h^m$. On noting \ref{A3}, we calculate with Hölder's and Young's inequalities.
\begin{align*}
    \abs{\frac12 \skp{\kappa(\phi_h^m) \trace\bbB_h^m}{1}_h }
    \leq \frac{1}{16} \norm{\bbB_h^m}_h^2 + \norm{\kappa(\phi_h^m)}_h^2
    \leq \frac{1}{16} \norm{\bbB_h^m}_h^2 + C_\kappa^2 \big( \norm{\phi_h^m}_h^2 + 1 \big),
\end{align*}
and, for any $n\in\{0,...,m\}$, we have with \ref{A4}, that
\begin{align*}
    \norm{\phi_h^n}_h^2 \leq  C_\psi^{-1} \big( \skp{\psi(\phi_h^n)}{1}_h + 1 \big).
\end{align*}
Using this and putting the terms in \eqref{eq:energy_estimate_dissipative2} with index $m$ to the left-hand side, we obtain
\begin{align*}
    &\Big(\frac{\beta\varepsilon}{2} - C \Delta t (1 + c_\infty^4) \Big) \norm{\nabla\phi_h^m}_{L^2}^2 
    + \big( \frac\beta\varepsilon - C_\kappa^2 C_\psi^{-1} - C \Delta t (1 + c_\infty^4) \big) \skp{\psi(\phi_h^m)}{1}_h
    \\
    &\quad 
    + \Big( \frac{1}{16} - C \Delta t (1 + c_\infty^4 )\Big) \norm{\bbB_h^m}_h^2
    + \skp{\frac18 \abs{\bbB_h^m}^2 - \frac12 \trace g_\delta(\bbB_h^m)}{1}_h
    \\
    &\quad 
    +  C \sum_{n=1}^{m} \Big( \norm{\nabla\phi_h^n - \nabla\phi_h^{n-1}}_{L^2}^2
    + \norm{\bbB_h^n- \bbB_h^{n-1}}_{h}^2 \Big)
    + C \Delta t \sum_{n=1}^{m} \Big( \norm{\nabla\mu_h^n}^2_{L^2}
    + \norm{\mu_h^n}_h^2
    + \norm{ \D(\bz_h^n)}^2_{L^2}
    \Big)
    \\
    &\quad 
    + C \Delta t \sum_{n=1}^{m} \Big( 
    \norm{\bbT_{e,h}^{\delta,n} \beta_\delta^{-1/2}(\bbB_h^n) }_h^2
    + \norm{\nabla \bbB_h^n}_{L^2}^2
    + \norm{\nabla \calI_h \trace \ln \beta_\delta(\bbB_h^n)}_{L^2}^2 \Big)
    \\
    &\leq
    \abs{\calF_h^\delta(\phi_h^0, \bbB_h^0) }
    + C \Delta t \sum_{n=0}^{m-1} \Big( 
    1 + \norm{\bbB_h^n}_h^2
    + \skp{\psi(\phi_h^n)}{1}_h
    + \norm{\nabla\phi_h^n}_{L^2}^2
    \Big) (1 + c_\infty^4),
\end{align*}
where 
\begin{align*}
    \calF_h^\delta(\phi_h^0, \bbB_h^0) 
    \coloneqq 
    \frac{\beta\varepsilon}{2} \norm{\nabla\phi_h^0}_{L^2}^2  
    + \frac\beta\varepsilon \skp{\psi(\phi_h^0)}{1}_h
    + \skp{\frac14 \abs{\bbB_h^0}^2 + \frac12 \kappa(\phi_h^0) \trace\bbB_h^0 - \frac12 \trace g_\delta(\bbB_h^0)}{1}_h.
\end{align*}
By the assumption \ref{A4}, it is $\frac\beta\varepsilon> C_\kappa^2 C_\psi^{-1}$. 
Moreover, it follows from \eqref{eq:reg_c} that $\frac18 \abs{\bbB_h^m}^2 - \frac12 \trace g_\delta(\bbB_h^m)$ is non-negative.
So, the coefficients on the left-hand side are non-negative supposed that $\Delta t>0$ is small enough. 
This, in particular, leads to the constraint $\Delta t < \Delta t_*$ for a (possibly small) constant $\Delta t_*>0$ that only depends on the model parameters and on $c_\infty = \max_{n\in\{1,...,N_T\}} \norm{\sigma_{\infty,h}^n}_{L^2(\partial\Omega)}$, but not on $\delta, h>0$.
Applying the discrete Gronwall inequality \eqref{eq:gronwall_discrete} and taking the maximum over $m\in\{1,...,N_T\}$ lead to
\begin{align}
    \label{eq:energy_estimate_dissipative3} \nonumber
    & \max_{m\in\{1,...,N_T\}} \left( \norm{\nabla\phi_h^m}_{L^2}^2 
    + \skp{\psi(\phi_h^m)}{1}_h
    + \norm{\bbB_h^m}_h^2
    + \skp{\frac18 \abs{\bbB_h^m}^2 - \frac12 \trace g_\delta(\bbB_h^m)}{1}_h \right)
    \\ \nonumber
    &+  \sum_{n=1}^{N_T} \Big( \norm{\nabla\phi_h^n - \nabla\phi_h^{n-1}}_{L^2}^2
    + \norm{\bbB_h^n- \bbB_h^{n-1}}_{h}^2 \Big)
    + \Delta t \sum_{n=1}^{N_T} \Big( \norm{\nabla\mu_h^n}^2_{L^2}
    + \norm{\mu_h^n}_h^2
    + \norm{ \D(\bz_h^n)}^2_{L^2}
    \Big)
    \\ \nonumber
    &\quad 
    + \Delta t \sum_{n=1}^{N_T} \Big( 
    \norm{\bbT_{e,h}^{\delta,n} \, \beta_\delta^{-1/2}(\bbB_h^n) }_h^2
    + \norm{\nabla \bbB_h^n}_{L^2}^2
    + \norm{\nabla \calI_h \trace \ln \beta_\delta(\bbB_h^n)}_{L^2}^2 \Big)
    \\
    &\leq C(T, c_\infty) 
    \big(1 + \abs{\calF_h^\delta(\phi_h^0, \bbB_h^0) } \big),
\end{align}
for a constant $C(T, c_\infty)>0$ that depends exponentially on $T$ and $c_\infty$.

From \eqref{eq:energy_estimate_dissipative3}, we also get with \ref{A4}, \eqref{eq:norm_equiv}, \eqref{eq:reg_c} and Korn's inequality, that
\begin{align}
    \label{eq:energy_estimate_dissipative4}
    \max_{m\in\{1,...,N_T\}} 
    \Big( \norm{\phi_h^m}_{L^2}^2 + \frac1\delta \norm{ \calI_h \abs{ [\bbB_h^m]_-} }_{L^1} \Big)
    + \Delta t \sum_{n=1}^{N_T} \norm{\bz_h^n}_{H^1}^2
    \leq C(T, c_\infty) \big(1 + \abs{\calF_h^\delta(\phi_h^0, \bbB_h^0) } \big) .
\end{align}

\subsubsection{Existence for the reduced problem}
The next step is to combine the stability estimates \eqref{eq:energy_estimate_dissipative3}--\eqref{eq:energy_estimate_dissipative4} and a fixed-point strategy to prove existence of a solution of the reduced problem \eqref{eq:phi_reduced}--\eqref{eq:B_reduced}. This will be shown with a proof by contradiction.

First, we need to find some specific mappings and a convex and compact finite dimensional set such that we can apply Brouwer's fixed point theorem. We define the following inner product
\begin{align*}
    \skpp{(\phi_h,\mu_h,\bz_h,\bbB_h)}{(\zeta_h,\rho_h,\bw_h,\bbG_h)}
    &\coloneqq 
    \skp{\phi_h}{\zeta_h}_h
    + \skp{\mu_h}{\rho_h}_h
    + \skp{\bz_h}{\bw_h}_{L^2}
    + \skp{\bbB_h}{\bbG_h}_h,
\end{align*}
for any $(\phi_h,\mu_h,\bz_h,\bbB_h), (\zeta_h,\rho_h,\bw_h,\bbG_h) \in (\calS_h)^2 \times \calV_{h,\text{div}} \times \calW_h$ on the Hilbert space $(\calS_h)^2 \times \calV_{h,\text{div}} \times \calW_h$.

For some given $(\phi_h^{n-1},\bv_h^{n-1}, \bbB_h^{n-1})\in \calS_h \times \calV_{h,\text{div}} \times \calW_h$ from the previous time step and $\sigma_h^n \in \calS_h$, $\bu_h^n \in \calV_h$ as above, let the mapping
\begin{align*}
    \calH\colon\, (\calS_h)^2 \times \calV_{h,\text{div}} \times \calW_h \longrightarrow (\calS_h)^2 \times \calV_{h,\text{div}} \times \calW_h
\end{align*}
be such that, for any $(\phi_h,\mu_h,\bz_h,\bbB_h) \in (\calS_h)^2 \times \calV_{h,\text{div}} \times \calW_h$,
\begin{align*}
    &\skpp{\calH (\phi_h,\mu_h,\bz_h,\bbB_h)}{(\zeta_h,\rho_h,\bw_h,\bbG_h)}
    \\
    &\coloneqq 
    \frac{1}{\Delta t} \skp{\phi_h - \phi_h^{n-1}}{\zeta_h}_h
    + \skp{\calI_h [m_{\phi,h}^{n-1}] \nabla\mu_h}{\nabla\zeta_h}_{L^2}
    + \skp{(\bz_h + \bu_h^n) \cdot \nabla\phi_h^{n-1}}{\zeta_h}_{L^2}
    + \skp{\phi_h^{n-1} \Gamma_{\bv,h}^{n-1} - \Gamma_{\phi,h}^{n-1}}{\zeta_h}_h
    \\
    &\quad 
    + \skp{- \mu_h + \frac\beta\varepsilon \psi_h'(\phi_h,\phi_h^{n-1}) - \chi_\phi \sigma_h^n + \frac12 \kappa_h'(\phi_h, \phi_h^{n-1}) \trace \bbB_h} {\rho_h}_h
    + \beta\varepsilon \skp{\nabla\phi_h}{\nabla\rho_h}_{L^2}
    \\
    &\quad 
    + \skp{2 \calI_h[\eta(\phi_h^{n-1})] \D(\bz_h)}{\D(\bw_h)}_{L^2}
    + \skp{2 \calI_h[\eta(\phi_h^{n-1})] \D(\bu_h^n)}{\D(\bw_h)}_{L^2}
    \\
    &\quad
    + \skp{\calI_h \bbT_{e}^\delta(\phi_h^{n-1}, \bbB_h)
    }{\nabla\bw_h}_{L^2}
    - \skp{(\mu_h + \chi_\phi \sigma_h^n) \nabla\phi_h^{n-1}}{\bw_h}_{L^2}
    \\
    &\quad
    - \frac12 \skp{\nabla \calI_h\left[ \kappa(\phi_h^{n-1}) \trace\beta_\delta(\bbB_h) \right]}{\bw_h}_{L^2}
    + \frac12 \sum_{i,j=1}^d \skp{\partial_{x_j} \calI_h\left[ \kappa(\phi_h^{n-1}) \right]  \trace \mathbf{\Lambda}^\delta_{i,j}(\bbB_h)}{(\bw_h)_i}_{L^2}
    \\
    &\quad 
    + \frac{1}{\Delta t} \skp{\bbB_h - \bbB_h^{n-1}}{\bbG_h}_h
    + \skp{\bz_h + \bu_h^n}{ \nabla\calI_h\left[\beta_\delta(\bbB_h^n) : \bbG_h \right]}_{L^2}
    - \sum_{i,j=1}^d \skp{(\bz_h + \bu_h^n)_i \, \mathbf{\Lambda}^\delta_{i,j}(\bbB_h)}{\partial_{x_j} \bbG_h}_{L^2}
    \\
    &\quad
    + \skp{\frac{1}{\tau(\phi_h^{n-1})} \bbT_{e}^\delta(\phi_h^{n-1}, \bbB_h) }{\bbG_h}_h 
    + \skp{\Gamma_{\bbB,h}^{n-1} \beta_\delta(\bbB_h)}{\bbG_h}_h
    \\
    &\quad 
    - \skp{2 \nabla(\bz_h + \bu_h^n)}{\calI_h\left[ \bbG_h \beta_\delta(\bbB_h) \right]}_{L^2}
    + \alpha \skp{\nabla\bbB_h}{\nabla\bbG_h}_{L^2},
\end{align*}
for all $(\zeta_h,\rho_h,\bw_h,\bbG_h) \in (\calS_h)^2 \times \calV_{h,\text{div}} \times \calW_h $. Here, we use the notation
\begin{align*}
    \bbT_{e}^\delta(\phi_h^{n-1}, \bbB_h)
    \coloneqq \bbB_h \beta_\delta(\bbB_h) + \kappa(\phi_h^{n-1}) \beta_\delta(\bbB_h) - \bbI 
    \qquad \forall\, \bbB_h \in \calW_h.
\end{align*}
Note that the mapping $\calH$ is continuous. A wanted solution $(\phi_h,\mu_h,\bz_h,\bbB_h)$, if it exists, corresponds to a zero of $\calH$, i.e.,
\begin{align*}
    \skpp{\calH (\phi_h,\mu_h,\bz_h,\bbB_h)}{(\zeta_h,\rho_h,\bw_h,\bbG_h)}=0 \qquad \forall\, (\zeta_h,\rho_h,\bw_h,\bbG_h) \in (\calS_h)^2 \times \calV_{h,\text{div}} \times \calW_h .
\end{align*}
The testing procedure of the stability proof is very helpful to prove the existence of a zero of $\calH$. To mimic the testing procedure, we introduce a linear transformation $f\colon \, (\calS_h)^2 \times \calV_{h,\text{div}} \times \calW_h \longrightarrow (\calS_h)^2 \times \calV_{h,\text{div}} \times \calW_h$ and its inverse by
\begin{align*}
    f\colon\, (\phi_h, \mu_h, \bz_h, \bbB_h) 
    \mapsto 
    \Big( \mu_h, \,
    \frac{\phi_h}{\Delta t} - 2 \mu_h, \,
    \bz_h, \, 
    \bbB_h \Big),
    \qquad
    f^{-1}\colon\, (\zeta_h, \rho_h, \bw_h, \bbG_h) 
    \mapsto 
    \Big(\Delta t (\rho_h + 2 \zeta_h), \,
    \zeta_h, \,
    \bw_h, \, 
    \bbG_h \Big).
\end{align*}
Now let $R>0$ be given. The goal is to prove the following assumption wrong. Let us assume that the continuous mapping $\calH \circ f^{-1}$ has no zero in the closed ball 
\begin{align*}
    \calB_R \coloneqq \left\{ (\zeta_h,\rho_h,\bw_h,\bbG_h) \in (\calS_h)^2 \times \calV_{h,\text{div}} \times \calW_h
    \ \mid \
    ||| (\zeta_h,\rho_h,\bw_h,\bbG_h)||| \leq R \right\},
\end{align*}
where
\begin{align*}
    |||(\zeta_h,\rho_h,\bw_h,\bbG_h)|||^2 
    \coloneqq 
    \skpp{(\zeta_h,\rho_h,\bw_h,\bbG_h)}{(\zeta_h,\rho_h,\bw_h,\bbG_h)}.
\end{align*}
Then, for such $R>0$, we define a continuous mapping $\calG_R: \calB_R \to \partial \calB_R$ by
\begin{align*}
    \calG_R (\zeta_h,\rho_h,\bw_h,\bbG_h)
    \coloneqq 
    -R \, \frac{(\calH  \circ f^{-1}) (\zeta_h,\rho_h,\bw_h,\bbG_h)}{ 
    ||| (\calH \circ f^{-1}) (\zeta_h,\rho_h,\bw_h,\bbG_h) ||| }
    \quad\quad \forall \,  (\zeta_h,\rho_h,\bw_h,\bbG_h)\in \calB_R.
\end{align*}
We deduce from Brouwer's fixed point theorem \cite[Chap.~8.1.4, Thm.~3]{evans_2010} that there exists at least one fixed point $(\zeta_h^R, \rho_h^R, \bw_h^R, \bbG_h^R) \in \calB_R$ of the mapping $\calG_R$, and we define $(\phi_h^R, \mu_h^R,\bz_h^R, \bbB_h^R) \coloneqq f^{-1}(\zeta_h^R, \rho_h^R, \bw_h^R, \bbG_h^R)$. Then, it holds $f(\phi_h^R, \mu_h^R, \bz_h^R, \bbB_h^R) = \left(\calG_R\circ f\right)(\phi_h^R, \mu_h^R,\bz_h^R, \bbB_h^R)$ and, as $\calG_R\colon \, \calB_R\to \partial\calB_R$, 
\begin{align}
\label{eq:existence_R}
    ||| f(\phi_h^R, \mu_h^R, \bz_h^R, \bbB_h^R) |||
    = ||| \left(\calG_R \circ f\right) (\phi_h^R, \mu_h^R, \bz_h^R, \bbB_h^R) |||
    = R.
\end{align}
Next, the strategy is to show
\begin{align*}
    0 < \skpp{\calH(\phi_h^R, \mu_h^R, \bz_h^R, \bbB_h^R)}{(\zeta_h, \rho_h, \bw_h, \bbG_h)} < 0
\end{align*}
for one specific tuple of test functions $(\zeta_h, \rho_h, \bw_h, \bbG_h)$ and for $R>0$ large enough. This will disprove the assumption that $\calH\circ f^{-1}$ has no zero in $\calB_R$. The first inequality will be shown with the strategy from the proof of the stability estimates \eqref{eq:energy_estimate_dissipative3}--\eqref{eq:energy_estimate_dissipative4}, and for the second inequality, we will use that $f(\phi_h^R, \mu_h^R, \bv_h^R, \bbB_h^R)$ is a fixed-point of $\calG_R$.
Choosing the tuple of test functions as 
\begin{align*}
    (\zeta_h, \rho_h, \bw_h, \bbG_h)
    \coloneqq \left( 
    \mu_h^R + \chi_\phi \sigma_h^n, \,
    \frac{\phi_h^R - \phi_h^{n-1}}{\Delta t} - 2\mu_h^R, \,
    \bz_h^R, \,
    \frac12\big( \bbB_h^R 
    + \calI_h[\kappa(\phi_h^{n-1})] \bbI 
    - \calI_h[\beta_\delta^{-1}(\bbB_h^R)]\big) 
    \right),
\end{align*}
we obtain analogously to the stability proof with a straightforward computation
\begin{align*}
    &\skpp{\calH (\phi_h^R,\mu_h^R,\bz_h^R,\bbB_h^R)}{(\zeta_h,\rho_h,\bw_h,\bbG_h)} 
    \\
    &= \skpp{\calH (\phi_h^R,\mu_h^R,\bz_h^R,\bbB_h^R)}{\left( 
    \mu_h^R + \chi_\phi \sigma_h^n, \,
    \frac{\phi_h^R - \phi_h^{n-1}}{\Delta t} - 2\mu_h^R, \,
    \bz_h^R, \,
    \frac12\big( \bbB_h^R 
    + \calI_h[\kappa(\phi_h^{n-1})] \bbI 
    - \calI_h[\beta_\delta^{-1}(\bbB_h^R)]\big) 
    \right)} 
    \\
    &\geq
    C \left( \frac{1}{\Delta t} - \frac{1}{\Delta t_*} \right)
    \left( \norm{\phi_h^R}_h^2
    + \norm{\nabla\phi_h^R}_{L^2}^2
    + \norm{\bbB_h^R}_h^2 \right)
    + \frac{1}{\Delta t} \left(
    \norm{\nabla\phi_h^R - \nabla\phi_h^{n-1}}_{L^2}^2 
    + \norm{\bbB_h^R - \bbB_h^{n-1}}_h^2 \right)
    \\
    %
    %
    &\quad 
    + C \left( \norm{\mu_h^R}_h^2
    + \norm{\nabla\mu_h^R}^2_{L^2}
    + \norm{\D(\bz_h^R)}_{L^2}^2
    + \norm{\nabla\bbB_h^R}_{L^2}^2 
    + \norm{\nabla\calI_h \trace \ln \beta_\delta(\bbB_h^R)}_{L^2}^2
    \right)
    \\
    &\quad 
    + C \norm{\bbT_{e}^\delta(\phi_h^{n-1},\bbB_h^R) \beta_\delta^{-1/2}(\bbB_h^R)}_h^2
    - C(\phi_h^{n-1}, \bbB_h^{n-1}, \bu_h^n, \sigma_h^n)
    \\
    &\geq C \left( \frac{1}{\Delta t} - \frac{1}{\Delta t_*} \right) 
    ||| (\phi_h^R, \mu_h^R, \bz_h^R, \bbB_h^R) |||^2 
    - C(\phi_h^{n-1}, \bbB_h^{n-1}, \bu_h^n, \sigma_h^n).
\end{align*}
Here, we have used the same assumptions as in the stability proof such that the coefficient of the first term is non-negative, i.e., $\frac\beta\varepsilon > C_\kappa^2 C_\psi^{-1}$, see \ref{A4}, and $\Delta t < \Delta t_*$,
where the constant $\Delta t_*>0$ only depends on the model parameters and $c_\infty = \max_{n\in\{1,...,N_T\}} \norm{\sigma_{\infty,h}^n}_{L^2(\partial\Omega)}$.
Due to norm equivalence in finite dimensions, we have for a constant $C(h)$ depending on $h$, $\phi_h^{n-1}$, $\bbB_h^{n-1}$, $\bu_h^n$, $\sigma_h^{n-1}$, $\sigma_h^n$ but not on $\phi_h^R$, $\mu_h^R$, $\bz_h^R$, $\bbB_h^R$, $\delta$, $\Delta t$, such that
\begin{align*}
    &\skpp{\calH (\phi_h^R,\mu_h^R,\bz_h^R,\bbB_h^R)}{(\zeta_h,\rho_h,\bw_h,\bbG_h)} 
    \geq C(h) \left(  ||| f(\phi_h^R, \mu_h^R, \bz_h^R ,\bbB_h^R)|||^2 - 1 \right)
    = C(h) \left( R^2 - 1 \right),
\end{align*}
which is greater than zero, if $R>0$ is large enough.

On the other side, as $f(\phi_h^R, \mu_h^R,\bv_h^R, \bbB_h^R)$ is a fixed-point of $\calG_R$, we have
\begin{align*}
    &\skpp{\calH (\phi_h^R,\mu_h^R,\bz_h^R,\bbB_h^R)}{(\zeta_h,\rho_h,\bw_h,\bbG_h)}
    \\
    &= - \frac{||| \calH^h (\phi_h^R,\mu_h^R, \bz_h^R, \bbB_h^R)  ||| }{R} 
    \skpp{(\calG_R \circ f)(\phi_h^R,\mu_h^R, \bz_h^R, \bbB_h^R)}{(\zeta_h,\rho_h,\bw_h,\bbG_h)}
    \\
    &= - \frac{||| \calH^h (\phi_h^R,\mu_h^R, \bz_h^R, \bbB_h^R)  ||| }{R} 
    \skpp{f(\phi_h^R,\mu_h^R, \bz_h^R, \bbB_h^R)}{(\zeta_h,\rho_h,\bw_h,\bbG_h)}.
\end{align*}
It suffices to show
\begin{align*}
    \skpp{f(\phi_h^R,\mu_h^R, \bz_h^R, \bbB_h^R)}{(\zeta_h,\rho_h,\bw_h,\bbG_h)}
    > 0,
\end{align*}
for $R>0$ large enough. Here we use \eqref{eq:reg_e}, \eqref{eq:existence_R} and Hölder's and Young's inequalities. It holds
\begin{align*}
    & \skpp{f(\phi_h^R,\mu_h^R, \bz_h^R, \bbB_h^R)}{(\zeta_h,\rho_h,\bw_h,\bbG_h)}
    \\
    &= 
    \skpp{\left(\mu_h^R, \frac{\phi_h^R}{\Delta t} - 2 \mu_h^R, \bz_h^R, \bbB_h^R  \right)}{ \left( 
    \mu_h^R + \chi_\phi \sigma_h^n, \,
    \frac{\phi_h^R - \phi_h^{n-1}}{\Delta t} - 2\mu_h^R, \,
    \bz_h^R, \,
    \frac12\big( \bbB_h^R 
    + \calI_h[\kappa(\phi_h^{n-1})] \bbI 
    - \calI_h[\beta_\delta^{-1}(\bbB_h^R)]\big) 
    \right) }
    \\
    &= 
    \norm{\mu_h^R}_h^2 
    + \skp{\mu_h^R}{\chi_\phi \sigma_h^n}_h
    + \nnorm{ \frac{\phi_h^R}{\Delta t} - 2 \mu_h^R }_h^2
    - \skp{\frac{\phi_h^R}{\Delta t}-2 \mu_h^R}{ \frac{\phi_h^{n-1}}{\Delta t}} 
    + \norm{\bz_h^R}_{L^2}^2
    \\
    &\quad+ \frac12 \skp{\bbB_h^R}{\bbB_h^R 
    + \kappa(\phi_h^{n-1}) \bbI 
    - \beta_\delta^{-1}(\bbB_h^R)}_h
    \\
    &\geq \frac12 \norm{\mu_h^R}_h^2 
    + \frac12 \nnorm{ \frac{\phi_h^R}{\Delta t} - 2 \mu_h^R }_h^2
    + \norm{\bz_h^R}_{L^2}^2
    + \frac18 \norm{\bbB_h^R}_h^2 
    - C(\sigma_h^n, \phi_h^{n-1})
    \\
    & \geq \frac18 ||| f(\phi_h^R, \mu_h^R, \bz_h^R, \bbB_h^R) |||^2 
    - C(\sigma_h^n, \phi_h^{n-1})
    \\
    &= \frac18 R^2 - C(\sigma_h^n, \phi_h^{n-1}),
\end{align*}
which is greater than zero for $R>0$ large enough, and hence,
\begin{align*}
    &\skpp{\calH (\phi_h^R,\mu_h^R,\bz_h^R,\bbB_h^R)}{(\zeta_h,\rho_h,\bw_h,\bbG_h)} < 0,
\end{align*}
for $R>0$ large enough. Summarized, we have shown 
\begin{align*}
    0 < \skpp{\calH (\phi_h^R,\mu_h^R,\bz_h^R,\bbB_h^R)}{(\zeta_h,\rho_h,\bw_h,\bbG_h)} < 0,
\end{align*}
for $R>0$ large enough, which yields a contradition. Hence, supposed that $R>0$ is large enough, the assumption that $\calH$ has no zero in $\calB_R$ has been proved wrong. Therefore, there exists a zero of $\calH$ in $\calB_R$ which corresponds to a solution $(\phi_h^n, \mu_h^n, \bz_h^n, \bbB_h^n)$ of the reduced problem \eqref{eq:phi_reduced}--\eqref{eq:B_reduced}.
In addition, the stability estimates \eqref{eq:energy_estimate_dissipative3}--\eqref{eq:energy_estimate_dissipative4} hold true for $(\phi_h^n, \mu_h^n, \bz_h^n, \bbB_h^n)$, as we used the same assumptions for that proof.

\subsubsection{Reconstruction of the pressure}
From now on, we set $\bv_h^n\coloneqq \bu_h^n + \bz_h^n$. Note that \eqref{eq:energy_estimate_dissipative4} and \eqref{eq:bound_uh} imply
\begin{align}
    \label{eq:bound_v}
    \Delta t \sum_{n=1}^{N_T} \norm{\bv_h^n}_{H^1}^2 
    \leq C(T,c_\infty) \big( 1 + \abs{\calF_h^\delta(\phi_h^0, \bbB_h^0)} \big).
\end{align}
Now we reconstruct the pressure $p_h^n \in \calS_h$ and hence justify the solvability of \eqref{eq:phi_reg}--\eqref{eq:B_reg}. We define a linear functional $\calF_h \colon\, \calV_h \to \bbR$ by
\begin{align*}
    \calF_h(\bw_h) 
    &\coloneqq
    \skp{2 \calI_h[\eta(\phi_h^{n-1})] \D(\bv_h^n)}{\D(\bw_h)}_{L^2}
    + \skp{\calI_h \bbT_{e,h}^{\delta,n}}{\nabla\bw_h}_{L^2}
    - \skp{(\mu_h^n + \chi_\phi \sigma_h^n) \nabla\phi_h^{n-1}}{\bw_h}_{L^2}
    \\
    &\quad
    - \frac12 \skp{\nabla \calI_h\left[ \kappa(\phi_h^{n-1}) \trace\beta_\delta(\bbB_h^n) \right]}{\bw_h}_{L^2}
    + \frac12 \sum_{i,j=1}^d \skp{\partial_{x_j} \calI_h\left[ \kappa(\phi_h^{n-1}) \right]  \trace \mathbf{\Lambda}^\delta_{i,j}(\bbB_h^n)}{(\bw_h)_i}_{L^2} ,
\end{align*}
for any $\bw_h\in\calV_h$. As \eqref{eq:v_reduced} is fulfilled for all $\bw_h \in \calV_{h,\mathrm{div}}$, it follows that the linear functional $\calF_h$ vanishes on $\calV_{h,\mathrm{div}}$, i.e., $\calF_h(\bw_h) = 0$ for all $\bw_h \in \calV_{h,\mathrm{div}}$. With \eqref{eq:LBB_p}--\eqref{eq:LBB_p_bound}, we get the existence of a unique pressure $p_h^n \in \calS_h$ satisfying 
\begin{align*}
    \skp{p_h^n}{\Div\bw_h}_{L^2} &= \calF_h(\bw_h) \quad \forall\, \bw_h\in\calV_h,
    \\
    \norm{p_h^n}_{L^2} &\leq C 
    \sup_{\bw_h\in\calV_h \setminus\{0\}} \frac{\calF_h(\bw_h)}{\norm{\bw_h}_{H^1}}.
\end{align*}
This justifies the existence of a solution $(\phi_h^{n},\mu_h^n,\sigma_{h}^{n},p_h^n,\bv_{h}^{n},\bbB_{h}^{n}) \in (\calS_h)^4 \times \calV_h \times \calW_h$ of \eqref{eq:phi_reg}--\eqref{eq:B_reg} with $\bv_h^n \coloneqq \bz_h^n + \bu_h^n$.

Now, we derive a uniform estimate for the pressure. 
Using \ref{A2} and Hölder's inequality, we have
\begin{align*}
    \abs{\calF_h(\bw_h)}
    &\leq C \Big(
    \norm{\bv_h^n}_{H^1} \norm{\bw_h}_{H^1}
    + \norm{\calI_h \bbT_{e,h}^{\delta,n}}_{L^2} 
    \norm{\bw_h}_{H^1} 
    + \norm{\mu_h^n}_{L^3} \norm{\nabla\phi_h^{n-1}}_{L^2} \norm{\bw_h}_{L^6}
    + \norm{\sigma_h^n}_{L^3} \norm{\nabla\phi_h^{n-1}}_{L^2} \norm{\bw_h}_{L^6} \Big)
    \\
    &\quad
    + \frac12 \abs{ \skp{\nabla \calI_h\left[ \kappa(\phi_h^{n-1}) \trace\beta_\delta(\bbB_h^n) \right]}{\bw_h}_{L^2}}
    + \frac12 \abs{\sum_{i,j=1}^d \skp{\partial_{x_j} \calI_h\left[ \kappa(\phi_h^{n-1}) \right]  \trace \mathbf{\Lambda}^\delta_{i,j}(\bbB_h^n)}{(\bw_h)_i}_{L^2} }.
\end{align*}
Recalling \eqref{eq:bound_Te_delta}, we have
\begin{align*}
    \norm{\calI_h \bbT_{e,h}^{\delta,n}}_{L^2} 
    &\leq C \left( 1 + \left( 1 + \norm{\bbB_h^n}_{L^2}^{1/2} \norm{\bbB_h^n}_{H^1}^{1/2} \right)
    \left( \norm{\bbB_h^n}_{H^1} 
    + \norm{\phi_h^{n-1}}_{H^1} + 1 \right) \right).
\end{align*}
Also, we have with \eqref{eq:askp}, \eqref{eq:bound_k_h} and \eqref{eq:bound_grad_k_h}, that
\begin{align*}
    \abs{ \skp{\nabla \calI_h\left[ \kappa(\phi_h^{n-1}) \trace\beta_\delta(\bbB_h^n) \right]}{\bw_h}_{L^2}}
    &\leq
    C \norm{\bw_h}_{H^1} 
    \norm{\calI_h \beta_\delta(\bbB_h^n)}_{H^1}
    \norm{\calI_h \kappa(\phi_h^{n-1})}_{H^1} 
    \\
    &\leq 
    C \norm{\bw_h}_{H^1} 
    (1 + \norm{\bbB_h^n}_{H^1} )
    ( 1 + \norm{\phi_h^{n-1}}_{H^1} ).
\end{align*}
Moreover, it holds with \eqref{eq:bound_grad_k_h}, \eqref{eq:Lambda_delta}, \eqref{eq:bound_beta_delta}, \eqref{eq:inverse_estimate} and Hölder's inequality, that
\begin{align*}
    \abs{\sum_{i,j=1}^d \skp{\partial_{x_j} \calI_h\left[ \kappa(\phi_h^{n-1}) \right]  \trace \mathbf{\Lambda}^\delta_{i,j}(\bbB_h^n)}{(\bw_h)_i}_{L^2} }
    &\leq 
    C \norm{\nabla\calI_h \kappa(\phi_h^{n-1})}_{L^2} 
    \max_{i,j} \, \norm{\mathbf{\Lambda}_{i,j}^\delta(\bbB_h^n)}_{L^3}
    \norm{\bw_h}_{L^6}
    \\
    &\leq 
    C \norm{\nabla\phi_h^{n-1}}_{L^2} (1 + \norm{\bbB_h^n}_{L^3}) \norm{\bw_h}_{H^1}.
\end{align*}
Combining these estimates with \eqref{eq:LBB_p_bound} with the Sobolev embedding $H^1(\Omega) \hookrightarrow L^6(\Omega)$ for dimension $d\in\{2,3\}$, we obtain
\begin{align*}
    \norm{p_h^n}_{L^2}
    &\leq C \Big(
    1 + \norm{\bv_h^n}_{H^1} 
    + \norm{\bbB_h^n}_{L^2}^{1/2} \norm{\bbB_h^n}_{H^1}^{3/2}
    + \norm{\bbB_h^n}_{L^2}^{1/2} \norm{\bbB_h^n}_{H^1}^{1/2}
    \norm{\phi_h^{n-1}}_{H^1}
    + \norm{\bbB_h^n}_{H^1}
    + \norm{\phi_h^{n-1}}_{H^1} 
    \Big) 
    \\
    &\quad
    + C \Big( 
    \norm{\mu_h^n}_{H^1} \norm{\phi_h^{n-1}}_{H^1} 
    + \norm{\sigma_h^n}_{H^1} \norm{\phi_h^{n-1}}_{H^1} 
    + \norm{\bbB_h^n}_{H^1} \norm{\phi_h^{n-1}}_{H^1}
    \Big) .
\end{align*}
Taking the power of $4/3$ on both sides, multiplying both sides by $\Delta t$, summing from $n=1,...,N_T$ and using Hölder's inequality, \eqref{eq:bound_sigma}, \eqref{eq:bound_uh}, \eqref{eq:energy_estimate_dissipative3} and Young's inequality, we get
\begin{align*}
    \Delta t \sum_{n=1}^{N_T} \norm{p_h^n}_{L^2}^{4/3}
    &\leq 
    CT + C T^{1/3} \left( \Delta t \sum_{n=1}^{N_T} \Big( 
    \norm{\bv_h^n}_{H^1}^2  
    + \norm{\bbB_h^n}_{H^1}^2 + \norm{\phi_h^{n-1}}_{H^1}^2\Big) \right)^{2/3}
    \\
    &\quad
    + C 
    \left( \max_{n\in\{1,...,N_T\}} \norm{\bbB_h^n}_{L^2}^2 \right)^{1/3}
    \left( \Delta t \sum_{n=1}^{N_T} 
    \norm{\bbB_h^n}_{H^1}^2 \right)
    \\
    &\quad 
    + C \left( \max_{n\in\{1,...,N_T\}} \norm{\bbB_h^n}_{L^2}^2 \right)^{1/3}
    \left( \Delta t \sum_{n=1}^{N_T} \norm{\bbB_h^n}_{H^1}^2 \right)^{1/3} 
    \left( \Delta t \sum_{n=1}^{N_T} \norm{\phi_h^{n-1}}_{H^1}^2 \right)^{2/3} 
    \\
    &\quad
    + C T^{1/3} \left( \max_{n\in\{1,...,N_T\}} \norm{\phi_h^{n-1}}_{H^1}^2 \right)^{2/3} 
    \left( \Delta t \sum_{n=1}^{N_T} \Big( \norm{\mu_h^n}_{H^1}^2 + \norm{\sigma_h^n}_{H^1}^2 + \norm{\bbB_h^n}_{H^1}^2 \Big) \right)^{2/3}
    \\
    &\leq 
    C(T,c_\infty) \big(1 + \abs{\calF_h^\delta(\phi_h^0, \bbB_h^0)}^{4/3} \big),
\end{align*}
for a constant $C(T,c_\infty)>0$ that depends exponentially on $T$ and $c_\infty = \max_{n\in\{1,...,N_T\}} \norm{\sigma_{\infty,h}^n}_{L^2(\partial\Omega)}$.
This implies
\begin{align}
    \label{eq:bound_p}
    \bigg( \Delta t \sum_{n=1}^{N_T} \norm{p_h^n}_{L^2}^{4/3} \bigg)^{3/4}
    \leq C(T,c_\infty) \big(1 + \abs{\calF_h^\delta(\phi_h^0, \bbB_h^0)} \big).
\end{align}

So far, we have shown the existence of solutions for the regularized discrete system \eqref{eq:phi_reg}--\eqref{eq:B_reg}, which satisfy the stability estimates \eqref{eq:energy_estimate_dissipative3}, \eqref{eq:energy_estimate_dissipative4}, \eqref{eq:bound_v}, \eqref{eq:bound_p}. We summarize this in the following lemma.

\begin{lemma}
Let $\delta\in(0,1)$, and let \ref{A1}--\ref{A6} and \ref{S} hold true. Let the discrete initial data $\phi_h^0\in\calS_h$ and $\bbB_h^0\in\calW_h$ be given. Moreover, for any $n\in\{1,...,N_T\}$, let the discrete boundary data $\sigma_{\infty,h}^n \in\calS_h$ be given.
Besides, assume that $\Delta t < \Delta t_*$,
where the constant $\Delta t_*>0$ depends only on the model parameters and on $c_\infty \coloneqq \max_{n\in\{1,...,N_T\}} \norm{\sigma_{\infty,h}^n}_{L^2(\partial\Omega)}$. 
Then, for all $n\in\{1,...,N_T\}$, there exists at least one solution tuple $(\phi_h^{n},\mu_h^n,\sigma_{h}^{n},p_h^n,\bv_{h}^{n},\bbB_{h}^{n}) \in (\calS_h)^4 \times \calV_h \times \calW_h$ to the regularized discrete problem \ref{P_delta_FE}.
Moreover, all solutions of \ref{P_delta_FE} are stable in the sense that
\begin{align}
\label{eq:stability_sigma_delta}
    &\max_{n=1,...,N_T} \norm{\sigma_h^n}_{H^1}^2
    \leq C \max_{n=1,...,N_T} \norm{\sigma_{\infty,h}^n}_{L^2(\partial\Omega)}^2,
\end{align}
and
\begin{align}
    \nonumber
    \label{eq:stability_FE_delta}
    &  \max_{n=1,...,N_T} \Big( 
    \norm{\phi_h^n}_{H^1}^2
    +  \norm{ \bbB_h^n }_{L^2}^2
    +  \norm{\calI_h \trace g_\delta(\bbB_h^n)}_{L^1}
    + \frac1\delta \norm{\calI_h \abs{ [\bbB_h^n]_-}}_{L^1}
    \Big)
    \\
    \nonumber
    &\quad
    +  \sum_{n=1}^{N_T} \Big(
    \norm{\nabla\phi_h^n - \nabla\phi_h^{n-1}}_{L^2}^2 
    +  \norm{\bbB_h^n - \bbB_h^{n-1}}_{L^2}^2 \Big)
    + \Delta t \sum_{n=1}^{N_T} \Big( 
    \norm{\mu_h^n}_{H^1}^2 
    + \norm{\bv_h^n}_{H^1}^2  
    + \norm{\bbB_h^n}_{H^1}^2  \Big)
    \\
    \nonumber
    &\quad
    + \Delta t \sum_{n=1}^{N_T} \Big( 
    \norm{\nabla\calI_h \trace \ln \beta_\delta(\bbB_h^n)}_{L^2}^2
    + \norm{\bbT_{e,h}^{\delta,n} \beta_\delta^{-1/2}(\bbB_h^n)}_h^2
    \Big)
    + \left( \Delta t \sum_{n=1}^{N_T} 
    \norm{p_h^n}_{L^2}^{4/3} 
    \right)^{3/4}
    \\
    &\leq C(T, c_\infty) \big(1 + \abs{\calF_h^\delta(\phi_h^0, \bbB_h^0)} \big),
\end{align}
where the constants $C, C(T, c_\infty)>0$ are independent of $\delta, h, \Delta t$, but $C(T,c_\infty)$ depends exponentially on $T$ and $c_\infty$.
\end{lemma}

\subsubsection{Limit passing in the regularization parameter}
We are now in the position to prove Theorem \ref{theorem:existence_FE} under the respective assumptions.
In particular, we pass to the limit $\delta\to 0$ and extract a converging subsequence of solutions of the $\delta$-regularized scheme \eqref{eq:phi_reg}--\eqref{eq:B_reg}. 

For any $n\in\{1,...,N_T\}$, we now write $\bbB_{h,\delta}^n$, etc., (with subscript $\delta$) for solutions of the regularized scheme \eqref{eq:phi_reg}--\eqref{eq:B_reg}.
As $(\phi_{h,\delta}^{n}, \mu_{h,\delta}^n, \sigma_{h,\delta}^{n}, p_{h,\delta}^n, \bv_{h,\delta}^{n}, \bbB_{h,\delta}^{n})$ are finite dimensional for fixed $(\calS_h)^4 \times \calV_h \times \calW_h$, the following subsequence convergence result follows directly from \eqref{eq:stability_sigma_delta}--\eqref{eq:stability_FE_delta}:
\begin{align*}
    f_{h,\delta}^n \to f_h^n \quad \text{pointwise on } \overline\Omega,
    \quad \forall\, f\in\{ \phi, \mu, \sigma, p, \bv, \bbB\}, \ \forall\, n\in\{1,...,N_T\},
\end{align*}
as $\delta \to 0$.
It follows from the uniform stability bounds \eqref{eq:stability_FE_delta}, that $\calI_h \abs{ [\bbB_{h,\delta}^n]_- }$ vanishes on $\overline\Omega$ in the limit $\delta\to 0$, so that $\bbB_h^n$ is positive semi-definite.
As $\beta_\delta$ is Lipschitz continuous with Lipschitz constant $L_{\beta_\delta}=1$ and, as $\bbB_h^n$ is positive semi-definite, we conclude
\begin{align*}
    \abs{ \beta_\delta(\bbB_{h,\delta}^n) - \bbB_h^n }
    \leq 
    \abs{ \beta_\delta(\bbB_{h,\delta}^n) - \beta_\delta(\bbB_h^n) }
    + \abs{ \beta_\delta(\bbB_h^n) - \bbB_h^n }
    \leq 
    \abs{ \bbB_{h,\delta}^n - \bbB_h^n }
    + \abs{ \beta_\delta(\bbB_h^n) - \bbB_h^n }
    \to 0,
\end{align*}
as $\delta \to 0$. Moreover, from the stability bounds we have for some $\widetilde{\bbC}_h^n \in \calW_h$ that
\begin{align*}
    \calI_h \big[ \bbT_{e,h}^{\delta,n} \beta_\delta^{-1/2}(\bbB_{h,\delta}^n) \big] \to \widetilde{\bbC}_h^n,
\end{align*}
as $\delta \to 0$. To identify $\widetilde{\bbC}_h^n$ with $\calI_h \big[ \bbT_{e,h}^{n} \, [\bbB_{h}^n]^{-1/2} \big]$, we use that $\beta_\delta(\bbB_{h,\delta}^n) \beta_\delta^{-1}(\bbB_{h,\delta}^n) = \bbI$ and that $\beta_\delta(\bbB_{h,\delta}^n) \to \bbB_h^n$, as $\delta\to 0$. Therefore, $\beta_\delta^{-1}(\bbB_{h,\delta}^n) \to [\bbB_h^n]^{-1}$, as $\delta\to 0$, and so $\bbB_h^n$ is positive definite.
Moreover, sending $\delta\to 0$, we have $\trace g_\delta(\bbB_{h,\delta}^n) \to \trace \ln(\bbB_h^n)$ and $\calI_h \trace \ln \beta_\delta(\bbB_{h,\delta}^n) \to \calI_h \trace\ln(\bbB_h^n)$ with similar arguments, as $\bbB_h^n$ is positive definite. Using these results and the concrete definition of $\mathbf\Lambda_{i,j}^\delta$, we have $\mathbf\Lambda_{i,j}^\delta(\bbB_{h,\delta}^n) \to \mathbf\Lambda_{i,j}(\bbB_n^n)$, $i,j\in\{1,...,d\}$, as $\delta\to 0$.
These have been the most relevant terms in \eqref{eq:phi_reg}--\eqref{eq:B_reg}. So, we can now pass to the limit in the system \eqref{eq:phi_reg}--\eqref{eq:B_reg} and obtain the existence of a solution to the unregularized system.
Using Fatou's lemma and, as the discrete initial datum $\bbB_h^0\in\calW_h$ is positive definite,
we can also pass to the limit $\delta\to 0$ in \eqref{eq:stability_sigma_delta}--\eqref{eq:stability_FE_delta} to justify the stability estimates \eqref{eq:stability_FE_sigma}--\eqref{eq:stability_FE}.


%
%
%
%
%
\section{Convergence of the numerical scheme}
\label{sec:convergence}
In this section, we improve the results from Section \ref{sec:approximation} and we show that there exists a subsequence of discrete solutions that converges to a weak solution of \eqref{eq:phi}--\eqref{eq:B} in the sense of Definition \ref{def:weak_solution}. This will then prove Theorem \ref{theorem:weak_solution}.

\subsection{Higher order estimates}

First, we derive higher order estimates that hold true uniformly in $h,\Delta t$. This helps us to apply compactness results.


\begin{lemma}~\\
Let \ref{A1}--\ref{A6} and \ref{S} hold true. Let the discrete initial data $\phi_h^0\in\calS_h$ and $\bbB_h^0\in\calW_h$ be given with $\bbB_h^0$ being  positive definite. Moreover, for any $n\in\{1,...,N_T\}$, let the discrete boundary data $\sigma_{\infty,h}^n \in\calS_h$ be given. 
Besides, assume that $\Delta t < \Delta t_*$,
where the constant $\Delta t_*>0$ depends only on the model parameters and on $c_\infty \coloneqq \max_{n\in\{1,...,N_T\}} \norm{\sigma_{\infty,h}^n}_{L^2(\partial\Omega)}$.
Then, in addition to \eqref{eq:stability_FE_sigma}--\eqref{eq:stability_FE}, all solutions of \ref{P_FE} fulfill
\begin{align}
    \nonumber
    \label{eq:higher_order1}
    & \left(\Delta t \sum\limits_{n=1}^{N_T} 
    \nnorm{\frac{\phi_h^n - \phi_h^{n-1}}{\Delta t}  }_{(H^1)'}^2 \right)^{1/2}
    + \left( \Delta t \sum\limits_{n=1}^{N_T} \nnorm{\frac{\bbB_h^n - \bbB_h^{n-1}}{\Delta t}  }_{(H^1)'}^{4/3} \right)^{3/4}
    + \left( \Delta t \sum\limits_{n=1}^{N_T} \nnorm{(\bbB_h^n)^{-1}}_{L^1} \right)^{2/3}
    \\
    &\leq  C(T, c_\infty) \big(
    1 + \norm{\bbB_h^0}_{L^2}^2 + \norm{\calI_h\trace\ln\bbB_h^0}_{L^1}
    + \norm{\nabla\phi_h^0}_{L^2}^2
    + \norm{\calI_h \psi(\phi_h^0)}_{L^1} \big) ,
\end{align}
where the constant $C(T, c_\infty)>0$ is independent of $h,\Delta t$, but depends exponentially on $T$ and $c_\infty$.
\end{lemma}

\begin{proof} 
We prove the estimate for the second summand of \eqref{eq:higher_order1}. The estimate for the first one follows with similar arguments. After that, we show the estimate for the third summand. Let $\bbG\in H^1(\Omega;\bbR^{d\times d}_\mathrm{S})$. Then, with $\bbG_h \coloneqq \calQ_h \bbG \in \calW_h$ in \eqref{eq:B_FE}, where $\calQ_h$ is the \textit{lumped} $L^2$ projector defined in \eqref{eq:proj_def}, it holds
\begin{align*}
    \frac{1}{\Delta t} \skp{\bbB_h^n - \bbB_h^{n-1}}{\bbG}_{L^2}
    &=
    \frac{1}{\Delta t} \skp{\bbB_h^n - \bbB_h^{n-1}}{\bbG_h}_h
    \\
    &=
    - \skp{\bv_h^n}{ \nabla\calI_h\left[\bbB_h^n : \bbG_h \right]}_{L^2}
    + \sum_{i,j=1}^d \skp{(\bv_h^n)_i \mathbf{\Lambda}_{i,j}(\bbB_h^n)}{\partial_{x_j} \bbG_h}_{L^2}
    \\
    &\quad 
    - \skp{\tfrac{1}{\tau(\phi_h^{n-1})} \bbT_{e,h}^n }{\bbG_h}_h 
    - \skp{\Gamma_{\bbB,h}^{n-1} \bbB_h^n}{\bbG_h}_h
    + \skp{2 \nabla\bv_h^n}{\calI_h\left[ \bbG_h \bbB_h^n \right]}_{L^2}
    - \alpha \skp{\nabla\bbB_h^n}{\nabla\bbG_h}_{L^2}.
\end{align*}
Using \eqref{eq:askp}, \eqref{eq:Lambda}, \eqref{eq:inverse_estimate}, \eqref{eq:bound_Te}, \eqref{eq:bound_kh'_tau}, \eqref{eq:bound_source}, \eqref{eq:Gagliardo}, \eqref{eq:norm_equiv}, Hölder's inequality and \eqref{eq:bound_sigma}, we get
\begin{align*}
    &\abs{\frac{1}{\Delta t} \skp{\bbB_h^n - \bbB_h^{n-1}}{\bbG}_{L^2}}
    \\
    &\leq 
    C \norm{\bv_h^n}_{H^1} \norm{\bbB_h^n}_{L^2}^{1/2} \norm{\bbB_h^n}_{H^1}^{1/2} \norm{\bbG_h}_{H^1}
    + C \left( 1 + \left( \norm{\bbB_h^n}_{L^2}^{1/2} \norm{\bbB_h^n}_{H^1}^{1/2} \right)
    \left( \norm{\bbB_h^n}_{H^1} 
    + \norm{\phi_h^{n-1}}_{H^1} + 1 \right) \right) \norm{\bbG_h}_{L^2}
    \\
    &\quad
    + C (1 + \norm{\sigma_h^n}_{L^6}) \norm{\bbB_h^n}_{L^2} \norm{\bbG_h}_{L^3}
    + C \norm{\nabla \bv_h^n}_{L^2} \norm{\bbB_h^n}_{L^3}   \norm{\bbG_h}_{L^6}
    + C \norm{\bbB_h^n}_{H^1} \norm{\bbG_h}_{H^1}
    \\
    &\leq
    C \norm{\bbG_h}_{H^1}
    + C \big( \norm{\bv_h^n}_{H^1} + \norm{\bbB_h^n}_{H^1} + \norm{\phi_h^{n-1}}_{H^1} \big)
    \norm{\bbB_h^n}_{L^2}^{1/2} \norm{\bbB_h^n}_{H^1}^{1/2} \norm{\bbG_h}_{H^1}
    \\
    &\quad
    + C(c_\infty) \norm{\bbB_h^n}_{H^1} \norm{\bbG_h}_{H^1}
    + C \norm{\phi_h^{n-1}}_{H^1} \norm{\bbG_h}_{H^1},
\end{align*}
which gives
\begin{align*}
    \nnorm{\frac{\bbB_h^n - \bbB_h^{n-1}}{\Delta t}}_{(H^1)'}
    &\leq
    C 
    + C \big( \norm{\bv_h^n}_{H^1} + \norm{\bbB_h^n}_{H^1} + \norm{\phi_h^{n-1}}_{H^1} \big)
    \norm{\bbB_h^n}_{L^2}^{1/2} \norm{\bbB_h^n}_{H^1}^{1/2} 
    + C(c_\infty) \norm{\bbB_h^n}_{H^1} 
    + C \norm{\phi_h^{n-1}}_{H^1} .
\end{align*}
Taking the power of $4/3$ on both sides, multiplying both sides by $\Delta t$ and summing from $n=1,...,N_T$, using Hölder's and Young's inequalities and \eqref{eq:stability_FE} yields
\begin{align*}
    &\Delta t \sum_{n=1}^{N_T} \nnorm{\frac{\bbB_h^n - \bbB_h^{n-1}}{\Delta t}}_{(H^1)'}^{4/3}
    \\
    &\leq C T
    +
    C \left( \Delta t \sum_{n=1}^{N_T} \big( \norm{\bv_h^n}_{H^1}^2 + \norm{\bbB_h^n}_{H^1}^2 + \norm{\phi_h^{n-1}}_{H^1}^2 \big)
    \right)^{2/3}
    \left( \max_{n\in\{1,...,N_T\}} \norm{\bbB_h^n}_{L^2}^2
    \right)^{1/3}
    \left( \Delta t \sum_{n=1}^{N_T}  \norm{\bbB_h^n}_{H^1}^2 
    \right)^{1/3}
    \\
    &\quad
    + C(c_\infty) T^{1/3} \left( \Delta t \sum_{n=1}^{N_T}  \norm{\bbB_h^n}_{H^1}^2 
    \right)^{2/3}
    + C T^{1/3} \left( \Delta t \sum_{n=1}^{N_T}  \norm{\phi_h^{n-1}}_{H^1}^2 
    \right)^{2/3}
    \\
    &\leq C(T, c_\infty) \big(
    1 + \norm{\bbB_h^0}_{L^2}^2 + \norm{\calI_h\trace\ln\bbB_h^0}_{L^1}
    + \norm{\nabla\phi_h^0}_{L^2}^2
    + \norm{\calI_h \psi(\phi_h^0)}_{L^1} \big)^{4/3}.
\end{align*}
Taking the power of $3/4$ on both sides shows the estimate of the second summand in \eqref{eq:higher_order1}.

We now show the estimate for the third summand in \eqref{eq:higher_order1}. On noting \cite[(6.39)]{barrett_boyaval_2009} and the positive definiteness of $\bbB_h^n$, it holds
\begin{align*}
    \norm{(\bbB_h^n)^{-1}}_{L^1} 
    \leq C \nnorm{\calI_h \abs{(\bbB_h^n)^{-1}} }_{L^1} 
    \leq C \norm{\calI_h \trace ((\bbB_h^n)^{-1})}_{L^1} .
\end{align*}
Combining this with the definition of $\bbT_{e,h}^n$, Hölder's and Young's inequalities, \eqref{eq:norm_equiv}, \eqref{eq:bound_k_h} and \eqref{eq:Gagliardo}, we have
\begin{align*}
    \norm{(\bbB_h^n)^{-1}}_{L^1} 
    &\leq 
    C \norm{\bbT_{e,h}^n \, (\bbB_h^n)^{-1/2}}_h^2
    + C \big( 1 + \norm{\bbB_h^n}_{L^3}^3 + \norm{\phi_h^{n-1}}_{L^3}^3 \big)
    \\
    &\leq C \norm{\bbT_{e,h}^n \, (\bbB_h^n)^{-1/2}}_h^2
    + C \big( 1 + \norm{\bbB_h^n}_{L^2}^{3/2} \norm{\bbB_h^n}_{H^1}^{3/2} 
    + \norm{\phi_h^{n-1}}_{H^1}^3 \big).
\end{align*}
Multiplying both sides by $\Delta t$, summing over all $n\in\{1,...,N_T\}$ and noting Hölder's and Young's inequalities and \eqref{eq:stability_FE}, we deduce
\begin{align*}
    \Delta t \sum_{n=1}^{N_T} \norm{(\bbB_h^n)^{-1}}_{L^1} 
    &\leq C \Delta t \sum_{n=1}^{N_T}
    \norm{\bbT_{e,h}^n \, (\bbB_h^n)^{-1/2}}_h^2
    + CT  + C T^{1/4} \Big( \max_{n=1,...,N_T} \norm{\bbB_h^n}_{L^2}^2 \Big)^{3/4} \Big( \Delta t \sum_{n=1}^{N_T}\norm{\bbB_h^n}_{H^1}^2 \Big)^{3/4} 
    \\
    &\quad 
    + CT \Big( \max_{n=1,...,N_T} \norm{\phi_h^{n-1}}_{H^1}^2 \Big)^{3/2}
    \\
    & \leq C(T, c_\infty) \big(
    1 + \norm{\bbB_h^0}_{L^2}^2 + \norm{\calI_h\trace\ln\bbB_h^0}_{L^1}
    + \norm{\nabla\phi_h^0}_{L^2}^2
    + \norm{\calI_h \psi(\phi_h^0)}_{L^1} \big)^{3/2}.
\end{align*}
Taking the power of $3/2$ on both sides finishes the proof.

\end{proof}

\subsection{Compactness and subsequence convergence}

For future reference, we recall the following compactness results from \cite[Sec.~8]{simon_1986}. 
Let $X, Y, Z$ be Banach spaces with a compact embedding $X \hookrightarrow \hookrightarrow Y$ and a continuous embedding $Y \hookrightarrow Z$. Let $1\leq p < \infty$ and $r>1$. Then, the following embedding are compact:
\begin{subequations}
\begin{alignat}{3}
    \label{eq:compact_Lp}
    &\{\eta \in L^p(0,T;X) 
    &&\mid \partial_t\eta \in L^1(0,T;Z) \}
    &&\hookrightarrow \hookrightarrow L^p(0,T;Y),
    \\
    \label{eq:compact_C}
    &\{\eta \in L^\infty(0,T;X)  
    &&\mid \partial_t\eta \in L^r(0,T;Z) \}
    &&\hookrightarrow \hookrightarrow  C([0,T];Y).
\end{alignat}


\end{subequations}

We introduce the following notation for affine and piecewise constant extensions of time discrete functions $a^n(\cdot)$, $n=0,...,N_T$:
\begin{subequations}
\begin{alignat}{2}
    \label{def:fun_Delta_t}
    a^{\Delta t}(\cdot, t) 
    &\coloneqq 
    \frac{t - t^{n-1}}{\Delta t} a^n(\cdot)
    + \frac{t^n - t}{\Delta t} a^{n-1}(\cdot)
    \quad\quad 
    && \forall\, t\in [t^{n-1},t^n], \ n\in \{1,...,N_T\},
    \\
    \label{def:fun_Delta_t_pm}
    a^{\Delta t,+}(\cdot, t) 
    &\coloneqq a^n(\cdot),
    \quad\quad 
    a^{\Delta t,-}(\cdot, t) 
    \coloneqq a^{n-1}(\cdot)
    \quad\quad 
    && \forall\, t\in (t^{n-1},t^n], \ n\in \{1,...,N_T\}.
\end{alignat}
\end{subequations}
In the following, we write $a^{\Delta t,\pm}$ for results that hold true for both $a^{\Delta t,+}$ and $a^{\Delta t,-}$, and we write $a^{\Delta t(,\pm)}$ for results that hold true for $a^{\Delta t}$, $a^{\Delta t,+}$ and $a^{\Delta t,-}$, respectively.

Using this notation, we reformulate the problem \ref{P_FE} continuously in time. Multiplying each equation by $\Delta t$ and summing from $n=1,...,N_T$, we obtain for any test functions $(\zeta_h$, $\rho_h$, $ \xi_h$, $ q_h$, $ \bw_h$, $ \bbG_h) \in (L^2(0,T;\calS_h))^4 \times L^2(0,T;\calV_h) \times L^2(0,T;\calW_h)$ that
\begin{subequations}
\begin{align}
    \label{eq:phi_FE_time}\nonumber
    0 &= 
    \int_0^T  \skp{ \partial_t \phi_h^{\Delta t}}{\zeta_h}_h \dt
    + \int_0^T \skp{\calI_h[m_{\phi,h}^{\Delta t,-}] \nabla\mu_h^{\Delta t,+}}{\nabla\zeta_h}_{L^2} \dt
    \\
    &\quad
    + \int_0^T \skp{\bv_h^{\Delta t,+} \cdot \nabla\phi_h^{\Delta t,-}}{\zeta_h}_{L^2} \dt
    + \int_0^T \skp{\phi_h^{\Delta t,-} \Gamma_{\bv,h}^{\Delta t,-} - \Gamma_{\phi,h}^{\Delta t,-}}{\zeta_h}_h \dt,
    \\[5pt]
    \label{eq:mu_FE_time}\nonumber
    0 &= \int_0^T
    \skp{- \mu_h^{\Delta t,+} + \frac\beta\varepsilon \psi_h'(\phi_h^{\Delta t,+},\phi_h^{\Delta t,-}) - \chi_\phi \sigma_h^{\Delta t,+} + \frac12 \kappa_h'(\phi_h^{\Delta t,+}, \phi_h^{\Delta t,-}) \trace\bbB_h^{\Delta t,+}} {\rho_h}_h \dt
    \\
    &\quad
    + \int_0^T \beta\varepsilon \skp{\nabla\phi_h^{\Delta t,+}}{\nabla\rho_h}_{L^2} \dt,
    \\[5pt] 
    \label{eq:sigma_FE_time}
    0 &= 
    \int_0^T\skp{\calI_h[m_{\sigma,h}^{\Delta t,-}] \nabla\sigma_h^{\Delta t,+}}{\nabla\xi_h}_{L^2} \dt
    + \int_0^T \skp{\sigma_h^{\Delta t,+} \Gamma_{\sigma,h}^{\Delta t,-}}{\xi_h}_h \dt
    + \int_0^T K \skp{\sigma_h^{\Delta t,+} - \sigma_{\infty,h}^{\Delta t,+}}{\xi_h}_{L^2(\partial\Omega)} \dt,
    \\[5pt]
    \label{eq:div_FE_time}
    0 &= \int_0^T \skp{\Div \bv_h^{\Delta t,+}}{q_h}_{L^2} \dt 
    - \int_0^T \skp{\Gamma_{\bv,h}^{\Delta t,-}}{q_h}_h \dt,
    \\[5pt]
    \nonumber
    \label{eq:v_FE_time}
    0 &= \int_0^T
    \skp{2 \calI_h[\eta(\phi_h^{\Delta t,-})] \D(\bv_h^{\Delta t,+})}{\D(\bw_h)}_{L^2} \dt
    - \int_0^T \skp{p_h^{\Delta t,+}}{\Div \bw_h}_{L^2} \dt
    + \int_0^T \skp{\calI_h \bbT_{e,h}^{\Delta t,+}}{\nabla\bw_h}_{L^2} \dt
    \\
    \nonumber
    &\quad 
    - \int_0^T \skp{(\mu_h^{\Delta t,+} + \chi_\phi \sigma_h^{\Delta t,+}) \nabla\phi_h^{\Delta t,-}}{\bw_h}_{L^2} \dt 
    - \int_0^T \frac12 \skp{\nabla \calI_h\left[ \kappa(\phi_h^{\Delta t,-}) \trace\bbB_h^{\Delta t,+} \right]}{\bw_h}_{L^2} \dt
    \\
    &\quad
    + \int_0^T \frac12 \sum_{i,j=1}^d \skp{\partial_{x_j} \calI_h\left[ \kappa(\phi_h^{\Delta t,-}) \right]  \trace \mathbf{\Lambda}_{i,j}(\bbB_h^{\Delta t,+})}{(\bw_h)_i}_{L^2} \dt,
    \\[5pt]
    \nonumber
    \label{eq:B_FE_time}
    0 &= \int_0^T
    \skp{ \partial_t \bbB_h^{\Delta t}}{\bbG_h}_h \dt
    + \int_0^T \skp{\bv_h^{\Delta t,+}}{ \nabla\calI_h\left[\bbB_h^{\Delta t,+} : \bbG_h \right]}_{L^2} \dt
    \\
    \nonumber
    &\quad
    - \int_0^T \sum_{i,j=1}^d \skp{(\bv_h^{\Delta t,+})_i \, \mathbf{\Lambda}_{i,j}(\bbB_h^{\Delta t,+})}{\partial_{x_j} \bbG_h}_{L^2} \dt
    + \int_0^T \skp{\frac{1}{\tau(\phi_h^{\Delta t,-})} \bbT_{e,h}^{\Delta t,+} }{\bbG_h}_h \dt
    \\
    &\quad 
    + \int_0^T \skp{\Gamma_{\bbB,h}^{\Delta t,-} \bbB_h^{\Delta t,+}}{\bbG_h}_h \dt
    - \int_0^T \skp{2 \nabla\bv_h^{\Delta t,+}}{\calI_h\left[ \bbG_h \bbB_h^{\Delta t,+} \right]}_{L^2} \dt
    + \int_0^T \alpha \skp{\nabla\bbB_h^{\Delta t,+}}{\nabla\bbG_h}_{L^2} \dt,
\end{align}
\end{subequations}
subject to the initial conditions $\phi_h^{\Delta t}(0) = \phi_h^0$ and $\bbB_h^{\Delta t}(0) = \bbB_h^0$. Here, we use the notation
\begin{align*}
    \bbT_{e,h}^{\Delta t,+} 
    &\coloneqq  (\bbB_h^{\Delta t,+})^2 + \kappa(\phi_h^{\Delta t,-}) \bbB_h^{\Delta t,+} - \bbI ,
    \\
    \Gamma_{f,h}^{\Delta t,-} &\coloneqq \Gamma_{f}(\phi_h^{\Delta t,-}, \sigma_h^{\Delta t,+}, \bbB_h^{\Delta t,-}) , \quad \text{where } f\in\{\phi,\bv,\bbB\},
    \\
    \Gamma_{\sigma,h}^{\Delta t,-} &\coloneqq \Gamma_{\sigma}(\phi_h^{\Delta t,-}, \bbB_h^{\Delta t,-}),
    \\
    m_{f,h}^{\Delta t,-} &\coloneqq m_{f}(\phi_h^{\Delta t,-}, \bbB_h^{\Delta t,-}) , \quad \text{where } f\in\{\phi,\sigma\}.
\end{align*}
Then, the following result is a direct consequence of Theorem \ref{theorem:existence_FE}, \eqref{eq:higher_order1} and \eqref{def:fun_Delta_t}--\eqref{def:fun_Delta_t_pm}.

\begin{corollary}~\\
Let \ref{A1}--\ref{A6} and \ref{S} hold true. 
Let the discrete initial data $\phi_h^0\in\calS_h$ and $\bbB_h^0\in\calW_h$ be given with $\bbB_h^0$ being  positive definite. Moreover, for any $n\in\{1,...,N_T\}$, let the discrete boundary data $\sigma_{\infty,h}^n \in\calS_h$ be given. Besides, assume that $\Delta t < \Delta t_*$,
where the constant $\Delta t_*>0$ depends only on the model parameters and on $c_\infty \coloneqq \norm{\sigma_{\infty,h}^{\Delta t,+}}_{L^\infty(0,T;L^2(\partial\Omega))}$. Then, there exist functions $\phi_h^{\Delta t(,\pm)}$, 
$\mu_h^{\Delta t,+}$, 
$\sigma_h^{\Delta t,+}$, 
$p_h^{\Delta t,+}$, 
$\bv_h^{\Delta t,+}$, 
$\bbB_h^{\Delta t(,\pm)}$ solving \eqref{eq:phi_FE_time}--\eqref{eq:B_FE_time} and constants $C>0$ and $C(T, c_\infty)>0$, both independent of $h, \Delta t$ and with $C(T,c_\infty)$ depending exponentially on $T, c_\infty$, such that
\begin{subequations}
\begin{align}
    \label{eq:bounds_time_sigma}
    \norm{\sigma_h^{\Delta t,+}}_{L^\infty(0,T;H^1)} 
    \leq C \norm{\sigma_{\infty,h}^{\Delta t,+}}_{L^\infty(0,T;L^2(\partial\Omega))},
\end{align}
and
\begin{align}
    \label{eq:bounds_time}
    \nonumber
    &\norm{\phi_h^{\Delta t(,\pm)}}_{L^\infty(0,T;H^1)}^2
    + \norm{\partial_t \phi_h^{\Delta t} }_{L^2(0,T;(H^1)')}
    + \tfrac{1}{\Delta t} \norm{\nabla \phi_h^{\Delta t} - \nabla \phi_h^{\Delta t,\pm}}_{L^2(0,T;L^2)}^2
    \\
    \nonumber
    &\quad
    + \norm{\mu_h^{\Delta t,+}}_{L^2(0,T;H^1)}^2
    + \norm{\bv_h^{\Delta t,+}}_{L^2(0,T;H^1)}^2
    + \norm{p_h^{\Delta t,+}}_{L^{4/3}(0,T;L^2)}
    \\
    \nonumber
    &\quad
    + \norm{\bbB_h^{\Delta t(,\pm)}}_{L^\infty(0,T;L^2)}^2 
    + \norm{\bbB_h^{\Delta t(,\pm)}}_{L^2(0,T;H^1)}^2 
    + \norm{\partial_t \bbB_h^{\Delta t} }_{L^{4/3}(0,T;(H^1)')}
    + \tfrac{1}{\sqrt{\Delta t}} \norm{\bbB_h^{\Delta t} - \bbB_h^{\Delta t,\pm}}_{L^2(0,T;L^2)}^2
    \\
    \nonumber
    &\quad 
    + \norm{\calI_h \trace\ln(\bbB_h^{\Delta t,+})}_{L^\infty(0,T;L^1)}
    + \norm{\calI_h \trace\ln(\bbB_h^{\Delta t,+})}_{L^2(0,T;H^1)}^2
    + \norm{(\bbB_h^{\Delta t,+})^{-1} }_{L^1(0,T;L^1)}^{2/3}
    \\
    &\leq  C(T, c_\infty) \Big(
    1 + \norm{\bbB_h^0}_{L^2}^2 
    + \Delta t \norm{\nabla \bbB_h^0}_{L^2}^2
    + \norm{\calI_h\trace\ln\bbB_h^0}_{L^1}
    + \norm{\nabla\phi_h^0}_{L^2}^2
    + \norm{\calI_h \psi(\phi_h^0)}_{L^1} \Big).
\end{align}
\end{subequations}
\end{corollary}
Here, we note the additional term $\Delta t \norm{\nabla \bbB_h^0}_{L^2}^2$ on the right-hand side of \eqref{eq:bounds_time}, as, compared to \eqref{eq:stability_FE}, we also have $\norm{\bbB_h^{\Delta t}}_{L^2(0,T;H^1)}^2$ and $\norm{\bbB_h^{\Delta t,-}}_{L^2(0,T;H^1)}^2$ on the left-hand side of \eqref{eq:bounds_time}. For a possible construction of initial data where the right-hand side of \eqref{eq:bounds_time} is bounded uniformly in $h,\Delta t>0$, we refer to Remark \ref{remark:initial_data}.


We now show that there exists at least one subsequence of discrete solutions which converges to some limit functions, as $(h,\Delta t)\to (0,0)$. 
For ease of presentation, we always write $\{h,\Delta t>0\}$ for any sequence $\{(h_k, \Delta t_k)\}_{k\in\bbN}$ with $(h_k, \Delta t_k)\to (0,0)$, as $k\to \infty$, and we do not relabel any further subsequence.

\begin{lemma}[Converging subsequences]
\label{lemma:convergence}~\\
Let \ref{A1}--\ref{A6} and \ref{S} hold true. 
Let the discrete initial data $\phi_h^0\in\calS_h$ and $\bbB_h^0\in\calW_h$ be given with $\bbB_h^0$ being  positive definite. Moreover, for any $n\in\{1,...,N_T\}$, let the discrete boundary data $\sigma_{\infty,h}^n \in\calS_h$ be given.
Assume that
\begin{subequations}
\begin{alignat}{3}
    \label{eq:phi0}
    \phi_h^0 &\to \phi_0  \quad && \text{strongly in } L^2(\Omega)
    \quad  &&\text{with } 
    \norm{\phi_h^0}_{H^1} \leq C \norm{\phi_0}_{H^1} , 
    \\
    \label{eq:B0} \nonumber
    \bbB_h^0 &\to \bbB_0 \quad && \text{weakly in } L^2(\Omega;\bbR^{d\times d}_{\mathrm{S}})\quad  &&\text{with } 
    \norm{\bbB_h^0}_{L^2}^2 + \Delta t \norm{\nabla \bbB_h^0}_{L^2}^2 
    + \norm{\calI_h \trace\ln\bbB_h^0}_{L^1} 
    \\
    & && && \quad \quad\leq C(1 + \norm{\bbB_0}_{L^2}^2 + \abs{\ln b_0}), 
    \\
    \label{eq:sigma_infty}
    \sigma_{\infty,h}^{\Delta t,+}|_{\partial\Omega} &\to \sigma_\infty \quad &&\text{strongly in } L^2(0,T;L^2(\partial\Omega))
    \quad  &&\text{with } 
    \norm{\sigma_{\infty,h}^{\Delta t,+}}_{L^\infty(0,T;L^2(\partial\Omega))}
    \leq c_\infty
    \coloneqq \norm{\sigma_\infty}_{L^\infty(0,T;L^2(\partial\Omega))} .
\end{alignat}
\end{subequations}
Suppose that $\Delta t < \Delta t_*$,
where the constant $\Delta t_*>0$ depends only on the model parameters and on $c_\infty$.

Then, there exists a (non-relabeled) subsequence of $$\left\{  \phi_h^{\Delta t(,\pm)}, 
\mu_h^{\Delta t,+}, 
\sigma_h^{\Delta t,+}, 
p_h^{\Delta t,+}, 
\bv_h^{\Delta t,+}, 
\bbB_h^{\Delta t(,\pm)} \right\}_{h,\Delta t>0},$$ 
where $\left\{  \phi_h^{\Delta t(,\pm)}, \mu_h^{\Delta t,+}, \sigma_h^{\Delta t,+}, p_h^{\Delta t,+}, \bv_h^{\Delta t,+}, \bbB_h^{\Delta t(,\pm)} \right\}$ solves the system \eqref{eq:phi_FE_time}--\eqref{eq:B_FE_time}, and limit functions 
\begin{align}
    \label{eq:limit_functions}
    \nonumber
    \phi &\in L^\infty(0,T;H^1) 
    \cap H^1(0,T; (H^1)'),
    \quad 
    \mu \in L^2(0,T; H^1),
    \\
    \nonumber
    \sigma &\in L^\infty(0,T; H^1) ,
    \quad 
    \bv \in L^2(0,T; H^1_{\mathrm{D}}(\Omega;\bbR^d)),
    \quad  p \in L^{4/3}(0,T;L^2),
    \\
    \bbB &\in L^\infty\left(0,T; L^2(\Omega; \bbR^{d\times d}_{\mathrm{S}})\right)
    \cap L^2\left(0,T; H^1(\Omega; \bbR^{d\times d}_{\mathrm{S}})\right)
    \cap W^{1,\frac{4}{3}}\left(0,T; (H^1(\Omega; \bbR^{d\times d}_{\mathrm{S}}))'\right)
\end{align}
with $\bbB$ positive definite a.e.~in $\Omega\times(0,T)$, subject to the initial conditions $\phi(0) = \phi_0 \text{ in } L^2(\Omega)$ and $\bbB(0) = \bbB_0 \text{ in } L^2(\Omega;\bbR^{d\times d}_\mathrm{S})$,
%
%
such that, as $(h,\Delta t) \to (0,0)$,
\begin{alignat}{3}
    \newsubeqblock
    \mysubeq
    \label{eq:conv_phi_Linf_H1}
    \phi_h^{\Delta t(,\pm)} &\to \phi \quad &&\text{weakly-$*$} \quad && \text{in } L^\infty(0,T;H^1),
    \\
    \mysubeq
    \label{eq:conv_dtphi_L2_H1'}
    \partial_t \phi_h^{\Delta t} &\to \partial_t\phi \quad
    &&\text{weakly} \quad && \text{in } L^2(0,T;(H^1)'),
    \\
    \label{eq:conv_mu_L2_H1}
    \mu_h^{\Delta t,+} &\to \mu \quad &&\text{weakly} \quad && \text{in } L^2(0,T;H^1),
    \\
    \label{eq:conv_sigma_Linf_H1}
    \sigma_h^{\Delta t,+} &\to \sigma \quad &&\text{weakly-$*$} \quad && \text{in } L^\infty(0,T; H^1),
    \\
    \label{eq:conv_p_L4/3_L2}
    p_h^{\Delta t,+} &\to p  \quad
    && \text{weakly} \quad &&\text{in } 
    L^{4/3}(0,T;L^2), 
    \\
    \label{eq:conv_v_L2_H1}
    \bv_h^{\Delta t,+} &\to \bv \quad
    && \text{weakly} \quad &&\text{in } 
    L^2(0,T;H^1_{\mathrm{D}}(\Omega;\bbR^d)), 
    \\
    \newsubeqblock
    \mysubeq
    \label{eq:conv_B_Linf_L2}
    \bbB_h^{\Delta t(,\pm)} &\to \bbB \quad
    && \text{weakly-$*$} \quad &&\text{in } 
    L^\infty \left(0,T;L^2(\Omega; \bbR^{d\times d}_{\mathrm{S}})\right), 
    \\
    \mysubeq
    \label{eq:conv_B_L2_H1}
    \bbB_h^{\Delta t(,\pm)} &\to \bbB \quad
    && \text{weakly} \quad &&\text{in } 
    L^2 \left(0,T;H^1(\Omega;\bbR^{d\times d}_{\mathrm{S}})\right), 
    \\
    \mysubeq
    \label{eq:conv_dtB_L4/3_H1'}
    \partial_t \bbB_h^{\Delta t} &\to \partial_t \bbB \quad
    && \text{weakly} \quad &&\text{in } 
    L^{4/3}\left(0,T;(H^1(\Omega;\bbR^{d\times d}_{\mathrm{S}}))'\right).
\end{alignat}
Moreover, it holds
\begin{alignat}{3}
    \newsubeqblock
    \mysubeq
    \label{eq:conv_phi_strong}
    \phi_h^{\Delta t(,\pm)} &\to \phi \quad &&\text{strongly} \quad && \text{in } 
    L^2(0,T;L^q)
    \\
    \mysubeq
    \label{eq:conv_B_strong}
    \bbB_h^{\Delta t(,\pm)} &\to \bbB \quad
    && \text{strongly} \quad &&\text{in } 
    L^2\left(0,T;L^q(\Omega;\bbR^{d\times d}_{\mathrm{S}})\right),
\end{alignat}
where $q\in \left[1,\frac{2d}{d-2}\right)$.

\end{lemma}


\begin{proof}
The subsequence convergences \eqref{eq:conv_phi_Linf_H1}--\eqref{eq:conv_dtB_L4/3_H1'} follow from \eqref{eq:bounds_time_sigma}--\eqref{eq:bounds_time} using standard weak($-*$) compactness results. 
To justify that $\bbB_h^{\Delta t(,\pm)}$ have the same limit function, we note that
\begin{align}
    \label{eq:B_same_limit}
    \norm{\bbB_h^{\Delta t} - \bbB_h^{\Delta t,\pm}}_{L^2(0,T;L^2)} 
    \leq C \sqrt{\Delta t} \to 0,
\end{align}
as $(h,\Delta t)\to (0,0)$, which follows from \eqref{eq:bounds_time}.
The identification of the limit functions of $\phi_h^{\Delta t(,\pm)}$ needs an additional argument, as \eqref{eq:bounds_time} only implies that their gradients have the same limit function. Here, one can derive the estimate
\begin{align}
    \label{eq:phi_same_limit}
    \norm{\phi_h^{\Delta t} - \phi_h^{\Delta t,\pm}}_{L^2(0,T;L^2)} 
    \leq C \sqrt{\Delta t} \to 0,
\end{align}
with, e.g., another testing procedure, taking $\zeta_h = \Delta t (\phi_h^n - \phi_h^{n-1})$ in \eqref{eq:phi_FE} and using \eqref{eq:stability_FE_sigma}, \eqref{eq:stability_FE}, see also, e.g., \cite[Lem.~4.4]{GKT_2022_viscoelastic} for a more general estimate.

The strong convergence results \eqref{eq:conv_phi_strong}--\eqref{eq:conv_B_strong} for $\phi_h^{\Delta t}$, $\bbB_h^{\Delta t}$ follow from \eqref{eq:conv_phi_Linf_H1}--\eqref{eq:conv_dtphi_L2_H1'}, \eqref{eq:conv_B_Linf_L2}--\eqref{eq:conv_dtB_L4/3_H1'} together with \eqref{eq:compact_Lp}, as the embedding $H^1 \hookrightarrow\hookrightarrow L^q$, $q\in[1,\frac{2d}{d-2})$, is compact. Note that \eqref{eq:conv_phi_strong}--\eqref{eq:conv_B_strong} also hold true for $\phi_h^{\Delta t,\pm}$, $\bbB_h^{\Delta t,\pm}$, which is due to \eqref{eq:B_same_limit}, \eqref{eq:phi_same_limit}, \eqref{eq:bounds_time} and the Gagliardo--Nirenberg inequality \eqref{eq:Gagliardo}.

To justify that the initial conditions are satisfied, we argue as follows. It follows from \eqref{eq:limit_functions} and \eqref{eq:compact_C}, that $\phi\in C([0,T]; L^2)$. Together with \eqref{eq:phi0}, we can conclude $\phi(0)=\phi_0$ in $L^2(\Omega)$. Similarly, we have from \eqref{eq:limit_functions} and \eqref{eq:compact_C} that $\bbB \in C([0,T]; (H^1(\Omega); \bbR^{d\times d}_\mathrm{S})')$. As $(H^1)' \hookrightarrow L^2$ is a continuous injection, it follows from, e.g., \cite[Chap.~3, Lem.~1.4]{temam_2001}, that $\bbB\colon\, [0,T]\to L^2(\Omega;\bbR^{d\times d}_{\mathrm{S}})$ is weakly continuous. Together with \eqref{eq:B0}, we can deduce $\bbB(0) = \bbB_0$ in $L^2(\Omega;\bbR^{d\times d}_{\mathrm{S}})$.

We now prove that $\bbB$ is positive definite a.e.~in $\Omega_T$. Using \eqref{eq:conv_B_strong}, we find a further subsequence and a set $N\subset\Omega_T$ with zero Lebesgue measure such that $\bbB_h^{\Delta t,+}(\mathbf{x}, t) \to \bbB(\mathbf{x}, t)$ for all $(\mathbf{x},t) \in N^c \coloneqq \Omega_T \setminus N$, as $(h,\Delta t)\to (0,0)$.
The goal is to show that the complement of $S \coloneqq \{(\mathbf{x},t)\in\Omega_T \mid \det \bbB(\mathbf{x},t) > 0 \}$
is a Lebesgue null set, i.e., $\abs{S^c}= 0$ where $S^c \coloneqq \Omega_T \setminus S$.
By the positive definiteness of $\bbB_h^{\Delta t,+}$, we have $\det(\bbB_h^{\Delta t,+})^{-1} > 0$ and $(\bbB_h^{\Delta t,+})^{-1} \to \bbB^{-1}$ on $S\cap N^c$, as $(h,\Delta t)\to (0,0)$.
We introduce the function
\begin{align*}
    f(\mathbf{x},t) \coloneqq 
    \begin{cases}
    \abs{\bbB^{-1}(\mathbf{x},t)} & \text{ if } (\mathbf{x},t)\in S, \\
    + \infty &\text{ if } (\mathbf{x},t) \not\in S.
    \end{cases}
\end{align*}
For $(\mathbf{x},t) \in S^c \cap N^c$, it holds $\bbB_h^{\Delta t,+}(\mathbf{x},t) \to \bbB(\mathbf{x},t)$, as $(h,\Delta t)\to (0,0)$, and $\det \bbB(\mathbf{x},t)= 0$, which implies $\det(\bbB_h^{\Delta t,+})^{-1}(\mathbf{x},t) \to +\infty$ and $\abs{(\bbB_h^{\Delta t,+})^{-1}(\mathbf{x},t)} \to +\infty$, as $(h,\Delta t)\to (0,0)$.
Thus, we obtain $\abs{(\bbB_h^{\Delta t,+})^{-1}(\mathbf{x},t)} \to f(\mathbf{x},t)$ for a.e.~$(\mathbf{x},t) \in \Omega_T$. On noting the Fatou lemma, \eqref{eq:bounds_time} and \eqref{eq:phi0}--\eqref{eq:sigma_infty}, we have
\begin{align*}
    \int_{\Omega_T} f(\mathbf{x},t) \dv(\mathbf{x},t)
    \leq \liminf_{(h,\Delta t)\to(0,0)} \int_{\Omega_T} \abs{(\bbB_h^{\Delta t,+})^{-1}(\mathbf{x},t)} \dv(\mathbf{x},t)
    \leq C.
\end{align*}
This shows that $f$ is integrable and therefore it holds $\abs{S^c}=0$ and $\bbB$ is positive definite a.e.~in $\Omega_T$.
This completes the proof.
\end{proof}

\subsection{Limit passing}
Now, we use the subsequence convergence result from Lemma \ref{lemma:convergence} to prove Theorem \ref{theorem:weak_solution}.
The limit passing in most of the terms in \eqref{eq:phi_FE_time}--\eqref{eq:B_FE_time} can be done with straight-forward arguments as in, e.g., \cite{barrett_boyaval_2009, GKT_2022_viscoelastic, metzger_2018}. We present only the most technical limit passings that hold true for at least one (non-relabeled) subsequence of $\{h,\Delta t >0\}$. 
The strategy is based on the subsequence convergence results \eqref{eq:conv_phi_Linf_H1}--\eqref{eq:conv_B_strong} together with technical estimates for the numerical errors due to the nodal interpolation operator $\calI_h$.

First, we pass to the limit in \eqref{eq:phi_FE_time}--\eqref{eq:B_FE_time} to show that a (time-)integrated version of \eqref{eq:div_weak}--\eqref{eq:B_weak} is fulfilled with smooth test functions in space and time. For instance, considering \eqref{eq:sigma_weak}, 
we justify the identity
\begin{align*}
    0 &= \int_0^T \delta(t) \skp{ m_\sigma(\phi, \bbB) \nabla \sigma}{\nabla \xi}_{L^2} \dt
    + \int_0^T \delta(t) \skp{\sigma \Gamma_\sigma(\phi,\bbB)}{\xi}_{L^2} \dt
    \\
    &\quad
    + \int_0^T \delta(t) K \skp{\sigma-\sigma_\infty}{\xi}_{L^2(\partial\Omega)} \dt
\end{align*}
for arbitrary $\delta \in C^\infty([0,T])$ and $\xi\in C^\infty(\overline{\Omega})$, and similarly for the other equations in \eqref{eq:div_weak}--\eqref{eq:B_weak}.
Then, the precise form of \eqref{eq:div_weak}--\eqref{eq:B_weak} follows with a density argument and the fundamental theorem of calculus of variations.

First, we pass to the limit in \eqref{eq:sigma_FE_time}. 
Here we take $\xi_h(\cdot,t) = \delta(t) \calI_h \xi$ for all $t\in[0,T]$, with arbitrary $\delta \in C^\infty([0,T])$ and $\xi\in C^\infty(\overline{\Omega})$. 
For ease of presentation, we denote the $j$-th term in \eqref{eq:sigma_FE_time} by $\eqref{eq:sigma_FE_time}_j$, $j\in\bbN$, and similarly for the other equations.
Recalling \eqref{eq:interp_H2}, we note that $\xi_h(\cdot,t)\to \delta(t) \xi$ strongly in $W^{1,\infty}(\Omega)$ for all $t\in[0,T]$, as $h\to 0$.
We then consider a non-relabeled subsequence such that $\xi_h\to \delta \xi$ and $\nabla\xi_h \to \delta\nabla\xi$ a.e.~in $\Omega_T$, as $h\to 0$.
For the limit passing in $\eqref{eq:sigma_FE_time}_1$, the key strategy is to control the numerical error of the nodal interpolation operator $\calI_h\colon\, C(\overline\Omega)\to \calS_h$ with the help of the Lipschitz assumption on $m_\sigma(\cdot,\cdot)$.
It follows from \ref{A2}, \eqref{eq:bounds_time}, \eqref{eq:conv_phi_strong}, \eqref{eq:conv_B_strong} and \eqref{eq:interp_lipschitz_multivar}, that
\begin{align*}
    &\norm{\calI_h[m_{\sigma,h}^{\Delta t,-}] - m_\sigma(\phi,\bbB)}_{L^2(\Omega_T)}
    \\
    &\leq \norm{\calI_h[m_{\sigma,h}^{\Delta t,-}] - m_{\sigma,h}^{\Delta t,-}}_{L^2(\Omega_T)} 
    + \norm{m_{\sigma,h}^{\Delta t,-} - m_\sigma(\phi,\bbB)}_{L^2(\Omega_T)} 
    \\
    &\leq C h \big(\norm{\nabla\phi_h^{\Delta t,-}}_{L^2(\Omega_T)} + \norm{\nabla\bbB_h^{\Delta t,-}}_{L^2(\Omega_T)}\big)
    + C \big( \norm{\phi_h^{\Delta t,-} - \phi}_{L^2(\Omega_T)}
    + \norm{\bbB_h^{\Delta t,-} - \bbB}_{L^2(\Omega_T)} \big)
    \to 0,
\end{align*}
as $(h,\Delta t)\to (0,0)$. Hence, we can extract a further (non-relabeled) subsequence such that $\calI_h m_{\sigma,h}^{\Delta t,-} \to m_\sigma(\phi,\bbB)$ a.e.~in $\Omega_T$. Using the generalized Lebesgue dominated convergence theorem, one gets $\calI_h[m_{\sigma,h}^{\Delta t,-}] \nabla\xi_h \to m_\sigma(\phi,\bbB) \nabla\xi \, \delta$ strongly in $L^2(\Omega_T)$, as $(h,\Delta t)\to (0,0)$. Combining this with the weak convergence result \eqref{eq:conv_sigma_Linf_H1}, we can pass to the limit in $\eqref{eq:sigma_FE_time}_1$.
Next, we deal with $\eqref{eq:sigma_FE_time}_2$. 
Here, we note
\begin{align*}
    &\abs{\int_0^T \skp{\sigma_h^{\Delta t,+} \Gamma_{\sigma,h}^{\Delta t,-} }{\xi_h}_h
    - \skp{\sigma \Gamma_{\sigma}(\phi,\bbB)}{\delta(t) \xi }_{L^2} \dt}
    \\
    &\leq 
    \abs{\int_0^T \skp{\sigma_h^{\Delta t,+} \Gamma_{\sigma,h}^{\Delta t,-} }{\xi_h}_h
    - \skp{ \calI_h \Gamma_{\sigma,h}^{\Delta t,-} }{\calI_h\big[\sigma_h^{\Delta t,+} \xi_h \big]}_{L^2} \dt}
    \\
    &\quad 
    + \abs{\int_0^T 
    \skp{ \calI_h \Gamma_{\sigma,h}^{\Delta t,-} }{\calI_h\big[\sigma_h^{\Delta t,+} \xi_h \big]}_{L^2}
    - \skp{ \calI_h \Gamma_{\sigma,h}^{\Delta t,-} }{\sigma_h^{\Delta t,+} \xi_h }_{L^2}
    \dt}
    \\
    &\quad 
    + \abs{\int_0^T \skp{ \calI_h \Gamma_{\sigma,h}^{\Delta t,-} }{\sigma_h^{\Delta t,+} \xi_h }_{L^2}
    - \skp{\sigma \Gamma_{\sigma}(\phi,\bbB)}{\delta(t) \xi}_{L^2} \dt}
    \eqqcolon I + II + III.
\end{align*}
Using \eqref{eq:lump_Sh_Sh}, \eqref{eq:bounds_time_sigma}, \eqref{eq:bound_source} and \eqref{eq:interp_stability}, we have
\begin{align*}
    I &\leq C h \int_0^T \norm{\calI_h \Gamma_{\sigma,h}^{\Delta t,-}}_{L^2}
    \norm{\calI_h\big[\sigma_h^{\Delta t,+} \xi_h \big] }_{H^1} \dt
    \\
    &\leq C h \big( 1 + \norm{\sigma_h^{\Delta t,+}}_{L^2(\Omega_T)} \big)
    \norm{\sigma_h^{\Delta t,+}}_{L^2(0,T;H^1)} \norm{\delta}_{L^\infty(0,T)} \norm{\xi}_{W^{1,\infty}(\Omega)} 
    \to 0,
\end{align*}
as $(h,\Delta t)\to (0,0)$. We obtain with the same arguments that $II \to 0$,
as $(h,\Delta t)\to (0,0)$. Moreover, we can apply the same strategy as for $\eqref{eq:sigma_FE_time}_1$ to show $III\to 0$, as $(h,\Delta t)\to (0,0)$. This justifies the limit passing in $\eqref{eq:sigma_FE_time}_2$. The limit passage in $\eqref{eq:sigma_FE_time}_3$ follows from \eqref{eq:conv_sigma_Linf_H1}, \eqref{eq:sigma_infty} and \eqref{eq:interp_H2}.

The limit passing in \eqref{eq:mu_FE_time} follows similarly to the limit passing in $\eqref{eq:sigma_FE_time}$, and so we continue with \eqref{eq:phi_FE_time}. First, we show, that (possibly for a further subsequence)
\begin{subequations}
\begin{alignat}{3}
    \label{eq:conv_sigma_strong}
    \sigma_h^{\Delta t,+} &\to \sigma \quad &&\text{strongly} \quad && \text{in } L^2(0,T; H^1),
    \\
    \label{eq:conv_nabla_phi_strong}
    \nabla\phi_h^{\Delta t(,\pm)} &\to \nabla\phi
    \quad
    &&\text{strongly} \quad && \text{in } L^2(0,T;L^2) .
\end{alignat}
\end{subequations}
As $\bigcup_{h>0} \calS_h$ is dense in $H^1(\Omega)$, we find a sequence of functions $\tilde\sigma_h\in L^2(0,T; \calS_h)$ such that $\tilde\sigma_h \to \sigma$ strongly in $L^2(0,T;H^1)$ and a.e.~in $\Omega_T$, and $\nabla\tilde\sigma_h\to \nabla\sigma$ a.e.~in $\Omega_T$, as $h\to 0$. Then, on noting \ref{A1}--\ref{A2} and, as $\tilde\sigma_h \in L^2(0,T;\calS_h)$ is an admissible test function in \eqref{eq:sigma_FE_time}, we have
\begin{align*}
    &C \norm{\sigma_h^{\Delta t,+} - \tilde\sigma_h}_{L^2(0,T;H^1)}^2
    \\
    &\leq 
    \int_0^T \skp{\calI_h[m_{\sigma,h}^{\Delta t,-}] \nabla(\sigma_h^{\Delta t,+}-\tilde\sigma_h)}{\nabla(\sigma_h^{\Delta t,+}-\tilde\sigma_h)}_{L^2} \dt
    \\
    &\quad
    + \int_0^T \skp{ (\sigma_h^{\Delta t,+}-\tilde\sigma_h) \Gamma_{\sigma,h}^{\Delta t,-}}{\sigma_h^{\Delta t,+}-\tilde\sigma_h}_h \dt
    + \int_0^T K \skp{\sigma_h^{\Delta t,+}-\tilde\sigma_h}{\sigma_h^{\Delta t,+}-\tilde\sigma_h}_{L^2(\partial\Omega)} \dt
    \\
    &=
    \int_0^T K \skp{\sigma_{\infty,h}^{\Delta t,+} - \tilde\sigma_h}{\sigma_h^{\Delta t,+}-\tilde\sigma_h}_{L^2(\partial\Omega)} \dt
    - \int_0^T \skp{\calI_h[m_{\sigma,h}^{\Delta t,-}] \nabla\tilde\sigma_h}{\nabla(\sigma_h^{\Delta t,+}-\tilde\sigma_h)}_{L^2} \dt
    \\
    &\quad
    - \int_0^T \skp{ \tilde\sigma_h \Gamma_{\sigma,h}^{\Delta t,-}}{\sigma_h^{\Delta t,+}-\tilde\sigma_h}_h \dt.
\end{align*}
Here, one can use the same arguments as for the limit passage in \eqref{eq:sigma_FE_time} to show that the right-hand side vanishes, as $(h,\Delta t)\to (0,0)$.
This gives rise to
\begin{align*}
    \norm{\sigma_h^{\Delta t,+} - \sigma}_{L^2(0,T;H^1)}
    \leq \norm{\sigma_h^{\Delta t,+} - \tilde\sigma_h}_{L^2(0,T;H^1)}
    + \norm{\tilde\sigma_h - \sigma}_{L^2(0,T;H^1)}
    \to 0,
\end{align*}
as $(h,\Delta t)\to (0,0)$, which proves \eqref{eq:conv_sigma_strong}.
In order to justify \eqref{eq:conv_nabla_phi_strong} for $\nabla\phi_h^{\Delta t,+}$, we choose $\rho_h = \phi_h^{\Delta t,+}$ in \eqref{eq:mu_FE_time} and use similar arguments as for the limit passing in \eqref{eq:sigma_FE_time} to show that
\begin{align*}
    \beta\varepsilon \norm{\nabla\phi_h^{\Delta t,+}}_{L^2(0,T;L^2)}^2
     &= \int_0^T
    \skp{ \mu_h^{\Delta t,+} - \frac\beta\varepsilon \psi_h'(\phi_h^{\Delta t,+},\phi_h^{\Delta t,-}) + \chi_\phi \sigma_h^{\Delta t,+} - \frac12 \kappa_h'(\phi_h^{\Delta t,+}, \phi_h^{\Delta t,-}) \trace\bbB_h^{\Delta t,+}} {\phi_h^{\Delta t,+}}_h \dt
    \\
    &\to \int_0^T \skp{ \mu - \frac\beta\varepsilon \psi'(\phi) + \chi_\phi \sigma - \frac12 \kappa'(\phi) \trace\bbB} {\phi}_{L^2} \dt,
\end{align*}
as $(h,\Delta t)\to (0,0)$.
Then, using the weak formulation \eqref{eq:mu_weak}, we can deduce that $\norm{\nabla\phi_h^{\Delta t,+}}_{L^2(0,T;L^2)} \to \norm{\nabla\phi}_{L^2(0,T;L^2)}$, as $(h,\Delta t)\to (0,0)$, which together with \eqref{eq:conv_phi_Linf_H1} implies the strong convergence result \eqref{eq:conv_nabla_phi_strong} for $\nabla\phi_h^{\Delta t,+}$.
%
%
%
%
%
%
On noting \eqref{eq:phi_same_limit}, we also obtain the strong convergence result \eqref{eq:conv_nabla_phi_strong} for $\nabla\phi_h^{\Delta t}$ and $\nabla\phi_h^{\Delta t,-}$. 
Now we consider the limit passing in \eqref{eq:phi_FE_time}.
Choosing the test function $\zeta_h (\cdot,t) = \delta(t) \calI_h \zeta$ in \eqref{eq:phi_FE_time} for all $t\in[0,T]$ with $\delta\in C^\infty([0,T])$ and $\zeta\in C^\infty(\overline{\Omega})$, and noting the strong convergence results \eqref{eq:conv_sigma_strong}--\eqref{eq:conv_nabla_phi_strong}, we can use similar arguments as for \eqref{eq:sigma_FE_time} to pass to the limit in \eqref{eq:phi_FE_time} and to justify \eqref{eq:phi_weak}. 
%
%

The limit passing in \eqref{eq:div_FE_time} follows with similar arguments as for \eqref{eq:sigma_FE_time}.
We now proceed with \eqref{eq:B_FE_time}. Here, we only consider $\eqref{eq:B_FE_time}_2$, $\eqref{eq:B_FE_time}_3$ and $\eqref{eq:B_FE_time}_4$, as the other terms can be dealt with similarly to \eqref{eq:phi_FE_time}--\eqref{eq:div_FE_time}.
Let $\bbG_h(\cdot, t) = \delta(t) \calI_h \bbG$ for all $t\in[0,T]$ with $\delta\in C^\infty([0,T])$ and $\bbG \in C^\infty(\overline{\Omega};\bbR^{d\times d}_{\mathrm{S}})$.
For $\eqref{eq:B_FE_time}_3$, we have with \eqref{eq:conv_v_L2_H1}, \eqref{eq:conv_B_strong}, \eqref{eq:error_Lambda} and \eqref{eq:interp_H2}, that
\begin{align*}
    - \int_0^T \sum_{i,j=1}^d \skp{(\bv_h^{\Delta t,+})_i \mathbf{\Lambda}_{i,j}(\bbB_h^{\Delta t,+})}{\partial_{x_j} \bbG_h}_{L^2} \dt
    \to - \int_0^T \delta(t) \skp{\bbB}{(\bv\cdot\nabla)\bbG}_{L^2} \dt,
\end{align*}
as $(h,\Delta t)\to (0,0)$. 
For $\eqref{eq:B_FE_time}_2$, we note
\begin{align*}
    &\int_0^T \skp{\bv_h^{\Delta t,+}}{ \nabla\calI_h\left[\bbB_h^{\Delta t,+} : \bbG_h \right]}_{L^2} \dt
    \\
    &= 
    \int_0^T \skp{\bv_h^{\Delta t,+}}{ \nabla\calI_h\left[\bbB_h^{\Delta t,+} : \bbG_h \right] - \nabla \left( \bbB_h^{\Delta t,+} : \bbG_h \right)}_{L^2} \dt
    \\
    &\quad
    - \int_0^T \skp{\Div \bv_h^{\Delta t,+}}{ \bbB_h^{\Delta t,+} : \bbG_h }_{L^2} \dt
    + \int_0^T \skp{\bv_h^{\Delta t,+} \cdot \bn}{\bbB_h^{\Delta t,+} : \bbG_h }_{L^2(\partial\Omega)} \dt,
\end{align*}
where we used integration by parts over $\Omega$. The first term on the right-hand side goes to zero, as $(h,\Delta t)\to (0,0)$, using Hölder's inequality, \eqref{eq:bounds_time}, \eqref{eq:interp_H2} and the local error estimate
\begin{align*}
    \nnorm{\nabla\calI_h \left[ \bbB_h^{\Delta t,+} : \bbG_h \right] - \nabla \left( \bbB_h^{\Delta t,+} : \bbG_h \right)}_{L^2(K)} 
    \leq C h_K \norm{\nabla \bbB_h^{\Delta t,+}}_{L^2(K)} \norm{\nabla\bbG_h}_{L^\infty(K)}
    \quad \forall\, K\in\calT_h,
\end{align*}
which follows from \eqref{eq:interp_H2} and the fact that $\bbB_h^{\Delta t,+}$ and $\bbG_h$ are affine functions on $K$.
Moreover, we have with \eqref{eq:conv_v_L2_H1}, \eqref{eq:conv_B_strong}, \eqref{eq:Gagliardo_boundary}, \eqref{eq:interp_H2} and integration by parts over $\Omega$, that
\begin{align*}
    &- \int_0^T \skp{\Div \bv_h^{\Delta t,+}}{ \bbB_h^{\Delta t,+} : \bbG_h }_{L^2} \dt
    + \int_0^T \skp{\bv_h^{\Delta t,+} \cdot \bn}{\bbB_h^{\Delta t,+} : \bbG_h }_{L^2(\partial\Omega)} \dt
    \\
    &\to 
    - \int_0^T \delta(t) \skp{\Div \bv}{ \bbB : \bbG }_{L^2} \dt
    + \int_0^T \delta(t) \skp{\bv \cdot \bn}{\bbB : \bbG }_{L^2(\partial\Omega)} \dt
    \\
    &= \int_0^T \delta(t) \skp{\bv}{\nabla (\bbB:\bbG)}_{L^2} \dt.
\end{align*}
Altogether, we can deduce
\begin{align*}
    \eqref{eq:B_FE_time}_2 + \eqref{eq:B_FE_time}_3 
    \to \int_0^T \delta(t) \skp{(\bv \cdot \nabla) \bbB}{\bbG}_{L^2} \dt,
\end{align*}
as $(h,\Delta t)\to (0,0)$.
Now, it remains to show the limit passing in $\eqref{eq:B_FE_time}_4$.
Here we note
\begin{align*}
    &\abs{ \int_0^T \skp{\frac{1}{\tau(\phi_h^{\Delta t,-})} \bbT_{e,h}^{\Delta t,+} }{\bbG_h}_h 
    - \skp{\frac{1}{\tau(\phi)} \bbT_e}{\bbG}_{L^2} \dt}
    \\
    &\leq 
    \abs{\int_0^T \skp{\frac{1}{\tau(\phi_h^{\Delta t,-})} \bbT_{e,h}^{\Delta t,+} }{\bbG_h}_h 
    -  \skp{ \calI_h\bbT_{e,h}^{\Delta t,+} }{\calI_h\bigg[ \frac{1}{\tau(\phi_h^{\Delta t,-})} \bbG_h \bigg]}_{L^2}
    \dt }
    \\
    &\quad
    + \abs{
    \int_0^T \skp{ \calI_h\bbT_{e,h}^{\Delta t,+} }{\calI_h\bigg[ \frac{1}{\tau(\phi_h^{\Delta t,-})} \bbG_h \bigg]}_{L^2}
    - \skp{ \calI_h\bbT_{e,h}^{\Delta t,+} }{\calI_h\bigg[ \frac{1}{\tau(\phi_h^{\Delta t,-})}\bigg] \bbG_h }_{L^2} \dt
    }
    \\
    &\quad
    + \abs{
    \int_0^T \skp{ \calI_h\bbT_{e,h}^{\Delta t,+} }{\calI_h\bigg[ \frac{1}{\tau(\phi_h^{\Delta t,-})}\bigg] \bbG_h }_{L^2} 
    - \delta(t) \skp{\frac{1}{\tau(\phi)} \bbT_e}{\bbG}_{L^2}\dt
    } 
    \eqqcolon I + II + III.
\end{align*}
Using \eqref{eq:bound_k_h}, \eqref{eq:bound_grad_k_h} and standard techniques for $\calI_h$ based on Hölder's inequality, \eqref{eq:interp_H2}, \eqref{eq:interp_stability} and \eqref{eq:inverse_estimate}, we have similarly to \eqref{eq:bound_Te}, that 
\begin{align*}
    \norm{\calI_h\bbT_{e,h}^{\Delta t,+}}_{W^{1,1}} 
    \leq C \big(1 + \norm{\bbB_h^{\Delta t,+}}_{L^2(0,T;H^1)} 
    + \norm{\phi_h^{\Delta t,-}}_{L^2(0,T;H^1)} \big).
\end{align*}
Combining this with \eqref{eq:interp_error_Sh_Sh}, \eqref{eq:inverse_estimate}, \eqref{eq:bound_kh'_tau} and \eqref{eq:bounds_time}, we have, as $(h,\Delta t)\to (0,0)$,
\begin{align*}
    I &\leq C h \int_0^T 
    \norm{\calI_h\bbT_{e,h}^{\Delta t,+}}_{W^{1,1}}
    \nnorm{\calI_h\bigg[ \frac{1}{\tau(\phi_h^{\Delta t,-})} \bbG_h \bigg] }_{L^\infty} \dt
    \\
    &\leq C h \big(1 + \norm{\bbB_h^{\Delta t,+}}_{L^2(0,T;H^1)} 
    + \norm{\phi_h^{\Delta t,-}}_{L^2(0,T;H^1)} \big)
    \norm{\delta}_{L^\infty(0,T)} \norm{\bbG}_{L^\infty(\Omega)} \to 0.
\end{align*}
For the second term, we have with Hölder's inequality, \ref{A2}, \eqref{eq:interp_error_Sh_Sh}, \eqref{eq:bound_Te}, \eqref{eq:bounds_time}, that
\begin{align*}
    II 
    &\leq C \int_0^T 
    \norm{\calI_h\bbT_{e,h}^{\Delta t,+}}_{L^2}
    \nnorm{\calI_h\bigg[ \frac{1}{\tau(\phi_h^{\Delta t,-})} \bbG_h \bigg]
    - \calI_h\bigg[ \frac{1}{\tau(\phi_h^{\Delta t,-})}\bigg] \bbG_h}_{L^2}
    \dt
    \\
    &\leq C h 
    \norm{\calI_h\bbT_{e,h}^{\Delta t,+}}_{L^{4/3}(0,T;L^2)} 
    \nnorm{\calI_h \bigg[ \frac{1}{\tau(\phi_h^{\Delta t,-})} \bigg]}_{L^\infty(\Omega_T)}
    \norm{\bbG_h}_{L^4(0,T;H^1)}
    \\
    &\leq C h 
    \norm{\delta}_{L^4(0,T)} \norm{\bbG}_{W^{1,\infty}(\Omega)} \to 0,
\end{align*}
as $(h,\Delta t)\to (0,0)$. 
On noting \ref{A2}, \eqref{eq:bounds_time}, \eqref{eq:conv_phi_strong}, \eqref{eq:interp_lipschitz} and the generalized Lebesgue dominated convergence theorem, we have similarly to the limit passage in $\eqref{eq:sigma_FE_time}_1$, that $\calI_h\Big[ \frac{1}{\tau(\phi_h^{\Delta t,-})}\Big] \bbG_h \to \frac{1}{\tau(\phi)} \delta \bbG$ strongly in $L^4(0,T;L^2(\Omega;\bbR^{d\times d}_{\mathrm{S}}))$. 
Moreover, we have with \ref{A3}, \eqref{eq:interp_lipschitz} and \eqref{eq:conv_phi_strong}, that $\calI_h\kappa(\phi_h^{\Delta t,-}) \to \kappa(\phi)$ strongly in $L^2(\Omega_T)$.
Combining this with \eqref{eq:bounds_time}, \ref{A3}, \eqref{eq:interp_error_Sh_Sh} and \eqref{eq:conv_B_strong}, we can extract a further (non-relabeled) subsequence such that $\calI_h \bbT_{e,h}^{\Delta t,+} \to \bbT_e$ strongly in $L^1(\Omega_T)$ and a.e.~in $\Omega_T$.
As $\calI_h \bbT_{e,h}^{\Delta t,+}$ is uniformly bounded in $L^{4/3}(0,T; L^2(\Omega;\bbR^{d\times d}_{\mathrm{S}}))$, there exists a limit function $\overline{\bbT_e} \in L^{4/3}(0,T;L^2(\Omega;\bbR^{d\times d}_{\mathrm{S}}))$ and a (non-relabeled) subsequence, such that $\calI_h \bbT_{e,h}^{\Delta t,+} \to \overline{\bbT_e}$ weakly in $L^{4/3}(0,T;L^2(\Omega;\bbR^{d\times d}_{\mathrm{S}}))$. As $\calI_h \bbT_{e,h}^{\Delta t,+} \to \bbT_e$ a.e.~in $\Omega_T$, we can identify $\overline{\bbT_e}$ with $\bbT_e$.
This allows to deduce that $III\to 0$, as $(h,\Delta t)\to (0,0)$.

The limit passing in \eqref{eq:v_FE_time} can be established similarly to \eqref{eq:phi_FE_time}--\eqref{eq:sigma_FE_time} and \eqref{eq:B_FE_time}. 
For $\eqref{eq:v_FE_time}_5$ and $\eqref{eq:v_FE_time}_6$, we only remark, that we additionally need to show that
$\nabla \calI_h \kappa(\phi_h^{\Delta t,-}) \to \nabla\kappa(\phi)$ weakly in $L^2(\Omega_T;\bbR^d)$. Here we can argue that (a non-relabeled subsequence of) $\nabla \calI_h \kappa(\phi_h^{\Delta t,-})$ converges weakly in $L^2(\Omega_T;\bbR^d)$ to a function $\gamma\in L^2(\Omega_T;\bbR^d)$, which is due to uniform boundedness and weak compactness. To identify $\gamma$ with $\nabla\kappa(\phi)$, we note that for any test function $\zeta\in C_0^\infty(\Omega_T;\bbR^d)$, it holds with integration by parts over $\Omega$, that
\begin{align*}
    \int_0^T \skp{\nabla \calI_h \kappa(\phi_h^{\Delta t,-})}{\xi}_{L^2} \dt
    = - \int_0^T \skp{\calI_h \kappa(\phi_h^{\Delta t,-})}{ \Div \xi}_{L^2} \dt
    \to - \int_0^T \skp{\kappa(\phi)}{ \Div \xi}_{L^2} \dt,
\end{align*}
as $(h,\Delta t)\to (0,0)$, where we used $\calI_h \kappa(\phi_h^{\Delta t,-}) \to \kappa(\phi)$ strongly in $L^2(\Omega_T)$, as $(h,\Delta t)\to (0,0)$.

The last step is devoted to recover \eqref{eq:stability_sigma}--\eqref{eq:stability} from \eqref{eq:bounds_time_sigma}--\eqref{eq:bounds_time}.
Using \eqref{eq:phi0}--\eqref{eq:sigma_infty} and \ref{A4}, we estimate
\begin{align*}
    \norm{\bbB_h^0}_{L^2}^2 
    + \Delta t \norm{\nabla \bbB_h^0}_{L^2}^2
    + \norm{\calI_h\trace\ln\bbB_h^0}_{L^1}
    + \norm{\nabla\phi_h^0}_{L^2}^2
    + \norm{\calI_h \psi(\phi_h^0)}_{L^1} 
    \leq C \big( 1 + \norm{\bbB_0}_{L^2}^2
    + \abs{\ln b_0} + \norm{\phi_0}_{H^1}^2 \big).
\end{align*}
As $\bbB$ is positive definite a.e.~in $\Omega_T$, and, as $\bbB_h^{\Delta t,+}\to \bbB$ a.e.~in $\Omega_T$, as $(h,\Delta t)\to (0,0)$, we can deduce from $\bbI = (\bbB_h^{\Delta t,+})^{-1} \bbB_h^{\Delta t,+}$, that $(\bbB_h^{\Delta t,+})^{-1} \to \bbB^{-1}$ a.e.~in $\Omega_T$, as $(h,\Delta t)\to (0,0)$.
Next, we show
\begin{align}
    \label{eq:conv_Ih_trlnB}
    \calI_h \trace\ln\bbB_h^{\Delta t,+} 
    = \big(\calI_h \trace\ln\bbB_h^{\Delta t,+} - \trace\ln\bbB_h^{\Delta t,+}\big) 
    + \trace\ln\bbB_h^{\Delta t,+} 
    \to \trace\ln\bbB \quad \text{a.e.~in } \Omega_T,
\end{align}
as $(h,\Delta t)\to (0,0)$.
On noting $\bbB_h^{\Delta t,+} \to \bbB$ a.e.~in $\Omega_T$, as $(h,\Delta t)\to (0,0)$, the continuity of the logarithm and the positive definiteness of $\bbB$, we get $\trace\ln\bbB_h^{\Delta t,+} \to \trace\ln\bbB$ a.e.~in $\Omega_T$ in the limit $(h,\Delta t)\to (0,0)$. 
Let $K\in\calT_h$ be an arbitrary simplex with vertices $P_0^K, ..., P_d^K$.
Using Hölder's inequality, the monotony of the logarithm and \eqref{eq:inverse_estimate}, we compute 
\begin{align*}
    \norm{\calI_h \trace\ln\bbB_h^{\Delta t,+} - \trace\ln\bbB_h^{\Delta t,+}}_{L^2(K)}^2
    &\leq C \abs{K} \norm{\calI_h \trace\ln\bbB_h^{\Delta t,+} - \trace\ln\bbB_h^{\Delta t,+}}_{L^\infty(K)}^2
    \\
    &\leq C \abs{K} \max_{i,j\in\{0,...,d\}} \norm{\trace\ln\bbB_h^{\Delta t,+}(P_i^K) - \trace\ln\bbB_h^{\Delta t,+}(P_j^K)}_{L^\infty(K)}^2
    \\
    &\leq C h_K^2 \abs{K} \norm{\nabla \calI_h \trace\ln\bbB_h^{\Delta t,+}}_{L^\infty(K)}^2
    \\
    &\leq C h_K^2 \norm{\nabla \calI_h \trace\ln\bbB_h^{\Delta t,+}}_{L^2(K)}^2.
\end{align*}
Summing over all $K\in \calT_h$, integrating over the time interval $(0,T)$ and noting \eqref{eq:bounds_time} and \eqref{eq:phi0}--\eqref{eq:sigma_infty}, we deduce that
\begin{align*}
    \norm{\calI_h \trace\ln\bbB_h^{\Delta t,+} - \trace\ln\bbB_h^{\Delta t,+}}_{L^2(\Omega_T)}^2
    \leq C h^2 \norm{\nabla \calI_h \trace\ln\bbB_h^{\Delta t,+}}_{L^2(\Omega_T)}^2
    \to 0,
\end{align*}
as $(h,\Delta t)\to (0,0)$. So, we can extract a further non-relabeled subsequence such that $\calI_h \trace\ln\bbB_h^{\Delta t,+} - \trace\ln\bbB_h^{\Delta t,+} \to 0$ a.e.~in $\Omega_T$, as $(h,\Delta t)\to (0,0)$. 
This shows \eqref{eq:conv_Ih_trlnB}.
Using the weak($-*$) lower semicontinuity of the $L^p$ norms, the Fatou lemma, \eqref{eq:conv_phi_Linf_H1}--\eqref{eq:conv_B_strong} and the subsequence convergence results established in this subsection, we can then recover \eqref{eq:stability_sigma}--\eqref{eq:stability} from \eqref{eq:bounds_time_sigma}--\eqref{eq:bounds_time}. This completes the proof of Theorem~\ref{theorem:weak_solution}.

\section{Numerical results}
\label{sec:numeric}
In this section, we present the results of several numerical computations in two and three space dimensions to show that the numerical scheme \ref{P_FE} is fully-practical. We first specify the model functions and parameters. Moreover, we comment on the solving strategy for the nonlinear scheme \ref{P_FE}. After that, we describe the construction of the initial and boundary data for the cases $d \in \{2,3\}$ separately. 
In all cases, we use the model functions
\begin{align*}
    \Gamma_\phi(\phi,\sigma,\bbB) = \calP \sigma(1+\phi), 
    \quad \Gamma_{\bv}(\phi,\sigma,\bbB) = \frac14 \calP \sigma(1+\phi),
    \quad \Gamma_{\bbB}(\phi,\sigma,\bbB) = \frac12 \calG \calP \sigma(1+\phi),
    \\
    \Gamma_\sigma(\phi,\bbB) = \frac12 \calC (1+\phi),
    \quad m_\phi(\phi,\bbB) = \frac12 (1+\phi)^2 + 10^{-12}, 
    \quad m_\sigma(\phi,\bbB) = 1,
    \\
    \quad \eta(\phi) = \overline\eta, 
    \quad \tau(\phi) = \overline\tau, 
    \quad \kappa(\phi) = \frac12 (1+\phi) \kappa_t, 
    \quad \psi(\phi) = \frac14 (1-\phi^2)^2,
\end{align*}
and, unless otherwise stated, the parameters are set as
\begin{align}
\label{eq:parameters}
    \nonumber
    \Delta t = 0.005, \quad
    \calP = 0.2, \quad
    \varepsilon = 0.02, \quad
    \beta = 0.1, \quad
    \chi_\phi = 4, \quad
    \kappa_t = 0, 
    \\
    \calC = 2, \quad
    K = 10, \quad
    \overline\eta = 10, \quad
    \calG = 0, \quad
    \overline\tau = 100, \quad 
    \alpha = 0.001.
\end{align}
For $\psi_h'(\cdot,\cdot)$, we use the convex-concave splitting of $\psi(\phi) = \psi_1(\phi) + \psi_2(\phi)$ with $\psi_1$ convex and $\psi_2$ concave. In particular, we set $\psi_h'(a,b) \coloneqq \psi_1'(a) + \psi_2'(b) = a^3 - b$ for all $a,b\in\bbR$.
These choices are motivated by, e.g., \cite{ebenbeck_garcke_nurnberg_2020, GKT_2022_viscoelastic, GarckeLSS_2016}. In order to fulfill the assumptions \ref{A1}--\ref{A4}, one can prescribe cut-offs for the model functions if $\abs{\phi}$ is large, which however is not necessary in practice as the phase-field variable $\phi$ usually stays close to the interval $[-1,1]$.

\begin{remark}
One can validate the numerical results with the Cahn--Hilliard--Stokes model  \cite{ebenbeck_garcke_nurnberg_2020} by choosing $\kappa_t = \calG = 0$ and $\overline\tau \approx 0$, as the influence of the viscoelasticity then can be neglected due to $\bbB \approx \bbI$.
\end{remark}

The following algorithm is implemented and solved with the linear algebra package PETSc \cite{petsc-user-ref_2021}. For the assembling, we use the finite element toolbox FEniCS \cite{fenics_book_2012}.
Let $n\in\{1,...,N_T\}$. Given $\phi_h^{n-1} \in\calS_h$, $\bbB_h^{n-1}\in\calW_h$ and $\sigma_{\infty,h}^n \in\calS_h$, we compute a solution of the nonlinear system \ref{P_FE} with the following strategy at each time step. First, we compute the nutrient $\sigma_h^n \in\calS_h$, as \eqref{eq:sig_FE} is decoupled from the other equations. More precisely, we find $\sigma_h^n\in\calS_h$ such that, for all $\xi_h\in\calS_h$,
\begin{align}
    \label{eq:sig_iterative}
    \skp{\nabla\sigma_h^n}{\nabla\xi_h}_{L^2}
    + \skp{\sigma_h^n \Gamma_{\sigma,h}^{n-1}}{\xi_h}_h
    + K \skp{\sigma_h^n}{\xi_h}_{L^2(\partial\Omega)}
    = K \skp{\sigma_{\infty,h}^n}{\xi_h}_{L^2(\partial\Omega)}.
\end{align}
Here one can apply a preconditioned cg-method as the system matrix associated with \eqref{eq:sig_iterative} is symmetric and positive definite.
Then, we use the following iterative scheme to approximate $\phi_h^n, \mu_h^n, p_h^n \in \calS_h$, $\bv_h^n \in \calV_h$ and $\bbB_h^n\in\calW_h$. Set $\ell = 0$ and $f_h^{n,0} \coloneqq f_h^{n-1}$ for all $f\in\{\phi,\mu, p, \bv, \bbB\}$ and repeat the following for $\ell = \{1, ..., \ell_{\mathrm{max}} \}$ until the approximate solution is ``good'' enough:

\begin{enumerate}[(i)]
\item Compute $(\delta\phi)_h^{n,\ell} \in \calS_h$ and $(\delta\mu)_h^{n,\ell} \in \calS_h$ such that, for all $\zeta_h, \rho_h \in\calS_h$,
\begin{subequations}
\begin{align}
    \label{eq:phi_iterative}
    \nonumber
    & \skp{(\delta\phi)_h^{n,\ell}}{\zeta_h}_h 
    + \Delta t \skp{\calI_h [m_{\phi,h}^{n-1}]  \nabla(\delta\mu)_h^{n,\ell}}{\nabla\zeta_h}_{L^2}
    \\ \nonumber
    &= 
    - \skp{\phi_h^{n,\ell-1} - \phi_h^{n-1}}{\zeta_h}_h
    - \Delta t \skp{\calI_h [m_{\phi,h}^{n-1}]  \nabla\mu_h^{n,\ell-1}}{\nabla\zeta_h}_{L^2}
    - \Delta t \skp{\bv_h^{n,\ell-1} \cdot \nabla\phi_h^{n-1}}{\zeta_h}_{L^2}
    \\
    &\quad
    - \Delta t \skp{\phi_h^{n-1} \Gamma_{\bv,h}^{n-1} - \Gamma_{\phi,h}^{n-1}}{\zeta_h}_h,
    \\ \nonumber
    \label{eq:mu_iterative}
    &\skp{- (\delta\mu)_h^{n,\ell} + \frac\beta\varepsilon \psi_1''(\phi_h^{n,\ell-1}) (\delta\phi)_h^{n,\ell}} {\rho_h}_h
    + \beta\varepsilon \skp{ \nabla(\delta\phi)_h^{n,\ell}}{\nabla\rho_h}_{L^2} 
    \\ \nonumber
    &= 
    \skp{\mu_h^{n,\ell-1} 
    - \frac\beta\varepsilon \psi_1'(\phi_h^{n,\ell-1}) 
    - \frac\beta\varepsilon \psi_2'(\phi_h^{n-1}) 
    + \chi_\phi \sigma_h^n 
    - \frac14 \kappa_t \trace\bbB_h^{n,\ell-1}} {\rho_h}_h
    \\
    &\quad - \beta\varepsilon \skp{ \nabla\phi_h^{n,\ell-1}}{\nabla\rho_h}_{L^2},
\end{align}
and update $\phi_h^{n,\ell} \coloneqq \phi_h^{n,\ell-1} + (\delta\phi)_h^{n,\ell}$ and $\mu_h^{n,\ell} \coloneqq \mu_h^{n,\ell-1} + (\delta\mu)_h^{n,\ell}$.

\item
Compute $p_h^{n,\ell} \in\calS_h$ and $\bv_h^{n,\ell} \in\calV_h$ such that, for all $q_h\in\calS_h$, $\bw_h\in\calV_h$,
\begin{align}
    \label{eq:div_iterative}
    &\skp{\Div\bv_h^{n,\ell}}{q_h}_{L^2} 
    = \skp{\Gamma_{\bv,h}^{n-1}}{q_h}_{h},
    \\ \nonumber
    \label{eq:v_iterative}
    & 2\overline\eta \skp{\D(\bv_h^{n,\ell})}{\D(\bw_h)}_{L^2}
    - \skp{p_h^{n,\ell}}{\Div\bw_h}_{L^2}
    \\ \nonumber
    &=
    - \skp{\calI_h \bbT_e(\phi_h^{n-1}, \bbB_h^{n,\ell-1})}{\nabla\bw_h}_{L^2}
    + \skp{(\mu_h^{n,\ell} + \chi_\phi \sigma_h^n) \nabla\phi_h^{n-1}}{\bw_h}_{L^2}
    \\
    &\quad
    + \frac12 \skp{\nabla \calI_h\left[ \kappa(\phi_h^{n-1}) \trace\bbB_h^{n,\ell-1} \right]}{\bw_h}_{L^2}
    - \frac12 \sum_{i,j=1}^d \skp{\partial_{x_j} \calI_h\left[ \kappa(\phi_h^{n-1}) \right]  \trace \mathbf{\Lambda}_{i,j}(\bbB_h^{n,\ell-1})}{(\bw_h)_i}_{L^2},
\end{align}
where $\bbT_e(\phi_h^{n-1}, \bbB_h^{n,\ell-1}) = (\bbB_h^{n,\ell-1})^2 + \kappa(\phi_h^{n-1}) \bbB_h^{n,\ell-1} - \bbI$.

\item Compute $\bbB_h^{n,\ell} \in \calW_h$ such that, for all $\bbG_h \in\calW_h$,
\begin{align}
    \nonumber
    \label{eq:B_iterative}
    &\skp{\bbB_h^n}{\bbG_h}_h
    + \Delta t  \alpha \skp{\nabla\bbB_h^n}{\nabla\bbG_h}_{L^2}
    \\ \nonumber
    &= \skp{\bbB_h^{n-1}}{\bbG_h}_h
    - \Delta t \skp{\bv_h^n}{ \nabla\calI_h\left[\bbB_h^{n,\ell-1} : \bbG_h \right]}_{L^2}
    + \Delta t  \sum_{i,j=1}^d \skp{(\bv_h^{n,\ell})_i \mathbf{\Lambda}_{i,j}(\bbB_h^{n,\ell-1})}{\partial_{x_j} \bbG_h}_{L^2}
    \\
    &\quad 
    - \Delta t  \frac{1}{\overline\tau}  \skp{\bbT_e(\phi_h^{n-1},\bbB_h^{n,\ell-1}) }{\bbG_h}_h 
    - \Delta t  \skp{\Gamma_{\bbB,h}^{n-1} \bbB_h^{n,\ell-1}}{\bbG_h}_h
    + \Delta t  \skp{2 \nabla\bv_h^{n,\ell}}{\calI_h\left[ \bbG_h \bbB_h^{n,\ell-1} \right]}_{L^2}.
\end{align}
\end{subequations}
\end{enumerate}
We use an incremental stopping criterion. In particular, the approximate solution is accepted if it satisfies $\norm{f_h^{n,\ell} - f_h^{n,\ell-1}}_{L^\infty} < \textit{tol }$ 
for all $f\in\{\phi,\mu, p, \bv, \bbB\}$.
Here, we use the tolerance $\textit{tol} = 10^{-7}$. 
Once the approximate solution is accepted, we set $f_h^n \coloneqq f_h^{n,\ell}$, $f\in\{\phi,\mu, p, \bv, \bbB\}$, and we continue with the next time step.

This iterative scheme corresponds to an approximative Newton iteration of the form
\begin{align*}
    \widetilde{\mathbf{DF}}(\mathbf{x}^{\ell-1}) (\delta\mathbf{x})^{\ell} 
    = - \mathbf{F}(\mathbf{x}^{\ell-1}),
    \quad
    \mathbf{x}^{\ell} \coloneqq (\delta\mathbf{x})^{\ell} + \mathbf{x}^{\ell-1},
\end{align*}
where $\mathbf{x}^{\ell-1}, \mathbf{x}^{\ell} \in \bbR^{N_h}$ denote the coefficient vectors of the previous iterate $(\phi_h^{n,\ell-1}, \mu_h^{n,\ell-1}, p_h^{n,\ell-1}, \bv_h^{n,\ell-1}, \bbB_h^{n,\ell-1})$ and the new iterate $(\phi_h^{n,\ell}, \mu_h^{n,\ell}, p_h^{n,\ell}, \bv_h^{n,\ell}, \bbB_h^{n,\ell})$, respectively. Here, $N_h \coloneqq \abs{(\calS_h)^3 \times \calV_h \times \calW_h} \in \bbN$ denotes the number of degrees of freedom for the finite dimensional system \eqref{eq:phi_FE}--\eqref{eq:mu_FE}, \eqref{eq:div_FE}--\eqref{eq:B_FE}.
Moreover, $\widetilde{\mathbf{DF}}\colon \, \bbR^{N_h}\to \bbR^{N_h\times N_h}$ is an approximation of the Jacobian of the map $\mathbf{F}\colon\, \bbR^{N_h}\to \bbR^{N_h}$ that is induced by \eqref{eq:phi_FE}--\eqref{eq:mu_FE}, \eqref{eq:div_FE}--\eqref{eq:B_FE}.
The idea is that one can first solve the linearized Cahn--Hilliard subsystem \eqref{eq:phi_iterative}--\eqref{eq:mu_iterative}, then continue with the Stokes subsystem \eqref{eq:div_iterative}--\eqref{eq:v_iterative} and thereafter compute a solution of the Oldroyd-B equation \eqref{eq:B_iterative}. This has the computational advantage that one can solve the decoupled linear subsystems \eqref{eq:phi_iterative}--\eqref{eq:mu_iterative}, \eqref{eq:div_iterative}--\eqref{eq:v_iterative} and \eqref{eq:B_iterative} successively instead of a fully-coupled linearized system at once. 
In practice, no more than three to five iteration steps were needed until the stopping criterion for \eqref{eq:phi_iterative}--\eqref{eq:B_iterative} was reached. 
Similarly to \eqref{eq:sig_iterative}, one can apply a preconditioned cg-method to solve \eqref{eq:B_iterative}. For the linearized Cahn--Hilliard subsystem \eqref{eq:phi_iterative}--\eqref{eq:mu_iterative}, we used a preconditioned bicgstab-method. The most expensive part in our computations was the Stokes subsystem \eqref{eq:div_iterative}--\eqref{eq:v_iterative} as it has the largest number of degrees of freedom. The Stokes system can be rewritten as a saddle point equation with a symmetric and indefinite system matrix, and thus we apply a Schur complement method with block diagonal preconditioning \cite[Chap.~4.2]{elman_silvester_wathen_2014}. The linear solvers are provided by the linear algebra package PETSc \cite{petsc-user-ref_2021}.

We use a mesh refinement strategy which is similar to \cite{gruen_klingbeil_2014}.
We refine the mesh close to the interfacial region, where we use $\abs{ \nabla\phi_h^n } > 0$ for its identification. 
Away from the interface, we use a coarse mesh corresponding to a uniform mesh with local cell diameter $h_c$. 
The local cell diameter in the interfacial region is $h_f$. 
The values of $h_f$, $h_c$ will be specified below for the cases $d\in\{2,3\}$. The benefit of the mesh refinement strategy is that we have good precision close to the interfacial region, while we reduce the total computational cost by using a local coarse mesh everywhere else.
The mesh always consists of up to $35000$ vertices in our computations for $d=2$ and up to $20000$ vertices for $d=3$, respectively.

The local error estimate \eqref{eq:error_Lambda} can be used to specify on which element $K\in\calT_h$ the quantity $\mathbf{\Lambda}_{i,j}(\bbB_h)|_K$ should be computed or can be replaced by $\delta_{i,j} \,\bbB_h|_{K}$. For simplicity, we compute $\mathbf{\Lambda}_{i,j}(\bbB_h)$ only in two space dimensions, and we use $\delta_{i,j} \,\bbB_h|_K$, $K\in\calT_h$, for the numerical tests in three space dimensions, which however is a good approximation due to \eqref{eq:error_Lambda}. Let us note that we observed no visuable difference in our numerical tests in two space dimensions if we replaced $\mathbf{\Lambda}_{i,j}(\bbB_h)$ by the naive choice $\delta_{i,j} \,\bbB_h|_{K}$. This suggests that the naive choice $\delta_{i,j} \,\bbB_h|_{K}$ is not a bad choice in practice, even though the theoretical analysis requires $\mathbf{\Lambda}_{i,j}(\bbB_h)|_K$ in the numerical scheme.

For the examples in two space dimensions, we fix $\Omega=[-5,5]^2 \subset \bbR^2$ and $\partial_{\mathrm{D}} \Omega = \{-5\}\times[-5,5]$, $\partial_{\mathrm{N}}\Omega = \partial\Omega \setminus \partial_{\mathrm{D}}\Omega$. 
Moreover, we use $h_c = 0.3125$ and $h_f \approx 0.01381$ for the largest and smallest cell diameters, respectively.
For the initial and boundary data, we take $\phi_h^0 = \calI_h \phi_0$, $\bbB_h^0 = \calI_h \bbB_0$ and $\sigma_{\infty,h}^n = \calI_h \sigma_\infty$, $n\in\bbN$, where $\calI_h$ denotes the standard nodal interpolation operator and
\begin{align*}
    \phi_0(\mathbf{x}) = - \tanh \Big( \frac{r(\mathbf{x})}{\sqrt{2} \varepsilon} \Big) \in (-1,1), 
    \quad 
    \bbB_0(\mathbf{x}) = \bbI,
    \quad
    \sigma_\infty(\mathbf{x}) = \sin\Big( \frac\pi2 \frac{\abs{\mathbf{x}}}{5 \sqrt{2}} \Big) \in [0,1],
\end{align*}
for all $\mathbf{x} = \abs{\mathbf{x}} ( \cos(\theta), \, \sin(\theta) )^\top \in \overline\Omega$, where $r(\mathbf{x}) = \abs{\mathbf{x}} - (1 + 0.1 \cos(2 \theta))$.





\begin{figure}[ht!]
\, \\[-4ex]
\centering
\subfloat
{\includegraphics[width=0.24\textwidth]{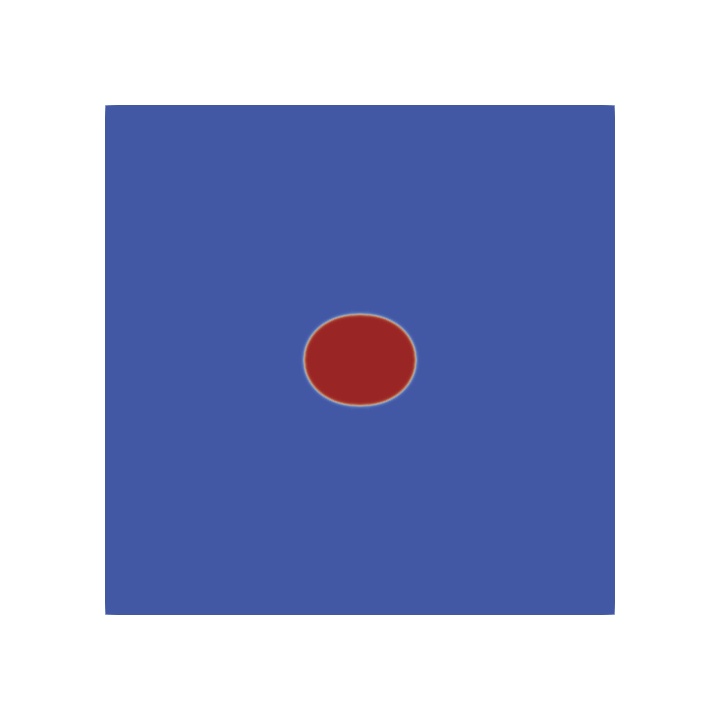}}
\hspace{-2em}
\subfloat
{\includegraphics[width=0.24\textwidth]{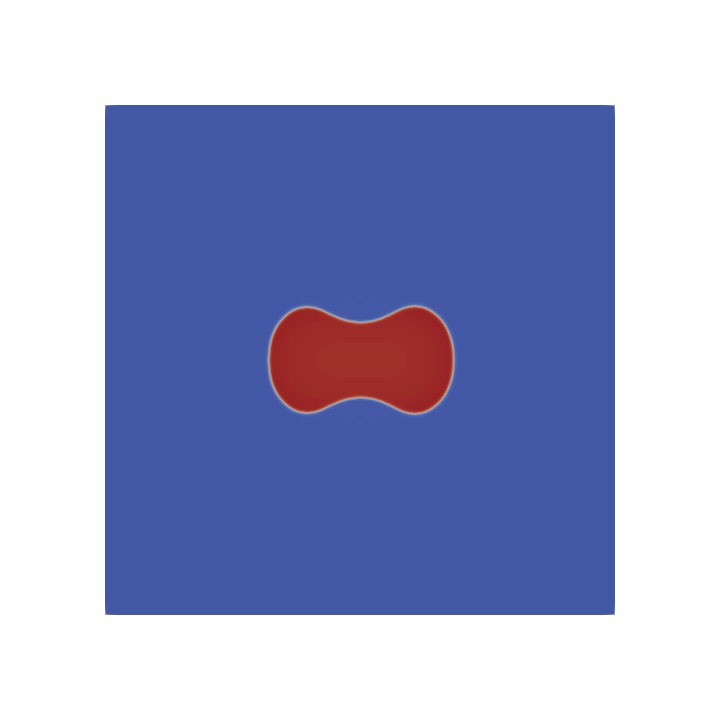}}
\hspace{-2em}
\subfloat
{\includegraphics[width=0.24\textwidth]{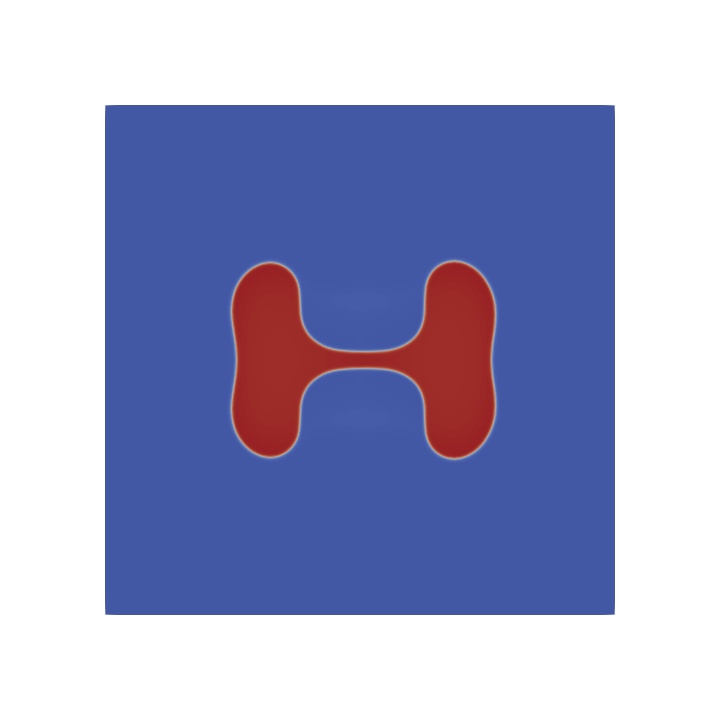}}
\hspace{-2em}
\subfloat
{\includegraphics[width=0.24\textwidth]{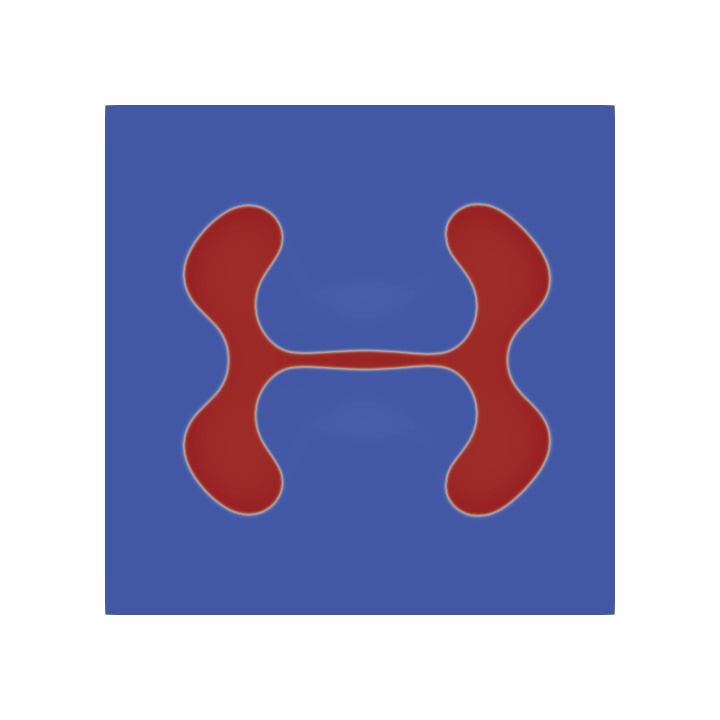}}
\\[-6ex]
\subfloat
{\includegraphics[width=0.24\textwidth]{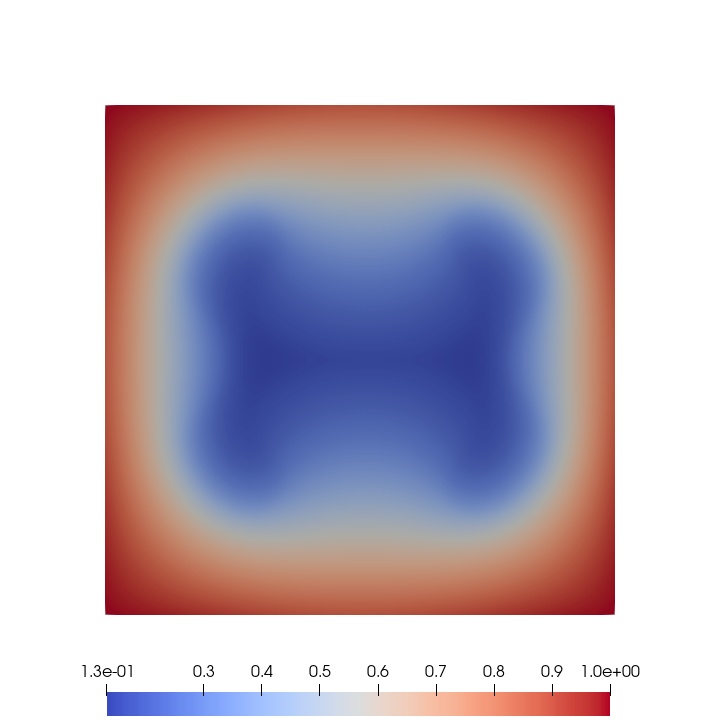}}
\hspace{-2em}
\subfloat
{\includegraphics[width=0.24\textwidth]{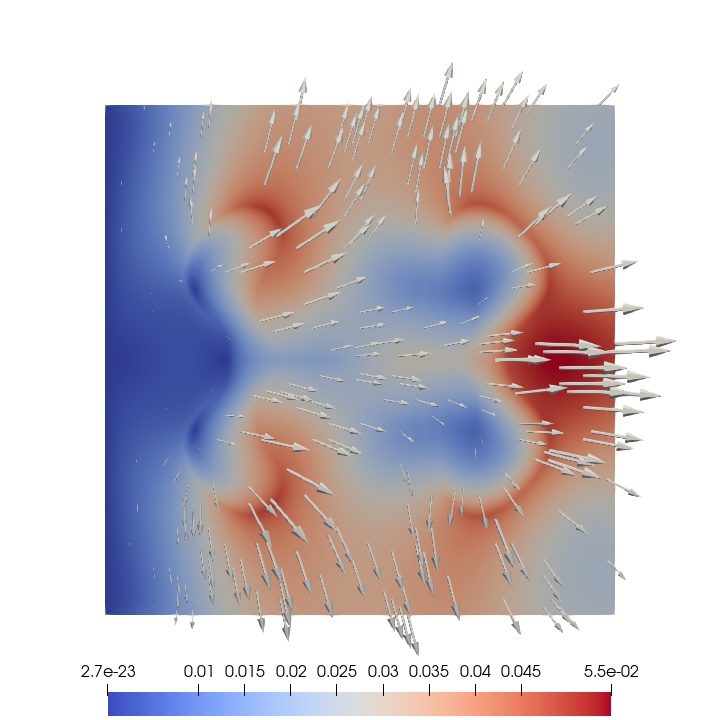}}
\hspace{-2em}
\subfloat
{\includegraphics[width=0.24\textwidth]{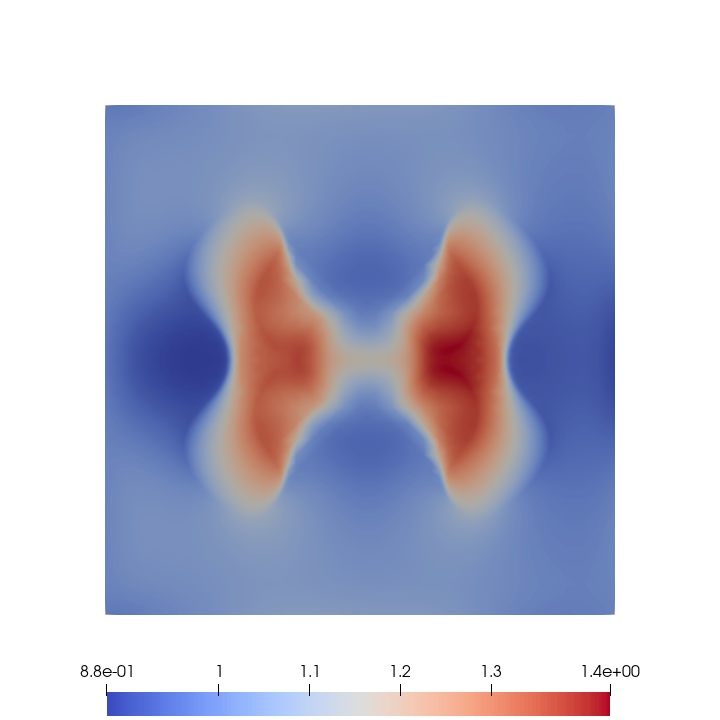}}
\hspace{-2em}
\subfloat
{\includegraphics[width=0.24\textwidth]{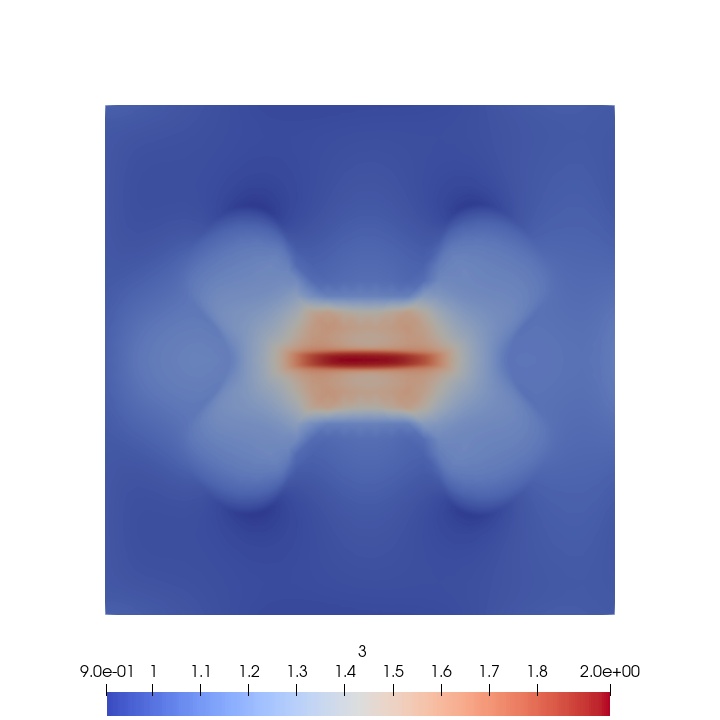}}
\caption{Snapshots of the first example with $\kappa_t=0$. First row: $\phi_h^n$ at $t=0, 5, 10, 14$, where $\phi_h^n=1$ (red) in the tumour tissue and $\phi_h^n=-1$ (blue) in the host tissue. Second row: the nutrient $\sigma_h^n$, $\abs{\bv_h^n}$ (with the velocity field $\bv_h^n$), and both eigenvalues of $\bbB_h^n$ at $t=14$.} 
\label{fig:2d_1a}
\end{figure}
\begin{figure}[ht!]
\, \\[-4ex]
\centering
\subfloat
{\includegraphics[width=0.24\textwidth]{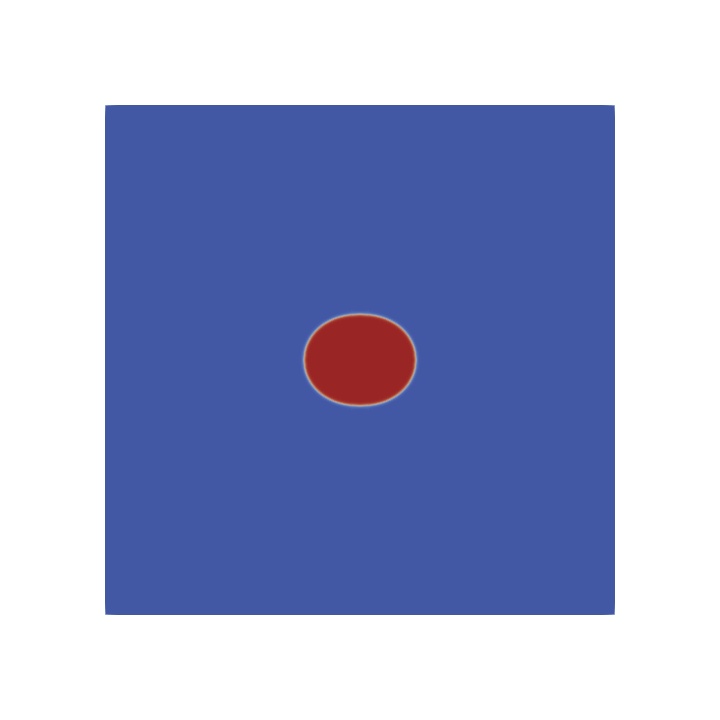}}
\hspace{-2em}
\subfloat
{\includegraphics[width=0.24\textwidth]{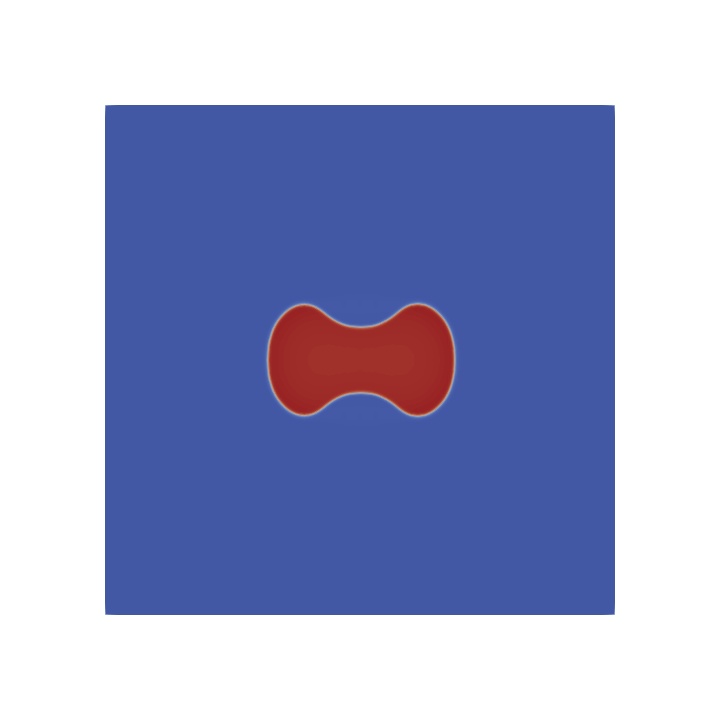}}
\hspace{-2em}
\subfloat
{\includegraphics[width=0.24\textwidth]{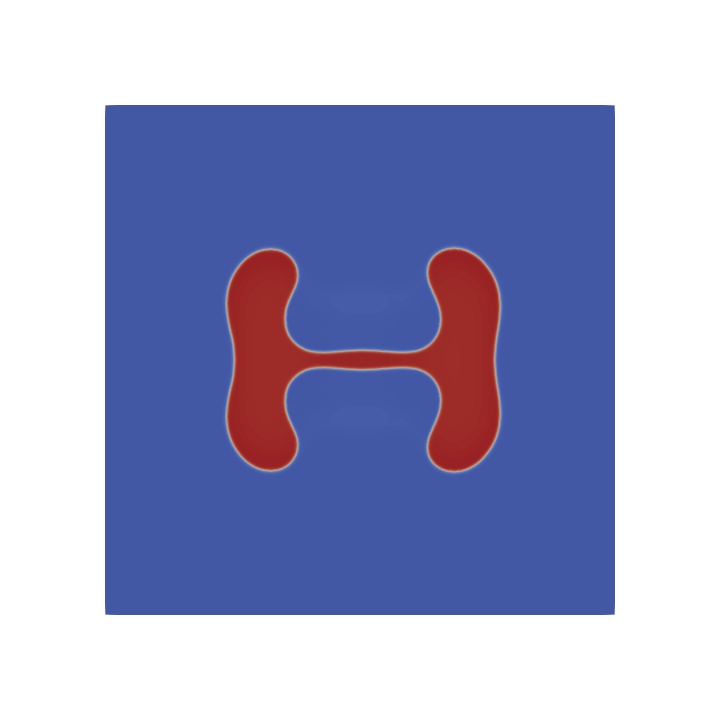}}
\hspace{-2em}
\subfloat
{\includegraphics[width=0.24\textwidth]{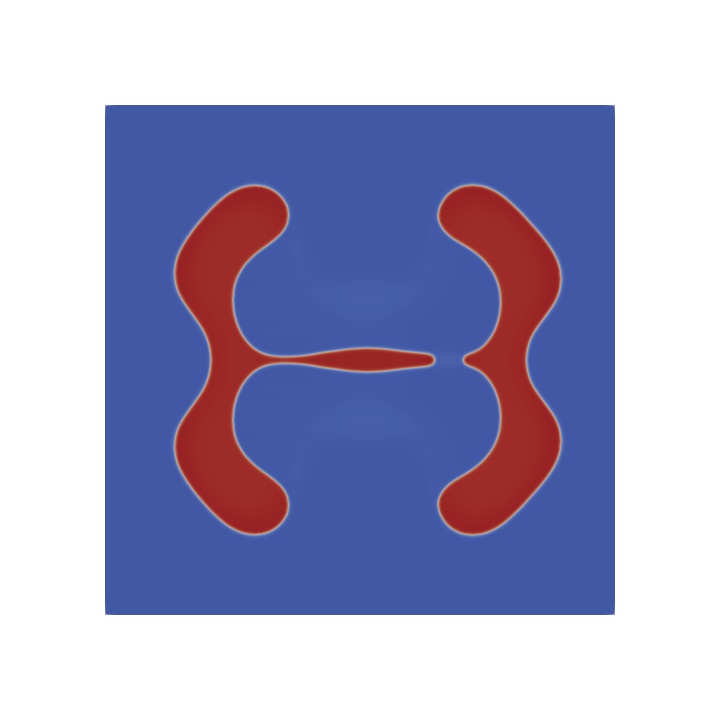}}
\\[-6ex]
\subfloat
{\includegraphics[width=0.24\textwidth]{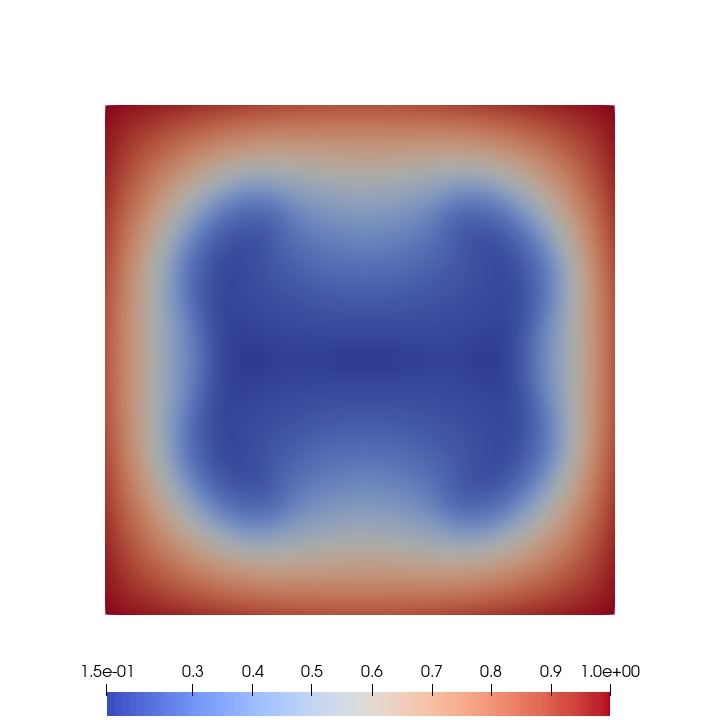}}
\hspace{-2em}
\subfloat
{\includegraphics[width=0.24\textwidth]{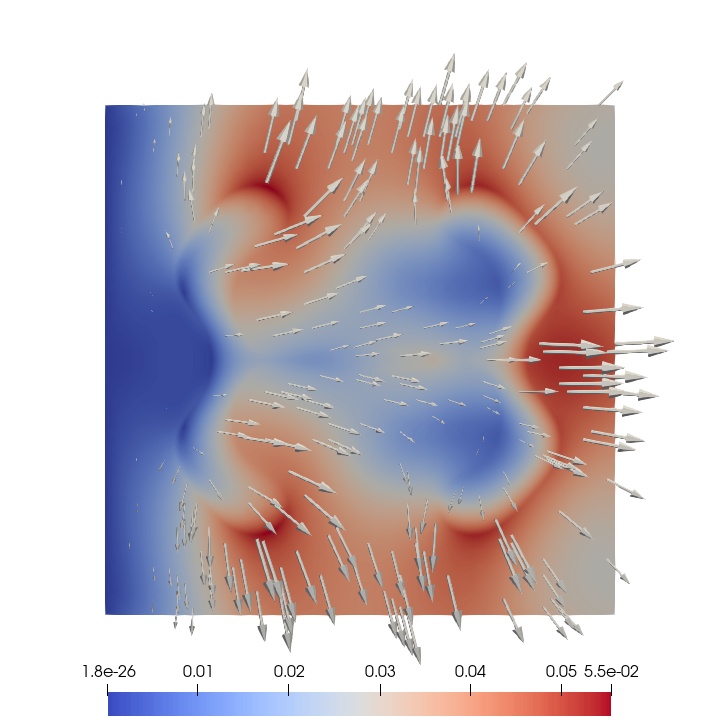}}
\hspace{-2em}
\subfloat
{\includegraphics[width=0.24\textwidth]{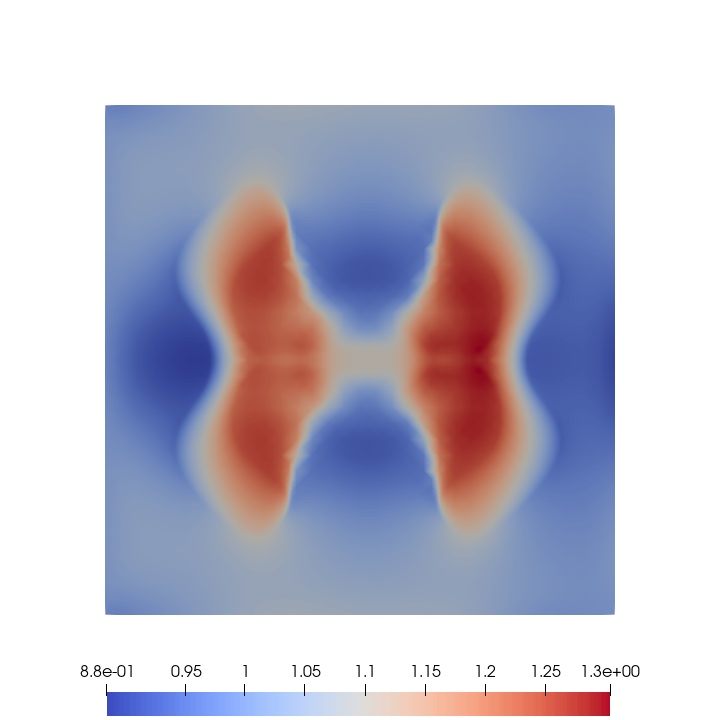}}
\hspace{-2em}
\subfloat
{\includegraphics[width=0.24\textwidth]{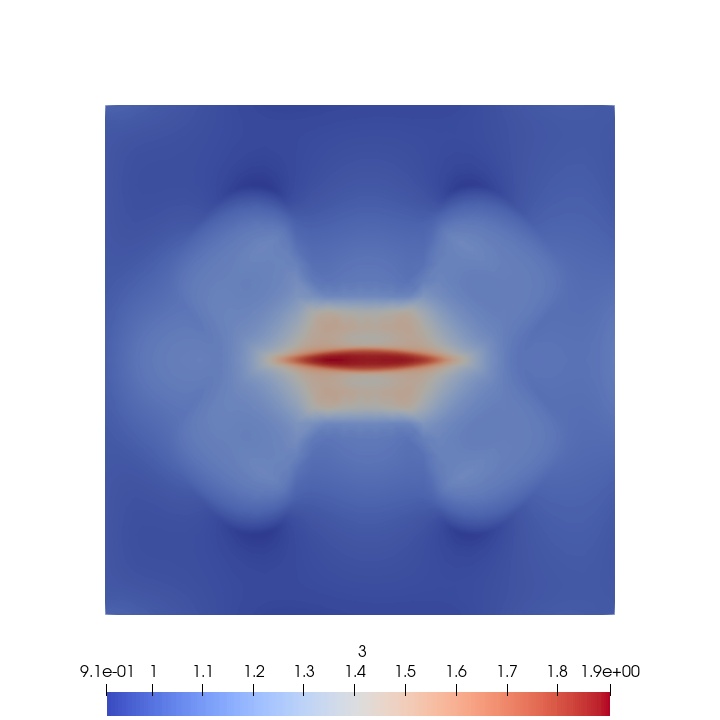}}
\caption{Snapshots of the first example with $\kappa_t=0.5$. First row: $\phi_h^n$ at $t=0, 5, 10, 14$, where $\phi_h^n=1$ (red) in the tumour tissue and $\phi_h^n=-1$ (blue) in the host tissue. Second row: the nutrient $\sigma_h^n$, $\abs{\bv_h^n}$ (with the velocity field $\bv_h^n$), and both eigenvalues of $\bbB_h^n$ at $t=14$.} 
\label{fig:2d_1b}
\end{figure}
\begin{figure}[ht!]
\, \\[-4ex]
\centering
\subfloat
{\includegraphics[width=0.24\textwidth]{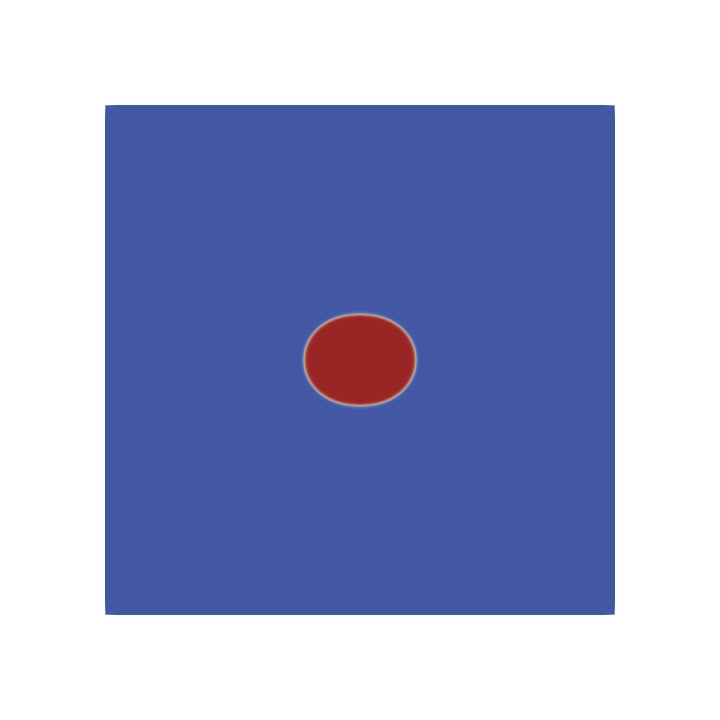}}
\hspace{-2em}
\subfloat
{\includegraphics[width=0.24\textwidth]{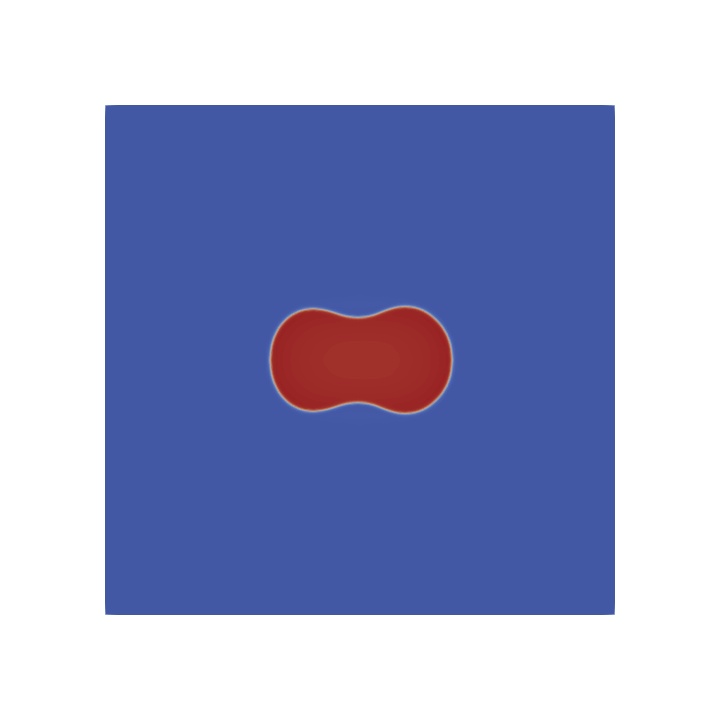}}
\hspace{-2em}
\subfloat
{\includegraphics[width=0.24\textwidth]{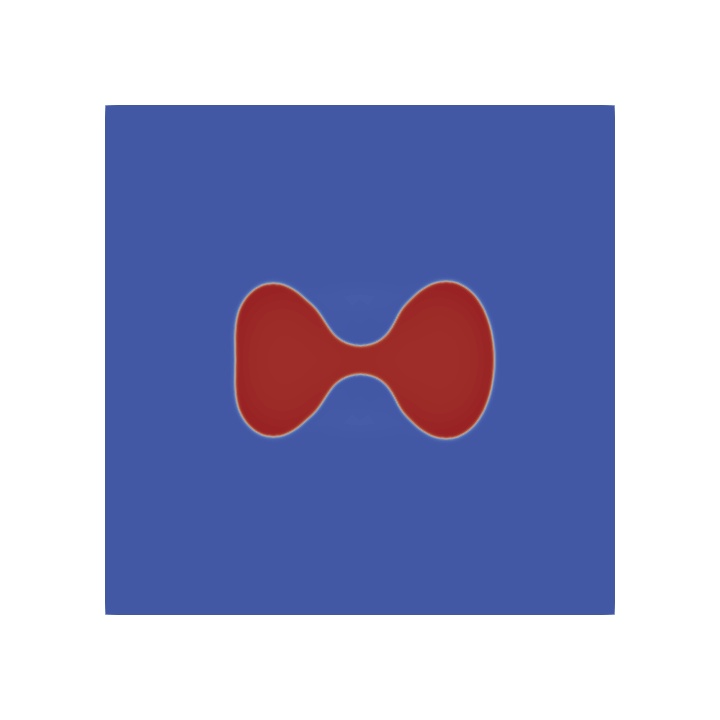}}
\hspace{-2em}
\subfloat
{\includegraphics[width=0.24\textwidth]{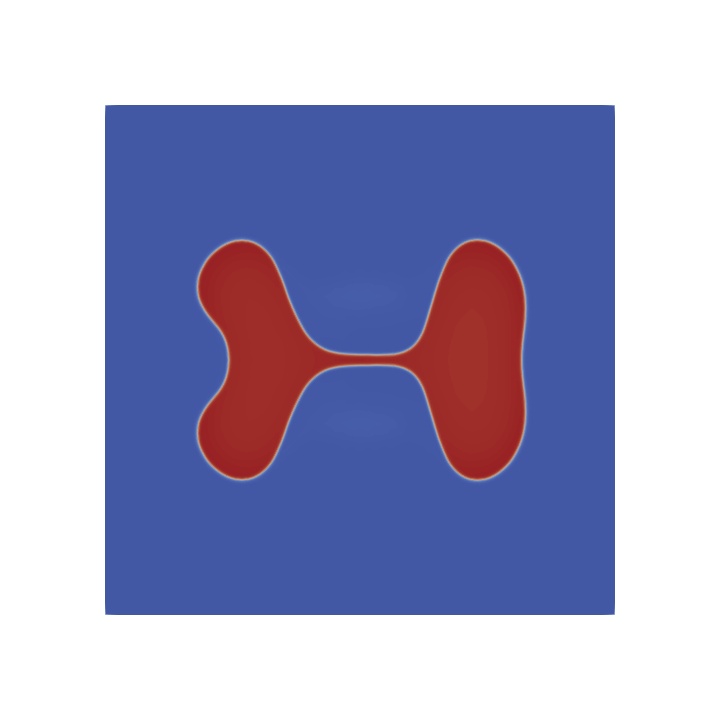}}
\\[-6ex]
\subfloat
{\includegraphics[width=0.24\textwidth]{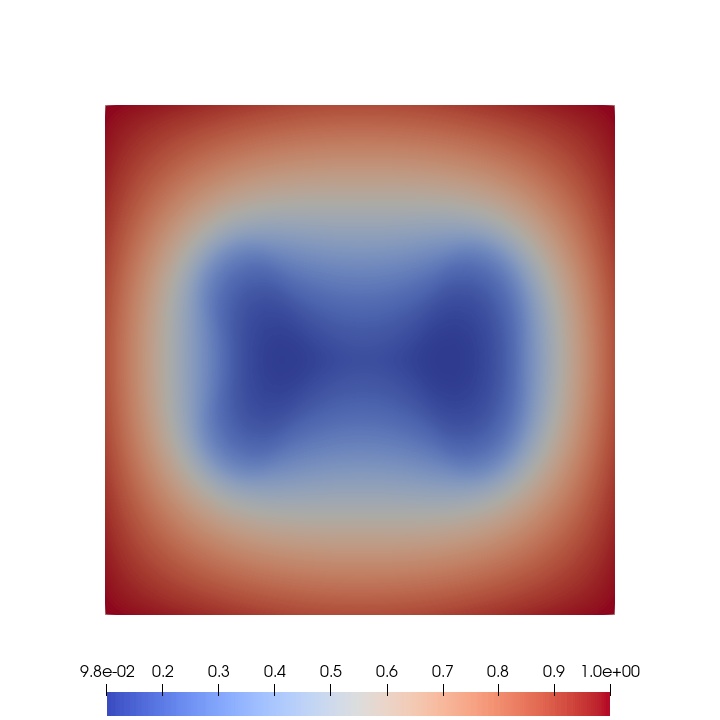}}
\hspace{-2em}
\subfloat
{\includegraphics[width=0.24\textwidth]{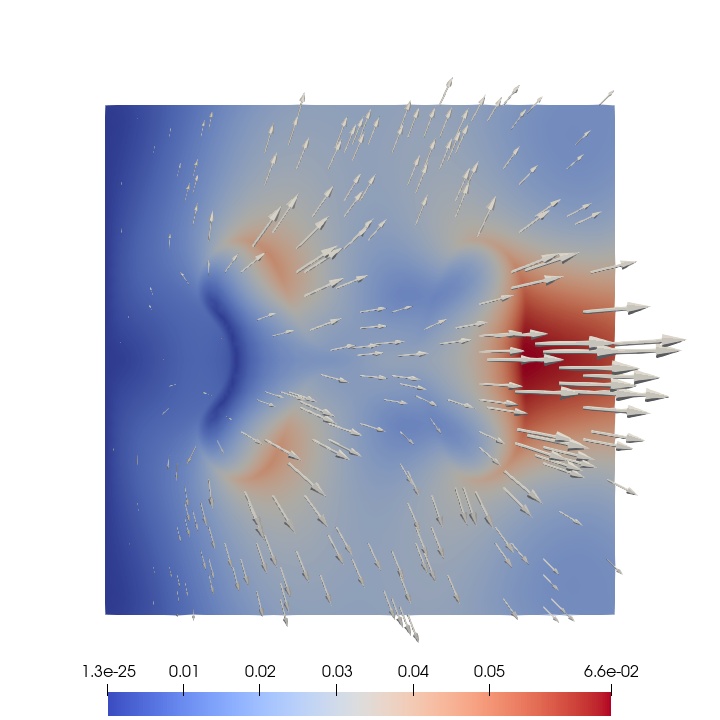}}
\hspace{-2em}
\subfloat
{\includegraphics[width=0.24\textwidth]{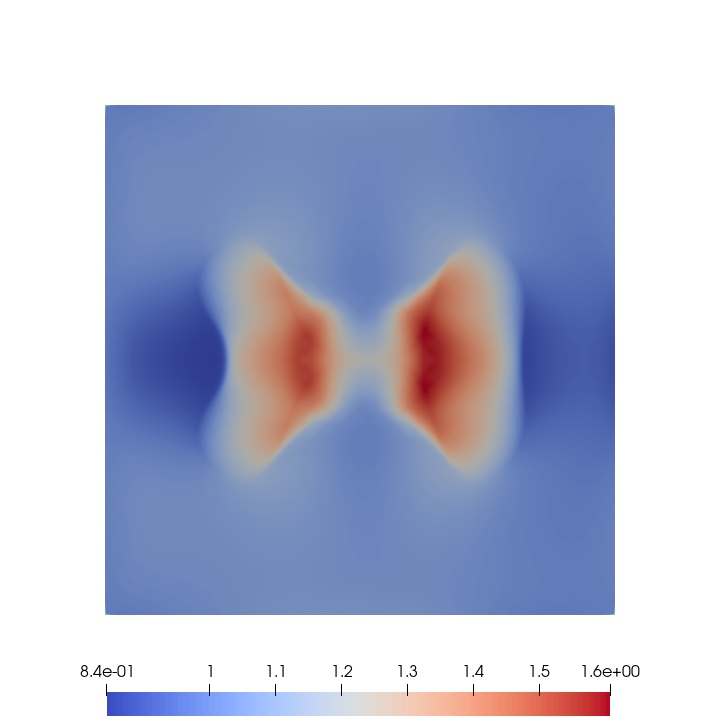}}
\hspace{-2em}
\subfloat
{\includegraphics[width=0.24\textwidth]{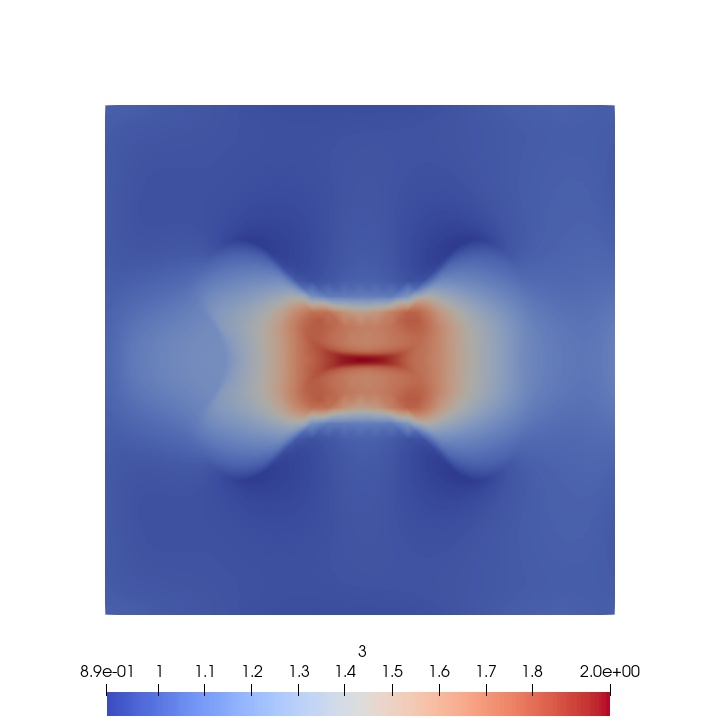}}
\caption{Snapshots of the first example with $\kappa_t=-0.5$. First row: $\phi_h^n$ at $t=0, 5, 10, 14$, where $\phi_h^n=1$ (red) in the tumour tissue and $\phi_h^n=-1$ (blue) in the host tissue. 
Second row: the nutrient $\sigma_h^n$, $\abs{\bv_h^n}$ (with the velocity field $\bv_h^n$), and both eigenvalues of $\bbB_h^n$ at $t=14$.} 
\label{fig:2d_1d}
\end{figure}

In the first example,
we vary $\kappa_t \in \{0, 0.5, -0.5\}$ while the other parameters are fixed as in \eqref{eq:parameters}. The numerical results are visualized in Figures \ref{fig:2d_1a}, \ref{fig:2d_1b} and \ref{fig:2d_1d} for the values $\kappa_t=0$, $0.5$ and $-0.5$, respectively.
In the first row of each figure, we show the order parameter $\phi_h^n$ at times $t\in\{0, 5, 10, 14\}$, and in the second row of each figure, we plot the nutrient $\sigma_h^n$, the velocity magnitude $\abs{\bv_h^n}$ with the corresponding velocity field and both eigenvalues of $\bbB_h^n$ at the final time $t=14$. 
Like for related tumour growth models \cite{ebenbeck_garcke_nurnberg_2020, GKT_2022_viscoelastic, GarckeLSS_2016}, we observe the building of fingers showing towards directions whith higher nutrient concentration, as time goes by, which can be interpreted as the chemotaxis effect. 
The fingers of the tumour in Figure \ref{fig:2d_1b} ($\kappa_t=0.5$) are more elongated compared to Figure \ref{fig:2d_1a} ($\kappa_t=0$). In contrast to that, the shape of the tumour in Figure \ref{fig:2d_1d} ($\kappa_t=-0.5$) is less advanced.
The reason for this behaviour is the term $\frac12 \kappa'(\phi) \trace\bbB = \frac14 \kappa_t \trace\bbB$ in the equation \eqref{eq:mu} for the chemical potential $\mu$, see also \eqref{eq:0_Jphi1} for the relation of $\mu$ with the diffusive flux which accounts for transport of the order parameter $\phi$. In particular, we can expect additional movement of the tumour in the direction of $-\nabla\trace\bbB$ for $\kappa_t>0$ or $+\nabla\trace\bbB$ for $\kappa_t<0$, respectively. Thus, the unstable growth is intensified for $\kappa_t>0$ and it is weakened for $\kappa_t<0$.
The asymmetry in Figures \ref{fig:2d_1b}--\ref{fig:2d_1d} can be explained with the form of the velocity field.
In addition, we even observe a topology change in Figure \ref{fig:2d_1b} at time $t=14$. Such a topology change does not cause any real problem for diffuse interface models, as the order parameter $\phi$ has a smooth transition between the phases $\pm1$. 

\begin{figure}[ht!]
\, \\[-4ex]
\centering
\subfloat
{\includegraphics[width=0.24\textwidth]{figures/example_1/1a_phi_0.jpeg}}
\hspace{-2em}
\subfloat
{\includegraphics[width=0.24\textwidth]{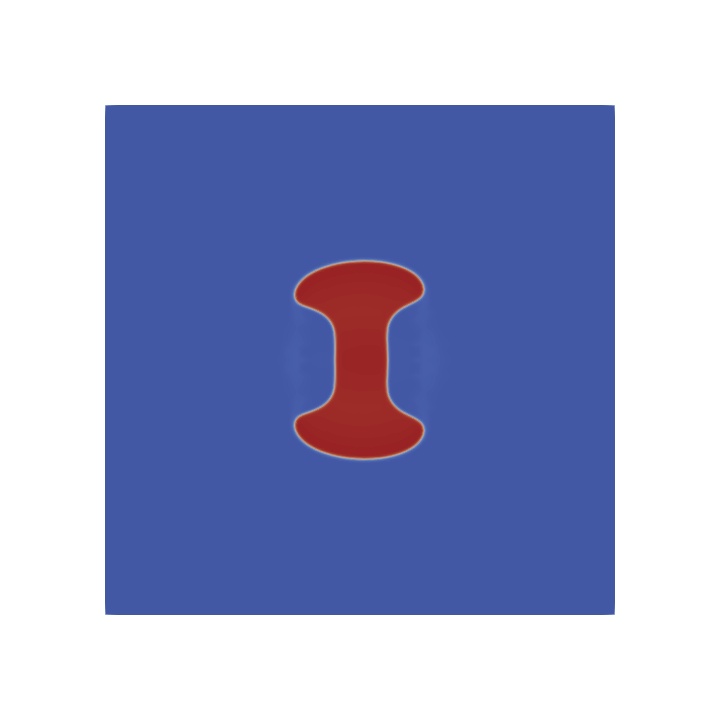}}
\hspace{-2em}
\subfloat
{\includegraphics[width=0.24\textwidth]{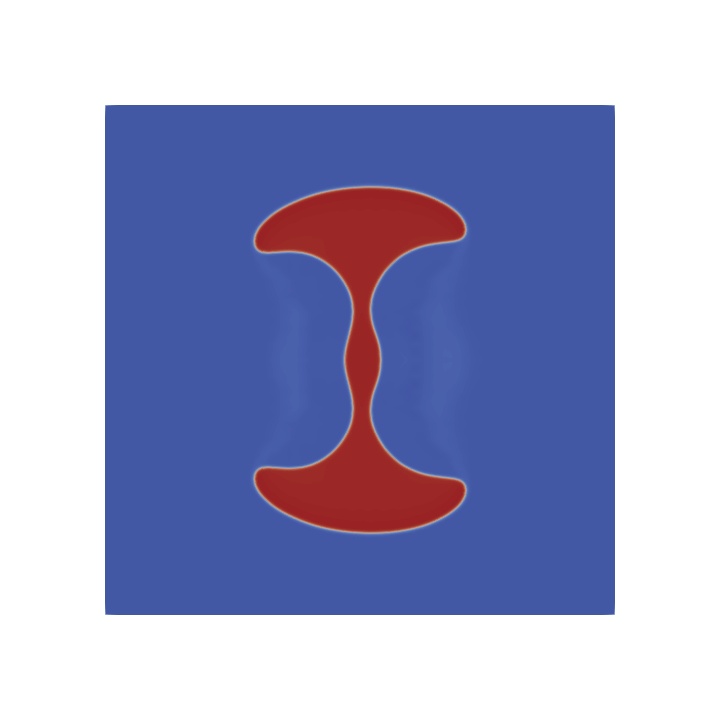}}
\hspace{-2em}
\subfloat
{\includegraphics[width=0.24\textwidth]{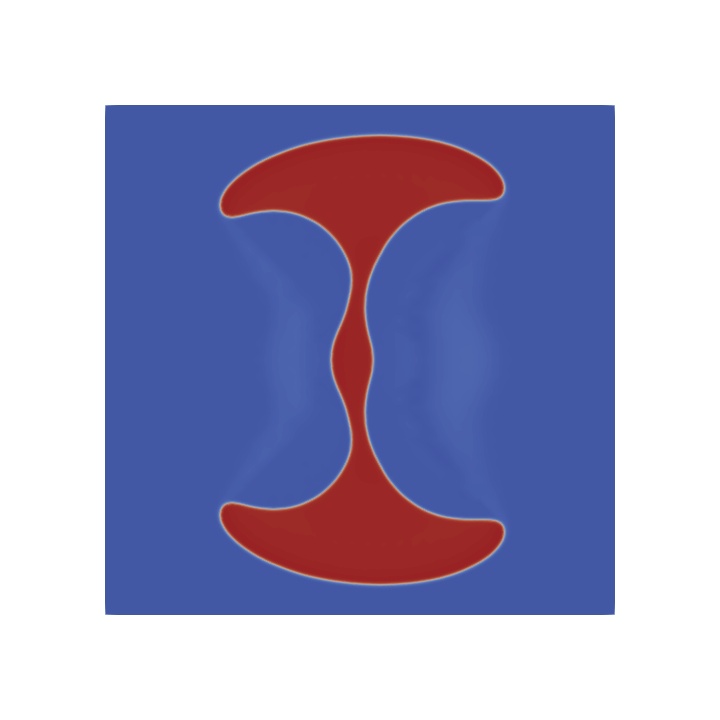}}
\\[-6ex]
\subfloat
{\includegraphics[width=0.24\textwidth]{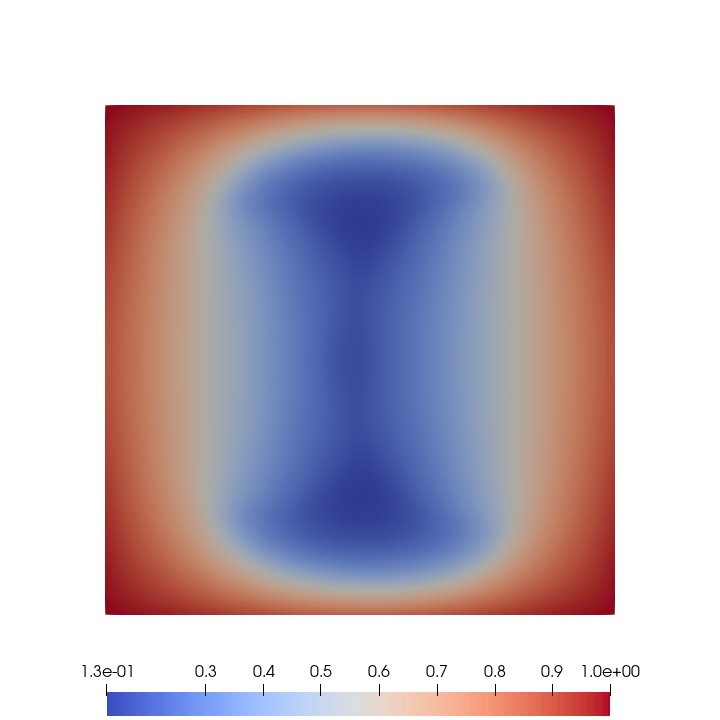}}
\hspace{-2em}
\subfloat
{\includegraphics[width=0.24\textwidth]{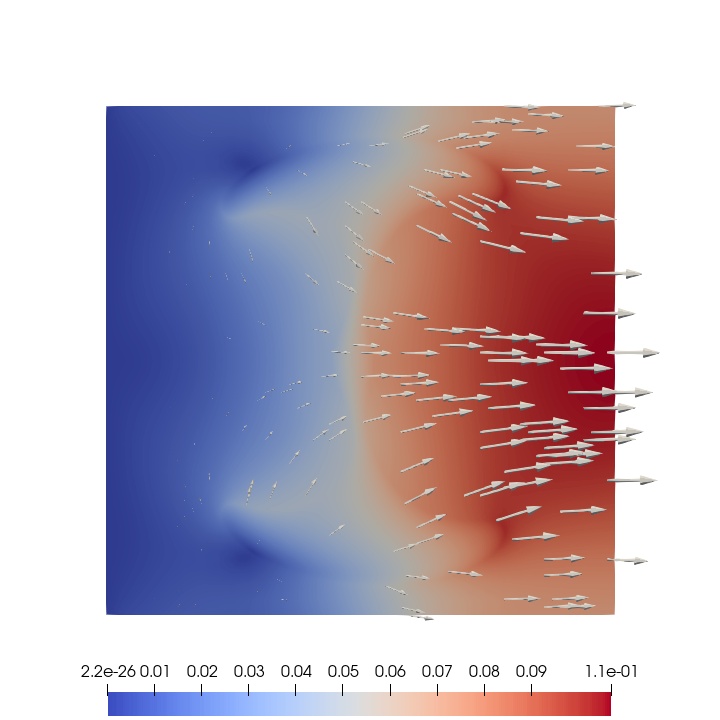}}
\hspace{-2em}
\subfloat
{\includegraphics[width=0.24\textwidth]{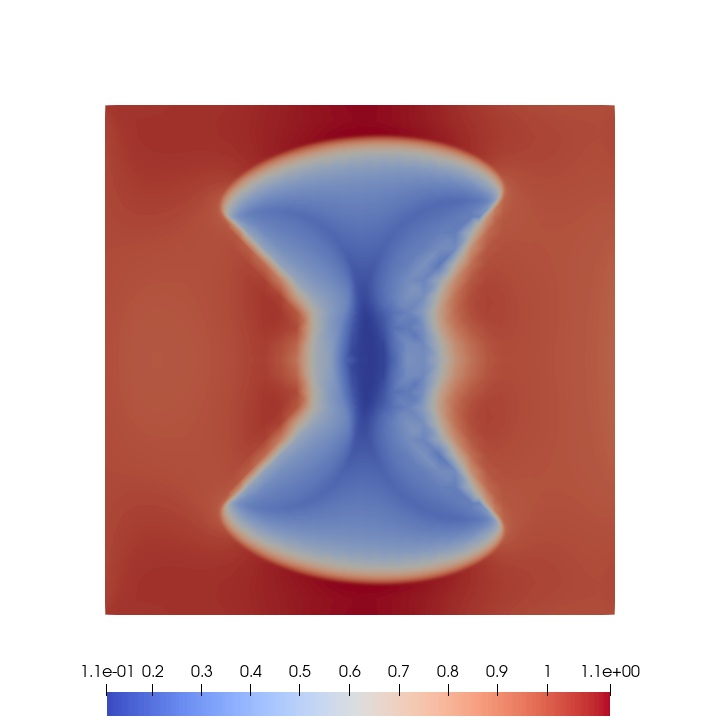}}
\hspace{-2em}
\subfloat
{\includegraphics[width=0.24\textwidth]{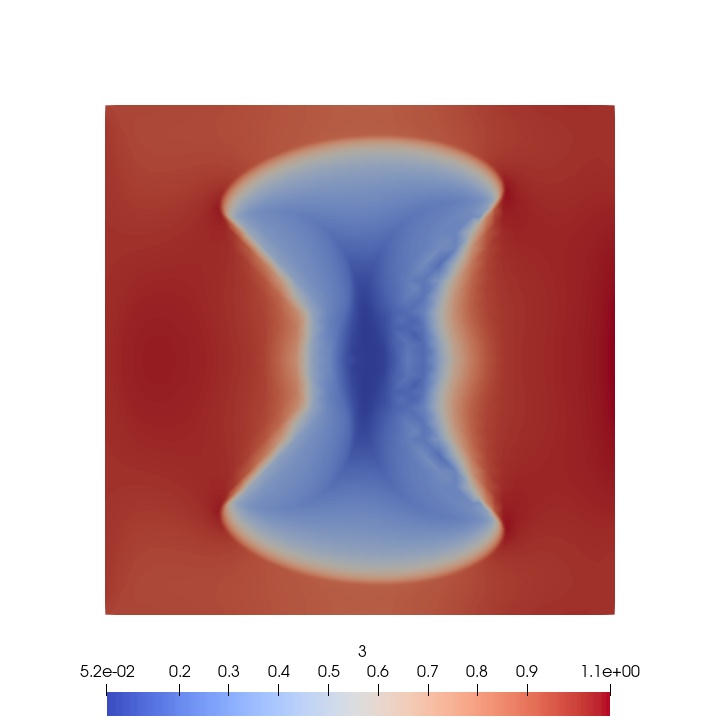}}
\caption{Snapshots of the second example with $\kappa_t=-2$. First row: $\phi_h^n$ at $t=0, 5, 10, 13.5$, where $\phi_h^n=1$ (red) in the tumour tissue and $\phi_h^n=-1$ (blue) in the host tissue.
Second row: the nutrient $\sigma_h^n$, $\abs{\bv_h^n}$ (with the velocity field $\bv_h^n$), and both eigenvalues of $\bbB_h^n$ at $t=13.5$.} 
\label{fig:2d_3b}
\end{figure}

\begin{figure}[ht!]
\, \\[-4ex]
\centering
\subfloat
{\includegraphics[width=0.24\textwidth]{figures/example_1/1a_phi_0.jpeg}}
\hspace{-2em}
\subfloat
{\includegraphics[width=0.24\textwidth]{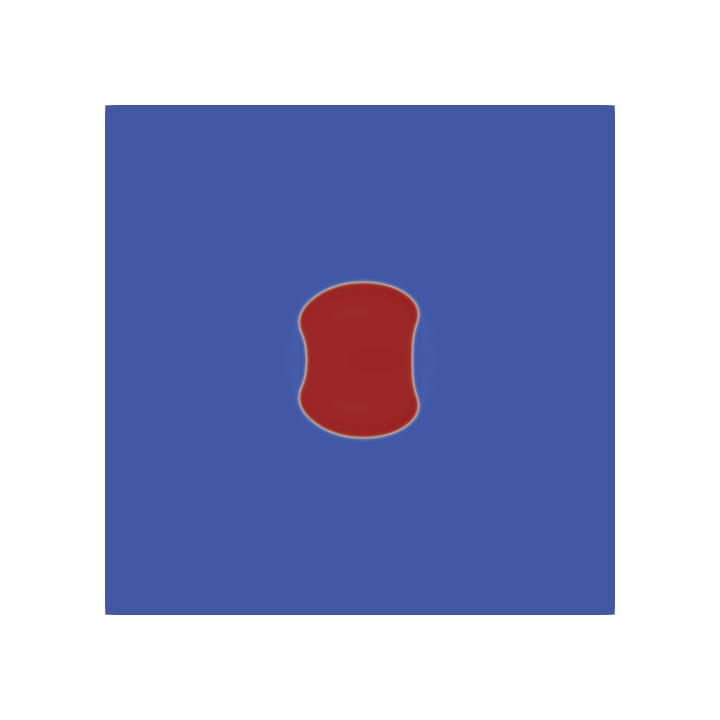}}
\hspace{-2em}
\subfloat
{\includegraphics[width=0.24\textwidth]{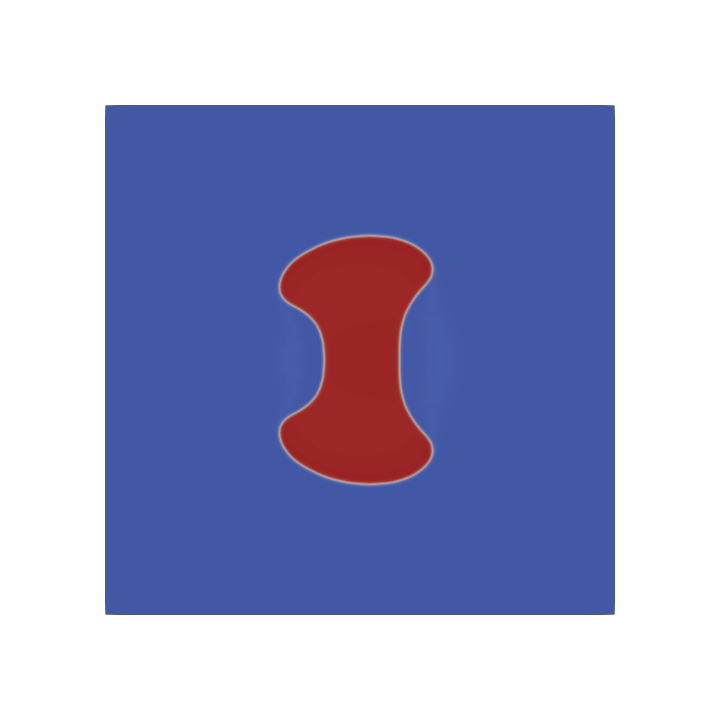}}
\hspace{-2em}
\subfloat
{\includegraphics[width=0.24\textwidth]{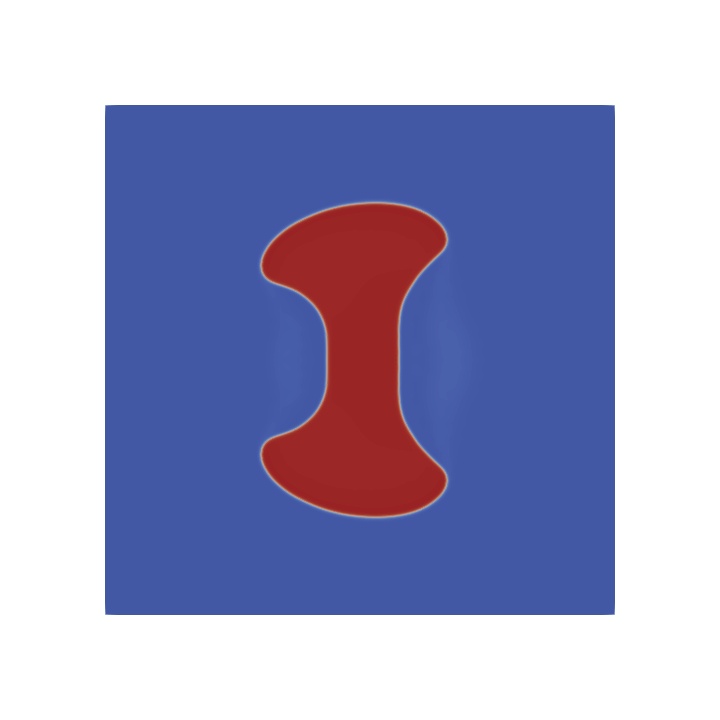}}
\\[-6ex]
\subfloat
{\includegraphics[width=0.24\textwidth]{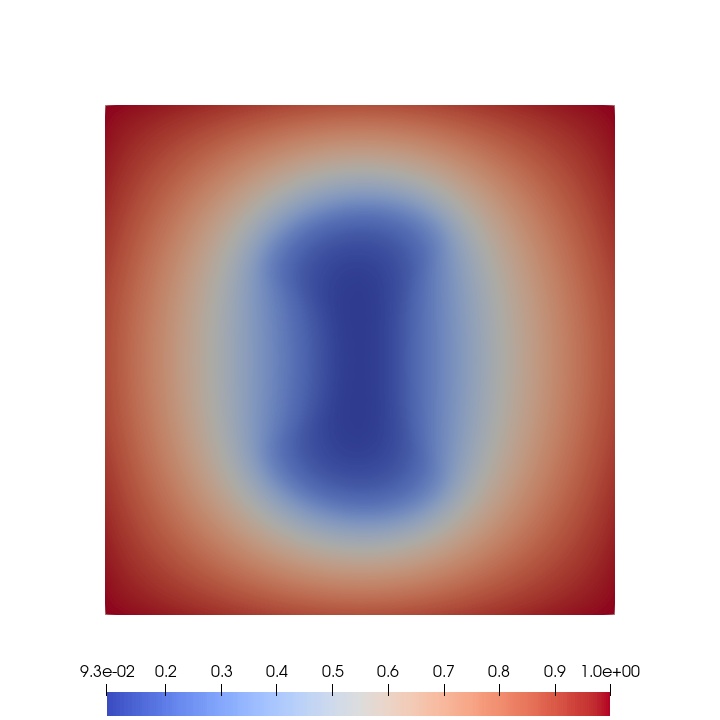}}
\hspace{-2em}
\subfloat
{\includegraphics[width=0.24\textwidth]{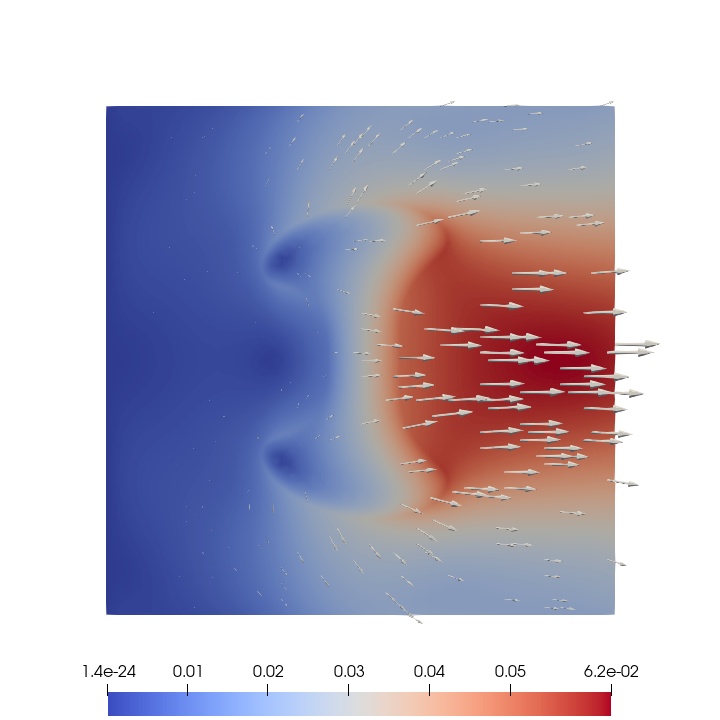}}
\hspace{-2em}
\subfloat
{\includegraphics[width=0.24\textwidth]{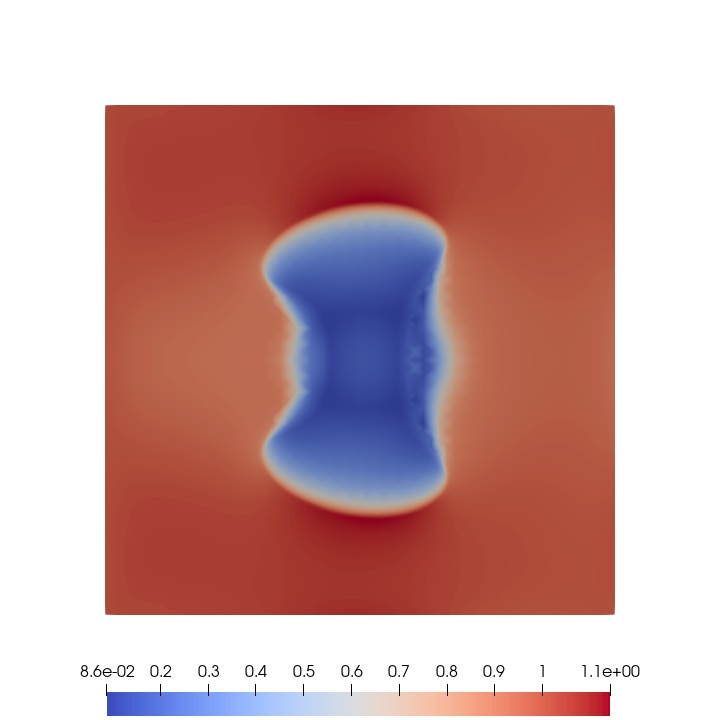}}
\hspace{-2em}
\subfloat
{\includegraphics[width=0.24\textwidth]{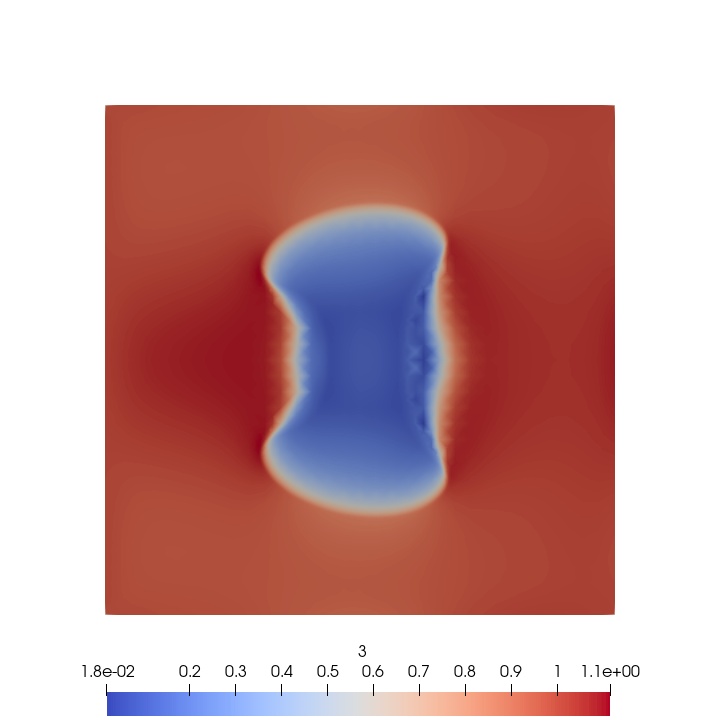}}
\caption{Snapshots of the second example with $\kappa_t=-1$. First row: $\phi_h^n$ at $t=0, 5, 10, 13.5$, where $\phi_h^n=1$ (red) in the tumour tissue and $\phi_h^n=-1$ (blue) in the host tissue.
Second row: the nutrient $\sigma_h^n$, $\abs{\bv_h^n}$ (with the velocity field $\bv_h^n$), and both eigenvalues of $\bbB_h^n$ at $t=13.5$.} 
\label{fig:2d_3e}
\end{figure}

\begin{figure}[ht!]
\, \\[-4ex]
\centering
\subfloat
{\includegraphics[width=0.24\textwidth]{figures/example_1/1a_phi_0.jpeg}}
\hspace{-2em}
\subfloat
{\includegraphics[width=0.24\textwidth]{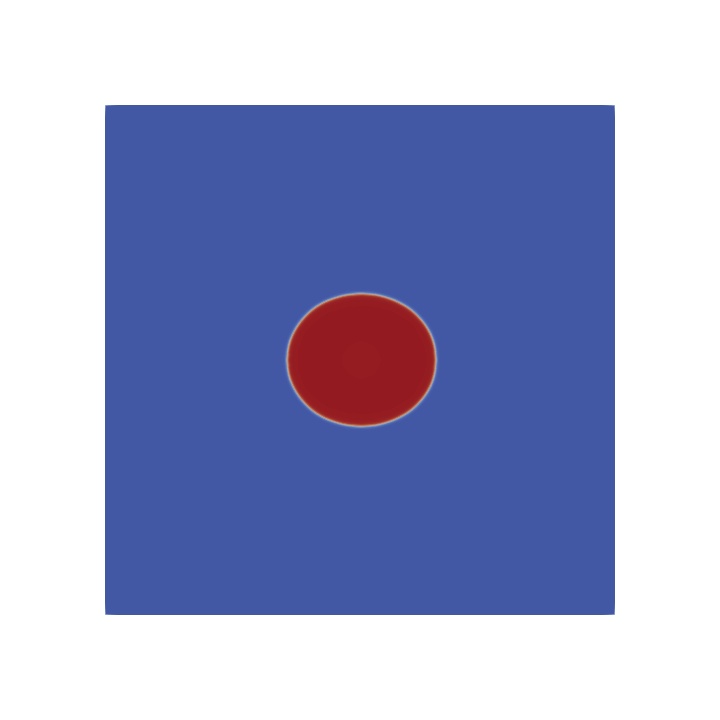}}
\hspace{-2em}
\subfloat
{\includegraphics[width=0.24\textwidth]{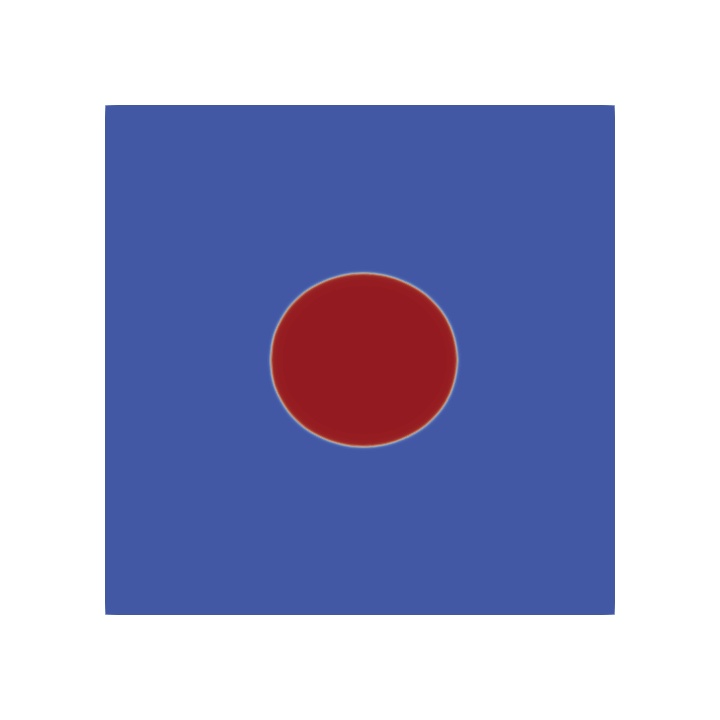}}
\hspace{-2em}
\subfloat
{\includegraphics[width=0.24\textwidth]{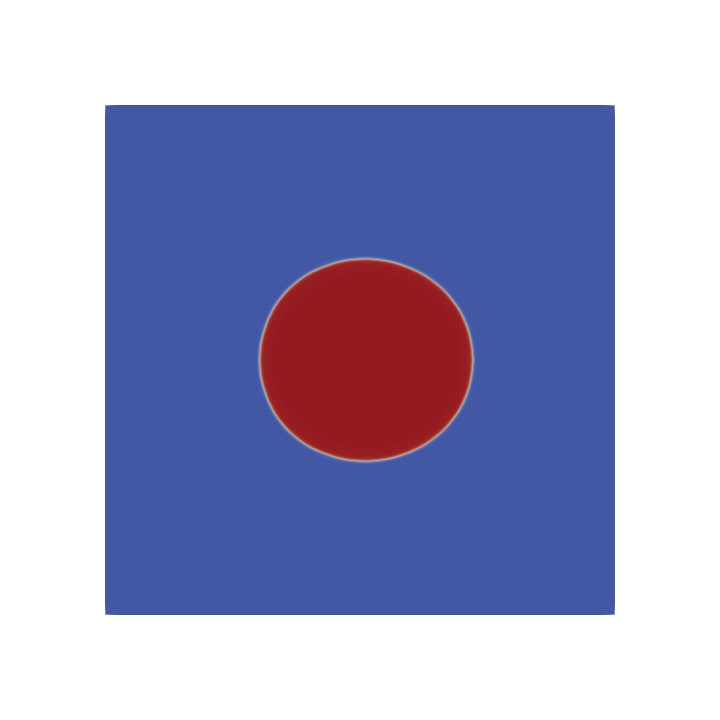}}
\\[-6ex]
\subfloat
{\includegraphics[width=0.24\textwidth]{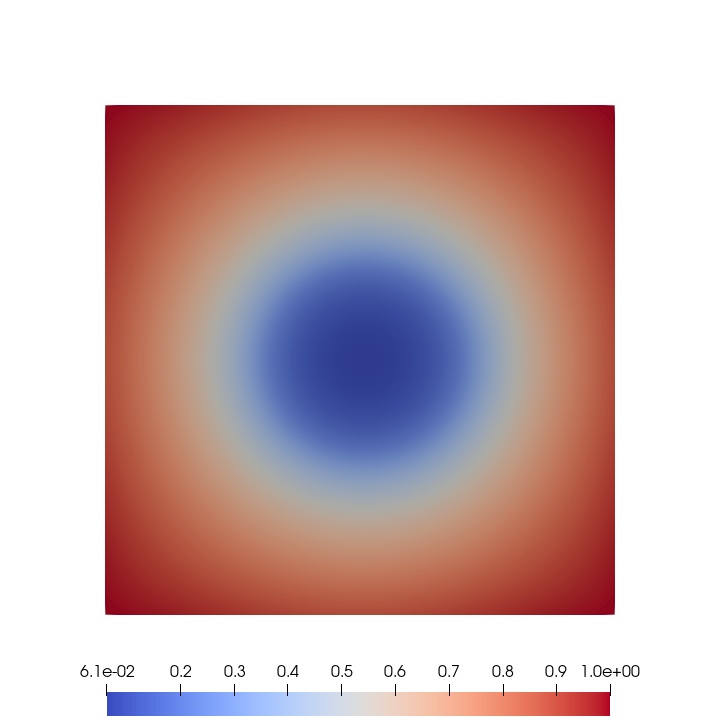}}
\hspace{-2em}
\subfloat
{\includegraphics[width=0.24\textwidth]{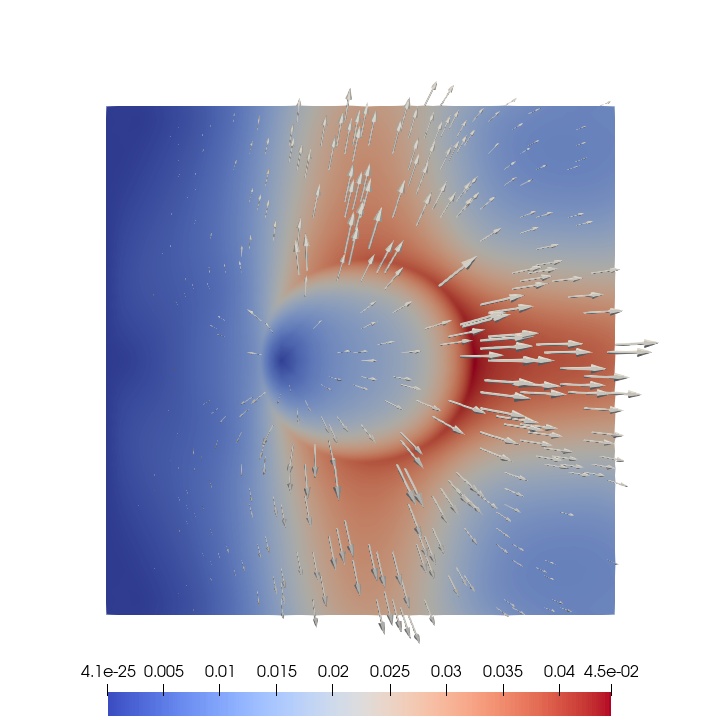}}
\hspace{-2em}
\subfloat
{\includegraphics[width=0.24\textwidth]{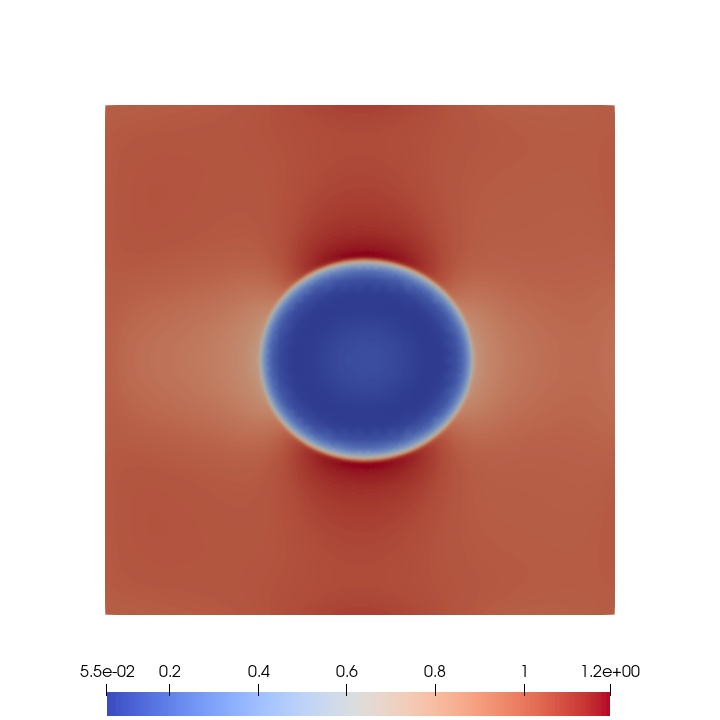}}
\hspace{-2em}
\subfloat
{\includegraphics[width=0.24\textwidth]{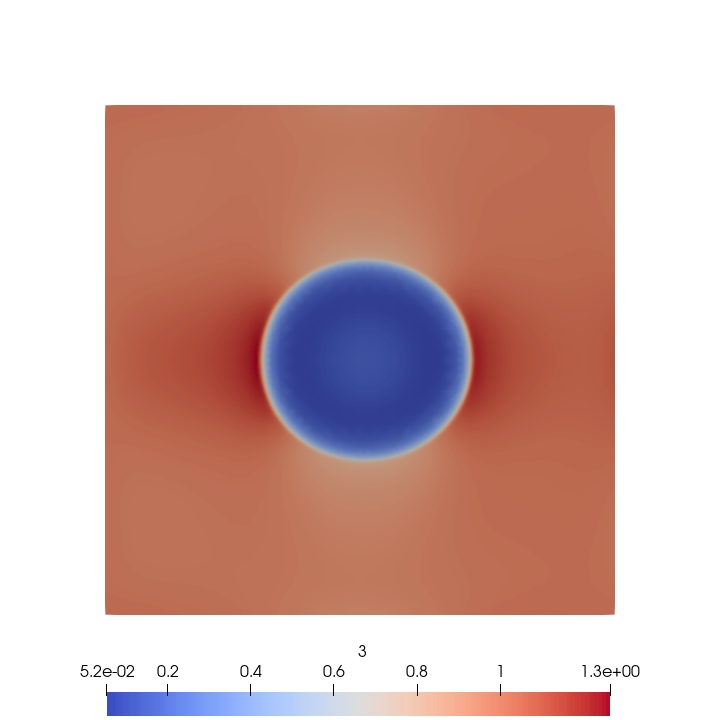}}
\caption{Snapshots of the second example with $\kappa_t=1$. First row: $\phi_h^n$ at $t=0, 5, 10, 13.5$, where $\phi_h^n=1$ (red) in the tumour tissue and $\phi_h^n=-1$ (blue) in the host tissue.
Second row: the nutrient $\sigma_h^n$, $\abs{\bv_h^n}$ (with the velocity field $\bv_h^n$), and both eigenvalues of $\bbB_h^n$ at $t=13.5$.} 
\label{fig:2d_3d}
\end{figure}

In the second example, we neglect the chemotaxis effect with $\chi_\phi=0$ and we include stress sources by growth with $\calG=4$. Moreover, we vary $\kappa_t\in\{-2, -1, 1\}$ and we choose the remaining parameters as in \eqref{eq:parameters}. 
%
%
The results are visualized in Figures \ref{fig:2d_3b}, \ref{fig:2d_3e} and \ref{fig:2d_3d} for the cases $\kappa_t=-2$, $-1$ and $1$, respectively. Compared to the first example, we now observe a different growth behaviour without chemotaxis ($\chi_\phi=0$). The shape of the tumour is now mainly influenced by the choice of $\kappa_t$. 
In the cases $\kappa_t<0$, the tumour grows along the positive gradient of $\trace\bbB$, so that its shape develops to a dumbbell which is more elongated for $\kappa_t=-2$ than for $\kappa_t=-1$.
Compared to before, spatial changes in the eigenvalues of $\bbB$ are intensified due to the additional stress sources (as $\calG>0$), and, as before, the asymmetric growth can be explained with the velocity field. For $\kappa_t>0$, the shape of the tumour is almost spherically symmetric, as now the tumour grows along the negative gradient of $\trace\bbB$ which points inside the tumour.

\begin{figure}[ht!]
\, \\[-4ex]
\centering
\subfloat
{\includegraphics[width=0.22\textwidth]{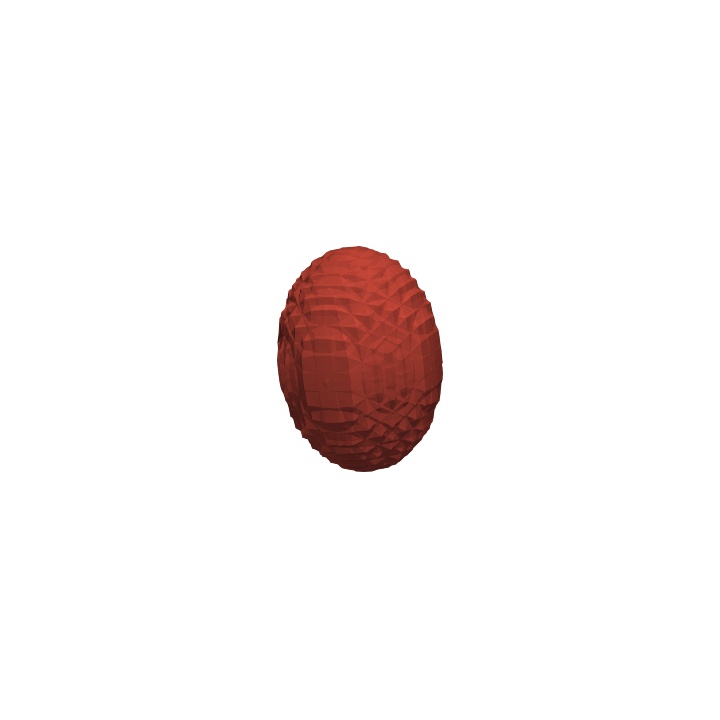}}
\hspace{-2em}
\subfloat
{\includegraphics[width=0.22\textwidth]{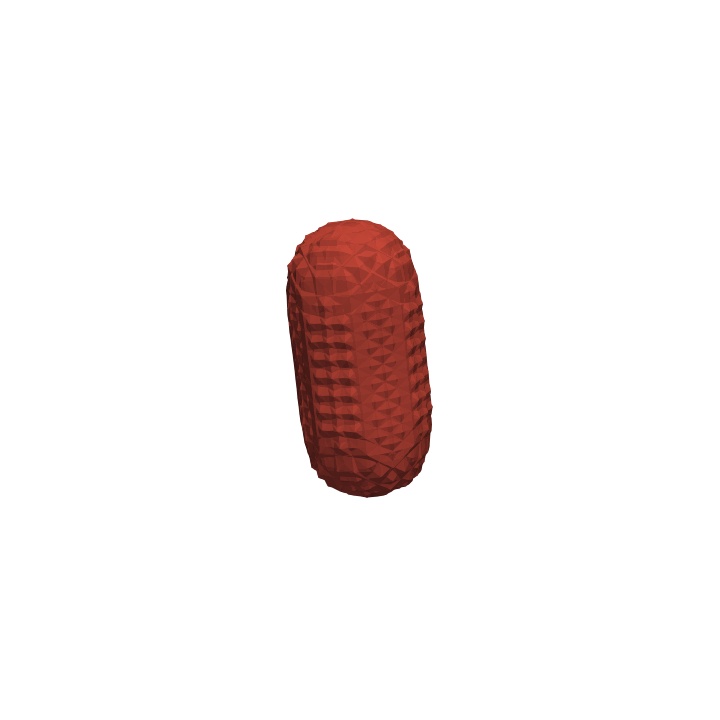}}
\hspace{-2em}
\subfloat
{\includegraphics[width=0.22\textwidth]{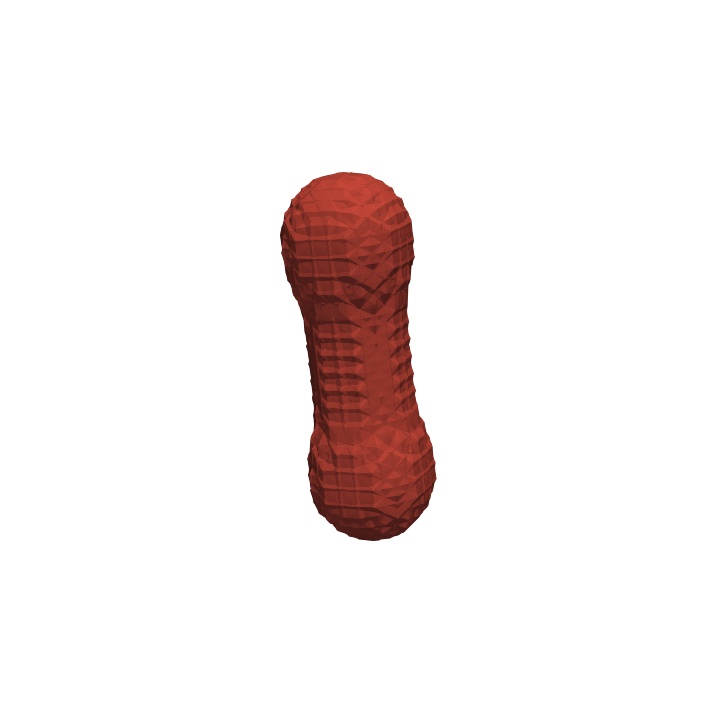}}
\hspace{-2em}
\subfloat
{\includegraphics[width=0.22\textwidth]{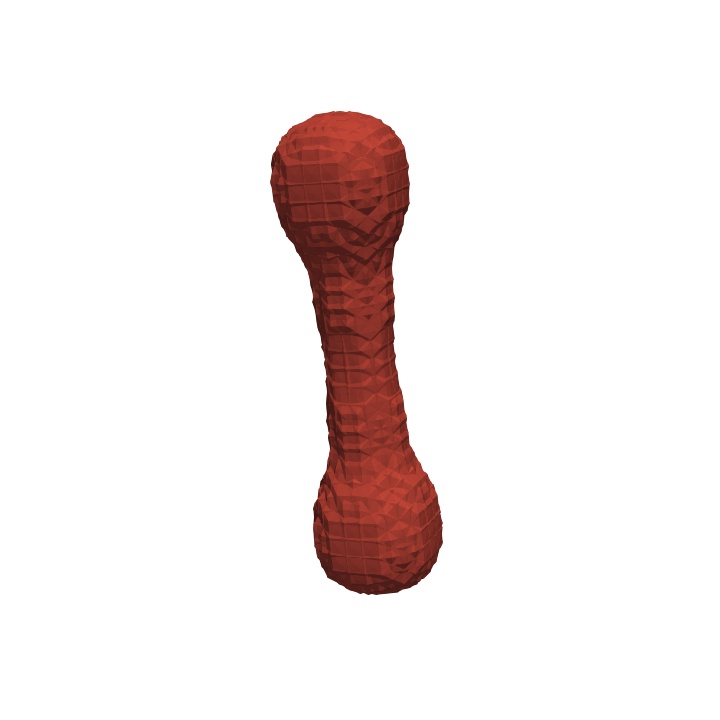}}
\hspace{-2em}
\subfloat
{\includegraphics[width=0.22\textwidth]{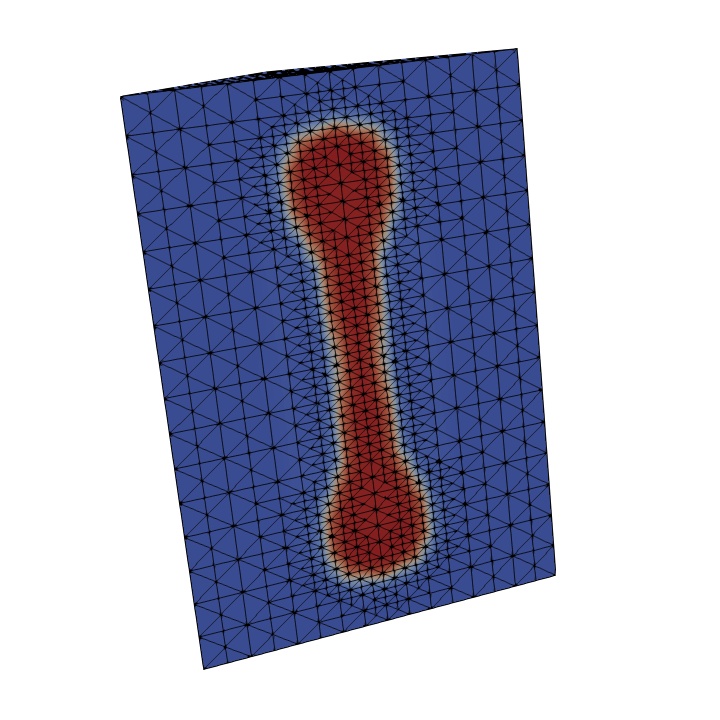}}
\\[-4ex]
\subfloat
{\includegraphics[width=0.22\textwidth]{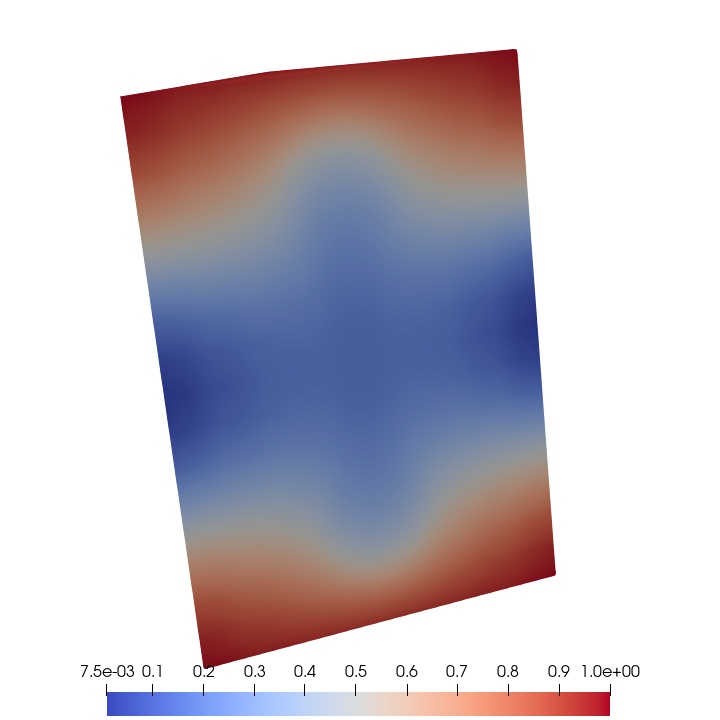}}
\hspace{-2em}
\subfloat
{\includegraphics[width=0.22\textwidth]{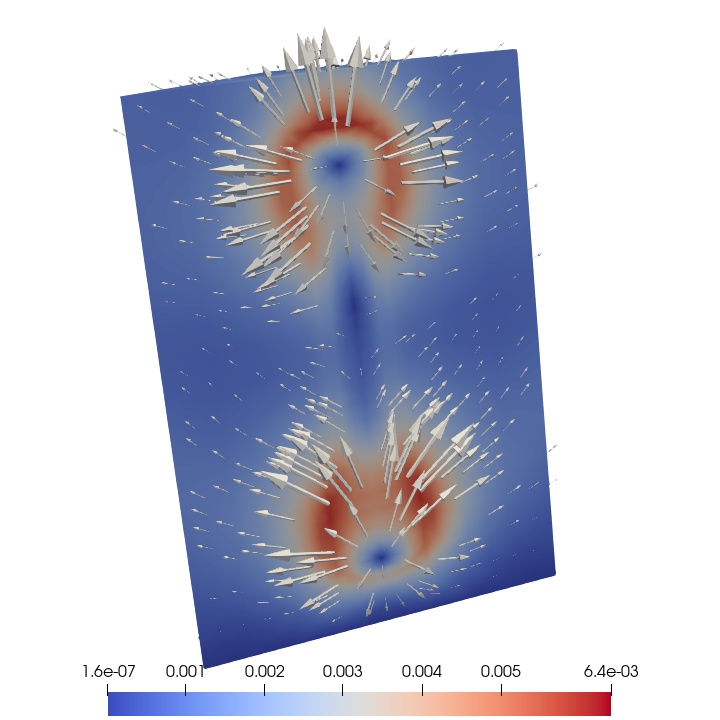}}
\hspace{-2em}
\subfloat
{\includegraphics[width=0.22\textwidth]{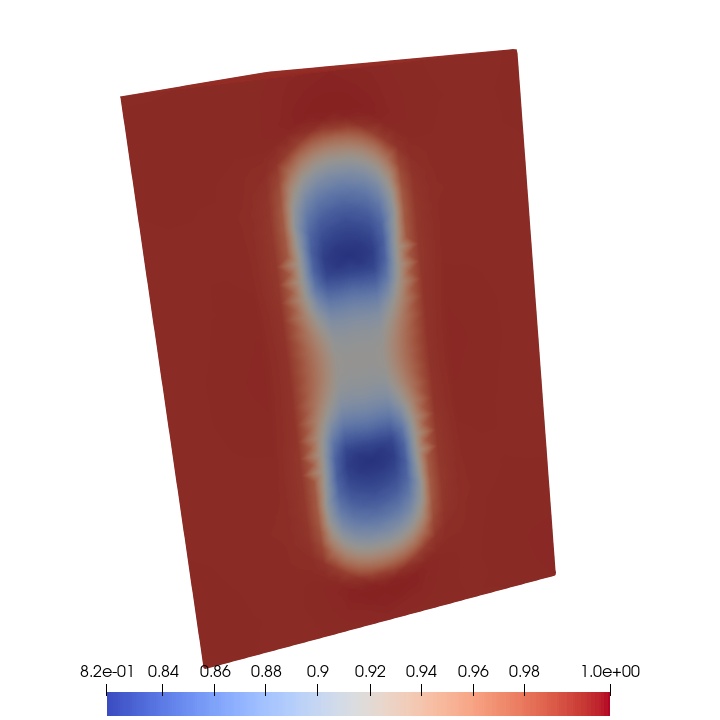}}
\hspace{-2em}
\subfloat
{\includegraphics[width=0.22\textwidth]{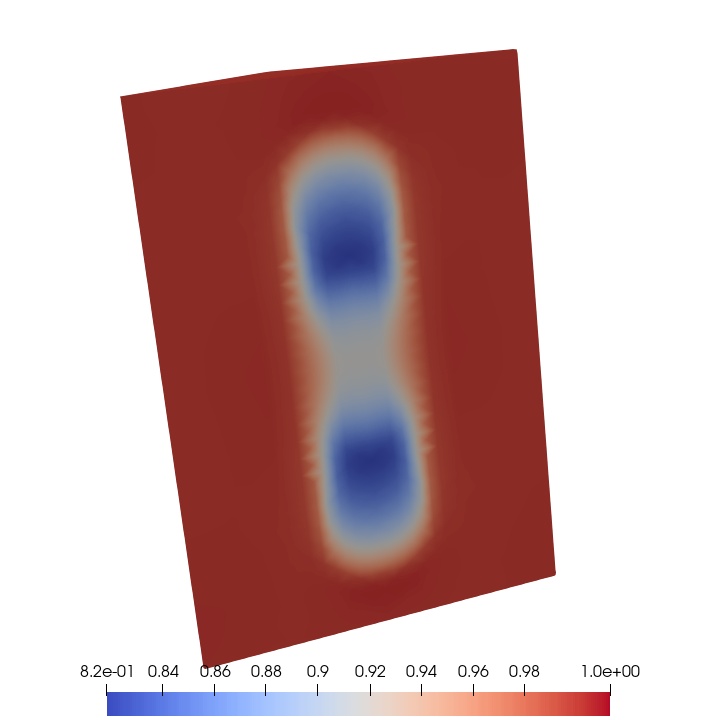}}
\hspace{-2em}
\subfloat
{\includegraphics[width=0.22\textwidth]{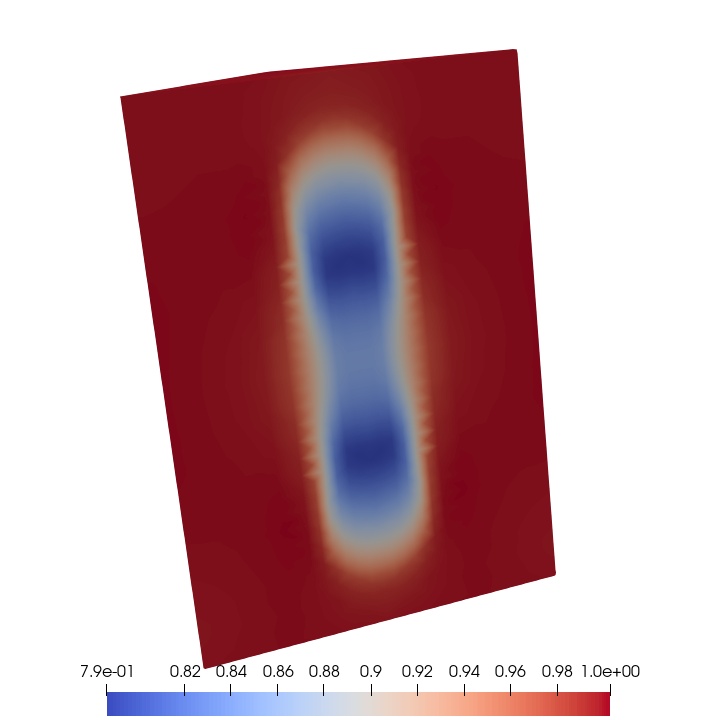}}
\caption{Snapshots of the example in 3d with $\kappa_t=0$. First row: time evolution of the interface of $\phi_h^n$ at $t=0, 1, 2, 3$, as well as $\phi_h^n$ with the corresponding mesh at $t=3$. Second row: the nutrient $\sigma_h^n$, the velocity magnitude $\abs{\bv_h^n}$ (with the velocity field $\bv_h^n$), and the three eigenvalues of $\bbB_h^n$ at $t=3$. A cut was taken along the plane with normal $(1,1,0)^\top$ and origin $(0,0,0)^\top$.} 
\label{fig:3dim}
\end{figure}

\begin{figure}[ht!]
\, \\[-4ex]
\centering
\subfloat
{\includegraphics[width=0.22\textwidth]{figures/3d/3d_1j_phi_0.jpeg}}
\hspace{-2em}
\subfloat
{\includegraphics[width=0.22\textwidth]{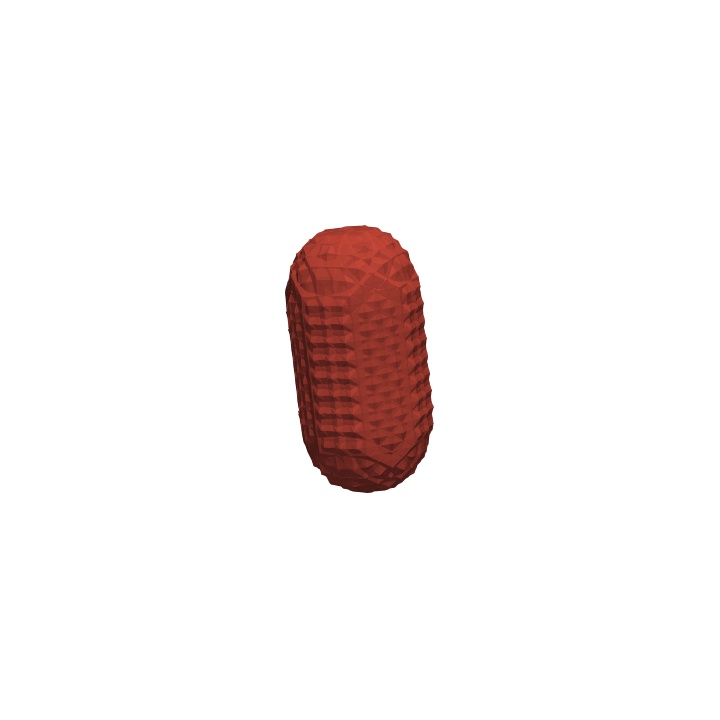}}
\hspace{-2em}
\subfloat
{\includegraphics[width=0.22\textwidth]{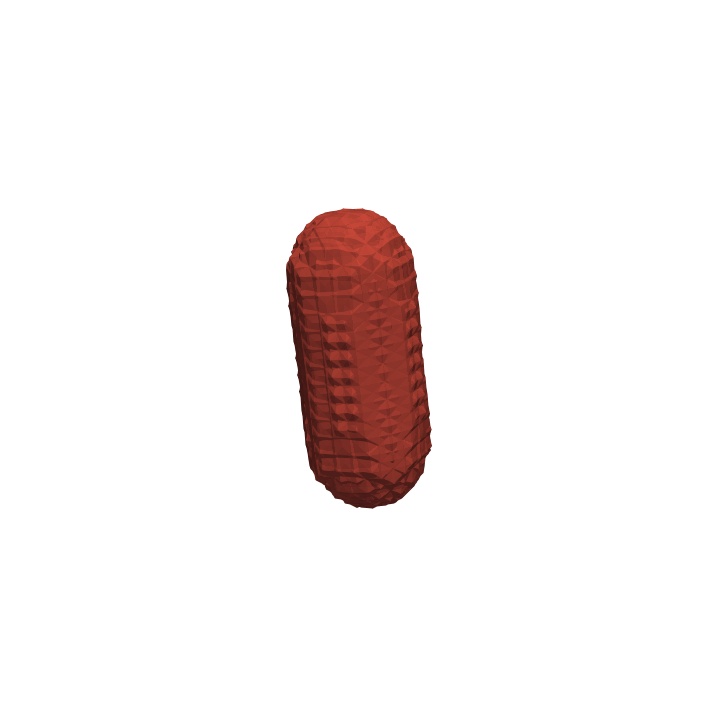}}
\hspace{-2em}
\subfloat
{\includegraphics[width=0.22\textwidth]{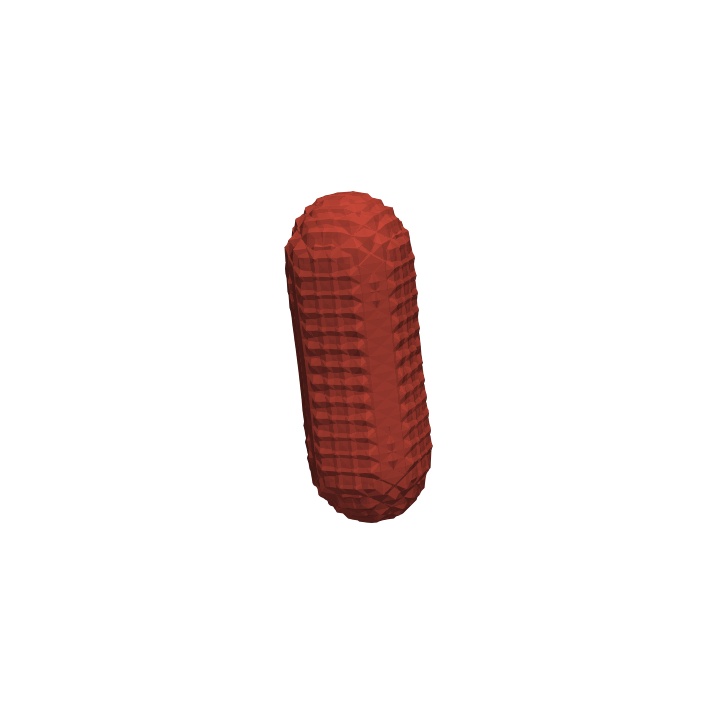}}
\hspace{-2em}
\subfloat
{\includegraphics[width=0.22\textwidth]{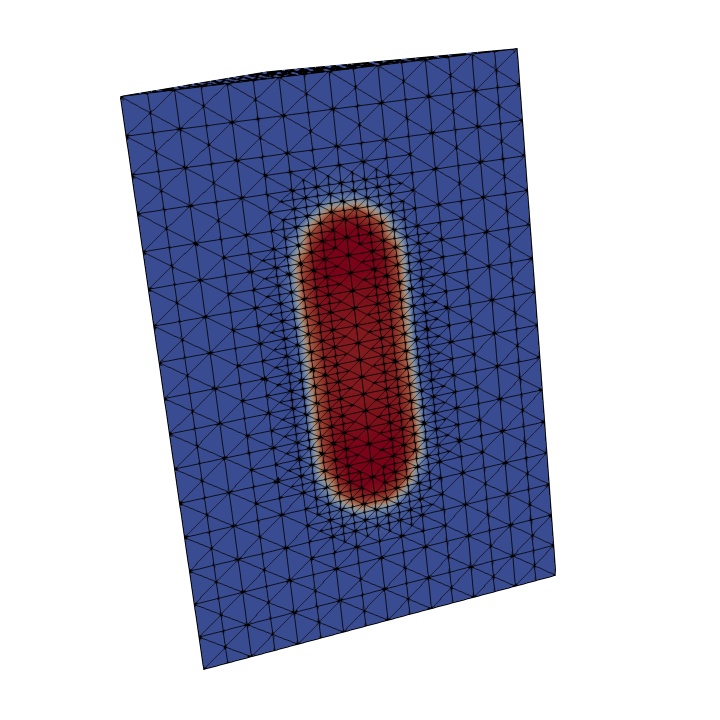}}
\\[-4ex]
\subfloat
{\includegraphics[width=0.22\textwidth]{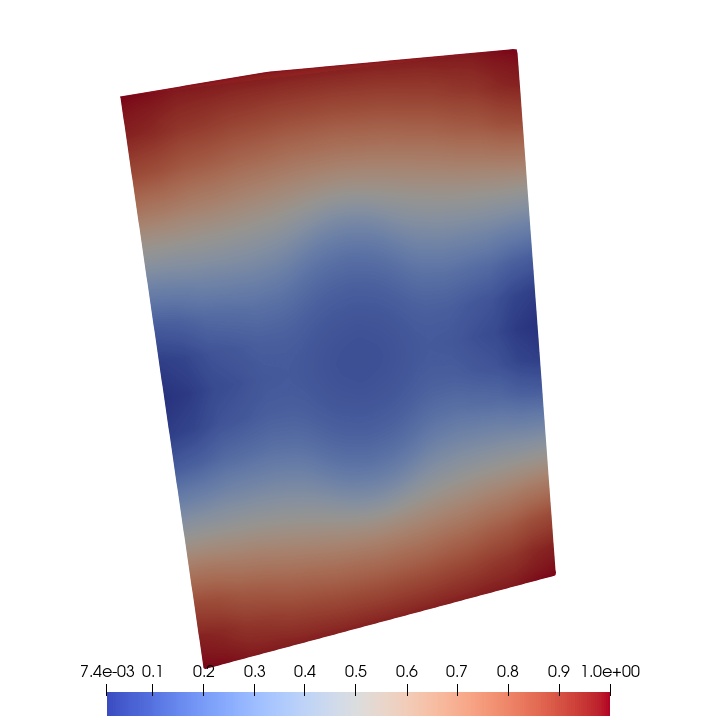}}
\hspace{-2em}
\subfloat
{\includegraphics[width=0.22\textwidth]{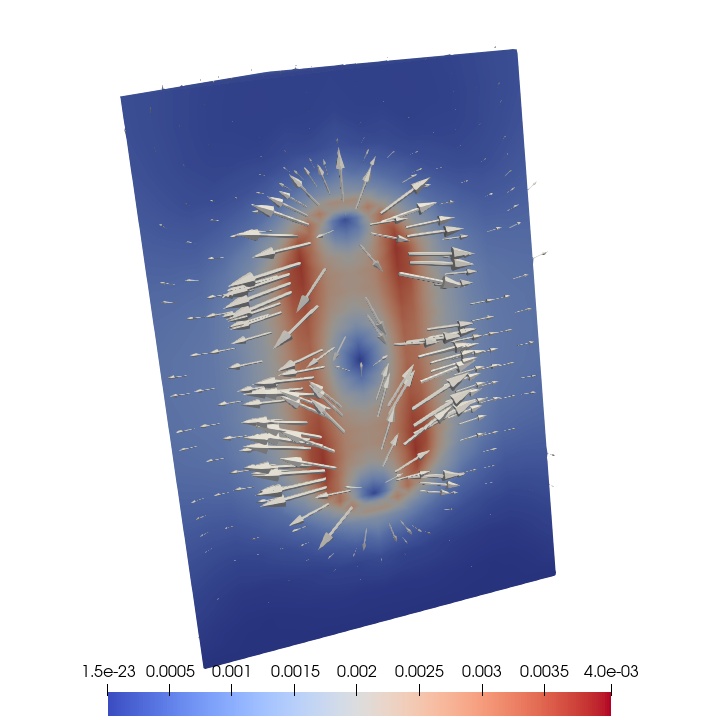}}
\hspace{-2em}
\subfloat
{\includegraphics[width=0.22\textwidth]{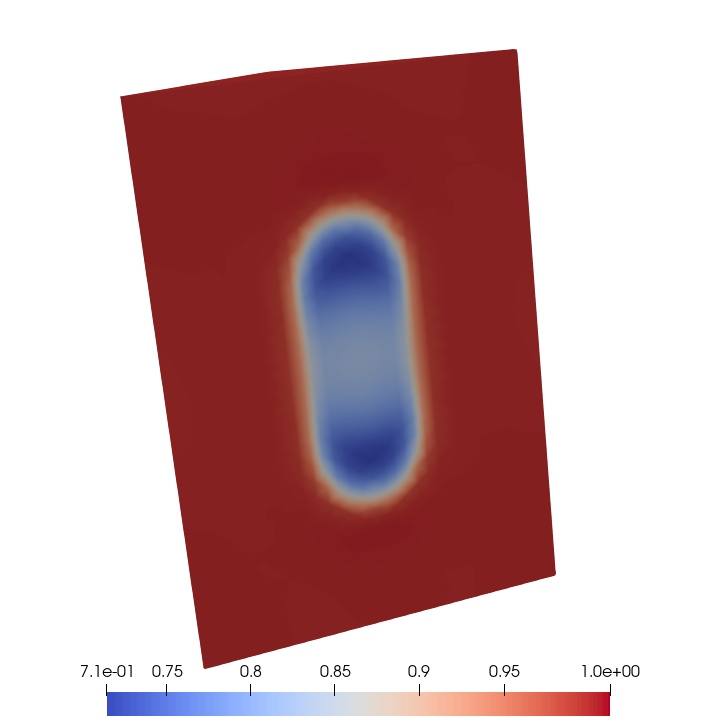}}
\hspace{-2em}
\subfloat
{\includegraphics[width=0.22\textwidth]{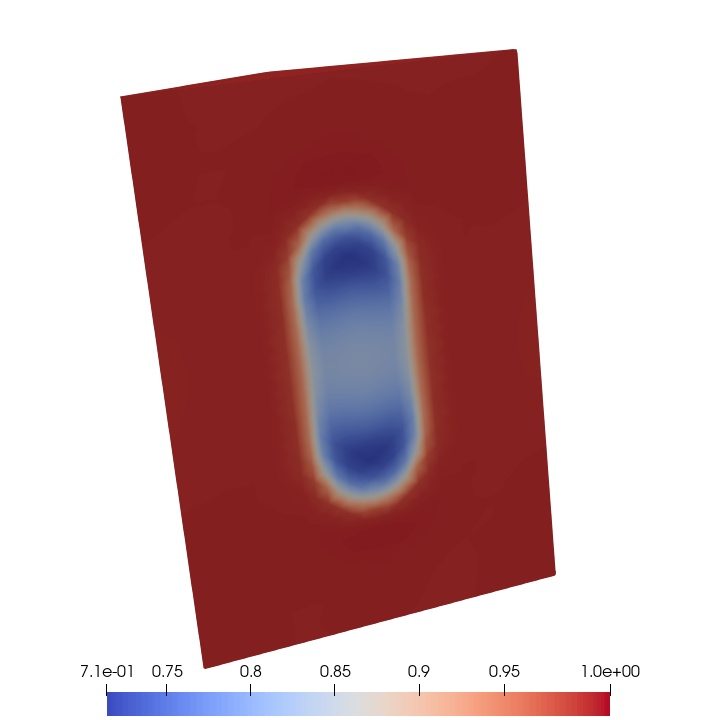}}
\hspace{-2em}
\subfloat
{\includegraphics[width=0.22\textwidth]{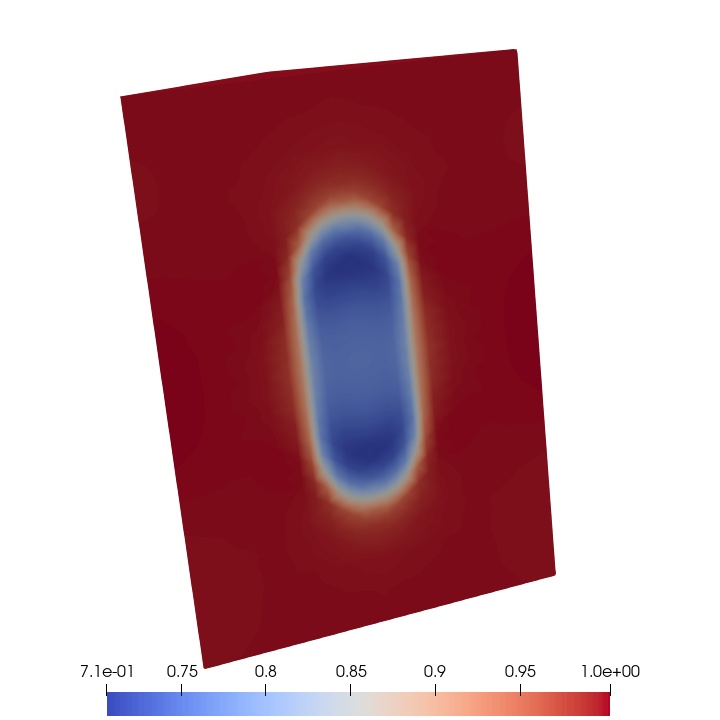}}
\caption{Snapshots of the example in 3d with $\kappa_t=4$. First row: time evolution of the interface of $\phi_h^n$ at $t=0, 1, 2, 3$, as well as $\phi_h^n$ with the corresponding mesh at $t=3$. Second row: the nutrient $\sigma_h^n$, the velocity magnitude $\abs{\bv_h^n}$ (with the velocity field $\bv_h^n$), and the three eigenvalues of $\bbB_h^n$ at $t=3$. A cut was taken along the plane with normal $(1,1,0)^\top$ and origin $(0,0,0)^\top$.} 
\label{fig:3dim_1g}
\end{figure}

\begin{figure}[ht!]
\, \\[-4ex]
\centering
\subfloat
{\includegraphics[width=0.22\textwidth]{figures/3d/3d_1j_phi_0.jpeg}}
\hspace{-2em}
\subfloat
{\includegraphics[width=0.22\textwidth]{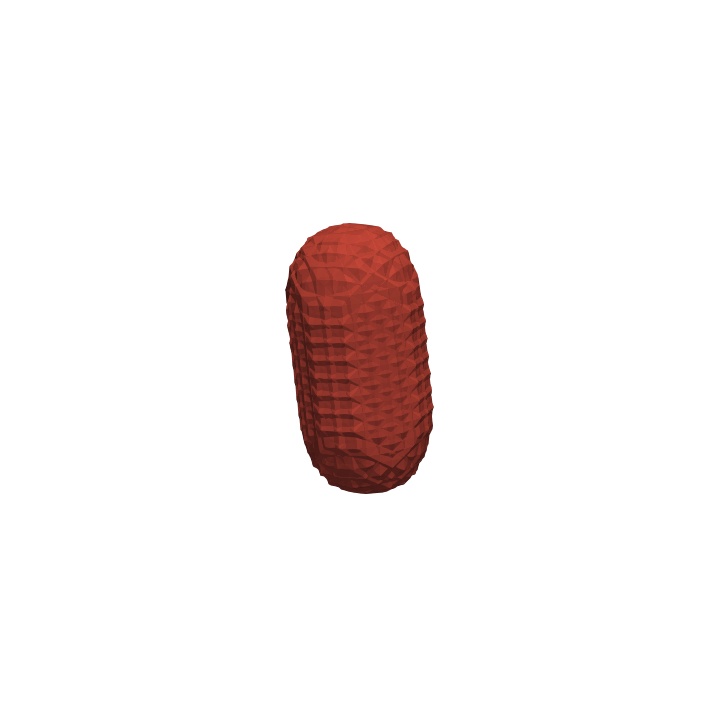}}
\hspace{-2em}
\subfloat
{\includegraphics[width=0.22\textwidth]{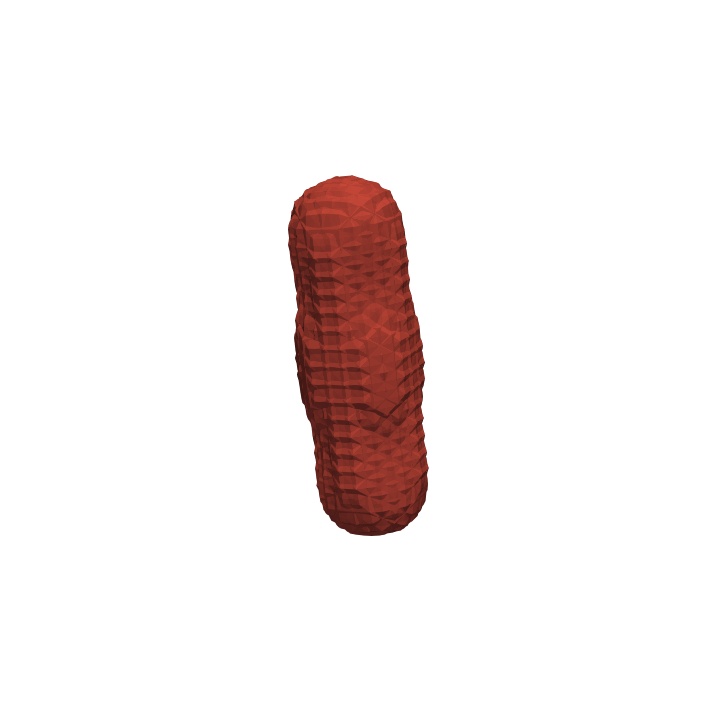}}
\hspace{-2em}
\subfloat
{\includegraphics[width=0.22\textwidth]{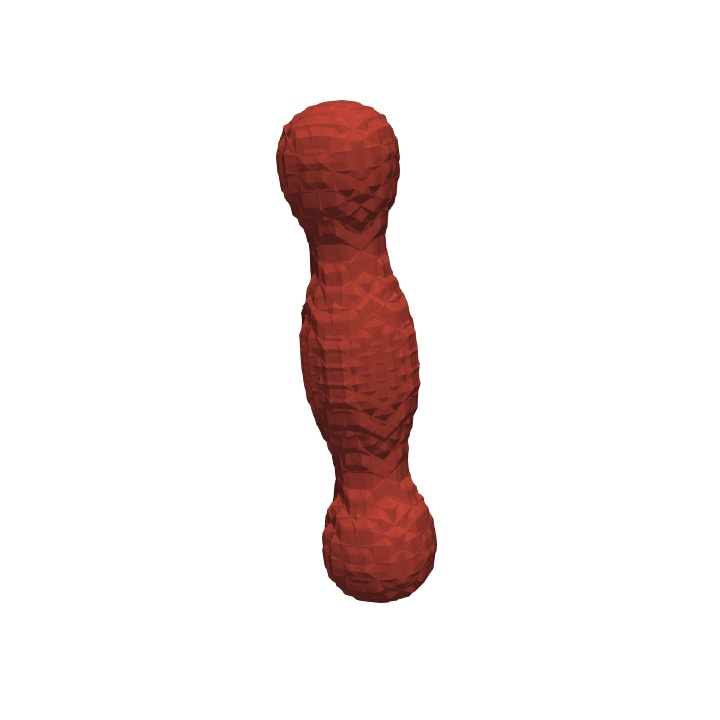}}
\hspace{-2em}
\subfloat
{\includegraphics[width=0.22\textwidth]{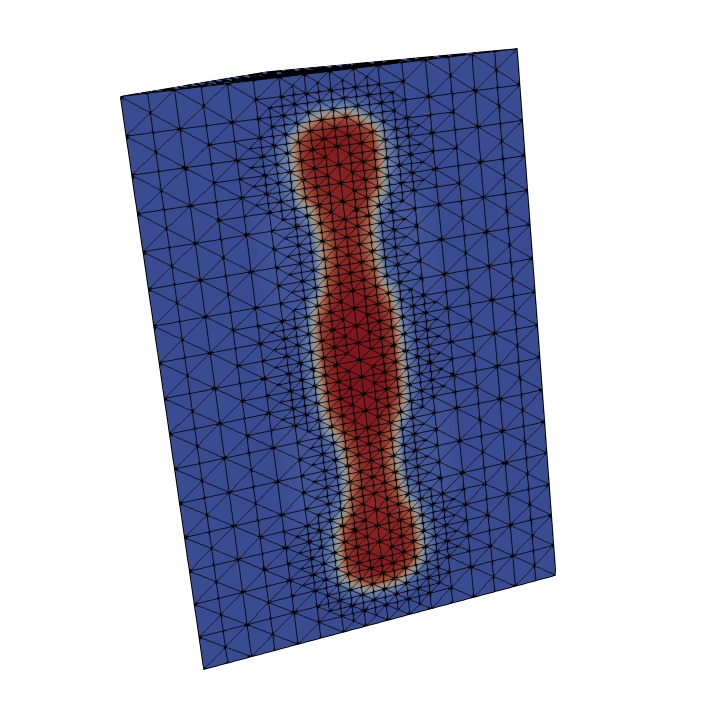}}
\\[-4ex]
\subfloat
{\includegraphics[width=0.22\textwidth]{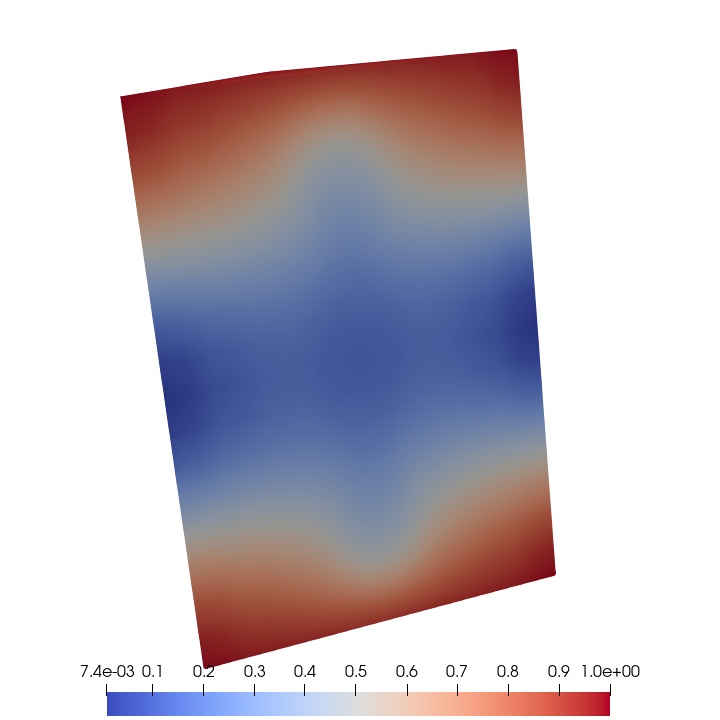}}
\hspace{-2em}
\subfloat
{\includegraphics[width=0.22\textwidth]{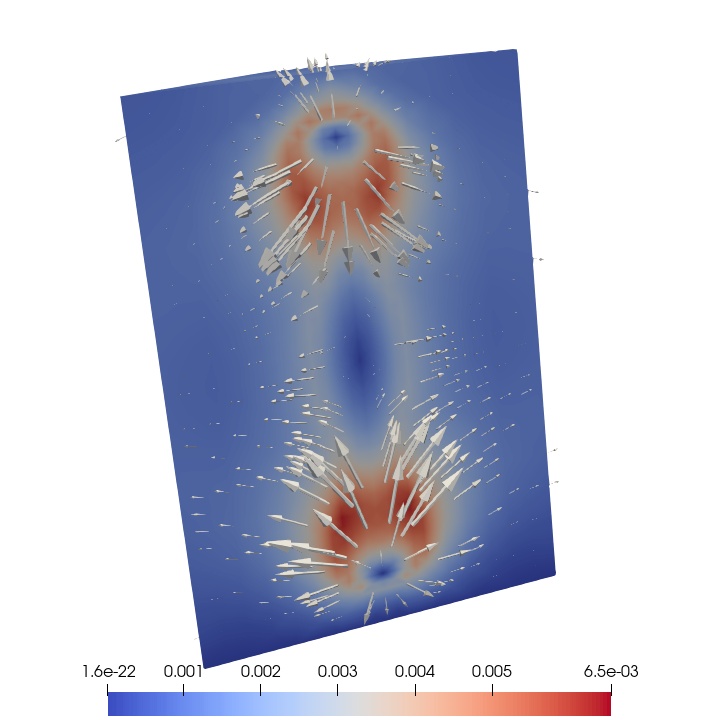}}
\hspace{-2em}
\subfloat
{\includegraphics[width=0.22\textwidth]{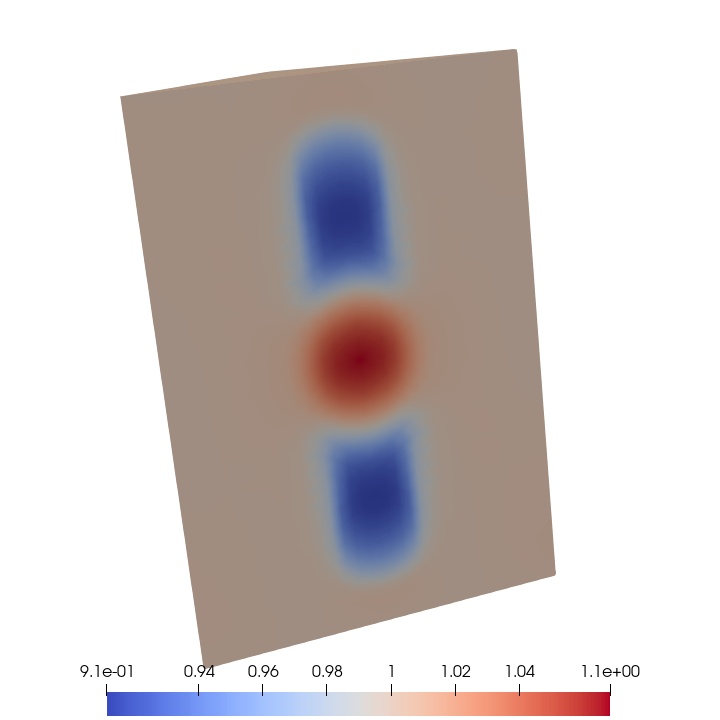}}
\hspace{-2em}
\subfloat
{\includegraphics[width=0.22\textwidth]{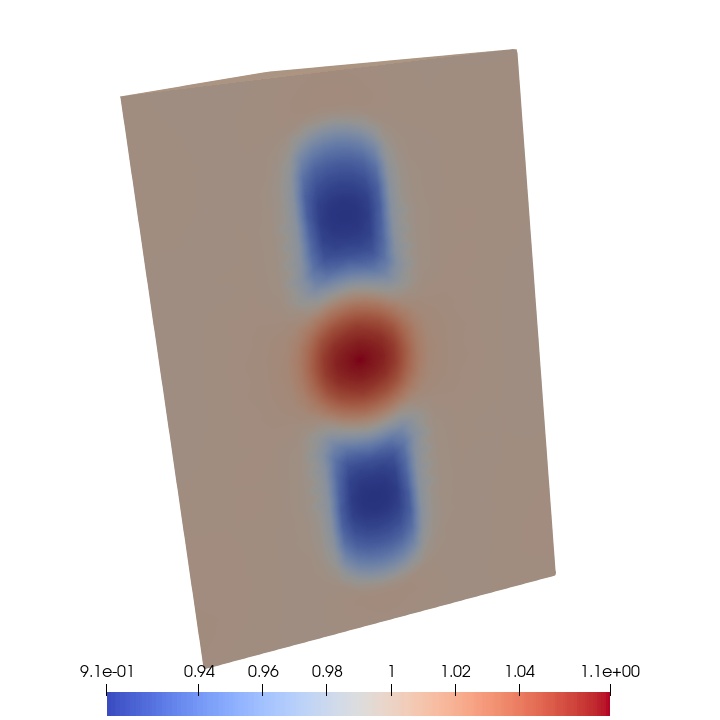}}
\hspace{-2em}
\subfloat
{\includegraphics[width=0.22\textwidth]{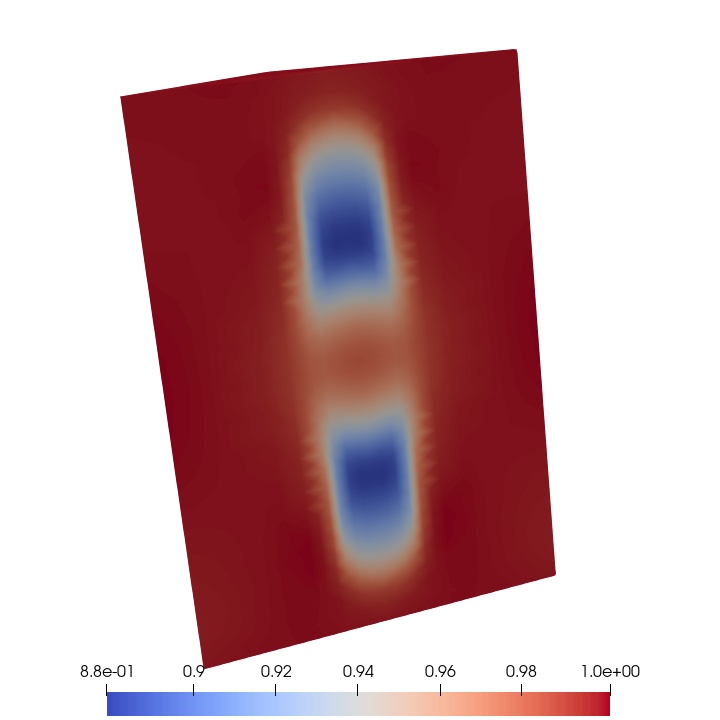}}
\caption{Snapshots of the example in 3d with $\kappa_t=-5$. First row: time evolution of the interface of $\phi_h^n$ at $t=0, 1, 2, 3$, as well as $\phi_h^n$ with the corresponding mesh at $t=3$. Second row: the nutrient $\sigma_h^n$, the velocity magnitude $\abs{\bv_h^n}$ (with the velocity field $\bv_h^n$), and the three eigenvalues of $\bbB_h^n$ at $t=3$. A cut was taken along the plane with normal $(1,1,0)^\top$ and origin $(0,0,0)^\top$.} 
\label{fig:3dim_1i}
\end{figure}

We now present the setting for the three-dimensional example.
We fix the domain as $\Omega=[-1,1]^2 \times [-2, 2] \subset \bbR^3$ and we set $\partial_{\mathrm{D}} \Omega = [-1,1]^2\times\{-2\}$, $\partial_{\mathrm{N}}\Omega = \partial\Omega \setminus \partial_{\mathrm{D}}\Omega$. 
We use $h_c \approx 0.4330$ and $h_f \approx 0.0541$ for the largest and smallest cell diameters, respectively.
We set $\varepsilon=0.04$, $\chi_\phi=2$, $\calP=10$, $\calG=2$, $\kappa_t\in\{0, 4, -5\}$ and we choose the remaining parameters as in \eqref{eq:parameters}.
For the initial and boundary data, we define $\phi_h^0 = \calI_h \phi_0$, $\bbB_h^0 = \calI_h \bbB_0$ and $\sigma_{\infty,h}^n = \calI_h \sigma_\infty$, $n\in\bbN$, where
\begin{align*}
    \phi_0(\mathbf{x}) = - \tanh \Big( \frac{r(\mathbf{x})}{\sqrt{2} \varepsilon} \Big) \in (-1,1), 
    \quad 
    \bbB_0(\mathbf{x}) = \bbI,
    \quad
    \sigma_\infty(\mathbf{x}) = -\frac12 \cos \big(\frac{\pi}{2} x_3 \big) + \frac12 \in [0,1],
\end{align*}
for all $\mathbf{x} = (x_1,x_2,x_3)^\top \in \overline\Omega$, where $r(\mathbf{x}) = \big( (2 x_1)^2 +  (2x_2)^2 + (1.4 x_3)^2 \big)^{1/2} - 1$.
To reduce the total number of degrees of freedom, we use the inf--sup stable mini-element \cite{girault_raviart_1986} for the velocity-pressure approximation in three space dimensions instead of the $\calP_2$/$\calP_1$-Taylor--Hood element $\calV_h \times \calS_h$. As remarked in Section \ref{sec:approximation}, the analysis of this work is also valid for the mini-element.
We visualize the numerical results of the example in three dimensions in Figures \ref{fig:3dim}, \ref{fig:3dim_1g} and \ref{fig:3dim_1i} for the cases $\kappa_t=0$, $4$ and $-5$, respectively. In the first row of each figure, we show the time evolution of the interface of $\phi_h^n$ at times $t\in\{0, 1, 2, 3\}$ as well as the order parameter $\phi_h^n$ at time $t=3$ with the corresponding mesh. For all cases, one observes that the tumour develops from an ellipsoid at $t=0$ to an elongated dumbbell at $t=3$. However, the shape of the dumbbell at the final time $t=3$ differs in each case. The tumour is either less elongated ($\kappa_t=4$) or more accumulated in the center ($\kappa_t=-5$), compared to the case $\kappa_t=0$.
In the second row in each figure, from left to right, we plot the nutrient $\sigma_h^n$ and the velocity magnitude $\abs{\bv_h^n}$ together with the corresponding velocity field $\bv_h^n$, and the spatial distribution of the three eigenvalues of $\bbB_h^n$ at time $t=3$. Note that a cut was taken through the domain along the plane with normal $(1,1,0)^\top$ and origin $(0,0,0)^\top$.

\section*{Acknowledgments}
The authors gratefully acknowledge the support by the Graduiertenkolleg 2339 IntComSin of the Deutsche Forschungsgemeinschaft (DFG, German Research Foundation) -- Project-ID 321821685. 
The authors also acknowledge discussions with Pierluigi Colli on viscoelastic effects in tumour growth modelling.


\addcontentsline{toc}{section}{References}
\printbibliography

@article{garcke_lam_2017, 
    title={{Well--posedness of a Cahn--Hilliard system modelling tumour growth with chemotaxis and active transport}}, 
    volume={28}, 
    DOI={10.1017/S0956792516000292}, 
    number={2}, 
    journal={European J. Appl. Math.},
    publisher={Cambridge University Press}, 
    author={Garcke, H. and Lam, K. F.}, 
    year={2017}, 
    pages={284–316}
}

@book{evans_2010,
  address = {Providence, R.I.},
  author = {Evans, L. C.},
  publisher = {American Mathematical Society},
  refid = {465190110},
  title = {Partial Differential Equations},
  year = {2010},
  edition = {2}
}

@article{GarckeLSS_2016,
author = {Garcke, H. and Lam, K. F. and Sitka, E. and Styles, V.},
title = {{A Cahn--Hilliard--Darcy model for tumour growth with chemotaxis and active transport}},
journal = {Math. Models Methods Appl. Sci.},
volume = {26},
number = {06},
pages = {1095-1148},
year = {2016},
doi = {10.1142/S0218202516500263},
}

@book{bartels_2016,
  author = {Bartels, S.},
  title = {Numerical Approximation of Partial Differential Equations},
  series = {Texts in Applied Mathematics},
  volume = {64},
  publisher = {Springer, [Cham]},
  year = {2016},
  pages = {xv+535},
  mrclass = {65-01 (65M60 65N30)},
  mrnumber = {3496530},
  doi = {10.1007/978-3-319-32354-1},
}

@inbook{byrne_tumour_2003,
author = {Byrne, H. M.},
year = {2003},
pages = {75-120},
title = {Modelling Avascular Tumour Growth},
booktitle = {Cancer Modelling and Simulation},
publisher = {CRC Press},
}

@book{fenics_book_2012,
  title = {Automated Solution of Differential Equations by the Finite Element Method},
  author = {A. Logg and K. A. Mardal and G. N. Wells and others},
  year = {2012},
  publisher = {Springer},
  doi = {10.1007/978-3-642-23099-8},
}

@article{simon_1986,
author = {Simon, J.},
year = {1986},
pages = {65--96},
title = {{Compact sets in the space $L^p(0,T; B)$}},
volume = {146},
JOURNAL = {Ann. Mat. Pura Appl. (4)},
FJOURNAL = {Annali di Matematica Pura ed Applicata. Serie Quarta},
doi = {10.1007/BF01762360}
}

@article{barrett_nurnberg_styles_2004,
author = {Barrett, J. W. and N{\"u}rnberg, R. and Styles, V.},
title = {Finite element approximation of a phase field model for void electromigration},
journal = {SIAM J. Num. Anal.},
volume = {42},
number = {2},
pages = {738-772},
year = {2004},
doi = {10.1137/S0036142902413421},
}

@article{barrett_nurnberg_2004,  
author={J. W. {Barrett} and R. {N{\"u}rnberg}},  
journal={IMA J. Numer. Anal.},
title={{Convergence of a finite‐element approximation of surfactant spreading on a thin film in the presence of van der Waals forces}},
year={2004},  
volume={24},  
number={2}, 
pages={323-363},  
doi={10.1093/imanum/24.2.323},  
}

@article{ebenbeck_garcke_nurnberg_2020,
title = {{Cahn--Hilliard--Brinkman systems for tumour growth}},
journal = {Discrete Contin. Dyn. Syst. Ser. S},
volume = {14},
number = {11},
pages = {3989--4033},
DOI={10.3934/dcdss.2021034},
year = {2021},
author = {M. Ebenbeck and H. Garcke and R. N{\"u}rnberg},
publisher={American Institute of Mathematical Sciences},
}

@article{GLNS_2018_multiphase_tumour_necrosis,
   title={{A multiphase Cahn--Hilliard--Darcy model for tumour growth with necrosis}},
   volume={28},
   DOI={10.1142/s0218202518500148},
   number={03},
   journal={Math. Models Methods Appl. Sci.},
   author={Garcke, H. and Lam, K. F. and N{\"u}rnberg, R. and Sitka, E.},
   year={2018},
   pages={525–577}
}

@article{AbelsGG_2012,
author = {Abels, H. and Garcke, H. and Gr{\"u}n, G.},
title = {Thermodynamically consistent, frame indifferent diffuse interface models for incompressible two-phase flows with different densities},
journal = {Math. Models Methods Appl. Sci.}, 
volume = {22},
number = {03},
pages = {1150013, 40},
year = {2012},
doi = {10.1142/S0218202511500138},
}

@article{bresch_2009,
author = {Bresch, D. and Colin, T. and Grenier, E. and Ribba, B. and Saut, O.},
title = {A viscoelastic model for avascular tumor growth},
JOURNAL = {Discrete Contin. Dyn. Syst.},
FJOURNAL = {Discrete and Continuous Dynamical Systems. Series A},
NUMBER = {Dynamical systems, differential equations and applications.
              7th AIMS Conference, suppl.},
PAGES = {101--108},
year = {2009},
doi={10.3934/proc.2009.2009.101}
}

@inbook{malek_prusa_2018,
author = {M{\'a}lek, J. and Pr{\r u}{\v s}a, V.},
year = {2018},
pages = {3--72},
title = {Derivation of equations for continuum mechanics and thermodynamics of fluids},
booktitle = {Handbook of Mathematical Analysis in Mechanics of Viscous Fluids},
publisher = {Springer International Publishing},
doi = {10.1007/978-3-319-13344-7_1}
}

@article{barrett_boyaval_2009,
author = {Barrett, J. W. and Boyaval, S.},
year = {2011},
pages = {1783-1837},
title = {Existence and approximation of a (regularized) {Oldroyd-B} model},
volume = {21},
number = {9},
journal = {Math. Models Methods Appl. Sci.},
doi = {10.1142/S0218202511005581}
}

@book{dahmen_reusken_numerik,
  author={Dahmen, W. and Reusken, A.},
  title={Numerik f{\"u}r Ingenieure und Naturwissenschaftler},
  edition={2}, 
  address={Berlin, Heidelberg},
  year={2008},
  publisher = {Springer},
  doi={10.1007/978-3-540-76493-9},
}

@book {girault_raviart_1986,
    AUTHOR = {Girault, V. and Raviart, P. A.},
     TITLE = {Finite element methods for {N}avier-{S}tokes equations},
    SERIES = {Springer Series in Computational Mathematics},
    VOLUME = {5},
      NOTE = {Theory and algorithms},
 PUBLISHER = {Springer-Verlag, Berlin},
      YEAR = {1986},
     PAGES = {x+374},
      ISBN = {3-540-15796-4},
   MRCLASS = {65N30 (65-02 76-08)},
  MRNUMBER = {851383},
MRREVIEWER = {Max D. Gunzburger},
       DOI = {10.1007/978-3-642-61623-5},
}

@article{barrett_lu_sueli_2017,
  title={{Existence of large-data finite-energy global weak solutions to a compressible Oldroyd-B model}},
  author={Barrett, J. W. and Lu, Y. and S{\"u}li, E.},
  journal={Commun. Math. Sci.}, 
  volume={15},
  number={5},
  year={2017},
  pages={1265--1323},
  publisher={International Press},
  doi={10.4310/CMS.2017.v15.n5.a5}
}

@article{sieber_2020,
author = {O. Sieber},
doi = {10.1515/jnma-2019-0019},
title = {{On convergent schemes for a two-phase Oldroyd-B type model with variable polymer density}},
fjournal = {Journal of Numerical Mathematics},
journal = {J. Numer. Math.},
number = {2},
volume = {28},
year = {2020},
pages = {99--129}
}

@article{barrett_2018_fene-p,
  title={{Finite element approximation of the FENE-P model}},
  author={Barrett, J. W. and Boyaval, S.},
  journal={IMA J. Numer. Anal.},
  volume={38},
  number={4},
  pages={1599--1660},
  year={2018},
  publisher={Oxford University Press},
  doi = {10.1093/imanum/drx061},
}

@incollection{frigeri_2017_tumour_degenerate,
  title={On a diffuse interface model for tumour growth with non-local interactions and degenerate mobilities},
  author={Frigeri, S. and Lam, K. F. and Rocca, E.},
  booktitle={Solvability, regularity, and optimal control of boundary value problems for PDEs},
  pages={217--254},
  year={2017},
  publisher={Springer}
}

@article{metzger_2018,
  title={On convergent schemes for two-phase flow of dilute polymeric solutions},
  author={Metzger, S.},
  journal={ESAIM: Math. Model. Numer. Anal.}, 
  volume={52},
  number={6},
  pages={2357--2408},
  year={2018},
  publisher={EDP Sciences},
  doi={10.1051/m2an/2018042}
}

@book{temam_2001,
  title={Navier--Stokes Equations: Theory and Numerical Analysis},
  author={Temam, R.},
  series={AMS/Chelsea publication},
  year={2001},
  publisher={AMS Chelsea Pub.}
}

@article{Lukacova_2017,
author = {Luk{\'a}{\v c}ov{\'a}-Medvid’ov{\'a}, M. and Mizerov{\'a}, H. and Ne{\v c}asov{\'a}, S. and Renardy, M.},
title = {{Global existence result for the generalized Peterlin viscoelastic model}},
journal = {SIAM J. Math. Anal.},
volume = {49},
number = {4},
pages = {2950-2964},
year = {2017},
doi = {10.1137/16M1068505},
}

@book{brenner_scott_2008,
  title={The Mathematical Theory of Finite Element Methods},
  author={Brenner, S. C. and Scott, L. R.},
  volume={3},
  year={2008},
  publisher={Springer}
}

@article{wise_lowengrub_2008,
title = {{Three-dimensional multispecies nonlinear tumor growth---I: Model and numerical method}},
journal = {J. Theoret. Biol.},
volume = {253},
number = {3},
pages = {524-543},
year = {2008},
doi = {10.1016/j.jtbi.2008.03.027},
author = {S. M. Wise and J. S. Lowengrub and H. B. Frieboes and V. Cristini},
}

@article{garcke_lam_signori_2021,
title = {{On a phase field model of Cahn--Hilliard type for tumour growth with mechanical effects}},
journal = {Nonlinear Anal. Real World Appl.},
volume = {57},
pages = {103192},
year = {2021},
doi = {https://doi.org/10.1016/j.nonrwa.2020.103192},
author = {H. Garcke and K. F. Lam and A. Signori},
}

@article{hawkins_2010,
author = {Oden, J. T. and Hawkins, A. and Prudhomme, S.},
title = {General diffuse-interface theories and an approach to predictive tumor growth modeling},
journal = {Math. Models Methods Appl. Sci.},
volume = {20},
number = {03},
pages = {477-517},
year = {2010},
doi = {10.1142/S0218202510004313},
}

@TechReport{petsc-user-ref_2021,
Author      = {S. Balay and others},
Title       = {{PETS}c Users Manual},
Number      = {ANL-95/11 - Revision 3.15},
Institution = {Argonne National Laboratory},
Year        = {2021}
}

@article{ebenbeck_2019_analysis,
title = {Analysis of a {Cahn--Hilliard--Brinkman} model for tumour growth with chemotaxis},
journal = {J. Differential Equations},
volume = {266},
number = {9},
pages = {5998-6036},
year = {2019},
doi = {10.1016/j.jde.2018.10.045},
author = {M. Ebenbeck and H. Garcke},
}

@article {knopf_2022_multiphase,
    AUTHOR = {Knopf, P. and Signori, A.},
     TITLE = {Existence of weak solutions to multiphase
              {C}ahn-{H}illiard-{D}arcy and {C}ahn-{H}illiard-{B}rinkman
              models for stratified tumor growth with chemotaxis and general
              source terms},
   JOURNAL = {Comm. Partial Differential Equations},
  FJOURNAL = {Communications in Partial Differential Equations},
    VOLUME = {47},
      YEAR = {2022},
    NUMBER = {2},
     PAGES = {233--278},
       DOI = {10.1080/03605302.2021.1966803},
}

@article{ambrosi_2009,
title={Cell adhesion mechanisms and stress relaxation in the mechanics of tumours},
author={D. Ambrosi and L. Preziosi},
year={2009},
journal={Biomech. Model. Mechanobiol.},
volume={8},
number={5},
pages={397--413},
doi={10.1007/s10237-008-0145-y}
}

@article{lowengrub_2021_viscoelastic,
  title={Stress generation, relaxation and size control in confined tumor growth},
  author={Yan, H. and Ramirez-Guerrero, D. and Lowengrub, J. and Wu, M.},
  journal = {PLoS. Comput. Biol},
  volume={17},
  number={12},
  year={2021},
  pages={e1009701},
  doi={10.1371/journal.pcbi.1009701}
}

@misc{garcke_2022_viscoelastic,
  title={A {Cahn--Hilliard} model coupled to viscoelasticity with large deformations},
  author={Agosti, A. and Colli, P. and Garcke, H. and Rocca, E.},
  year={2022},
  eprint={2204.04951},
  archivePrefix={arXiv},
  primaryClass={math.AP}
}

@misc{garcke_2023_viscoelastic,
  title={A {Cahn--Hilliard} phase field model coupled to an {Allen--Cahn} model of viscoelasticity at large strains},
  author={Agosti, A. and Colli, P. and Garcke, H. and Rocca, E.},
  year={2023},
  eprint={2301.08341},
  archivePrefix={arXiv},
  primaryClass={math.AP}
}

@article{trautwein_2021,
      title={Numerical analysis for a {Cahn--Hilliard} system modelling tumour growth with chemotaxis and active transport}, 
      author={H. Garcke and D. Trautwein},
      year={2022},
      journal={J. Numer. Math.},
      fjournal={Journal of Numerical Mathematics},
      volume={30},
      number={4},
      pages={295--324},
      doi={10.1515/jnma-2021-0094}
}

@article{GKT_2022_viscoelastic,
    title={Viscoelastic {Cahn--Hilliard} models for tumour growth},
    author={H. Garcke and B. Kov{\'a}cs and D. Trautwein},
    year={2022},
    journal={Math. Models Methods Appl. Sci.},
    fjournal={Mathematical Models and Methods in Applied Sciences},
    VOLUME = {32},
    YEAR = {2022},
    NUMBER = {13},
    PAGES = {2673--2758},
    doi={10.1142/S0218202522500634}
}

@article {bathory_2021_viscoelastic_3D,
    AUTHOR = {Bathory, M. and Bul\'{i}\v{c}ek, M. and M\'{a}lek, J.},
     TITLE = {Large data existence theory for three-dimensional unsteady
              flows of rate-type viscoelastic fluids with stress diffusion},
   JOURNAL = {Adv. Nonlinear Anal.},
  FJOURNAL = {Advances in Nonlinear Analysis},
    VOLUME = {10},
      YEAR = {2021},
    NUMBER = {1},
     PAGES = {501--521},
       DOI = {10.1515/anona-2020-0144},
}

@article {barrett_suli_2011,
    AUTHOR = {Barrett, J. W. and S\"{u}li, E.},
     TITLE = {Finite element approximation of kinetic dilute polymer models
              with microscopic cut-off},
   JOURNAL = {ESAIM Math. Model. Numer. Anal.},
  FJOURNAL = {ESAIM. Mathematical Modelling and Numerical Analysis},
    VOLUME = {45},
      YEAR = {2011},
    NUMBER = {1},
     PAGES = {39--89},
       DOI = {10.1051/m2an/2010030},
}

@book{galdi_2011,
    AUTHOR = {Galdi, G. P.},
     TITLE = {An introduction to the mathematical theory of the {N}avier-{S}tokes equations},
    SERIES = {Springer Monographs in Mathematics},
   EDITION = {2},
      NOTE = {Steady-state problems},
 PUBLISHER = {Springer, New York},
      YEAR = {2011},
     PAGES = {xiv+1018},
       DOI = {10.1007/978-0-387-09620-9},
}

@book {ern_guermond_2004,
    AUTHOR = {Ern, A. and Guermond, J. L.},
     TITLE = {Theory and practice of finite elements},
    SERIES = {Applied Mathematical Sciences},
    VOLUME = {159},
 PUBLISHER = {Springer-Verlag, New York},
      YEAR = {2004},
     PAGES = {xiv+524},
      ISBN = {0-387-20574-8},
   MRCLASS = {65-02 (65M60 65N30 74S05 76M10 78M10)},
  MRNUMBER = {2050138},
MRREVIEWER = {R. S. Anderssen},
       DOI = {10.1007/978-1-4757-4355-5},
}

@book {BBF2013_fem,
    AUTHOR = {Boffi, D. and Brezzi, F. and Fortin, M.},
     TITLE = {Mixed finite element methods and applications},
    SERIES = {Springer Series in Computational Mathematics},
    VOLUME = {44},
 PUBLISHER = {Springer, Heidelberg},
      YEAR = {2013},
     PAGES = {xiv+685},
   MRCLASS = {65-02 (65M60 65N30)},
  MRNUMBER = {3097958},
MRREVIEWER = {Beny Neta},
       DOI = {10.1007/978-3-642-36519-5},
}

@article{RAJAGOPAL2000207,
title = {A thermodynamic frame work for rate type fluid models},
journal = {Journal of Non-Newtonian Fluid Mechanics},
volume = {88},
number = {3},
pages = {207-227},
year = {2000},
issn = {0377-0257},
doi = {https://doi.org/10.1016/S0377-0257(99)00023-3},
author = {K.R. Rajagopal and A.R. Srinivasa},
}

@misc{LT_2022_arxiv,
    title={On a diffuse interface model for incompressible viscoelastic two-phase flows},
    author={Y. Liu and D. Trautwein},
    year={2022},
    eprint={2212.13507},
    archivePrefix={arXiv},
    primaryClass={math.AP}
}

@article{mao_2013,
author = {Mao, Y. and Tournier, A. L and Hoppe, A. and Kester, L. and Thompson, B. J and Tapon, N.},
title = {Differential proliferation rates generate patterns of mechanical tension that orient tissue growth},
journal = {The EMBO Journal},
volume = {32},
number = {21},
pages = {2790-2803},
doi = {https://doi.org/10.1038/emboj.2013.197},
year = {2013}
}

@article{voutouri_2014,
author = {Voutouri, C. and Mpekris, F. and Papageorgis, P. and Odysseos, A. D. and Stylianopoulos, T.},
title = {Role of constitutive behavior and tumor-host mechanical interactions in the state of stress and growth of solid tumors},
journal = {PLoS ONE},
volume = {9},
number = {8},
year = {2014},
pages = {e104717},
doi = {10.1371/journal.pone.0104717}
}

@article{nia_2016,
  title={Solid stress and elastic energy as measures of tumour mechanopathology},
  author={Nia, H. T. and Liu, H. and Seano, G. and Datta, M. and Jones, D. and Rahbari, N. and Incio, J. and Chauhan, V. P. and Jung, K. and Martin, J. D. and others},
  journal={Nature Biomedical Engineering},
  volume={1},
  number={1},
  pages={0004},
  year={2016},
  publisher={Nature Publishing Group UK London}
}

@article{delarue_2014,
title = {Compressive Stress Inhibits Proliferation in Tumor Spheroids through a Volume Limitation},
journal = {Biophysical Journal},
volume = {107},
number = {8},
pages = {1821-1828},
year = {2014},
issn = {0006-3495},
doi = {https://doi.org/10.1016/j.bpj.2014.08.031},
author = {M. Delarue and F. Montel and D. Vignjevic and J. Prost and J.-F. Joanny and G. Cappello},
}

@article{helmlinger_1997,
  title={Solid stress inhibits the growth of multicellular tumor spheroids},
  author={Helmlinger, G. and Netti, P. A. and Lichtenbeld, H. C. and Melder, R. J. and Jain, R. K.},
  fjournal={Nature Biotechnology},
  journal = {Nat Biotechnol},
  volume={15},
  number={8},
  pages={778--783},
  year={1997},
  publisher={Nature Publishing Group US New York}
}

@article{legoff_2016,
  title={Mechanical forces and growth in animal tissues},
  author={LeGoff, L. and Lecuit, T.},
  journal={Cold Spring Harbor perspectives in biology},
  volume={8},
  number={3},
  pages={a019232},
  year={2016},
  publisher={Cold Spring Harbor Lab}
}

@article{forgacs_1998,
title = {Viscoelastic Properties of Living Embryonic Tissues: a Quantitative Study},
journal = {Biophysical Journal},
volume = {74},
number = {5},
pages = {2227--2234},
year = {1998},
doi = {https://doi.org/10.1016/S0006-3495(98)77932-9},
author = {G. Forgacs and R. A. Foty and Y. Shafrir and M. S. Steinberg},
}

@article{doubrovinski_2017,
  title={Measurement of cortical elasticity in Drosophila melanogaster embryos using ferrofluids},
  author={Doubrovinski, K. and Swan, M. and Polyakov, O. and Wieschaus, E. F.},
  journal={Proceedings of the National Academy of Sciences},
  volume={114},
  number={5},
  pages={1051--1056},
  year={2017},
  publisher={National Acad Sciences}
}

@article{northcott_2018,
  title={Feeling stress: the mechanics of cancer progression and aggression},
  author={Northcott, J. M. and Dean, I. S. and Mouw, J. K. and Weaver, V. M.},
  journal={Frontiers in cell and developmental biology},
  volume={6},
  pages={17},
  year={2018},
  publisher={Frontiers Media SA}
}

@article{lucci_agosti_2022,
  title={Coupling solid and fluid stresses with brain tumour growth and white matter tract deformations in a neuroimaging-informed model},
  author={Lucci, G. and Agosti, A. and Ciarletta, P. and Giverso, C.},
  journal={Biomechanics and Modeling in Mechanobiology},
  volume={21},
  number={5},
  pages={1483--1509},
  year={2022},
  publisher={Springer}
}

@article {gruen_klingbeil_2014,
    AUTHOR = {Gr\"{u}n, G. and Klingbeil, F.},
     TITLE = {Two-phase flow with mass density contrast: stable schemes for
              a thermodynamic consistent and frame-indifferent
              diffuse-interface model},
   JOURNAL = {J. Comput. Phys.},
  FJOURNAL = {Journal of Computational Physics},
    VOLUME = {257},
      YEAR = {2014},
    NUMBER = {part A},
     PAGES = {708--725},
      ISSN = {0021-9991},
   MRCLASS = {65M60 (65M12 76M10 76T99)},
  MRNUMBER = {3129556},
MRREVIEWER = {Maurizio Brocchini},
       DOI = {10.1016/j.jcp.2013.10.028},
       URL = {https://doi.org/10.1016/j.jcp.2013.10.028},
}

@article {constantin_kliegl_2012,
    AUTHOR = {Constantin, P. and Kliegl, M.},
     TITLE = {Note on global regularity for two-dimensional {O}ldroyd-{B}
              fluids with diffusive stress},
   JOURNAL = {Arch. Ration. Mech. Anal.},
  FJOURNAL = {Archive for Rational Mechanics and Analysis},
    VOLUME = {206},
      YEAR = {2012},
    NUMBER = {3},
     PAGES = {725--740},
   MRCLASS = {76A10 (35B65 35Q84)},
  MRNUMBER = {2989441},
MRREVIEWER = {Thomas Hagen},
       DOI = {10.1007/s00205-012-0537-0},
}

@article{frigeri_grasselli_rocca_2015, 
title={On a diffuse interface model of tumour growth}, 
volume={26}, 
DOI={10.1017/S0956792514000436}, 
number={2}, 
fjournal={European Journal of Applied Mathematics}, 
journal={European J. Appl. Math.},
publisher={Cambridge University Press}, 
author={Frigeri, S. and Grasselli, M. and Rocca, E.}, 
year={2015}, 
pages={215--243}
}

@article{agosti_ACGV_2017,
author = {Agosti, A. and Antonietti, P. F. and Ciarletta, P. and Grasselli, M. and Verani, M.},
title = {A Cahn-Hilliard–type equation with application to tumor growth dynamics},
journal = {Mathematical Methods in the Applied Sciences},
volume = {40},
number = {18},
pages = {7598-7626},
doi = {https://doi.org/10.1002/mma.4548},
year = {2017}
}

@article{agosti_CGAC_2018,
author = {Agosti, A. and Cattaneo, C. and Giverso, C. and Ambrosi, D. and Ciarletta, P.},
title = {A computational framework for the personalized clinical treatment of glioblastoma multiforme},
journal = {ZAMM - Journal of Applied Mathematics and Mechanics / Zeitschrift für Angewandte Mathematik und Mechanik},
volume = {98},
number = {12},
pages = {2307-2327},
doi = {https://doi.org/10.1002/zamm.201700294},
year = {2018}
}

@book {elman_silvester_wathen_2014,
    AUTHOR = {Elman, H. C. and Silvester, D. J. and Wathen, A.
              J.},
     TITLE = {Finite elements and fast iterative solvers: with applications
              in incompressible fluid dynamics},
    SERIES = {Numerical Mathematics and Scientific Computation},
   EDITION = {Second},
 PUBLISHER = {Oxford University Press, Oxford},
      YEAR = {2014},
     PAGES = {xiv+479},
       DOI = {10.1093/acprof:oso/9780199678792.001.0001},
}

@article{Chen_Lowengrub_2014,
title = {Tumor growth in complex, evolving microenvironmental geometries: A diffuse domain approach},
journal = {J. Theoret. Biol.},
fjournal = {Journal of Theoretical Biology},
volume = {361},
pages = {14-30},
year = {2014},
issn = {0022-5193},
doi = {https://doi.org/10.1016/j.jtbi.2014.06.024},
author = {Y. Chen and J. S. Lowengrub},
}

@article{brunk_2022_3d,
  title={Existence, regularity and weak-strong uniqueness for three-dimensional Peterlin viscoelastic model},
  author={Brunk, A. and Lu, Y. and Luk{\'a}{\v{c}}ov{\'a}-Medvi{\v{d}}ov{\'a}, M.},
  journal={Communications in Mathematical Sciences},
  volume={20},
  number={1},
  pages={201--230},
  year={2022},
  publisher={International Press of Boston}
}

@article{dostalik_prusa_skrivan_2019,
    author = {Dostal{\'i}k, M. and Pr{\r u}{\v s}a, V. and Sk{\v r}ivan, T.},
    title = "{On diffusive variants of some classical viscoelastic rate-type models}",
    journal = {AIP Conference Proceedings},
    volume = {2107},
    number = {1},
    year = {2019},
    month = {05},
    doi = {10.1063/1.5109493},
    note = {020002},
}

@article{chupin_2018,
author = {Chupin, L.},
title = {Global Strong Solutions for Some Differential Viscoelastic Models},
journal = {SIAM Journal on Applied Mathematics},
volume = {78},
number = {6},
pages = {2919-2949},
year = {2018},
doi = {10.1137/18M1186873},
}

\end{document}